\documentclass[12pt,twoside]{amsart}
\usepackage{amssymb}
\usepackage{amscd}
\usepackage[abbrev,alphabetic]{amsrefs}
\usepackage{hyperref}
\usepackage{comment}
\usepackage{array,multirow,tabularx,longtable}
\usepackage{tikz}
\usetikzlibrary{cd}
\usepackage{here}
\usepackage{multirow}
\usepackage[margin=1.25in]{geometry}

\usepackage{fancybox}
\usepackage{ascmac}

\title[Fano threefolds in positive characteristic IV]
{Fano threefolds in positive characteristic IV} 
\author{Hiromu Tanaka} 
\address{Department of Mathematics, 
Graduate School of Science, 
Kyoto University, 
Kyoto 606-8502, JAPAN} 
\email{tanaka.hiromu.7z@kyoto-u.ac.jp}
\newcommand{\pr}[0]{{\operatorname{pr}}}
\newcommand{\Blowdown}[0]{{\operatorname{Blowdown}}}
\newcommand{\Blowup}[0]{{\operatorname{Blowup}}}
\newcommand{\Bl}[0]{{\operatorname{Bl}}}

\newcommand{\Br}[0]{{\operatorname{Br}}}

\newcommand{\NE}[0]{{\operatorname{NE}}}
\newcommand{\red}[0]{{\operatorname{red}}}

\newcommand{\Ker}[0]{{\operatorname{Ker}}}

\renewcommand{\Im}[0]{{\operatorname{Im}}}

\newcommand{\Spec}[0]{{\operatorname{Spec}}}

\newcommand{\Bs}[0]{{\operatorname{Bs}}}
\newcommand{\Supp}[0]{{\operatorname{Supp}}}
\newcommand{\Pic}[0]{{\operatorname{Pic}}}

\newcommand{\Ex}[0]{{\operatorname{Ex}}}

\newcommand{\Ext}[0]{{\operatorname{Ext}}}

\newtheorem{thm}{Theorem}[section]
\newtheorem{lem}[thm]{Lemma}
\newtheorem{cor}[thm]{Corollary}
\newtheorem{prop}[thm]{Proposition}
\newtheorem{claim}[thm]{Claim}    
    
\newtheorem*{claim*}{Claim}         
\newtheorem{step}{Step}

\theoremstyle{definition}

\newtheorem{dfn}[thm]{Definition}

\newtheorem{rem}[thm]{Remark}
      
\newtheorem{nota}[thm]{Notation}         
\newtheorem{nasi}[thm]{}

\makeatletter
  
  \@addtoreset{equation}{thm}
  \makeatother

\newcommand{\cred}{\color{red}}

\newcommand{\MO}{\mathcal{O}}

\newcommand{\R}{\mathbb{R}}
\newcommand{\Q}{\mathbb{Q}}
\newcommand{\Z}{\mathbb{Z}}
\newcommand{\F}{\mathbb{F}}

\renewcommand{\P}{\mathbb{P}}

\newcommand{\m}{\mathfrak{m}}

\newcommand{\wt}{\widetilde}
\begin{document}

\maketitle

\begin{abstract}
Based on 
the former parts, 
we classify smooth Fano threefolds of positive characteristic. 
\end{abstract}

\tableofcontents

\section{Introduction}

This article  is the last part of our series of papers. 
In the former parts \cite{TanI}, \cite{TanII}, \cite{ATIII}, 
we have settled the classification of 
Fano threefolds in positive characteristic 
for the following cases: 
\begin{itemize}
\item 
Fano threefolds with $\rho(X) \leq 2$. 
\item 
Primitive Fano threefolds. 
\end{itemize}
Moreover, we  have proved that 
\begin{itemize}
\item 
if $X$ is a Fano threefold with $\rho(X) = r_X =1$ and $|-K_X|$ is very ample, 
then there is no smooth curve $\Gamma$ on $X$ such that the blowup of $X$ along $\Gamma$ is Fano. 
\end{itemize}
Based on these results, we complete the classification of Fano threefolds in positive characteristic. 
More precisely, the main theorem is as follows. 



\begin{thm}
Let $k$ be an algebraically closed field of characteristic $p>0$ and 
let $X$ be a Fano threefold over $k$, i.e., $X$ is a three-dimensional smooth projective variety over $k$ such that $-K_X$ is ample. 
Then $X$ is isomorphic to one of threefolds listed in Section \ref{s-table}. 
For example, if $\rho(X)=3$, then $X$ is one of No.\ 3-1, 3-2, ...., 3-31 
in Subsection \ref{ss-table-pic3}. 
\end{thm}

The classification list is almost  identical to that of characteristic zero \cite{MM81}, \cite{MM03}, \cite[Chapter 12]{IP99}. 
We here point out the following differences. 

\begin{enumerate}
\renewcommand{\labelenumi}{(\roman{enumi})}
\item 
Given a Fano threefold $X$ and a double cover $f: X \to Y$, 
the classification tables in characteristic zero give the branched divisor \cite{MM81}. 
The branched divisor behaves differently in characteristic two, 
whilst the description of the invertible sheaf $\mathcal L := (f_*\MO_X/\MO_Y)^{-1}$ is of characteristic free. 
Thus, instead of the branched divisor, 
our classification tables (Section \ref{s-table}) include what $\mathcal L$ is. 
\item 
Let $X$ be a Fano threefold and let $f: X \to S$ be a conic bundle.  
There are the following two  phenomena which  occur only in characteristic two. 
\begin{enumerate}
\item $f$ might be a  wild conic bundle, i.e., no fibre of $f$ is smooth. 
Such a phenomenon happens only when $X$ is of No.\ 2-24 or 3-10 (Theorem \ref{t-wild-cb}). 
\item Even if $f$ is generically smooth, the discriminant divisor $\Delta_f$ might be non-reduced. 
In the end, this phenomenon has hardly  affected the proof of our classification result. 
\end{enumerate}
\end{enumerate}


\subsection{Overviews}

In this paper, we follow the original strategy by Mori-Mukai \cite{MM83}. 
Recall that the classification for the case when $\rho(X) \leq 2$ is carried out 
in the former parts \cite{TanI}, \cite{TanII}, \cite{ATIII}. 
Hence it is enough to treat the case when $\rho(X) \geq 3$. 
The key tool is the notion of conic bundles.  
The following results illustrate its significance. 
\begin{enumerate}
\item[(A)] A Fano threefold $X$ has a conic bundle structure if 
$\rho(X) \geq 3$ and $X$ is not of No.\ 3-18 (Theorem \ref{t-pic3-structure}, Theorem \ref{t-pic4-CB}). 
\item[(B)] 
If $X$ is a Fano threefold and $f: X \to S$ is a conic bundle, 
then $S$ is a smooth del Pezzo surface (Proposition \ref{p-FCB-dP}). 
Moreover, $X \simeq S \times \P^1$ if $K_S^2 \leq 7$ (Proposition \ref{p-FCB-triv}). 
\end{enumerate}
By using (B), we will see that $X \simeq S \times \P^1$ when $\rho(X) \geq 6$. 
In what follows, we mainly explain the proof of our classification for the case when $\rho(X)=3$, 
as the problem is simpler for the case when $\rho(X) \geq 4$. 

\medskip

Let $X$ be a Fano threefold with $\rho(X)=3$. 
The first step for the classification is to prove that 
one of (I)--(V) holds (Subsection \ref{ss-pic3-structure}). 
\begin{enumerate}
\item[(I)] $X$ has a conic bundle structure over $\P^2$. 
\item[(II)] $X$ has a conic bundle structure over $\F_1$. 
\item[(III)] $X$ is primitive. 
\item[(IV)] $X$ is of No.\ 3-18. 
\item[(V)] 
$X$ is of No.\ 3-6, 3-10, or 3-25. 
\end{enumerate}
Since the cases (III)--(V) are classified, 
it is enough to treat the remaining cases (I) and (II). 

(I) Assume that $X$ is a Fano threefold and $f: X \to S =\P^2$ is a conic bundle with $\rho(X)=3$. 
Then $X$ is obtained by taking a blowup of $Y$ along a regular subsection of a Fano conic bundle $g : Y \to S=\P^2$ with $\rho(Y)=2$ (Lemma \ref{l-FCB-P^2-1}). 
Here $B_Y$ is called a subsection $B_Y$ of $g$ if 
$B_Y$ is a smooth curve on $Y$ such that the induced morphism to $S$ is a 
closed immersion. 
In this case, we obtain the elementary transform $Y'$: 
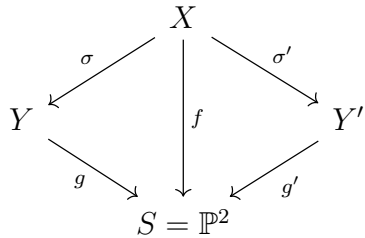
\[
\begin{tikzcd}
& X \arrow[rd, "\sigma'"] \arrow[ld, "\sigma"'] \arrow[dd, "f"]\\
Y \arrow[rd, "g"'] & & Y' \arrow[ld, "g'"]\\
&S=\P^2
\end{tikzcd}
\]
where $\sigma'$ is a blowup along a smooth curve and $g'$ is another conic bundle. 
There are two cases: either $Y'$ is Fano or not. 
If $Y'$ is Fano, then there are only finitely many possibilities for 
classes of $Y$ and $Y'$. 
If $Y'$ is not Fano, then we have a strong restriction: $(-K_{Y/S}) \cdot B_Y = 2(B^2+1)$. 
Mainly by using these facts, we get a   classification list of the triples $(X, Y, Y')$ (Theorem \ref{t-ele-tr-P2}). 
For more details, see Subsection \ref{ss-pic3-P2}.

(II) Assume that $X$ is a Fano threefold and $f: X \to S =\F_1$ is a conic bundle with $\rho(X)=3$. 
In this case, $X$ is obtained by applying to $(-) \times_{\P^2} \F_1$ 
to another Fano threefold $\widetilde X$ 
with a conic bundle structure $\widetilde f : \widetilde X \to \widetilde S = \P^2$, 
where the blowup centre of $\widetilde{X} \to X$ is a smooth fibre of $\widetilde{f}$ (Lemma \ref{l-F1-pic3}). 
By $\rho(\widetilde{X})=2$, such a conic bundle $\widetilde f : \widetilde X \to \widetilde S = \P^2$ has been already classified. 
Then it is not so hard to get a classification 
of the pairs $(X, X')$ (Theorem \ref{t-F1-pic3}). 

\medskip

By the argument as above, 
we obtain a list of possibilities for Fano threefolds $X$ with $\rho(X)=3$. 
However, we need to care the overlapping problem, e.g., 
two Fano threefolds $X_{{\rm (I)}}$ and $X_{{\rm (II)}}$ obtained from (I) and (II) might coincide. 
In order to settle this problem, 
we shall determine all the contractions for every case (Subsection \ref{ss-pic3-classify}). 
As a consequence,  it holds that $X$ has exactly three extremal rays unless $X$ is of No.\ 3-9, 3-14, 3-19.

\medskip

Let $X$ be a Fano threefold with $\rho(X)=4$. 
In this case, 
we first prove that there is a conic bundle structure $f: X \to S$ with $S \in \{ \P^1 \times \P^1, \F_1\}$ (Subsection \ref{ss-pic4-structure}). 
Similar to the case $\rho(X)=3$, 
we shall classify such conic bundles 
(Subsection \ref{ss-pic4-P1P1}, Subsection \ref{ss-pic4-F1}). 
Contrary to the case when $\rho(X)=3$, 
the overlapping problem is settled by introducing the set $\Blowdown(X)$, 
which consists of the No.\ 3-xx appearing as a blowdown from $X$. 

Given the classification of Fano threefolds with $\rho(X)=4$, 
it is not so hard to classify Fano threefolds $X$ with $\rho(X)=5$. 
Indeed, we can show that 
either $X \simeq S \times \P^1$ for a smooth del Pezzo surface $S$ with $K_S^2 = 6$ 
or $X \simeq Y \times_Z Y'$, where 
$Z \to \F_1$ is a Fano conic bundle and
each of $Y \to Z$ and $Y' \to Z$ is a  blowup along a regular subsection of $Z \to \F_1$. 
Moreover, the blowup centres of $Y \to Z$ and $Y' \to Z$ are mutually disjoint, 
and there appear only three possibilities for the triple $(Y, Y', Z)$  (Lemma \ref{l-pic5-F1}).

\begin{rem}
As mentioned above, we shall determine all the contraction morphisms 
for the case when $X$ is a Fano threefold with $\rho(X)=3$. 
Matsuki establishes the same result in characteristic zero 
even without any assumption on $\rho(X)$ \cite[The list starting from page 50]{Mat86}. 
\end{rem}

\subsection{Related questions}




\subsubsection{Mukai's description}

Let $X$ be a Fano threefold such that $\Pic\,X  = \Z K_X$ and $|-K_X|$ is very ample. 
Then we have the closed embedding $X \subset \P^{g+1}$ induced by $|-K_X|$, 
where  $g$ satisfies $3 \leq g \leq 12$ and $g \neq 11$. 
If $g \leq 5$, then $X$ is a complete intersection in $\P^{g+1}$ (Table \ref{table-pic1} in Subsection \ref{ss-table-pic1}). 
In characteristic zero, Mukai predicted concrete descriptions 
even for the case $g \geq 6$,  
which has been 
confirmed by Bayer-Kuznetsov-Macri \cite[Theorem 1.1]{BKM}. 
It is an open question whether this result is extended to the case of positive characteristic.

\subsubsection{Number of deformation families} 

In characteristic zero, it is known that Fano threefolds form 
exactly 105 deformation families. 
Our classification table in Section \ref{s-table} 
is identical to that of characteristic zero \cite{MM81}. 
On the other hand, our series of papers do not 
provide the number of deformation families in positive characteristic.

\vspace{5mm}

\textbf{Acknowledgements:} 
The author would like to thank Kento Fujita and 
Tatsuro Kawakami 
for answering questions. 
The author was funded by JSPS KAKENHI Grant numbers JP22H01112 and JP23K03028. 

\section{Preliminaries}

\subsection{Notation}\label{ss-notation}

In this subsection, we summarise notation used in this paper. 

\begin{enumerate}
\item We will freely use the notation and terminology in \cite{Har77} and \cite{KM98}. 
In particular, $D_1 \sim D_2$ means linear equivalence of Weil divisors. 
\item 
Throughout this paper, 
we work over an algebraically closed field $k$ 
of characteristic $p>0$ unless otherwise specified. 
\item For an integral scheme $X$, 
we define the {\em function field} $K(X)$ of $X$ 
as the local ring $\MO_{X, \xi}$ at the generic point $\xi$ of $X$. 
For an integral domain $A$, $K(A)$ denotes the function field of $\Spec\,A$. 

\item 
For a scheme $X$, its {\em reduced structure} $X_{\red}$ 
is the reduced closed subscheme of $X$ such that the induced closed immersion 
$X_{\red} \to X$ is surjective. 
\item 
Our notation will not distinguish between invertible sheaves and 
Cartier divisors. For example, we will write $L+D$ for 
an invertible sheaf $L$ and a Cartier divisor $D$. 
\item We say that $X$ is a {\em variety} (over $k$) if 
$X$ is a separated integral scheme which is of finite type over $k$. 
We say that $X$ is a {\em curve} (resp. a {\em surface}, resp. a {\em threefold})  
if $X$ is a variety over $k$ of dimension one (resp. {\em two}, resp. {\em three}). 
\item 
Given a variety $Y$ and a closed subscheme $Z$ of $Y$, 
$\Bl_Z\,Y$ denotes the blowup of $Y$ along $Z$. 
In this case, $Z$ is called {\em the (blowup) centre} of the induced blowup $\Bl_Z\,Y \to Y$. 
Let  $\Ex(f)$ be the exceptional divisor 
equipped with reduced scheme structure. 
In particular, if $Y$ is a smooth threefold and $Z$ is a smooth curve on $Y$, 
then we have $K_X \sim f^*K_Y +\Ex(f)$ for $X := \Bl_Z\,Y$.
\item We say that $f: X \to Y$ be a {\em contraction} if $f$ is a morphism of schemes satisfying $f_*\MO_X = \MO_Y$. 
Here the equality $f_*\MO_X = \MO_Y$ means that the induced ring homomorphism $\MO_Y \to f_*\MO_X$ is an isomorphism. 
\item 
 We say that $X$ is a {\em Fano threefold}  
if $X$ is a three-dimensional smooth projective variety over $k$ such that $-K_X$ is ample. 
A Fano threefold $X$ is {\em imprimitive} if 
there exists a Fano threefold $Y$ and a smooth curve $B$ on $Y$ such that 
$X$ is isomorphic to the blowup $\Bl_B Y$ of $Y$ along $B$. 
We say that a Fano threefold $X$ is {\em primitive} if $X$ is not imprimitive. 
\item For a Fano threefold $X$, the {\em index} $r_X$ of $X$ 
is the largest positive integer $r$ 
that divides $-K_X$ in $\Pic\,X$. 
For a curve $C$ on $X$, the positive integer $-\frac{1}{r_X}K_X \cdot C$ is called 
the {\em degree} of $C$ (on $X$). 
If $|H|$ is very ample for a Cartier divisor $H$ satisfying $-K_X \sim r_X H$, 
then a curve of degree one (resp. two) on $X$ is called a {\em line} (resp. {\em conic}). 
\item 
For the definition of types of extremal rays for smooth projective threefolds, 
we refer to \cite[Definition 3.3]{ATIII}. 
\item $Q$ denotes a smooth quadric hypersurface on $\P^4$. 
For $1 \leq d \leq 5$ and $d =7$, let $V_d$ be a Fano threefold of index two such that $(-K_{V_d}/2)^3 =d$. 
Let $W$ be a smooth divisor on $\P^2 \times \P^2$ of bidegree $(1, 1)$. 
Note that such a threefold $W$ is unique up to isomorphisms \cite[Lemma 5.16]{ATIII}. 
\item For a Fano threefold $Y$ with $\rho(Y)=1$, 
$\MO_Y(1)$ denotes an invertible sheaf which  generates $\Pic\,Y (\simeq \Z)$ and 
we set $\MO_Y(\ell) := \MO_Y(1)^{\otimes \ell}$. 
\item $\F_m := \P_{\P^1}(\MO_{\P^1} \oplus \MO_{\P^1}(m))$ for every $m \in \Z_{\geq 0}$. 
$\tau : \F_1 \to \P^2$ often denotes the blowdown of the $(-1)$-curve of $\F_1$. 

\item 
Given a closed subscheme $Z$ on $\P^N_k$, 
we define $\langle Z \rangle$ as the smallest linear subvariety of $\P^N_k$ containing $Z$. 
\item 
Each of $Y_{\text{a-bc}}, Y'_{\text{a-bc}}$, and $Z_{\text{a-bc}}$ 
denotes a Fano threefold of No.\ a-bc. 
The definitions of No.\ a-bc will be given in 
Definition \ref{d-pic2} ($\rho=2$), 
Definition \ref{d-pic3} ($\rho=3$), 
Definition \ref{d-pic4} ($\rho=4$), 
Definition \ref{d-pic5} ($\rho=5$).
In many theorems, we write something like \lq\lq ($X$ is a-bc)", e.g., 
\lq\lq ($X$ is 3-21)" appears in Lemma \ref{l-P2-nonFano}. 
This is not a rigorous statement but added just for convenience. 
For example, the rigorous definition of No.\ 3-21 (Definition \ref{d-pic3}) will be given in the latter part than Lemma \ref{l-P2-nonFano}. 
\item 
\begin{itemize}
\item We say that a divisor $D$ on $\P^a \times \P^b \times \P^c$ 
is of {\em tridegree} $(d_1, d_2, d_3)$ if 
$\MO_{\P^a \times \P^b \times \P^c}(D) \simeq \MO_{\P^a \times \P^b \times \P^c}(d_1, d_2, d_3)$. 
We define the {\em bidegree} of a divisor on $\P^a \times \P^b$ in a similar way. 
\item 
We say that a curve $B$ on $\P^a \times \P^b \times \P^c$ 
is of {\em tridegree} $(d_1, d_2, d_3)$ if 
$\pr^*_i\MO(1) \cdot B = d_i$ for every $i \in \{1, 2, 3\}$. 
When $(a, b) \neq (1, 1)$, 
we define the {\em bidegree} of a curve on $\P^a \times \P^b$ in a similar way.   
\end{itemize}
\end{enumerate}

\begin{rem}\label{r curve bidegree}
Let $f : \P^1_1 \times \P^1_2 \times \P^1_3 \to \P^1_1 \times \P^1_2$ 
be the projection onto the first and second direct product factors. 
Take a curve $B$ on $\P^1_1 \times \P^1_2 \times \P^1_3$ such that 
the induced morphism $B \to B' := f(B)$ is an isomorphism. 
If $B$ is of tridegree $(d_1, d_2, d_3)$, then $B'$ is of bidegree $(d_2, d_1)$. 
\end{rem}

\subsection{Case of relative Picard number two}

Let $\varphi : X \to Z$ be a contraction of projective normal varieties. 
Then $N_1(X/Z)$ is the $\R$-linear subspace of $N_1(X)$ defined by 
\[
N_1(X/Z) := \{ [a_1 C_1 + \cdots +a_mC_m] \in N_1(X) \,|\, 
\]
\[
m \in \Z_{>0}, 
a_i \in \R, C_i \text{ is a curve such that }\varphi(C_i)\text{ is a point for every }i\}. 
\]
We set $N^1(X/Z) := (\Pic\,X \otimes_{\Z} \R) /\equiv_Z$, where $\equiv_Z$ denotes the numerical equivalence over $Z$, 
i.e., $L \equiv_Z L'$ is defined by $L \cdot C = L' \cdot C$ for any $C \in N_1(X/Z)$. 
We then obtain $\rho(X/Z) := \dim_{\R} N_1(X/Z) = \dim_{\R} N^1(X/Z)$. 
Set 
\[
\NE(X/Z) := \{ [a_1 C_1 + \cdots +a_mC_m] \in N_1(X) \,|\, 
\]
\[
m \in \Z_{>0}, 
a_i \in \R_{\geq 0}, C_i \text{ is a curve such that }\varphi(C_i)\text{ is a point for every }i\}. 
\]

\begin{rem}\label{r-rel-pic}
We have a sequence 
\begin{equation}\label{e1-rel-pic}
0\to N^1(Z) \xrightarrow{\alpha} N^1(X) \xrightarrow{\beta} N^1(X/Z) \to 0
\end{equation}
of $\R$-linear maps of $\R$-vector spaces 
for $\alpha := \varphi^*$ and the natural surjection $\beta$. 
It is easy to see that $\alpha$ is injective and $\beta \circ \alpha =0$. 
In particular, we obtain a surjective $\R$-linear map $\overline{\beta} : N^1(X)/N^1(Z) \to N^1(X/Z)$. 
Hence 
\[
\rho(X) - \rho(Z) =\dim_{\R} (N^1(X) / N^1(Z)) \geq \dim N^1(X/Z) = \rho(X/Z).
\]
Then the following are equivalent. 
\begin{enumerate}
    \item The sequence (\ref{e1-rel-pic}) is exact. 
    \item $\Im (\alpha) \supset \Ker(\beta)$. 
    \item $\rho(X/Z) = \rho(X) - \rho(Z)$. 
\end{enumerate}

\end{rem}

\begin{lem}\label{l-rel-pic}
Let $X$ be a Fano threefold and let $\varphi : X \to Z$ be a contraction 
to a projective normal variety $Z$. 
Then the following hold. 
\begin{enumerate}
\item 
The natural surjection $\Pic\,Z \to (\Pic\,Z)/\equiv$ is an isomorphism. 
\item 
If $\varphi$ is birational or $\varphi$ is a conic bundle, then $\rho(X/Z) = \rho(X) - \rho(Z)$. 
\item 
If $\rho(X) = \rho(Z) +2$, then $\rho(X/Z)=2$. 
\end{enumerate}
\end{lem}

\begin{proof}
Let us show (1). 
Fix a Cartier divisor $L$ on $Z$ with $L \equiv 0$. 
Then $\varphi^*L \equiv 0$, which implies $\varphi^*L \sim 0$. 
Hence $H^0(Z, L) \simeq H^0(X, \varphi^*L) \neq 0$, which implies $L \sim 0$. 
Thus (1) holds.

Let us show (2). 
By Remark \ref{r-rel-pic}, 
it suffices to prove Remark \ref{r-rel-pic}(2). 
Fix a Cartier divisor $L_X$ on $X$ with $L_X \equiv_Z 0$. 
It is enough to find an integer $m>0$ and a Cartier divisor $L_Z$ 
satisfying $mL_X \sim \varphi^*L_Z$. 
For an ample Cartier divisor $A_Z$ on $Z$, 
$L_X + m \varphi^*A_Z$ is nef for $m \gg 0$ by the cone theorem. 
If $\varphi$ is a conic bundle, then 
we are done by 
\cite[Theorem 1.1]{CT20}. 
Hence we may assume that $\varphi$ is birational. 
Since $\varphi^*A_Z$ is big, 
we may assume that $L_X + m\varphi^*A_Z$ is nef and big. 
Then $L_X +m\varphi^*A_Z$ is semi-ample by \cite[Theorem H]{BMPSTWW}. 
By $L_X + m\varphi^*A_Z \equiv_Z 0$, we can find $n \in \Z_{>0}$ and a Cartier divisor $L_Z$ such that $n(L_X + m\varphi^*A_Z) \sim \varphi^*L_Z$. 
Thus Remark \ref{r-rel-pic}(2) holds. 
This completes the proof of (2).


Let us show (3). 
Assume $\rho(X) = \rho(Z) +2$. 
Then $\rho(X/Z) \leq \rho(X) -\rho(Z) =2$ (Remark \ref{r-rel-pic}). 
Hence it suffices to show $\rho(X/Z) >1$. 
For an ample Cartier divisor $A_Z$ on $Z$, 
we have $\NE(X) \cap (\varphi^*A_Z)^{\perp} \neq \{0\}$, 
because there exists a curve $C$ on $X$ contracted by $\varphi :X \to Z$. 
Take an extremal ray $R$ of 
$\NE(X) \cap (\varphi^*A_Z)^{\perp}$, 
which is automatically an extremal ray of $\NE(X)$. 
Let $f: X \to Y$ be the contraction of $R$. 
Then we get the facotorisation: $\varphi : X \xrightarrow{f} Y \xrightarrow{g} Z$. 
Hence we obtain another natural surjection $\gamma : N^1(X/Z) \to N^1(X/Y)$. 
This is not injective, because $\gamma(f^*A_Y)=0$ for an ample Cartier divisor $A_Y$ on $Y$.  
Hence $\rho(X/Z) = \dim_{\R} N^1(X/Z) > \dim_{\R} N^1(X/Y)=1$. 
Thus (3) holds. 
\qedhere





\end{proof}

\begin{prop}\label{p-2ray}
Let $X$ be a Fano threefold and let $\varphi : X \to Z$ be a contraction 
to a projective normal variety $Z$ such that $\rho(X/Z) = 2$ (e.g.,  $\rho(X) -\rho(Z)=2$). 
Then there exist two extremal rays $R$ and $R'$ of $\NE(X)$ and the following diagram consisting of contractions of projective normal varieties 
\begin{equation}\label{e1-2ray}
\begin{tikzcd}
    & X \arrow[ld, "f"'] \arrow[rd, "f'"] \arrow[dd, "\varphi"]\\
    Y \arrow[rd, "g"'] & & Y' \arrow[ld, "g'"]\\
    & Z
\end{tikzcd}
\end{equation}
such that 
\begin{enumerate}
\item $f : X \to Y$ and $f' : X \to Y'$ are the contractions of $R$ and $R'$ respectively, and 
\item $F_Z = R+R'$ for the extremal face $F_Z$ of $\NE(X)$ corresponding to $\varphi$.  
\end{enumerate}
\end{prop}


\begin{proof}
For an ample Cartier divisor $A_Z$ on $Z$, we have 
\[
F_Z = \NE(X) \cap (\varphi^* A_Z)^{\perp} = \NE(X/Z). 
\]
Then $F_Z$ is a rational polyhedral cone, because so is $\NE(X)$. 
Recall that $\NE(X/Z)$ generates $N_1(X/Z)$, which is of dimension $\rho(X/Z)=2$. 
Therefore, 
we obtain $F_Z = R +R'$ for distinct extremal rays $R$ and $R'$ 
of $\NE(X)$. 
Let $f:X \to Y$ and $f':X \to Y'$ be the contractions of $R$ and $R'$, respectively. 
By construction, we get the commutative diagram (\ref{e1-2ray}). 
\qedhere

\end{proof}

\subsection{Non-Fano blowdowns}

\begin{lem}\label{l-blowup-formula}
Let $Y$ be a smooth projective threefold. 
Take a smooth curve $C$ on $Y$ and let $\sigma : X \to Y$ be the blowup along $C$. Set $D := \Ex(\sigma). $
Then the following hold. 
\begin{enumerate}
\item 
$K_X \sim \sigma^*K_Y +D$. 
\item 
$(-K_X)^3 = (-K_Y)^3 -2 (-K_Y) \cdot C +2p_a(C) -2$. 
\item 
$(-K_X)^2 \cdot D = (-K_Y) \cdot C -2p_a(C) +2$. 
\item 
$(-K_X) \cdot D^2 =2p_a(C) -2$. 
\item 
$D^3 = -\deg N_{C/Y} = 
-(-K_Y) \cdot C -2p_a(C) +2$. 
\end{enumerate}
\end{lem}

\begin{proof}
See, for example, \cite[Lemma 3.21]{TanII}. 
\end{proof}

\begin{lem}\label{l-nonFano-blowdown}
Let $X$ be a Fano threefold and let $f : X \to Y$ be a blowup along a smooth curve $\Gamma$ on $Y$  such that $Y$ is a non-Fano smooth projective threefold. 
Set $D := \Ex(f)$. 
Then the following holds. 
\begin{enumerate}
\item $\Gamma \simeq \P^1$. 
\item $N_{\Gamma/Y} \simeq \MO_{\P^1}(-1) \oplus \MO_{\P^1}(-1)$. 
\item $D \simeq \P^1 \times \P^1$. 
\item $\MO_X(K_X)|_D \simeq \MO_{\P^1 \times \P^1}(-1, -1)$ and 
$\MO_X(D)|_D\simeq \MO_{\P^1 \times \P^1}(-1, -1)$. 
\item $-K_Y$ is semi-ample. 
\end{enumerate}
\end{lem}

\begin{proof}
The assertions (1)-(4) follow from \cite[Lemma 4.7]{ATIII}. 
Let us show (5). 
By $-K_X+D \sim f^*(-K_Y)$, 
it is enough to show that $-K_X+D$ is semi-ample. 
As (4) implies $(-K_X +D)|_D \sim 0$, 
$-K_X+D$ is semi-ample by \cite[Corollary 3.4]{CMM14}. 
Thus (5) holds. 
\end{proof}

\begin{prop}\label{p-nonFano-iff}
Let $Y$ be a smooth projective threefold and  
let $f : X \to Y$ be a blowup along a smooth curve $\Gamma$.  
Assume that $X$ is  Fano. 
Then the following hold. 
\begin{enumerate}
\item $-K_Y$ is semi-ample. 
\item If $C$ is a curve on $Y$ such that $C \neq \Gamma$, then $-K_Y \cdot C>0$. 
\item $-K_Y$ is ample if and only if $-K_Y \cdot \Gamma >0$. 
\end{enumerate}
\end{prop}

\begin{proof}
Let us show (1). 
If $-K_Y$ is ample, then there is nothing to show. 
If $-K_Y$ is not ample, then $-K_Y$ is semi-ample by Lemma \ref{l-nonFano-blowdown}. 
Thus (1) holds. 

Let us show (2). 
We have $-K_X +D \sim -f^*K_Y$ for $D:=\Ex(f)$. 
For the proper transform $C_X$ of $C$ on $X$, we obtain 
\[
-K_Y \cdot C = -f^*K_Y \cdot C_X = (-K_X +D) \cdot C_X = -K_X \cdot C_X + D \cdot C_X \geq -K_X \cdot C_X >0. 
\]
Thus (2) holds. 
The assertion (3) follows from (1) and (2). 
\end{proof}

\begin{cor}\label{c-disjoint-blowup}
Let $V$ be a Fano threefold. 
Take mutually disjoint smooth curves $C_1$ and $C_2$ on $V$. 
If $\Bl_{C_1 \amalg C_2}\,V$ is Fano, then 
both $\Bl_{C_1}\,V$ and $\Bl_{C_2}\,V$ are Fano. 
\end{cor}

\begin{proof}
The assertion immediately follows from Proposition \ref{p-nonFano-iff}(3). 
\end{proof}

\begin{prop}\label{p-nonFano-flop}
Let $X$ be a Fano threefold and let $f : X \to Y$ be a blowup along a smooth curve $B_Y$, where $Y$ is a non-Fano smooth projective threefold. 
Then there exists the following commutative diagram consisting of birational morphisms 
\begin{equation}\label{e1-nonFano-flop}
\begin{tikzcd}
    & X \arrow[ld, "f"'] \arrow[rd, "f'"] \arrow[dd, "\varphi"]\\
    Y \arrow[rd, "g"'] & & Y' \arrow[ld, "g'"]\\
    & Z
\end{tikzcd}
\end{equation}
such that the following hold. 
\begin{enumerate}
\item $Y'$ is a non-Fano smooth projective threefold. 
\item $Z$ is a projective normal Gorenstein threefold. 
\item $D:=\Ex(f) = \Ex(f') =\Ex(\varphi) \simeq \P^1 \times \P^1$ and $\varphi(D)$ is a point. 
\item $\rho(X/Z) = \rho(X) -\rho(Z)=2$. 
\item $B_{Y'} := f'(D)$ is a smooth curve on $Y'$, and $f' : X \to Y'$ is the blowup along $B_{Y'}$. 
\item $\Ex(g) = B_Y$ and $\Ex(g') = B_{Y'}$. 
\item $g^*K_Z \sim K_Y$ and $g'^*K_{Z'} \sim K_{Y'}$. 
\end{enumerate}
\end{prop}

\begin{proof}
Since $-K_Y$ is semi-ample (Lemma \ref{l-nonFano-blowdown}) and 
$B_Y$ is a unique curve on $Y$ satisfying $-K_Y \cdot B_Y=0$ (Proposition \ref{p-nonFano-iff}(2)), 
there exists a birational morphism $g : Y \to Z$ 
to a projective normal threefold $Z$ such that $\Ex(g) =B_Y$. 
Since $g: Y \to Z$ is a flopping contraction, 
$Z$ is Gorestenin and $K_Y \sim g^*K_Z$ \cite[Proposition 6.10]{TanI}. 
Thus (2) holds.



Let us show (4). 
Since $\varphi : X \to Z$ is birational, 
we have $\rho(X/Z) = \rho(X) -\rho(Z)$ (Lemma \ref{l-rel-pic}). 
It holds that $\rho(X/Z) > \rho(X/Y) =1$. 
Since $D \simeq \P^1 \times \P^1$ and $\varphi(D)$ is a point, 
we get $\rho(X/Z) \leq \rho(D) =2$, and hence $\rho(X/Z) = \rho(X) -\rho(Z)=2$. 
Thus (4) holds. 


By $\rho(X/Z)=2$, 
we get the commutative diagram (\ref{e1-nonFano-flop}), 
where $f' : X \to Y'$ is the contraction of the extremal ray of $\NE(X/Z)$ 
not corresponding to $f : X \to Y$ (Proposition \ref{p-2ray}). 
It is clear that $f'$ is of type $E$ and $D = \Ex(f')$ (see \cite[Definition 3.3]{TanII} for the definition of the type). 
If $f'$ is not type $E_1$, then $f'(D)$ would be a point, 
and hence $Y' \simeq Z$, which contradicts $\rho(Y') = \rho(X)-1 =\rho(Y) > \rho(Z)$. 
Hence $f'$ is of type $E_1$. 
Then $f'$ is the blowup along $B_{Y'}$ for $B_{Y'} :=f'(D)$. 
Hence (3) and (5) hold. 
Moreover, $\Ex(g') =B_{Y'}$. Thus (6) holds. 
What is remaining is to prove (1) and (7). 
Since $K_Z$ is Cartier and $\dim \Ex(g') = \dim B_{Y'}=1$, we obtain $K_{Y'} \sim g'^*K_Z$, and hence $K_{Y'} \cdot B_{Y'}=0$. 
In particular, $Y'$ is not Fano. 
Therefore, (1) and (7) hold. 
\end{proof}


\begin{lem}\label{l-blowup-formula2}
Let $Y$ be a smooth projective threefold. 
Take a smooth curve $C$ on $Y$ and let $\sigma : X \to Y$ be the blowup along $C$. Set $D := \Ex(\sigma). $
Assume that $X$ is Fano. 
Then the following hold. 
\begin{enumerate}
\item $(-K_Y) \cdot C > 2p_a(C) -2$. 
\item $(-K_X)^3 <(-K_Y)^3$. 
\item If 
$(-K_X)^3 = (-K_Y)^3-4$, then 
$(p_a(C), -K_Y \cdot C) \in \{ (0, 1), (1, 2), (2, 3)\}$. 
\end{enumerate}
\end{lem}

\begin{proof}
The assertion (1) follows from $(-K_X)^2 \cdot D >0$ and Lemma \ref{l-blowup-formula}(3). 
Let us show (2). 
By Lemma \ref{l-blowup-formula}(2), we have 
\begin{equation}\label{e1-blowup-formula2}
(-K_Y)^3 -(-K_X)^3 = 2 ( (-K_Y) \cdot C - (p_a(C) -1)). 
\end{equation}
It is enough to show that $(-K_Y) \cdot C > p_a(C) -1$. 
If $p_a(C) \geq 1$, then this holds by (1). 
If $p_a(C) =0$, then this follows from $-K_Y \cdot C \geq 0$ 
(Proposition \ref{p-nonFano-iff}(1)). 
Thus (2) holds. 

Let us show (3).  
Assume $(-K_X)^3 = (-K_Y)^3-4$. 
By (\ref{e1-blowup-formula2}), we get  $4 
= 2 (-K_Y) \cdot C -2p_a(C) +2$, i.e., 
$(-K_Y) \cdot C =p_a(C) +1$. 
By (1), it holds that $2p_a(C) -2 < (-K_Y) \cdot C =p_a(C) +1$, 
which implies $p_a(C) <3$. 
Hence $(p_a(C), -K_Y \cdot C) \in \{ (0, 1), (1, 2), (2, 3)\}$. 
Thus (3) holds. 
\qedhere


\end{proof}

\subsection{Fano threefolds with $\rho=2$}\label{ss-pic2-III}


\begin{prop}\label{p-pic2-Pic}
Let $X$ be a Fano threefold with $\rho(X)=2$. 
Let $R_1$ and $R_2$ be the extremal rays of $\NE(X)$. 
For each $i \in \{1, 2\}$, 
\begin{itemize}
\item let $f_i : X \to Y_i$ be the contraction of $R_i$, 
\item fix an ample Cartier divisor $H_i$ on $Y_i$ which generates $\Pic\,Y_i (\simeq \Z)$, 
\item take an extremal rational curve $\ell_i$ on $X$ with $R_i = \R_{\geq 0}[\ell_i]$, and 
\item set $\mu_i := -K_X \cdot \ell_i$ 
(which is called the length of $R_i$). 
\end{itemize}
Then the following hold 
\begin{enumerate}
\item $\Pic\,X = \Z H_1 \oplus \Z H_2$. 
\item $-K_X \sim \mu_2 H_1 + \mu_1 H_2$. 
\item $H_1 \cdot \ell_2 = H_2 \cdot \ell_1 =1$. 
\end{enumerate}
\end{prop}

\begin{proof}
See    \cite[Proposition 5.9]{ATIII}. 
\end{proof}

\noindent
We now recall the definition of No.\ 2-xx. 

\begin{dfn}\label{d-pic2}
We say that $X$ is 
{\em 2-xx} or 
{\em of No.\ 2-xx} if 
$X$ is a Fano threefold with $\rho(X)=2$  
such that 
$(-K_X)^3$ and 
the types of the extremal rays are as in Table \ref{table-pic2} in Subsection \ref{ss-table-pic2}. 
For example, No.\ 2-1 and No.\ 2-8 are defined as follows. 
\begin{itemize}
\item 
A Fano threefold $X$ is {\em 2-1} or of {\em No.\ 2-1} if 
$\rho(X)=2$, 
$(-K_X)^3 = 4$, 
one of the extremal rays is of type $D_1$, and the other is of type $E_1$. 
\item 
A Fano threefold $X$ is {\em 2-8} or of {\em  No.\ 2-8} if 
$\rho(X)=2$, 
$(-K_X)^3 = 14$, one of the extremal rays is of type $C_1$, and the other is of type $E_3$ or $E_4$. 
\end{itemize}
\end{dfn}

\noindent

By \cite[Section 9]{ATIII}, 
a Fano threefold $X$ with $\rho(X)=2$ satisfies one and only one of the possibilities listed in Table \ref{table-pic2} in Subsection \ref{ss-table-pic2},  except for the column \lq\lq blowups".

\subsection{Some non-Fano criteria}

\begin{lem}\label{l-line-meeting}
Let $Y$ be a smooth projective threefold and let $\sigma: X \to Y$ be a blowup along a smooth curve $\Gamma$. 
Assume that $X$ is Fano. 
Then the following hold. 
\begin{enumerate}
\item 
If $C$ is a smooth curve on $Y$ satisfying $\Gamma \neq C$, 
then 
$\dim_k \MO_{\Gamma \cap C} < -K_Y \cdot C$. 
\item 
If $L$ is a curve on $Y$ such that $-K_Y \cdot L=1$, 
then $\Gamma = L$ or $\Gamma \cap L =\emptyset$. 
\end{enumerate}
\end{lem}


\begin{proof}
Let us show (1). 
Fix a smooth curve $C$ on $Y$ such that $\Gamma \neq C$. 
Let $C_X$ be the proper transform of $C$ on $X$. 
For $E := \Ex(\sigma)$, we have $K_X = \sigma^*K_Y +E$. 
By the scheme-theoretic equality $E = \sigma^{-1}(\Gamma)$ 
and the induced isomorphism $C_X \xrightarrow{\sigma|_{C_X}, \simeq} C$, the following holds: 
\[
E \cdot C_X = \deg_{C_X}( E|_{C_X}) = 
\dim_k \MO_{E|_{C_X}}
= \dim_k \MO_{E \cap C_X} 
= \dim_k \MO_{\Gamma \cap C}. 
\]
Then 
\[
0> K_X \cdot C_X = (\sigma^*K_Y +E) \cdot C_X 
= K_Y \cdot C + \dim_k \MO_{\Gamma \cap C}. 
\]
Thus (1) holds. 

Let us show (2). 
Assume that $L$ is a curve on $Y$ such that $-K_Y \cdot L=1$ and  $\Gamma \neq L$. 
For the proper transform $L_X$ of $L$ on $X$, 
it holds that 
\[
0> K_X \cdot L_X = (\sigma^*K_Y +E) \cdot L_X 
= K_Y \cdot L + E \cdot L_X = -1 +E \cdot L_X.
\]
We then get $1>E \cdot L_X$, which implies $E \cap L_X = \emptyset$, i.e., $\Gamma \cap L = \emptyset$. 
Thus (2) holds. 
\end{proof}

\begin{lem}\label{l-line-blowup}
Let $Y$ be a smooth projective threefold and let $\sigma: X \to Y$ be a blowup along a smooth curve $\Gamma$. 
Assume that $X$ is Fano, $\Gamma \simeq \P^1$, and $-K_Y \cdot \Gamma =1$. 
Then $N_{\Gamma/Y} \simeq \MO_{\P^1} \oplus \MO_{\P^1}(-1)$. 
\end{lem}

\begin{proof}
We have $\deg N_{\Gamma/Y} = (-K_Y) \cdot \Gamma +2p_a(\Gamma) -2 =-1$ (Lemma \ref{l-blowup-formula}). 
We then get $N_{\Gamma/Y} \simeq \MO_{\P^1}(n) \oplus \MO_{\P^1}(-n-1)$ for some $n \in \Z_{\geq 0}$, because one of the direct product factors is of 
non-negative degree. 
Let $s$ be the section of the $\P^1$-bundle $D :=\Ex(\sigma) = \P_{\Gamma}(N_{\Gamma/Y}^*) \to \Gamma$ corresponding to the projection 
$N_{\Gamma/Y}^* \simeq \MO_{\P^1}(-n) \oplus \MO_{\P^1}(n+1) \to \MO_{\P^1}(-n)$. 
Since $\MO_D(-D)$ is isomorphic to the tautological bundle $\MO_D(1)$ of the $\P^1$-bundle $D \to \Gamma$, 
 we obtain $\MO_D(-D) \cdot s =\MO_D(1) \cdot s= \deg (\MO_D(1)|_s) 
= \deg \MO_{\P^1}(-n) =-n$. 
Hence 
$0<-K_X \cdot s = (-\sigma^*K_Y -D) \cdot s = -K_Y \cdot \Gamma + \MO_D(-D) \cdot s = 1-n$. 
Thus $n=0$. 
\end{proof}


    

\begin{lem}\label{l-DB-bound}
Let $\sigma: X \to Y$ be a blowup along a smooth curve $B_Y$, 
where $X$ is a Fano threefold and $Y$ is a smooth projective threefold. 
Let $D$ be a prime divisor on $Y$ such that $B_Y \not\subset D$. 
Then $D \cdot B_Y < (-K_Y)^2 \cdot D$. 
\end{lem}


\begin{proof}
For $E :=\Ex(\sigma)$ and the prime divisor $D_X :=\sigma^*D$, the following holds: 
\[
0 <(-K_X)^2 \cdot D_X = (\sigma^*K_{Y} +E)^2 \cdot D_X 
= (-K_{Y})^2 \cdot D +2 \sigma^*K_{Y}\cdot E \cdot D_X +E^2 \cdot D_X.
\]
We have  $E^2 \cdot D_X=(E|_E) \cdot (D_X|_E)  = - D \cdot B_Y$ 
and $\sigma^*K_{Y}\cdot E \cdot D_X = \sigma^*K_Y \cdot E \cdot \sigma^*D =0$. 
Hence $D \cdot B_Y < (-K_{Y})^2 \cdot D$. 
\end{proof}

\subsection{Wild conic bundles}

Mori-Saito have classified Fano threefolds which admitting wild conic bundle structures. 

\begin{thm}\label{t-wild-cb}
Let $X$ be a Fano threefold. 
Assume that there exists a conic bundle $f: X \to S$ such that no fibre of $f$ is smooth. 
Then the following hold. 
\begin{enumerate}
\item $p=2$. 
\item One of the following holds. 
\begin{enumerate}
    \item 
    $Y \simeq \P^2$ and $X$ is isomorphic to 
    a prime divisor on $\P^2 \times \P^2$ of bidegree $(1, 2)$. 
    Furthermore, $\rho(X) =2$, $(-K_X)^3 = 30$, and $X$ is primitive ($X$ is  2-24). 
    \item  
    $S \simeq \P^1 \times \P^1$ and $X$ is isomorphic to 
    a prime divisor on $P:=\P_{\P^1 \times \P^1}( \MO(0, 1) \oplus \MO(1, 0) \oplus \MO)$ 
    which is linearly equivalent to $\MO_P(1)^{\otimes 2}$, 
    where $\MO_P(1)$ denotes the tautological bundle with respect to the $\P^2$-bundle structure $P \to \P^1 \times \P^1$. 
    Furthermore, $\rho(X)=3$, $(-K_X)^3 =26$, $X$ is imprimitive, 
    and $X \simeq \Bl_{C \amalg C'}\,Q$ for two conics $C$ and $C'$ on $Q$ satisfying $C \cap C' = \emptyset$ 
     ($X$ is  3-10). 
\end{enumerate}
\end{enumerate}
\end{thm}

\begin{proof}
See \cite[Corollary 8 and Remark 10]{MS03}. 
\end{proof}

\section{Fano conic bundles}

In this section, we establish some foundational results on Fano conic bundles. 
We start by giving its definition. 
Although it is not natural to impose the generically smooth assumption, 
this restriction is harmless for our purpose, because 
Fano threefolds admitting wild conic bundle structures are classified  by Mori-Saito \cite{MS03} (cf. Theorem \ref{t-wild-cb}). 

\begin{dfn}\label{d-FCB}
\begin{enumerate}
\item We say that $f: X \to S$ is a {\em threefold conic bundle} 
(resp. {\em threefold $\P^1$-bundle}) if 
$X$ is a smooth projective threefold and $f: X \to S$ is a generically smooth conic bundle (resp. a $\P^1$-bundle). 
\item We say that $f: X \to S$ is a {\em Fano conic bundle} 
(resp. {\em Fano $\P^1$-bundle}) 
if $X$ is a Fano threefold and $f: X \to S$ is a 
generically smooth conic bundle (resp. a $\P^1$-bundle). 
\item Given a theefold conic bundle $f: X \to S$, 
we say that  $\Gamma \subset X$ is a {\em subsection} of $f$ if $\Gamma$ is a closed subscheme of $X$ and the induced composite morphism $\Gamma \hookrightarrow X \xrightarrow{f} S$ is a closed immersion. 
We say that $\Gamma \subset X$ is a {\em regular subsection} of $f$ 
if $\Gamma$ is a subsection of $f$, $\Gamma$ is a smooth curve, 
and $f(\Gamma) \cap \Delta_f = \emptyset$, 
where $\Delta_f$ denotes the discriminant divisor of $f$ 
\cite[Definition 3.4]{Tan-conic}.
\end{enumerate}
\noindent
For the definition of conic bundles and discriminant divisors $\Delta_f$ 
we refer to \cite[Definition 2.3]{Tan-conic}. 
\end{dfn}


\begin{lem}\label{l-CB-criterion}
Let $f: X \to S$ be a  smooth projective morphism, 
where $X$ is a smooth threefold and $S$ is a smooth surface. 
Assume that $f_*\MO_X=\MO_S$ and $-K_X$ is $f$-ample. 
Then the following hold. 
\begin{enumerate}
\item $f$ is a conic bundle. 
\item Take a regular subsection $\Gamma$ of $f$ and let $\sigma: X' \to X$ be the blowup along $\Gamma$. 
Then $f \circ \sigma: X' \to X \to S$ is a conic bundle. 
\end{enumerate}
\end{lem}

\begin{proof}
The assertion (1) follows from the fact that $X_s \simeq \P^1$ for any closed point $s \in S$. 
Let us show (2). 
Set $\Gamma_S := f(\Gamma)$. 
Fix a closed point $s \in \Gamma_S$. 
After replacing $S$ by  an open neighbourhood of $s \in S$, 
we can find a smooth curve $T$ on $S$ such that 
$s \in T$ and $T +\Gamma$ is simple normal crossing. 
Take the base change: $f_T : X_T \to T$. 
Then the base change $X'_T := X' \times_S T$ coincides with the blowup of $X_T$ at the point $\Gamma \cap f^{-1}(T)$. 
In this case, it is well known that $X'_s$ is a conic. 
\end{proof}

\subsection{Elementary transforms}

\begin{prop}\label{p-ele-tf}
Let $f : X \to S$ be a Fano conic bundle. 
Let $B$ be a curve on $S$ such that $f^{-1}(B)$ is not irreducible. 
Then the following hold. 
\begin{enumerate}
\item There exist prime divisors $D_1$ and $D_2$ on $X$ such that 
$D_1 \neq D_2$ and the set-theoretic equality $f^{-1}(B)_{\red} = D_1 \cup D_2$ holds. 
\item For each $i \in \{1, 2\}$, 
there exist a threefold conic bundle $g_i : Y_i \to S$ and a blowup 
$\sigma_i : X \to Y_i$ along a smooth curve $B_i$ satisfying $\Ex(\sigma_i)=D_i$. 
Moreover, the following diagram is commutative. 
\begin{equation}\label{e1-ele-tf}
\begin{tikzcd}
& X \arrow[rd, "\sigma_2"] \arrow[ld, "\sigma_1"'] \arrow[dd, "f"]\\
Y_1 \arrow[rd, "g_1"'] & & Y_2 \arrow[ld, "g_2"]\\
&S
\end{tikzcd}
\end{equation}
\item $\rho(X) = \rho(Y_1)+1 = \rho(Y_2)+1$. In particular, $\rho(Y_1) = \rho(Y_2)$. 
\item 
$B$ is a smooth connected component of $\Delta_f$, i.e., 
$B \subset \Delta_f$ and $\Delta_f$ is smooth around $B$. 
\item $\Delta_f = \Delta_{g_1} \amalg B = \Delta_{g_2} \amalg B$. 
In particular, $\Delta_{g_1} = \Delta_{g_2}$. 
\end{enumerate}
\end{prop}

\begin{proof}
Fix a general closed point $P$ on $B$. 
Take the irreducible decomposition $f^{-1}(B) = D_1 \cup D_2 \cup \cdots \cup D_n$, which is a set-theoretic equality with $n \geq 2$. 
For every $1 \leq i \leq n$, we set 
\[
\zeta_i := (f^{-1}(P) \cap D_i)_{\red}. 
\]

\setcounter{step}{0}

\begin{claim}
For every $1 \leq i \leq n$,  
\begin{enumerate}
\renewcommand{\labelenumi}{(\roman{enumi})}
\item $D_i$ is a prime divisor on $X$, 
\item $f(D_i) =B$, and 
\item $\zeta_i$ is a curve. 
\end{enumerate}
\end{claim}

\begin{proof}
Since $B$ is an effective Cartier divisor on $S$, its pullback $f^*(B) = f^{-1}(B)$ is an effective Cartier divisor on $X$. 
Thus (i) holds. Then (ii) holds, because every fibre of $f$ is one-dimensional. 
Let us show (iii). 
Since $\zeta_i$ is pure one-dimensional, 
it suffices to show that $\zeta_i$ is irreducible. 
Suppose that $\zeta_i$ is not irreducible. 
Since every fibre of $f$ has at most two irreducible components and 
$P$ is a general closed point on $B$, 
we would get a set-theoretic equality 
$D_i|_{f^{-1}(B^{\circ})} = f^{-1}(B)|_{f^{-1}(B^{\circ})} = f^{-1}(B^{\circ})$ 
for some non-empty open subset $B^{\circ} \subset B$, 
which contradicts the fact that $D_j \to B$ is surjective for $j \neq i$. 
Therefore, each $\zeta_i$ is a curve. 
Thus (iii) holds. 
This completes the proof of Claim. 
\end{proof}

Since $P$ is a general closed point on $B$, 
it holds that $\zeta_i \neq \zeta_{i'}$ 
for $i \neq i'$ (as otherwise, we would get $D_i = D_{i'}$).
Therefore, we obtain $n=2$, because $f^{-1}(P)$ has at most two irreducible components. 
Thus (1) holds. 

Moreover, we get $\zeta_1 \not\subset D_2$ and $\zeta_2 \not\subset D_1$. 
We have a cycle-theoretic equality $f^*P = a_1 \zeta_1 + a_2 \zeta_2$ for 
some $a_1, a_2 \in \Z_{>0}$.  
By $D_1 \cdot \zeta_2>0$ and 
\[
0 = D_1 \cdot f^*P = D_1 \cdot (a_1 \zeta_1 + a_2 \zeta_2), 
\]
we get $D_1 \cdot \zeta_1 <0$. By symmetry, we obtain $D_2 \cdot \zeta_2 <0$. 


\begin{claim}
For each $i \in \{1, 2\}$, 
there exists a curve $\ell_i$ on $X$ such that 
$D_i \cdot \ell_i<0$, 
$\ell_i$ is contracted by $f|_{D_i} : D_i \to B$, and $\R_{\geq 0}[\ell_i]$ is an extremal ray of $\NE(X)$. 
\end{claim}

\begin{proof}[Proof of Claim]
By symmetry, we may assume that $i=1$. 
Fix an ample Cartier divisor $A_S$ on $S$. 
Note that we have 
\[
[\zeta_1] \in \NE(X) \cap (f^*A_S)^{\perp} =\R_{\geq 0}[m_1] + \cdots + \R_{\geq 0}[m_s]
\]
where each $\R_{\geq 0}[m_j]$ is an extremal ray of $\NE(X)$ generated by a curve $m_j$. 
We can write $\zeta_1  \equiv a_1 m_1 + \cdots +a_s m_s$ 
for some $a_1, ..., a_s \in \R_{\geq 0}$. 
By $D_1 \cdot \zeta_1<0$, we get $D_1 \cdot m_j <0$ for some $j$. 
It follows from $D_1 \cdot m_j <0$ and $f^*A_S \cdot m_j=0$ 
that $\ell_1 := m_j$ is a curve on $D_1$ contracted by $f|_{D_1} : D_1 \to B$. 
This completes the proof of Claim. 
\qedhere

\end{proof}

Let us show (2) and (3). 
Fix $i \in \{1, 2\}$. 
Let $\sigma_i : X \to Y_i$ be the contraction of the extremal ray $\R_{\geq 0}[\ell_i]$. 
By $D_i \cdot \ell_i <0$, 
every curve contracted by $\sigma_i$ is contained in $D_i$, 
and hence $\sigma_i$ is a birational morphism satisfying $\Ex(\sigma_i) \subset D_i$. 
By the classification of extrmeal rays on smooth threefolds \cite[Theorem 1.1]{Kol91}, 
$\Ex(\sigma_i)$ is a prime divisor, which implies $\Ex(\sigma_i)= D_i$. 
On the other hand, we have $f^*A_S \cdot \ell_i=0$ for an ample divisor $A_S$ on $S$. 
Therefore, a curve contracted by $\sigma_i : X \to Y_i$ is contracted by $f: X \to S$, i.e., we get a factorisation: 
\[
f: X \xrightarrow{\sigma_i} Y_i \xrightarrow{g_i} S. 
\]
In particular, every fibre of $\sigma_i$ is of dimension $\leq 1$. 
Again by the classification of extrmeal rays on smooth threefolds \cite[Theorem 1.1]{Kol91}, 
$\sigma_i$ is of type $E_1$. 
In particular, each $Y_i$ is a smooth projective threefold. 
Moreover, each $g_i : Y_i \to S$ is a conic bundle 
by Lemma \ref{l-CB-criterion}. 
Then (2) and (3) hold. 

Let us show (4). 
Clearly, $B \subset \Delta_f$. 
Fix a closed point $Q \in B$. 
In order to prove (4), it suffices to show that $\Delta_f$ is smooth at $Q$. 
Since $f: X \to S$ is a generically smooth conic bundle between smooth varieties, 
it is enough to show that $f^{-1}(Q)$ is not irreducible 
\cite[Theorem 4.3]{Tan-conic}. 
For each $i \in \{1, 2\}$, 
we set $\xi_i := (f^{-1}(Q) \cap D_i)_{\red}$, 
which is non-empty and pure one-dimensional. 
Since we have  a set-theoretic equality $f^{-1}(Q) = \xi_1 \cup \xi_2$, 
it suffices to prove that $\xi_1 \cap \xi_2$ contains no curve. 
Suppose that there exists a curve $C \subset \xi_1 \cap \xi_2$. 
Then $f(C)$ is a point and 
we get $C \subset \xi_1 \cap \xi_2 \subset D_1 \cap D_2$. 
However, this leads to a contradiction 
\[
[C] \in \R_{\geq 0}[\ell_1] \cap \R_{\geq 0}[\ell_2] =\{0\}.
\]
Thus (4) holds. 
Since $g_i : Y_i \to S$ is a conic bundle by (2), 
the assertion (5) follows from (4) and the fact that a generically reduced 
irreducible conic is automatically smooth. 
\qedhere



\end{proof}

\begin{dfn}\label{d-ele-tf}
We use the same notation as in the statement of Proposition \ref{p-ele-tf}. 
In this case, the diagram (\ref{e1-ele-tf}) is called an 
{\em elementary transform} (over $S$). 
The threefold $Y_2$ is called the {\em elementary transform of $f: X \xrightarrow{\sigma} Y_1 \xrightarrow{g_1} S$} (or the {\em elementary transform}  of $Y_1$).  
\end{dfn}

\begin{prop}\label{p-ele-tf-numbers}
We use the same notation as in the statement of Proposition \ref{p-ele-tf}. 
Set $E_i := g_i^{-1}(B)$ for each $i \in \{1, 2\}$. 
Then the following hold. 
\begin{enumerate}
\item $K_{Y_i/S} \cdot B_i =-(B_i \cdot B_i)_{E_i}$ for each $i \in \{1, 2\}$. 
\item $(B_1 \cdot B_1)_{E_1} + (B_2 \cdot B_2)_{E_2} = B^2$. 
\item $(-K_{Y_2})^3 = (-K_{Y_1})^3 -4 (-K_{Y_1/S}) \cdot B_1 +2B^2 
= (-K_{Y_1})^3 -4 (B_1 \cdot B_1)_{E_1} +2B^2$. 
\item $-K_{Y_2} \cdot B_2 = B^2 +2 (- K_S \cdot B) -(-K_{Y_1} \cdot B_1)$. 
\end{enumerate}
\end{prop}


\begin{proof}
Let us show (1). 
Fix $i \in \{1, 2\}$. 
We have  the induced $\P^1$-bundle 
structure $g_i|_{E_i} : E_i \to B$. 
We then get 
\[
K_{Y_i/S} \cdot B_i 
= K_{E_i/B} \cdot B_i 
= (K_{E_i} +B_i -B_i - (g_i|_{E_i})^*K_B) \cdot B_i  
= \deg K_{B_i} -B_i^2 - \deg K_{B_i} = -B_i^2. 
\]
Thus (1) holds.

Let us show (2). 
Note that we have $\sigma_1|_{D_2} : D_2 \xrightarrow{\simeq} E_1$ and $(\sigma_1|_{D_2})(B_X) = B_1$ for $B_X := D_1 \cap D_2$ 
(because the induced morphism $\sigma_1|_{D_2} : D_2 \to E_1$ is a finite birational morphism between normal surfaces). 
Hence $(B_X^2)_{D_2} = (B_1^2)_{E_1}$. 
Then the assertion (2) follows from 
\[
(B_1^2)_{E_1} + (B_2^2)_{E_2} = 
(B_X^2)_{D_1} + (B_X^2)_{D_2} =D^2_1 \cdot D_2 + D_1 \cdot D_2^2 
\]
\[
= (D_1+D_2) \cdot D_1 \cdot D_2 = (f^*B) \cdot B_X = B^2, 
\]
where the last equality follows from $f|_{B_X} : B_X \xrightarrow{\simeq} B$. 

Let us show (3). The following hold (Lemma \ref{l-blowup-formula}(2)): 
\begin{eqnarray*}
(-K_X)^3 &=& (-K_{Y_1})^3 -2 (-K_{Y_1}) \cdot B_1 +2p_a(B_1) -2\\
(-K_X)^3 &=& (-K_{Y_2})^3 -2 (-K_{Y_2}) \cdot B_2 +2p_a(B_2) -2. 
\end{eqnarray*}
By $B_1 \simeq B \simeq B_2$, we obtain $p_a(B_1) =p_a(B) = p_a(B_2)$ and 
$(g_1^*K_S) \cdot B_1 = K_S \cdot B = (g_2^*K_S) \cdot B_2$. 
We then get 
\[
 (-K_{Y_1})^3 -2 (-K_{Y_1/S}) \cdot B_1 =  (-K_{Y_2})^3 -2 (-K_{Y_2/S}) \cdot B_2. 
\]
Hence 
\begin{eqnarray*}
 (-K_{Y_2})^3 - (-K_{Y_1})^3  
 &=& -2 (-K_{Y_1/S}) \cdot B_1 +2 (-K_{Y_2/S}) \cdot B_2\\
  &=& -4 (-K_{Y_1/S}) \cdot B_1 +
  2(-K_{Y_1/S}) \cdot B_1+ 2(-K_{Y_2/S}) \cdot B_2\\
 &\overset{{\rm (1)}}{=}& -4 (-K_{Y_1/S}) \cdot B_1 +2 (B_1^2)_{E_1}+2 (B_2^2)_{E_2}\\
 &\overset{{\rm (2)}}{=}& -4 (-K_{Y_1/S}) \cdot B_1 +2B^2.
\end{eqnarray*}
Thus (3) holds. 
The assertion (4) follows from 
\[
-K_{Y_1} \cdot B_1 -K_{Y_2} \cdot B_2 + 2 K_S \cdot B 
\overset{{\rm (1)}}{=} (B_1^2)_{E_1} +(B_2^2)_{E_2} 
\overset{{\rm (2)}}{=} B^2. 
\]
\end{proof}



\begin{prop}\label{p-ele-tf-1Fano}
We use the same notation as in the statement of Proposition \ref{p-ele-tf}. 
Then the following hold. 
\begin{enumerate}
\item One of $Y_1$ and $Y_2$ is Fano. 
\item If $Y_1$ is Fano and $Y_2$ is not Fano, 
then $B \simeq \P^1$, $g^{-1}_1(B) \simeq \P^1 \times \P^1$, $-K_{Y_2} \cdot B_2 =0$, and $-K_{Y_1/S} \cdot B_1 = 2 (B^2+1)$. 
\end{enumerate}
\end{prop}


\begin{proof}
Let us show (1). 
Suppose that none of $Y_1$ nor $Y_2$ is Fano. 
Then $D_1 \simeq D_2 \simeq \P^1 \times \P^1$ (Lemma \ref{l-nonFano-blowdown}). 
By $f^{-1}(B)_{\red} = D_1 \cup D_2$, we get $D_1 \cap D_2 \neq \emptyset$. 
Pick a curve $C$ on $X$ satisfying $C \subset D_1 \cap D_2$. 
Since $\MO_X(-D_1)|_{D_1}$ is ample  (Lemma \ref{l-nonFano-blowdown}), 
we get $D_1 \cdot C = (D_1|_{D_1}) \cdot C <0$. 
On the other hand, $D_1|_{D_2}$ is an effective Cartier divisors on $D_2 = \P^1 \times \P^1$, which is nef. Hence $D_1 \cdot C = (D_1|_{D_2}) \cdot C \geq 0$. 
This is absurd. 
Thus (1) holds.

Let us show (2). 
By Lemma \ref{l-nonFano-blowdown} and Proposition \ref{p-nonFano-iff}, we have that 
$-K_{Y_2} \cdot B_2 =0$, 
$B \simeq B_2 \simeq \P^1$, and 
$g^{-1}_1(B) \simeq \Ex(\sigma_2) \simeq \P^1 \times \P^1$. 
It follows from Proposition \ref{p-ele-tf-numbers}(4) that 
$-K_{Y_1/S} -K_{Y_2/S} =B^2$. 
Then 
\[
-K_{Y_1/S} \cdot B_1 = B^2 - (-K_{Y_2/S} \cdot B_2) = 
B^2 -( -K_{Y_2} \cdot B_2) + (-K_S \cdot B)
\]
\[
=B^2 -0  -(K_S +B) \cdot B +B^2=2(B^2 +1). 
\]
Thus (2) holds. 
\end{proof}

\begin{lem}\label{l-ele-tf-K-relation}
We use the same notation as in the statement of Proposition \ref{p-ele-tf}. 
Then 
\[
-2K_X \sim -\sigma_1^*K_{Y_1} - \sigma_2^*K_{Y_2}-f^*B. 
\]
\end{lem}

\begin{proof}
It holds that 
\begin{eqnarray*}
K_X &\sim & \sigma_1^*K_X + E_{X/Y_1}\\
K_X &\sim& \sigma_2^*K_X + E_{X/Y_2},
\end{eqnarray*}
where each $E_{X/Y_i}$ denotes the exceptional prime divsior of $\sigma_i : X \to Y_i$. 
By $f^*B = E_{X/Y_1} +E_{X/Y_2}$, we obtain 
\[
2K_X \sim  
(\sigma_1^*K_{Y_1} + E_{X/Y_1}) + (\sigma_2^*K_{Y_2} + E_{X/Y_2})
=\sigma_1^*K_{Y_1} + \sigma_2^*K_{Y_2} +f^*B,  
\]
as required. 
\end{proof}

\subsection{Reduction to smaller Fano conic bundles}

Let $f: X \to S$ be a Fano conic bundle. 
The purpose of this subsection is to construct a smaller Fano conic bundle when 
$f$ or $S$ is not \lq\lq minimal". 
More specifically, such a construction is applicable 
when $\rho(X) > \rho(S) +1$ (Proposition \ref{p-smaller-same-base}) or $S$ has a curve $E$ with $E^2 <0$ (Proposition \ref{p-FCB-(-1)}).  
To this end, we will establish some structural results on Fano conic bundles 
(Proposition \ref{p-FCB-centre}, Lemma \ref{l-FCB-pic-irre}). 
We start with the following auxiliary result.

\begin{lem}\label{l-double-conic}
Let $f: X \to S$ be a conic bundle, 
where $X$ is a smooth threefold and $S$ is a smooth surface. 
Fix a closed point $s \in S$ whose  fibre $f^{-1}(s)$ is not reduced. 
Set $C := f^{-1}(s)_{\red} (\simeq \P^1)$. 
Then $N_{C/X} \simeq \MO_{\P^1}(a) \oplus \MO_{\P^1}(-a-1)$ for some positive integer $a$. 
\end{lem}

The following argument is based on \cite[Lemma 3.25]{Mor82}. 

\begin{proof}
Set $Z := f^{-1}(s)$. Let $I_C$ and $I_Z$ be the coherent ideal 
sheaves on $X$ corresponding to $C$ and $Z$, respectively. 

\begin{claim*}
$I_C^2 \subsetneq I_Z \subsetneq I_C$.     
\end{claim*}

\begin{proof}[Proof of Claim]
The inclusion $I_Z \subsetneq I_C$ follows from 
$C\neq Z$ and
$C  = f^{-1}(s)_{\red} =Z_{\red}$. 
The remaining proof consists of two parts. 
\begin{enumerate}
\item $I_C^2 \subset I_Z$. 
\item $I_C^2 \neq I_Z$. 
\end{enumerate}

Let us  show (1). 
Possibly after replacing $S$ by an open neighbourhood of $s\in S$, 
we have a scheme-theoretic equality $D_1 \cap D_2 = s$ for some 
simple normal crossing divisor $D_1 +D_2$, 
where each $D_i$ is a smooth prime divisor on $S$. 
In particular, we get $Z =f^{-1}(s) =f^{-1}(D_1) \cap f^{-1}(D_2)$. 
Hence  $Z$ is Cohen-Macaulay, and hence $I_Z$ has no embedded points. 
Thus $I_Z$ is a primary ideal. 

Therefore, it is enough to show $I_C^2 \MO_{X, \xi} \subset I_Z\MO_{X, \xi}$ for the generic point $\xi$ of $C$. 
It follows from $-K_X \cdot Z =2$ that 
$A := \MO_{X, \xi}/I_Z\MO_{X, \xi}$ is an artinian local ring with ${\rm length}_A\,A=2$. 
We have the sequence $A \supsetneq I_CA \supset I_C^2A$ of ideals. 
By ${\rm length}_A\,A=2$, we obtain $I_CA =I_C^2A$ or $I_C^2A=0$. 
Suppose $I_CA = I_C^2A$. Since $I_CA$ is the maximal ideal of the local ring $A$, 
it follows from Nakayama's lemma that $I_CA =0$, i.e., $I_C \MO_{X, \xi} = I_Z\MO_{X, \xi}$. 
This would  imply $I_C = I_Z$, which is absurd. 
Hence $I_C^2A =0$. We then get $I_C^2\MO_{X, \xi} \subset I_Z\MO_{X, \xi}$. 
This completes the proof of (1). 

Let us show (2). 
Suppose $I_C^2  = I_Z$. 
It suffices to derive a contradiction. 
We get $I_Z = I_C^2 = \m_s^2 \MO_X$ 
for the coherent ideal $\MO_S$-module $\m_s$ corresponding to 
the closed point $s \in S$. 
Since $f: X \to S$ is faithfully flat, 
so is the base change $f' : X' := X \times_S S' \to S'$ 
for $S' := \Spec\,(\MO_S/\m_s^2)$. 
Then the induced ring homomorphism 
\[
\MO_{S'} \to f'_*\MO_{X'}
\]
is injective (because $\MO_{S'} = \MO_{S', s} \to \MO_{X', x}$ 
is a local  flat ring homomorphism for a closed point $x$ lying over $s$, which is injective). 
However, 
we then get the following contradiction: 
\[
3= \dim_k(\MO_{S, s}/\m_s^2) = h^0(S', \MO_{S'})
\]
\[
 \leq h^0(X', \MO_{X'}) = h^0(X, \MO_X/\m_s^2\MO_X) = h^0(X, \MO_X/I_Z)
=h^0(Z, \MO_Z)=1. 
\]
Thus (2) holds. 
This completes the proof of Claim. 
\end{proof}

We have an exact sequence 
\[
0 \to I_C/I_Z \to \MO_Z \to \MO_C \to 0. 
\]
By $h^0(Z, \MO_Z ) =h^0(C, \MO_C) =k$, we obtain $H^0(Z, I_C/I_Z)=0$. 
We also have an exact sequence of coherent $\MO_C$-modules: 
\[
0 \to I_Z/I_C^2 \to I_C/I_C^2 \to I_C/I_Z \to 0. 
\]
Since $I_C/I_C^2$ is a locally free $\MO_C$-module of rank $2$, 
$I_Z/I_C^2$ is a locally free $\MO_C$-module of rank $\leq 2$. 
Suppose ${\rm rank}\,I_Z/I_C^2 =2$. 
Then the inclusion $I_Z/I_C^2 \hookrightarrow I_C/I_C^2$ is an isomorphism around the generic point $\xi$ of $C$. 
Hence we would get  $I_C\MO_{X, \xi} = I_Z\MO_{X, \xi}$, which leads to a contradiction:  $I_C =I_Z$. 
Thus $I_Z/I_C^2$ is an invertible $\MO_C$-module. 
Then $I_C/I_Z \simeq L \oplus T$ for an invertible $\MO_C$-module $L$ and a coherent $\MO_C$-module $T$ with $\dim (\Supp\,T) \leq 0$. 
By $H^0(Z, I_C/I_Z)=0$, we get $T=0$. 
Then also $I_C/I_Z$ is an invertible $\MO_C$-module. 
Again by $H^0(C, I_C/I_Z)=0$, we get $\deg_C (I_C/I_Z) <0$. 
Recall that $\deg_C (I_C/I_C^2) = -\deg N_{C/X} 
= -(-K_X) \cdot C -2p_a(C) +2 =1$ (Lemma \ref{l-blowup-formula}(5)). 
Hence $\deg_C (I_Z/I_C^2) = \deg_C(I_C/I_C^2) - \deg_C (I_C/I_Z) >1$. 
Therefore, we get 
$\Ext^1_{\MO_C}(I_C/I_Z, I_Z/I_C^2) \simeq H^1(C, (I_C/I_Z)^{-1} \otimes I_Z/I_C^2)=0$, and hence 
\[
I_C/I_C^2 \simeq I_Z/I_C^2 \oplus I_C/I_Z. 
\]
Since $I_C/I_C^2$ has a direct summand $I_C/I_Z$ of negative degree, 
its dual $N_{C/X}$ has a direct summand $(I_C/I_Z)^{-1}$ of positive  degree. 
We are done by $\deg N_{C/X} = -1$. 
\end{proof}

\begin{prop}\label{p-FCB-centre}
Let $g : Y \to S$ be a threefold conic bundle and let $\Gamma$ be a smooth curve on $Y$. 
For the blowup $\sigma : X \to Y$ along $\Gamma$, assume that $X$ is a Fano threefold. 
Set $f := g \circ \sigma$: 
\[
f : X \xrightarrow{\sigma} Y \xrightarrow{g} S. 
\]
Then the following hold. 
\begin{enumerate}
\item $\Gamma$ does not intersect any non-smooth fibres. 
\item One of {\rm (a)} and {\rm (b)} holds. 
\begin{enumerate}
    \item $\Gamma$ is a smooth fibre of $g$. Moreover, 
    there exists a commutative diagram 
    \begin{equation}\label{e1-FCB-structure}
    \begin{tikzcd}
    X \arrow[r, "\sigma"] \arrow[rd, "f"] \arrow[d, "h"] & Y \arrow[d, "g"]\\
    T \arrow[r, "\sigma'"] &S
    \end{tikzcd}
\end{equation}
    such that the square diagram in (\ref{e1-FCB-structure}) is carterian, 
\begin{enumerate}
\item 
$\sigma': T \to S$ is the blowup at the point $g(\Gamma)$, and 
\item $h: X \to T$  is a Fano conic bundle.  
\end{enumerate}
    \item $\Gamma$ is a regular subsection of $g$. Moreover, 
    $f : X \to S$ is a Fano conic bundle and $\Delta_f = \Delta_g \amalg g(\Gamma)$. 
\end{enumerate}
\end{enumerate}
\end{prop}

\begin{proof}
Let us show (1). 
Suppose that $\Gamma$ intersects a non-smooth fibre $g^{-1}(s)$ for some closed point $s \in S$. 
Note that $-K_X \cdot \ell =1$ for every curve $\ell$ satisfying $\ell \subset g^{-1}(s)$. 
If $\Gamma \not\subset g^{-1}(s)$, then we can find a curve $\ell \subset g^{-1}(s)$ such that 
$\Gamma \cap \ell \neq \emptyset$, $\Gamma \neq \ell$, and $-K_X \cdot \ell =1$, which contradict Lemma \ref{l-line-meeting}. 
Hence $\Gamma \subset g^{-1}(s)$. 
If $g^{-1}(s)$ is reducible, then the irreducible component $\ell$ of $g^{-1}(s)$ other than $\Gamma$ satisfies 
$\Gamma \cap \ell \neq \emptyset$, $\Gamma \neq \ell$, and $-K_X \cdot \ell =1$. 
This is absurd again by Lemma \ref{l-line-meeting}. 
Thus $g^{-1}(s)$ is irreducible and non-reduced. 
In this case, we have $N_{\Gamma/X} \simeq \MO_{\Gamma}(a) \oplus \MO_{\Gamma}(-a-1)$ for some $a>0$ 
(Lemma \ref{l-double-conic}), 
which contradicts Lemma \ref{l-line-blowup}. 
Thus (1) holds. 

We now show (2). 
Assume that $g(\Gamma)$ is a point. 
Let us show (a). 
By (1), $\Gamma$ is a smooth fibre of $g$. 
We then obtain the carterian diagram (\ref{e1-FCB-structure}) satisfying (i). 
Since a base change of a conic bundle is a conic bundle, 
also (ii) holds. 
Thus (a) holds. 

Assume that $g(\Gamma)$ is a curve. 
Let us show (b). 
Set $E := \Ex(\sigma)$. 
Fix a fibre $Z$ of $g: Y \to S$ intersecting $\Gamma$ and let $Z_X$  be the  proper transform of $Z$ on $X$. 
By (1), we have $Z_X \xrightarrow{\simeq} Z \simeq \P^1$. 
Then 
\[
0< E \cdot Z_X = (K_X -\sigma^*K_Y) \cdot Z_X = K_X \cdot Z_X - K_Y \cdot Z  = K_X \cdot Z_X +2 <2. 
\]
Therefore, 
we get $(\Gamma \cdot Z)_{g^{-1}(g(\Gamma))} =1$, and hence 
$\Gamma$ is a section of the $\P^1$-bundle 
$g^{-1}(g(\Gamma)) \to g(\Gamma)$. 
Then $\Gamma$ is a regular subsection of $g$. 
By Lemma \ref{l-CB-criterion},  $f:X \to S$ is a Fano conic bundle. 
Thus (b) holds. 
This completes the proof of (2). 
\end{proof}

Let $f: X \to S$ be a Fano conic bundle. 
If there is a curve $B$ on $S$ with $f^{-1}(B)$ reducible, 
then we can find a smaller Fano conic bundle $g: Y \to S$, i.e., we have the following  factorisation for some blowdown $\sigma: X \to Y$ (Proposition \ref{p-ele-tf}): 
\[
f: X \xrightarrow{\sigma} Y \xrightarrow{g} S
\]
This result enables us to reduce the classification of Fano conic bundles to the case when $f^{-1}(B)$ is  irreducible for every curve $B$ on $S$. 
This situation is characterised by $\rho(X) = \rho(S)+1$ as the  following lemma shows. 

\begin{lem}\label{l-FCB-pic-irre}
Let $f : X \to S$ be a Fano conic bundle. 
Then the following are equivalent. 
\begin{enumerate}
\item $\rho(X) = \rho(S)+1$. 
\item For every curve $B$ on $S$, $f^{-1}(B)$ is irreducible. 
\end{enumerate}
\end{lem}

\begin{proof}
Note that a curve $C$ on $S$ is an effective Cartier divisor, and hence its pullback $f^*(C) =f^{-1}(C)$ is an effective Cartier divisor on $X$. 
Since $f: X \to S$ is a Fano conic bundle, 
we get $\rho(X/S) = \rho(X) -\rho(S)$ (Lemma \ref{l-rel-pic}). 


Let us show $(1) \Rightarrow (2)$. 
Assume that (2) does not hold. 
Then there exists a curve $B$ on $S$ such that $f^{-1}(B)$ is not irreducible. 
By Proposition \ref{p-ele-tf}, $\rho(X) = \rho(Y_1)+1 > \rho(S)+1$ 
under the same notation as in Proposition \ref{p-ele-tf}. 
Thus (1) does not hold. This completes the proof of $(1) \Rightarrow (2)$.


Let us show $(2) \Rightarrow (1)$. 
Assume that (1) does not hold, i.e., $\rho(X/S) = \rho(X) -\rho(S) \geq 2$. 
Fix an ample Cartier divisor $A_S$ on $S$. 
Then we can find an extremal ray $R$ of $\NE(X)$ with 
$R \subsetneq \NE(X/S) = \NE(X) \cap (f^*A_S)^{\perp}$. 
Let $\varphi : X \to Z$ be the contraction of $R$. 
We have the factorisation 
\[
f: X \xrightarrow{\varphi} Z \xrightarrow{\psi} S,  
\]
where both $\varphi$ and $\psi$ are contractions. 
By $\dim X \geq \dim Z \geq \dim S$, 
we get $\dim Z = 3$ or $\dim Z=2$. 
Assume that $\dim Z =3$, i.e., $\varphi : X \to Z$ is birational. 
Since $\varphi(\Ex(\varphi))$ is irreducible and 
of dimension $\leq 1$, 
we can find a curve $B$ on $S$ such that $f(\Ex(\varphi)) = \psi(\varphi(\Ex(\varphi))) \subset B$. 
Then it is easy to see that $f^{-1}(B)$ is not irreducible, i.e., (2) does not  hold. 
Hence we may assume that $\dim Z =2$. 
By $\rho(Z) > \rho(S)$, there exists a curve $C$ on $Z$ 
such that $\psi(C)$ is a point.  
Pick a curve $B$ passing through $\psi(C)$. 
Then $\psi^{-1}(B)$ is not irreducible, and hence $f^{-1}(B)$ is not irreducible. 
This completes the proof of $(2) \Rightarrow (1)$. 
\end{proof}

\begin{prop}\label{p-smaller-same-base}
Let $f : X \to S$ be a Fano conic bundle such that $\rho(X) >\rho(S)+1$. 
Then there exist morphisms 
\[
f : X \xrightarrow{\sigma} Y \xrightarrow{g} S
\]
such that $g: Y \to S$ is a Fano conic bundle and $\sigma : X \to Y$ 
is a blowup along a regular subsection $B_Y$ of $g$. 
\end{prop}

\begin{proof}
Lemma \ref{l-FCB-pic-irre} enables us to find a curve $B$ such that $f^{-1}(B)$ is not irreducible. 
By Proposition \ref{p-ele-tf} and Proposition \ref{p-ele-tf-1Fano}, 
we get morphisms 
$f : X \xrightarrow{\sigma} Y \xrightarrow{g} S$ such that  $g: Y \to S$ is a Fano conic bundle and $\sigma : X \to Y$ 
is a blowup along smooth curve $B_Y$. 
Automatically, $B_Y$ is a regular subsection of $g$ (Proposition \ref{p-FCB-centre}). 
\end{proof}


\begin{prop}\label{p-FCB-(-1)}
Let $f : X \to S$ be a Fano conic bundle. 
Let $E$ be a curve on $S$ such that $E^2 <0$ and $f^{-1}(E)$ is irreducible. 
Then $E$ is a $(-1)$-curve on $S$ and there exists a cartesian diagram 
\begin{equation}\label{e1-FCB-(-1)}
\begin{tikzcd}
X \arrow[r, "\sigma"] \arrow[d, "f"] & X'\arrow[d, "f'"]\\
S \arrow[r, "\tau"] & S'
\end{tikzcd}
\end{equation}
such that 
\begin{enumerate}
\item $\tau: S \to S'$ is the blowdown of $E$ (so, $S'$ is a smooth projective surface), 
    \item 
$f' : X' \to S'$ is a Fano conic bundle, 
    \item the scheme-theoretic fibre $\Gamma := f'^{-1}(\tau(E))$  over the point $\tau(E)$ is smooth, and 
    \item $\sigma$ is  the blowup along $\Gamma$. 
\end{enumerate}
\end{prop}


\begin{proof}
Set $D := f^{-1}(E)_{\red}$, which is a prime divisor. 
As $f^*(E)$ is an irreducible effective Cartier divisor, 
we have $f^*(E) = c D$ for some $c \in \Z_{>0}$.

\setcounter{step}{0}

\begin{step}\label{s1-FCB-(-1)}
There exists a birational morphism $\sigma: X \to X'$ such that 
\begin{enumerate}
\renewcommand{\labelenumi}{(\roman{enumi})}
\item $X'$ is a smooth projective threefold, 
\item $\Gamma := \sigma(D)$ is a smooth curve on $X'$, and 
\item $\sigma$ coincides with the blowup along $\Gamma$. 
\end{enumerate}
\end{step}

\begin{proof}[Proof of Step \ref{s1-FCB-(-1)}]
Fix a curve $Z$ on $X$ such that $f(Z)=E$. 
By $[Z] \in \NE(X) = \sum_{i=1}^n \R_{\geq 0}[\ell_i]$, 
we can write $[Z] \equiv \sum_{i=1}^n a_i [\ell_i]$, 
where $\R_{\geq 0}[\ell_i]$ is an extremal ray of $\NE(X)$, $\ell_i$ is a   curve on $X$, 
and $a_i \in \R_{\geq 0}$ for every $1 \leq i \leq n$. 
By 
\[
0> f_*Z \cdot E = \sum_{i=1}^n a_i E \cdot f_*(\ell_i), 
\]
we can find a curve $\ell := \ell_i$ on $X$ such that 
$\R_{\geq 0}[\ell]$ is an extremal ray of $\NE(X)$ and 
$E \cdot f_*(\ell) <0$, and hence $f(\ell) =E$. 
Let $\sigma : X \to X'$ be the contraction of the extremal ray $\R_{\geq 0}[\ell]$. 
By 
\[
\ell \cdot f^*E = f_*\ell \cdot E <0, 
\]
any curve $C$ satisfying $[C] \in \R_{\geq 0}[\ell]$ is contained in $D =(f^*E)_{\red}$. 
Hence $\sigma : X \to X'$ is a birational morphism (i.e, of type $E$) such that 
$\Ex(\sigma) = D$. 
Since there is a morphism $f : X \to S$ such that $f(D)$ is a curve, 
$\sigma$ is of type $E_1$, i.e., 
$X'$ is a smooth projective threefold and $\sigma$ is the blowup along $\Gamma :=\sigma(D)$. 
This completes the proof of Step \ref{s1-FCB-(-1)}. 
\end{proof}

\begin{step}\label{s2-FCB-(-1)}
$D \simeq \P^1 \times \P^1$. 
Moreover, $D \to \Gamma$ and $D \to \widetilde E$ are the projections, where 
$D \to \widetilde E$ is the Stein factorisation of $f|_D : D \to E$. 
\end{step}

\begin{proof}[Proof of Step \ref{s2-FCB-(-1)}]
By Step \ref{s1-FCB-(-1)}, $D$ is a $\P^1$-bundle over $\Gamma$. 
Note that $D$ has another surjection $f|_D : D \to E$ to a curve $E$. 
Since a fibre of $f|_D :D \to E$ consists of rational curves, 
we get $\Gamma \simeq \P^1$. 
Then $D$ is a $\P^1$-bundle over $\P^1$ which has two morphisms to curves. 
Hence $D \simeq \P^1 \times \P^1$. 
This completes the proof of Step \ref{s2-FCB-(-1)}. 
\end{proof}

\begin{step}\label{s3-FCB-(-1)}
$X'$ is a Fano threefold. 
\end{step}

\begin{proof}[Proof of Step \ref{s3-FCB-(-1)}]
Since $X'$ is a smooth projective threefold, 
it is enough to show that $-K_{X'}$ is ample. 
By Proposition \ref{p-nonFano-iff}(3), 
it suffices to show that $-K_{X'} \cdot \Gamma >0$. 
Pick a curve $\Gamma_X$ contained in a fibre of $f|_D : D \to E$. 
By $f^*E = cD$, we obtain $D \cdot \Gamma_X =0$. 
Therefore, 
\[
0> K_X \cdot \Gamma_X = (\sigma^*K_{X'}+D)\cdot \Gamma_X = \sigma^*K_{X'} \cdot \Gamma_X 
= bK_{X'} \cdot \Gamma
\]
for some $b \in \Z_{>0}$. 
This completes the proof of Step \ref{s3-FCB-(-1)}. 
\end{proof}

\begin{step}\label{s4-FCB-(-1)}
There exists a birational morphism $\tau : S \to S'$ 
such that $S'$ is a projective normal surface, $\Ex(\tau) = E$, and $\rho(S) = \rho(S') +1$. 
\end{step}

\begin{proof}[Proof of Step \ref{s4-FCB-(-1)}]
Although the proof is very similar to the one of \cite[Theorem 3.21]{Tan14}, we here give a proof for the sake of completeness. 

First of all, we construct a Cartier divisor $N_S$ on $S$ such that (a)-(c) hold. 
\begin{enumerate}
\item[(a)] $N_S = H_S+ nE$ for some ample divisor $H_S$ and $n \in \Z_{>0}$. 
\item[(b)] For a curve $B$ on $S$, $N_S \cdot B =0$ if and only if $B = E$. 
\item[(c)] $N_S$ is semi-ample. 
\end{enumerate}
Fix an ample Cartier divisor $A_S$ on $S$. 
We define $\lambda \in \Q_{>0}$ by $(A_S + \lambda E) \cdot E =0$, 
i.e., $\lambda = \frac{A \cdot E}{-E^2}$. 
For $m := -E^2$, we set $N_S := m(A+\lambda E)$. 
Then (a) holds for $H_S :=mA_S$ and $n :=m\lambda$. 
By (a), $N_S$ is big. 
We have $N_S \cdot E =0$ by construction. 
For a curve $B$ on $S$ satisfying $B \neq E$, we obtain $N_S \cdot B = m(A_S+\lambda E) \cdot B \geq mA_S \cdot E >0$. 
Thus $N_S$ is nef and (b) holds. 
Let us show that $N_S$ is semi-ample. 
By \cite[Proposition 1.6]{Kee99} or \cite[Corollary 3.4]{CMM14}, 
it is enough to show that $N_S|_E$ is semi-ample. 
Since $f|_D : D \to E$ has connected fibres, 
it suffices to prove that $f^*N_S|_D$ is semi-ample \cite[Lemma 2.11(3)]{CT20}. 
This holds, because $D \simeq \P^1 \times \P^1$ and $f^*N_S|_D$ is nef. 
This completes the proof of (a)-(c). 

By (b) and (c), we obtain a birational morphism $\tau :S \to S'$ to a projective normal surface $S'$ such that $\Ex(\tau)=E$. 
It suffices to prove $\rho(S) = \rho(S')+1$. 
It is enough to show that the sequence 
\[
0 \to \Pic\,S' \otimes_{\Z} \Q  \xrightarrow{\tau^*} \Pic\,S \otimes_{\Z} \Q \xrightarrow{\cdot E}  \Q \to 0 
\]
is exact. 
Take a Cartier divisor $M$ on $S$ such that $M \cdot E =0$. 
It suffices to show that $r M \sim \tau^*M_{S'}$ for some $r \in \Z_{>0}$ and 
Cartier divisor $M_{S'}$ on $S'$. 
We set 
\[
\widetilde{N}_S := \ell N_S + M 
\]
for sufficiently large integer $\ell \gg 0$. 
Then we have $\widetilde{N}_S = \ell N_S+M = \ell( H_S + nE) +M = (\ell H_S + M) + \ell nE$. 
The following properties, corresponding to (a)-(c), hold: 
\begin{enumerate}
\item[(A)] $\widetilde{N}_S = \widetilde{H}_S+ \widetilde{n}E$ for some ample divisor $\widetilde{H}_S$ and $\widetilde n \in \Z_{>0}$. 
\item[(B)] For a curve $B$ on $S$, $\widetilde{N}_S \cdot B =0$ if and only if $B = E$.
\item[(C)] $\widetilde{N}_S$ is semi-ample. 
\end{enumerate}
Indeed, (A)-(C) hold for $\widetilde H_S := (\ell H_S + M)$ and $\widetilde n := \ell n$ by applying the  same proof as that of (a)-(c). 
Hence $|rN_S|$ and $|r \widetilde{N}_S|$ are base point free for some $r \in \Z_{>0}$. 
For each of the morphisms induced by $|rN_S|$ and $|r \widetilde{N}_S|$, 
its Stein factorisation coincides with $\tau : S \to S'$. 
Thus 
$r N_S \sim \tau^*D$ and 
$r\widetilde{N}_S \sim \tau^*\widetilde{D}$ for some Cartier divisors $D$ and $\widetilde D$ on $S'$. 
Then $rM = r\widetilde{N}_S -r\ell N_S \sim \tau^*(\widetilde D -\ell D)$.  
This completes the proof of Step \ref{s4-FCB-(-1)}. 
\end{proof}

By Step \ref{s1-FCB-(-1)} and Step \ref{s4-FCB-(-1)}, 
we obtain the commutative diagram (\ref{e1-FCB-(-1)}).

\begin{step}\label{s5-FCB-(-1)}
The following hold. 
\begin{enumerate}
    \item[(iv)] 
$S'$ is a smooth projective surface. 
\item[(v)] $f'$ is a conic bundle. 
\item[(vi)] $E$ is a $(-1)$-curve and $\tau : S \to S'$ is its blowdown. 
\item[(vii)] The diagram (\ref{e1-FCB-(-1)}) is cartesian. 
\end{enumerate} 
\end{step}

\begin{proof}[Proof of Step \ref{s5-FCB-(-1)}]
By Step \ref{s1-FCB-(-1)} and Step \ref{s4-FCB-(-1)}, we get  
\[
\rho(X') -\rho(S') = (\rho(X)-1) -(\rho(S')-1) =1. 
\]
Then $f'$ is a contraction of an extremal ray of type C. 
In particular, (iv) and (v) hold. 
Then (iv) and Step \ref{s4-FCB-(-1)} imply (vi). 
Since $f'$ is flat, 
(vi) and Step \ref{s1-FCB-(-1)} imply (vii) \cite[Section 8.1, Proposition 1.12(c)]{Liu02}.
This completes the proof of Step \ref{s5-FCB-(-1)}. 
\end{proof}
By Step \ref{s1-FCB-(-1)}, 
Step \ref{s3-FCB-(-1)}, and Step \ref{s5-FCB-(-1)}, 
we get the cartesian diagram (\ref{e1-FCB-(-1)}) satisfying (1), (2), and (4). 
The remaining one (3) follows from Proposition \ref{p-FCB-centre}. 
This completes the proof of Proposition \ref{p-FCB-(-1)}. 
\qedhere


\end{proof}


\begin{prop}\label{p-FCB-dP}
Let $f : X \to S$ be a Fano conic bundle. 
Then $-K_S$ is ample, i.e., $S$ is a smooth del Pezzo surface. 
\end{prop}

\begin{proof}
By Proposition \ref{p-ele-tf}, Proposition \ref{p-ele-tf-1Fano}, and Lemma \ref{l-FCB-pic-irre}, 
the problem is reduced to the case when $\rho(X) = \rho(S)+1$. 
Recall that $S$ is a smooth rational surface \cite[Lemma 4.3]{ATIII}. 
By \cite[Corollary 4.10]{Eji19}, $-K_S$ is big. 
We can write $-K_S \equiv A + \Delta$ for some ample $\Q$-divisor $A$ and an effective $\Q$-divisor $\Delta$. 
Let $\Delta =\sum_{i=1}^r c_i C_i$ be the irreducible decomposition, 
where $c_i \in \Q_{\geq 0}$ and $C_i$ is a prime divisor for every $1 \leq i \leq r$. 
Since the $\Q$-divisor $-(K_X + \Delta +\frac{1}{2}A ) \equiv \frac{1}{2}A$ is ample, 
it follows from  the cone theorem \cite[Theorem 3.13(2)]{Tan14} that 
\[
\overline{\NE}(S) = \overline{\NE}(S)_{ K_S + \Delta +\frac{1}{2}A \geq 0} + \sum_{j=1}^s \R_{\geq 0}[\ell_i] = \sum_{j=1}^s \R_{\geq 0}[\ell_j], 
\]
where $\ell_1, ..., \ell_s$ are curves on $S$. 
Fix $1 \leq j \leq s$. 
By Kleiman's criterion for ampleness, 
it is enough to show that $K_S \cdot \ell_j <0$. 
By Proposition \ref{p-FCB-(-1)}, it suffices to show that $K_S \cdot \ell_j <0$ or $\ell_j^2 <0$. 
Suppose $K_S \cdot \ell_j \geq 0$ and $\ell_j^2 \geq 0$. 
The latter one implies $\Delta \cdot \ell_j \geq 0$. 
Then the numerical equivalence $K_S+ \Delta +A \equiv 0$ leads to  the following contradiction: 
\[
0 = (K_S + \Delta +A) \cdot \ell_j = 
K_S \cdot \ell_j + \Delta \cdot \ell_j +A\cdot \ell_j \geq 0+ 0+A \cdot \ell_j>0.
\]
\end{proof}

\begin{cor}\label{c-FCB-(-1)-2}
Let $f : X \to S$ be a Fano conic bundle. 
Take  a curve $E$ on $S$ satisfying $E^2 <0$. 
Then either 
\begin{enumerate}
    \item $\Delta_f \cap E = \emptyset$, or 
    \item $E$ is a connected component of $\Delta_f$ and $\Delta_f$ is smooth around $E$. 
\end{enumerate}
\end{cor}

\begin{proof}
If $f^{-1}(E)$ is irreducible, then (1) holds (Proposition \ref{p-FCB-(-1)}). 
Hence we may assume that $f^{-1}(E)$ is not irreducible. 
Then (2) holds by Proposition \ref{p-ele-tf}. 
\end{proof}

\subsection{Non-trivial Fano conic bundles}

\begin{dfn}
We say that $f: X \to S$ is a {\em trivial} conic bundle or a {\em trivial} $\P^1$-bundle if 
$X \simeq S \times \P^1$ and $f$ coincides with the first projection. 
We say that $f$ is {\em non-trivial} if $f$ is not trivial. 
\end{dfn}

Let $f: X \to S$ a non-trivial Fano conic bundle. 
The purpose of this subsection is to prove that $S$ is isomorphic to one of $\P^2, \P^1 \times \P^1$, and $\F_1$ (Proposition \ref{p-FCB-triv}). 
We start with the following auxiliary result. 

\begin{lem}\label{l-FCB-criterion}
Let $\varphi : V \to T$ and $\psi : T \to \P^1$ be $\P^1$-bundles: 
\[
V \xrightarrow{\varphi} T \xrightarrow{\psi} \P^1, 
\]
where $T$ is a smooth projective surface and $V$ is a smooth projective threefold. 
Then the following hold. 
\begin{enumerate}
\item If the $\P^1$-bundle $\varphi|_{\varphi^{-1}(T_z)} : \varphi^{-1}(T_z) \to T_z$ is trivial 
for every closed point $z \in \P^1$ and the fibre $T_z := \psi^{-1}(z)$ over $z$, then 
there exist a $\P^1$-bundle $\rho: V_0 \to \P^1$ and  a cartesian diagram: 
\[
\begin{tikzcd}
V \arrow[r] \arrow[d, "\varphi"'] & V_0\arrow[d, "\rho"]\\
T \arrow[r, "\psi"]& \P^1.
\end{tikzcd}
\]
\item 
If $V$ is a Fano threefold, then one of the following holds. 
\begin{enumerate}
    \item $\varphi^{-1}(T_z) \simeq \P^1 \times \P^1$ for every closed point $z \in Z$. 
    \item $\varphi^{-1}(T_z) \simeq \F_1$ for every closed point $z \in Z$. 
\end{enumerate}
\end{enumerate}
\end{lem}

\begin{proof}
By $\Br\,(T) =0$ \cite[Proposition 2.7(3)]{ATIII}, $\varphi: V \to T$ is a $\P^1$-bundle if and only if 
$\varphi : V \to T$ is a flat morphism from a smooth projective threefold such that any fibre of $\varphi$ is $\P^1$ \cite[Proposition 2.8]{ATIII}. 

Let us show (1). 
We can write $V = \P_T(E)$ for some locally free sheaf $E$ on $T$ of rank two. 
It is enough to find a locally free sheaf $E'$ on $\P^1$ such that  $E \simeq \psi^*E'$. 
Indeed, this implies $V \simeq V_0 \times_{\P^1} T$ for $V_0 := \P_{\P^1}(E')$. 

For a closed point $z \in \P^1$ and the fibre $T_z := \psi^{-1}(z)$ over $z$, 
we have the following cartesian diagrams: 
\[
\begin{tikzcd}
V\arrow[d, "\varphi"] & \varphi^{-1}(T_z) =T_z \times \P^1\arrow[l, hook'] \arrow[d, "\varphi|_{\varphi^{-1}(T_z)} = {\rm pr}_1"]\\
T\arrow[d, "\psi"] & T_z \arrow[l, hook'] \arrow[d]\\
\P^1 & z, \arrow[l, hook']
\end{tikzcd}
\]
where each horizontal arrow is a closed immersion and 
all the vertical arrows are $\P^1$-bundles. 
By $\varphi^{-1}(T_z) = \P_{T_z}(E|_{T_z})$, we have that $E|_{T_z} \simeq L_z^{\oplus 2}$ for some invertible sheaf $L_z$ on $T_z \simeq \P^1$. 
Since $\deg (E|_{T_z})$ is independent of $z \in \P^1$, 
also $d:=\deg L_z =  \frac{1}{2} \deg (E|_{T_z})$ is independent of $z \in \P^1$. 
We may replace $E$ by $E \otimes \MO_T(-d)$. 
Hence the problem is reduced to the case when $E|_{T_z} \simeq \MO_{T_z}^{\oplus 2}$. 
Then the function $z \mapsto h^0(T_z, E|_{T_z})$ is constant with 
$h^0(T_z, E|_{T_z})=2$. 
Hence $\psi_*E$ is a locally free sheaf of rank $2$, and $\psi_*E \otimes k(z) \xrightarrow{\simeq} H^0(T_z, E|_{T_z})$ \cite[Ch. III, Corollary 12.9]{Har77}. 
Then the induced homomorphism $\theta: \psi^*\psi_*E \to E$ is surjective. 
Both $\psi^*\psi_*E$ and $E$ are locally free sheaves of rank $2$, $\theta$ is an isomorphism. 
Hence we get $E \simeq \psi^*E'$ for $E' := \psi_*E$. 
Thus (1) holds.

Let us show (2). 
Fix a closed point $z \in \P^1$  and set $V_z$ to be the fibre of $V \to \P^1$ over $z$. 
We now show that $V_z \simeq \P^1 \times \P^1$ or $V_z \simeq \F_1$. 
Note that  $V_z$ is a $\P^1$-bundle over $T_z \simeq \P^1$. 
It suffices to show that $-K_{V_z}$ is ample, which follows from 
\[
K_{V_z} \sim (K_V + V_z)|_{V_z} \sim K_V|_{V_z}. 
\]
Thus $V_z \simeq \P^1 \times \P^1$ or $V_z \simeq \F_1$. 
Note that the same conclusion holds even when $z$ is the geometric generic point $\overline{\xi}$ of $\P^1$, i.e., $V_{\overline{\xi}} \simeq \P^1 \times \P^1$ or $V_{\overline{\xi}} \simeq \F_1$. 
Since $V$ is a Fano threefold, we have the following left square diagram (Proposition \ref{p-2ray})
\[
\begin{tikzcd}
& V \arrow[ld, "\varphi"'] \arrow[rd, "\alpha"]  &&&& V_z \arrow[ld, "\varphi_z"'] \arrow[rd, "\alpha_z"] \\
T \arrow[rd, "\psi"']& & U \arrow[ld, "\beta"] && T_z \simeq \P^1 \arrow[rd, "\psi_z"']& & U_z \arrow[ld, "\beta_z"]\\
& \P^1 &&&& z\arrow[llll] 
\end{tikzcd}
\]
where 
$\alpha : V \to U$ is the contraction of the extremal ray of $\NE(V/\P^1)$ not corresponding to $\varphi$, 
and $\beta$ is the induced contraction. 
If $z$ is a closed point or the geometric generic point $\overline{\xi}$ of $\P^1$, 
then we obtain the above right square by taking the base change $(-) \times_{\P^1} \{z\}$. 

\begin{claim*}
\begin{enumerate}
\item[(i)]  $\alpha_{\overline \xi}: V_{\overline{\xi}} \to U_{\overline{\xi}}$ is a contraction which is not an isomorphism. 
\item[(ii)] $\alpha_z : V_z \to U_z$ is not a finite morphism for every closed point $z \in \P^1$. 
\end{enumerate}
 
\end{claim*}

\begin{proof}[Proof of Claim]
Let us show (i). 
Note that $\alpha_{\overline \xi}$ is a contraction, 
because $\alpha_{\overline \xi}$ is obtained by a flat base change of the contraction $\alpha$. 
Hence it is enough to show that $\alpha_{\overline \xi}$ is not an isomorphism. 
Suppose that $\alpha_{\overline \xi}: V_{\overline{\xi}} \to U_{\overline{\xi}}$ is  an isomorphism. 
Then also $\alpha_{\xi} : V_{\xi} \to U_{\xi}$ is an isomorphism for the generic point $\xi$ of $\P^1$. 
Hence 
$\alpha : V \to U$ is of type E. 
We then get $\Ex(\alpha) = V_z$ for some closed point $z \in Z$, 
because both $\Ex(\alpha)$ and $V_z$ are prime divisors on $V$, and 
$\Ex(\alpha)$ is contained in some fibre $V_z$. 
Hence it holds that 
\[
\dim \beta^{-1}(z) =\dim \alpha(V_z) = \dim \alpha(\Ex(\alpha)) \leq 1, 
\]
which contradicts the fact that $\beta^{-1}(z)$ is a nonzero effective Cartier divisor on a threefold $U$. 
Thus (i) holds. 
If there exists a closed point $z \in \P^1$ such that $\alpha_z : V_z \to U_z$ is a finite morphism, 
then $\alpha : V \to U$ is birational and $\alpha_{\xi} : V_{\xi} \to U_{\xi}$ is an isomorphism, 
which contradicts (i). Hence (ii) holds. 
This completes the proof of Claim.     
\end{proof}


(a) Assume that $V_{\overline{\xi}} \simeq \P^1 \times \P^1$. 
By Claim(i), we get $\dim U_{\overline \xi} =1$, which implies $\dim U =2$. 
Hence $\alpha : V \to U$ is of type C. 
In particular, both $\alpha$ and $\beta$ are flat of relative dimension one. 
For any closed point $z \in \P^1$, 
$\alpha_z : V_z \to U_z$ is a flat surective morphism to a curve $U_z$ such that a fibre of $\varphi_z$ dominates $U_z$. 
Hence $V_z \simeq \P^1 \times \P^1$. 

(b) Assume that $V_{\overline{\xi}} \simeq \F_1$. 
Then $\dim U_{\overline \xi} =2$, which implies $\dim U =3$. 
Hence $\alpha : V \to U$ is of type E. 
By Claim(ii), 
$\alpha$ is of type $E_1$ and, for any closed point $z \in \P^1$, 
$\alpha_z : V_z \to U_z$ is a generically finite morphism 
which contracts at least one curve. 
Then the Stein factorisation $V_z \to S$ of $\alpha_z : V_z \to U_z$ is a non-trivial birational contraction. 
Hence $V_z \simeq \F_1$. 
\end{proof}

\begin{prop}\label{p-FCB-triv}
Let $f: X \to S$ be a Fano conic bundle. 
Assume that $S \not\simeq \P^2$, $S \not\simeq \P^1 \times \P^1$, and $S \not\simeq \F_1$. 
Then $f: X \to S$ is a trivial conic bundle, i.e., 
an isomorphism $X \simeq S \times \P^1$ holds and $f$ coincides with the first projection. 
\end{prop}

\begin{proof}
Since $S$ is a smooth del Pezzo surface (Proposition \ref{p-FCB-dP}), 
our assumption implies that there is a birational morphism $\alpha : S \to \P^1 \times \P^1$ which is the blowup along finitely many closed points: 
$P_1 \amalg \cdots \amalg P_n$. 
For every $1 \leq i \leq n$, 
let $L_i$ and $M_i$ be the prime divisors on $\P^1 \times \P^1$ 
passing through $P_i$ that are the fibres of the first and second projections, respectively. 
Let $L'_i \subset S$ and $M'_i \subset S$ be the proper transforms of $L_i$ and $M_i$, respectively. 
Set $E_i \subset S$ to be the $(-1)$-curve lying over $P_i$.

\begin{claim*}
$f: X \to S$ is a $\P^1$-bundle, i.e., there exists a locally free sheaf $E$ on $S$ of rank $2$ such that $X$ is $S$-isomorphic to $\P_S(E)$.     
\end{claim*}

\begin{proof}[Proof of Claim]
Suppose $\Delta_f \neq \emptyset$. 
By $\Br(S)=0$ \cite[Proposition 2.7(3)]{ATIII}, it is enough to derive a contradiction  \cite[Proposition 2.8]{ATIII}. 
We can find $a_i, b_i, c_i \in \Z_{>0}$ such that the divisor 
\[
A := \sum_{i=1}^n a_i L'_i + \sum_{i=1}^n b_i M'_i + \sum_{i=1}^n c_i E_i
\]
is ample (indeed, $ m \alpha^*(L_1 +M_1) - \sum_{i=1}^n E_i$ is ample for $m \gg 0$, 
and hence it is enough to add the nef divisor $\sum_{i=1}^n \alpha^*(L_i+M_i)$).  
Hence $\Delta_f \cap A \neq \emptyset$, 
which implies $\Delta_f \cap (L'_i + M'_i + E_i) \neq \emptyset$ for some $i$. 
For example, assume that $\Delta_f \cap L'_i \neq \emptyset$ (the rest of the argument is the same for the other cases). 
Corollary \ref{c-FCB-(-1)-2} implies that $L'_i$ is a connected component of $\Delta_f$. 
Then it follows from $L'_i \cap E_i \neq \emptyset$ that 
$\Delta_f \cap E_i \neq \emptyset$ and $E_i$ is not a connected component 
of $\Delta_f$, which contradicts Corollary \ref{c-FCB-(-1)-2}. 
This completes the proof of Claim. 
\end{proof}

Since $f : X \to S$ is a $\P^1$-bundle, 
we have a cartesian diagram (Proposition \ref{p-FCB-(-1)}):  
\begin{equation}\label{e1-FCB-triv}
\begin{tikzcd}
X \arrow[r, "\beta"] \arrow[d, "f"] & \widetilde{X}\arrow[d, "\widetilde{f}"]\\
S \arrow[r, "\alpha"] & \widetilde S :=\P^1 \times \P^1,  
\end{tikzcd}
\end{equation}
where $\widetilde f: \widetilde X \to \widetilde S= \P^1 \times \P^1$ is a Fano conic bundle. 
Since $f$ is a $\P^1$-bundle and the diagram (\ref{e1-FCB-triv}) is cartesian, 
any fibre of $\widetilde f$ is isomorphic to $\P^1$. 
By ${\rm Br}(\P^1 \times \P^1)=0$, we see that $\widetilde X \simeq \P_{\widetilde S}(\widetilde E)$ 
for some locally free sheaf $\widetilde E$ on $\widetilde S=\P^1 \times \P^1$ of rank $2$ \cite[Subsection 2.3]{ATIII}.

We now show that the $\P^1$-bundle $\widetilde{f}^{-1}( L_1) \to  L_1$ is trivial. 
By $L'_1 \xrightarrow{\simeq} L_1$ and the cartesian diagram (\ref{e1-FCB-triv}), 
it is enough to show that $f^{-1}(L'_1) \to L'_1$ is trivial. 
This follows from Proposition \ref{p-FCB-(-1)} and the fact that $L'_1$ is a $(-1)$-curve. 
Similarly, $\widetilde{f}^{-1}(M_1) \to M_1$ is a trivial $\P^1$-bundle.

Since $\widetilde X$ is a Fano threefold and 
the $\P^1$-bundle $\widetilde{f}^{-1}(L_1) \to  L_1$ is trivial, 
$\widetilde{f}^{-1}(L) \to  L$ is a trivial $\P^1$-bundle for every closed point $x \in \P^1$ and the fibre $L := \{x \} \times \P^1$ over $x$ (Lemma \ref{l-FCB-criterion}(2)). 
By Lemma \ref{l-FCB-criterion}(1), we get the right cartesian square  in the following diagram: 
\begin{equation}\label{e2-FCB-triv}
\begin{tikzcd}
\widetilde{f}^{-1}(M_1)=\widetilde X \times_{\widetilde S} M_1 \arrow[r] \arrow[d, "h"] \arrow[rr, bend left, "\simeq"]& 
\widetilde X \arrow[r] \arrow[d, "\widetilde f"] & W\arrow[d, "g"]\\
M_1 \arrow [r, hook]\arrow[rr, bend right, "\simeq"] & \widetilde S \arrow[r, "{\rm pr}_1"] & \P^1,  
\end{tikzcd}
\end{equation}
where $g: W \to \P^1$ is a $\P^1$-bundle. 
Note that the left square in (\ref{e2-FCB-triv}) is the cartesian diagram 
induced by the closed immersion $M_1 \hookrightarrow \P^1 \times \P^1 =\widetilde S$. 
The induced $\P^1$-bundle $h : \widetilde{f}^{-1}(M_1)=\widetilde X \times_{\widetilde S} M_1 \to M_1$ is trivial, and hence also $g : W \to \P^1$ is trivial. 
Since $\widetilde f : \widetilde X \to \widetilde S$ is a base change of a trivial $\P^1$-bundle $W \to \P^1$, 
also $\widetilde f : \widetilde X \to \widetilde S$ is trivial. 
Finally, by the cartesian diagram (\ref{e1-FCB-triv}), 
$f: X \to S$ is trivial. 
\end{proof}









\section{$\rho=3$}

The purpose of this subsection is to classify Fano threefolds of Picard number $3$. 
In what follows, we overview its proof and contents of this section. 

Let $X$ be a Fano threefold with $\rho(X)=3$. 
By case study depending on the type of extremal rays, 
we show that at least one of (I)-(V) holds. 
\begin{enumerate}
\item[(I)] $X$ has a conic bundle structure over $\P^2$. 
\item[(II)] $X$ has a conic bundle structure over $\F_1$. 
\item[(III)] $X$ is primitive. 
\item[(IV)] 
There exist a line $L$ and a conic $C$ on $\P^3$ such that $L \cap C = \emptyset$ 
and $X \simeq \Bl_{L \amalg C} \P^3$. 
\item[(V)] $X \simeq \Bl_{B_1 \amalg B_2}\,V$, where 
$B_1$ and $B_2$ are smooth curves on  $V$ satisfying one of the following. 
\begin{itemize}
\item $V = \P^3$ and both $B_1$ and $B_2$ are lines. 
\item $V = \P^3$, $B_1$ is a line, and $B_2$ is an elliptic curve of degree $4$.  
\item $V=Q$, and both $B_1$ and $B_2$ are conics. 
\end{itemize}
\end{enumerate}
This result is established in Subsection \ref{ss-pic3-structure} 
(Theorem \ref{t-pic3-structure}). 
Moreover, the cases (III)-(V) are explicitly classified. 
We then classify 
Fano conic bundles over $\P^2$ and $\F_1$ (i.e., (I) and (II)) 
in Subsection \ref{ss-pic3-P2} (Theorem \ref{t-ele-tr-P2}) and  Subsection \ref{ss-pic3-F1} (Theorem \ref{t-F1-pic3}), respectively. 
Finally, in order to check the overlapping, 
we will determine the number of extremal rays and their types 
in Subsection \ref{ss-pic3-classify}. 
For example, if $(-K_X)^3 =24$, then one and only one of the following holds. 
\begin{itemize}
\item Only (I) holds. In this case, $X$ has exactly three extremal rays and all of them are of type $E_1$ (No.\ 3-7, Proposition \ref{p-pic3-7}). 
\item Both (I) and (II) hold. 
 In this case, $X$ has exactly three extremal rays. 
 One of them is of type $C_1$, and the others are of type $E_1$ 
(No.\ 3-8, Proposition \ref{p-pic3-8}). 
\end{itemize}
For later usage, 
we give a classification of Fano conic bundles $X \to \P^1 \times \P^1$ with $\rho(X)=3$ (Subsection \ref{ss FCB P1P1 rho3}). 

\subsection{Reduction to conic bundle case}\label{ss-pic3-structure}

The purpose of this subsection is to show that, except for one case (No.\ 3-18),   
every Fano threefold $X$ with $\rho(X)=3$  has a conic bundle structure  (Theorem \ref{t-pic3-structure}). 
More precisely, Theorem \ref{t-pic3-structure} reduces the classification of Fano threefolds $X$ with $\rho(X)=3$ to the case when $X$ has a conic bundle structure over $\P^2$ or $\F_1$.

\begin{lem}\label{l-exactly3}
Let $X$ be a Fano threefold with $\rho(X)=3$. 
\begin{enumerate}
\item 
Let $\varphi : X \to Z$ be a contraction to a projective normal variety $Z$. 
Let $F$ be the extremal face corresponding to $\varphi$, i.e., 
$F := \NE(X) \cap (\varphi^*A_Z)^{\perp}$ for an ample Cartier divisor $A_Z$ on $Z$. 
Then the following hold. 
\begin{enumerate}
\item $\dim F=1$ if and only if $\rho(Z)=2$. 
\item $\dim F=2$ if and only if $\rho(Z)=1$. 
\end{enumerate}
Here we set $\dim F := \dim\,\langle F \rangle$ for the linear subspace $\langle F \rangle$ of $N_1(X)$ generated by $F$. 
\item 
For every $i \in \{1, 2, 3\}$, let $\varphi_i : X \to Z_i$ be a contraction to a projective normal variety $Z_i$ with $\rho(Z_i)=1$. 
Assume that 
the extremal faces corresponding to $\varphi_1, \varphi_2, \varphi_3$ are different 
and 
none of 
$\varphi_1 \times \varphi_2 : X \to Z_1 \times Z_2, \varphi_2 \times \varphi_3 : X \to Z_2 \times Z_3, \varphi_3 \times \varphi_1 : X \to Z_3 \times Z_1$ is a  finite morphism. 
Then $X$ has exactly three extremal rays. 
\item 
Assume that we have the following diagram 
\[
\begin{tikzcd}
& & X \arrow[ld, "f_1"'] \arrow[rd, "f_2"] \arrow[dd, "\varphi"] \arrow[lldd, "\varphi_1"', bend right=40] \arrow[rrdd, "\varphi_2", bend left=40] \\
& Y_1  \arrow[ld, "h_1"'] \arrow[rd, "g_1"]& & Y_2 \arrow[ld, "g_2"'] \arrow[rd, "h_2"]\\
Z_1 & & V & & Z_2
\end{tikzcd}
\]
where each $f_i : X \to Y_i$ is a contraction of an extremal ray and 
$\varphi_1 : X \to Z_1, \varphi: X \to V, \varphi_2 : X \to Z_2$ are contractions of two-dimensional extremal faces. 
We further assume that $f_1 : X \to Y_1$ and $f_2 : X\to Y_2$ correspond to distinct extremal rays and 
$\varphi_1 \times \varphi_2 : X \to Z_1 \times Z_2$ is not a finite morphism. 
Then $X$ has exactly three extremal rays and 
the Stein factorisation of $\varphi_1 \times \varphi_2 :X \to Z_1 \times Z_2$ is the contraction of 
the extremal ray not corresponding to $f_1: X \to Y_1$ nor $f_2 : X \to Y_2$. 
\end{enumerate}
\end{lem}

\begin{proof}
Let us show (1). 
The assertion (a) holds by \cite[Proposition 3.12]{TanII}. 
Then (a) implies (b) 
(for an extremal ray $R$ of $\NE(X)$ contained in $F$, 
its contraction $X \to Y$ 
satisfies $\rho(Y)=2$ and induces a factorisation $\varphi : X \to Y \to Z$).  
Thus (1) holds. 

Let us show (2). 
Given a curve $C$ on $X$, 
$(\varphi_1 \times \varphi_2)(C)$ is a point if and only if $\varphi_1(C)$ and $\varphi_2(C)$ is a point. 
For the Stein factorisation $\varphi_1 \times \varphi_2 : X \xrightarrow{\psi_{12}} Z_{12} \to Z_1 \times Z_2$ of $\varphi_1 \times \varphi_2$, 
$\psi_{12}$ is the contraction corresponds to the extremal face $F_1 \cap F_2$. 
It follows from (1)  that $\dim F_1 = \dim F_2 =2$. 
By $F_1 \neq F_2$, either $F_1 \cap F_2 = \{0\}$ or $F_1 \cap F_2$ is an extremal ray. 
Since $\varphi_1 \times \varphi_2 : X \to Z_1 \times Z_2$ is not a finite morphism, 
$F_1 \cap F_2$ is an extremal ray. 
By symmetry, also $F_2 \cap F_3$ and $F_3 \cap F_1$ are extremal rays. 
Thus (2) holds. 
The assertion (3) follows from (2). 
\end{proof}

\begin{lem}\label{l-exactly3 2P^1}
Let $X$ be a Fano threefold with $\rho(X)=3$. 
Assume that we have the following diagram 
\[
\begin{tikzcd}
& & X \arrow[ld, "f_1"'] \arrow[rd, "f_2"] \arrow[dd, "\varphi"] \arrow[lldd, "\varphi_1"', bend right=40] \arrow[rrdd, "\varphi_2", bend left=40] \\
& Y_1  \arrow[ld, "h_1"'] \arrow[rd, "g_1"]& & Y_2 \arrow[ld, "g_2"'] \arrow[rd, "h_2"]\\
Z_1:=\P^1 & & V & & Z_2 :=\P^1
\end{tikzcd}
\]
where each $f_i : X \to Y_i$ is a contraction of an extremal ray and 
$\varphi_1 : X \to Z_1, \varphi: X \to V, \varphi_2 : X \to Z_2$ are contractions of two-dimensional extremal faces. 
We further assume that $f_1 : X \to Y_1$ and $f_2 : X\to Y_2$ correspond to distinct extremal rays. 
Then $X$ has exactly three extremal rays and 
$\varphi_1 \times \varphi_2 :X \to Z_1 \times Z_2$ is the contraction of 
the extremal ray not corresponding to $f_1: X \to Y_1$ nor $f_2 : X \to Y_2$. 
\end{lem}

\begin{proof}
Since $\varphi := \varphi_1 \times \varphi_2 :X \to Z_1 \times Z_2 = \P^1 \times \P^1$ is not a finite morphism, 
$X$ has exactly three extremal rays (Lemma \ref{l-exactly3}). 
Take the Stein factorisation  of $\varphi$: 
\[
\varphi :X \xrightarrow{\psi} S \xrightarrow{\theta}  Z_1 \times Z_2
\]
It is enough to show that $\theta$ is an isomorphism (Lemma \ref{l-exactly3}). 

We now show that $S \simeq \P^1 \times \P^1$ and 
the induced composite morphism 
$\pi_i : S \to Z_1 \times Z_2 \xrightarrow{{\rm pr}_i} Z_i$ is a contraction for each $i \in \{1, 2\}$. 
By $(\pi_i)_*\MO_S = (\pi_i)_* \psi_*\MO_X 
= (\varphi_i)_*\MO_X = \MO_{Z_i}$, 
each $\pi_i : S \to Z_i$ is a contraction. 
By $\dim S = \dim (Z_1 \times Z_2)=2$ and $3 = \rho(X) > \rho(S) \geq \rho(Z_1\times Z_2) =\rho(\P^1 \times \P^1)=2$, 
$\psi : X \to S$ is a contraction of an extremal ray. 
In particular, 
$\psi: X \to S$ is a Fano conic bundle, and $S$ is a del Pezzo surface (Proposition \ref{p-FCB-dP}). 
As $S$ has two contractions to $\P^1$, we obtain $S \simeq \P^1 \times \P^1$.

By $(\pi_i)_*\MO_S = \MO_{Z_i}$, 
we may assume that $\pi_i = \pr_i$ for each $i \in \{1, 2\}$. 
Then each of $S$ and $Z_1 \times Z_2$ are the fibre product of 
$Z_1 \to \Spec\,k \leftarrow Z_2$. 
Hence we  get $\theta: S \xrightarrow{\simeq} Z_1 \times Z_2$ 
by the universal property of fibre products. 
\end{proof}

\begin{lem}\label{l-K-disjoint-blowup}
Let $V$ be a Fano threefold with $\rho(V)=1$. 
Let $B_1$ and $B_2$ be smooth curves on $V$ such that $B_1 \cap B_2 = \emptyset$ and $X := \Bl_{B_1 \amalg B_2} V$ is Fano. 
For each $i \in \{1, 2\}$, set $Y_i  := \Bl_{B_i} V$ and 
let $h_i : Y_i \to Z_i$ be the contraction of the extremal ray not corresponding to the blowup $g_i : \Bl_{B_i} V \to V$ (note that each $Y_i$ is a Fano threefold with $\rho(Y_i)=2$ by Corollary \ref{c-disjoint-blowup}). 
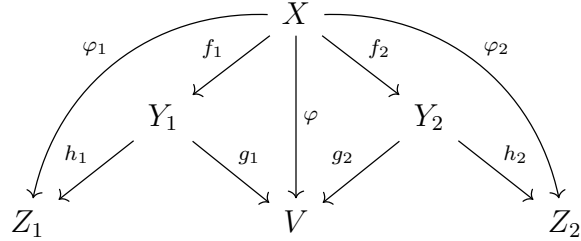
\[
\begin{tikzcd}
& & X \arrow[ld, "f_1"'] \arrow[rd, "f_2"] \arrow[dd, "\varphi"] \arrow[lldd, "\varphi_1"', bend right=40] \arrow[rrdd, "\varphi_2", bend left=40] \\
& Y_1  \arrow[ld, "h_1"'] \arrow[rd, "g_1"]& & Y_2 \arrow[ld, "g_2"'] \arrow[rd, "h_2"]\\
Z_1 & & V & & Z_2
\end{tikzcd}
\]
Let $H_V, H_{Z_1}, H_{Z_2}$ be the ample Cartier divisors which generate of $\Pic\,V, \Pic\,Z_1, \Pic\,Z_2$, respectively. 
Then the following holds. 
\begin{enumerate}
\item 
$-K_X 
\sim -f_1^*K_{Y_1} -f_2^*K_{Y_2} - (-\varphi^*K_V) \sim \varphi_1^*H_{Z_1} + \varphi_2^*H_{Z_2} + (\mu_1+\mu_2-r_V) \varphi^*H_V$, 
where $r_V$ denotes the index of $V$ and 
each $\mu_i$ is the length of the extremal ray corresponding to $h_i$. 
\item 
If $\varphi_1 \times \varphi_2 : X \to Z_1 \times Z_2$ is not a finite morphism, 
then $\mu_1 + \mu_2 > r_V$. 
\end{enumerate}
\end{lem}

\begin{proof}
Let us show (1).
We have 
\[
K_X \sim f_1^*K_{Y_1} +E_{X/Y_1}, \quad K_X \sim f_2^*K_{Y_2}+E_{X/Y_2}, 
\quad 
K_X \sim \varphi^*K_V + E_{X/Y_1}+E_{X/Y_2}
\]
for the exceptional prime divisors $E_{X/Y_1}$ and $E_{X/Y_2}$ of $f_1$ and $f_2$, respectively. 
Then we obtain 
\begin{eqnarray*}
-K_X &=& -K_X -K_X +K_X \\
&\sim& -(f_1^*K_{Y_1} +E_{X/Y_1}) -(f_2^*K_{Y_2}+E_{X/Y_2}) 
+(\varphi^*K_V + E_{X/Y_1}+E_{X/Y_2}) \\
&=&  -f_1^*K_{Y_1} -f_2^*K_{Y_2} - (-\varphi^*K_V).     
\end{eqnarray*}
It follows from Proposition \ref{p-pic2-Pic}(2) that 
$-K_{Y_1} \sim h_1^*H_{Z_1} + \mu_1 g_1^*H_V$ and $-K_{Y_2} \sim h_2^*H_{Z_2} + \mu_2 g_2^*H_V$. Hence 
\begin{eqnarray*}
-K_X &\sim & -f_1^*K_{Y_1} -f_2^*K_{Y_2} - (-\varphi^*K_V) \\
&\sim&  (\varphi_1^*H_{Z_1} + \mu_1 \varphi^*H_V) + 
(\varphi_2^*H_{Z_2} + \mu_2 \varphi^*H_V) - \varphi^*(r_VH_{V}) \\
&=&  \varphi_1^*H_{Z_1} + \varphi_2^*H_{Z_2} + (\mu_1+\mu_2-r_V) \varphi^*H_V. 
\end{eqnarray*}
Thus (1) holds.

Let us show (2). 
Assume that $\varphi_1 \times \varphi_2 : X \to Z_1 \times Z_2$ is not a finite morphism. 
Then we can find a curve $C$ on $X$ such that $\varphi_1(C)$ and $\varphi_2(C)$ are points. 
In particular, $\varphi_1^*H_1 \cdot C = \varphi_2^*H_2 \cdot C=0$. 
Hence 
\[
0< (-K_X) \cdot C  = 
(\varphi_1^*H_{Z_1} + \varphi_2^*H_{Z_2} + (\mu_1+\mu_2-r_V) \varphi^*H_V) \cdot C 
= (\mu_1+\mu_2-r_V) \varphi^*H_V \cdot C. 
\]
Since $\varphi^*H_V$ is nef, we obtain $\varphi^*H_V \cdot C \geq 0$. 
This, together with $(\mu_1+\mu_2-r_V) \varphi^*H_V \cdot C >0$, implies 
$\varphi^*H_V \cdot C >0$ and $\mu_1+\mu_2-r_V>0$. 
Thus (2) holds. 
\end{proof}

\begin{prop}\label{p-V-2dP}
We use the same notation as in the statement of Lemma \ref{l-K-disjoint-blowup}. 
Assume that $Z_1 \simeq Z_2 \simeq \P^1$. 
Then the following hold. 
\begin{enumerate}
    \item $X$ has exactly three extremal rays and $h := \varphi_1 \times \varphi_2 : X \to Z_1 \times Z_2 (\simeq \P^1 \times \P^1)$ is the contraction of the extremral ray not corresponding to $f_1$ nor $f_2$. 
    \item If $(-K_{Y_1})^3 \geq (-K_{Y_2})^3$, then one of the following holds. 
    \begin{enumerate}
        \item $(-K_X)^3 = 44$, $V=\P^3$, $Y_1$ and $Y_2$ are of No.\ 2-33, 
        $B_1$ and $B_2$ are lines, and $h$ is of type $C_2$ ($X$ is 3-25). 
        \item $(-K_X)^3 = 22$, $V=\P^3$, $Y_1$ is of No.\ 2-33, 
        $Y_2$ is 2-25, $B_1$ is a line, and $B_2$ is an elliptic curve of degree $4$, $h$ is of type $C_1$, and $\Delta_h$ is of bidegree $(2, 3)$ ($X$ is 3-6).
        \item $(-K_X)^3 = 26$, $V=Q$, $Y_1$ and $Y_2$ are  of No.\ 2-29, $B_1$ and $B_2$ are conics, $h$ is of type $C_1$, and $\Delta_h$ is of bidegree $(2, 2)$ ($X$ is 3-10). 
    \end{enumerate}
\end{enumerate}
\end{prop}

\begin{proof}
Fix a closed point $P_i$ of $Z_i = \P^1$ and 
let $(Y_i)_{P_i}$ be the fibre of $h_i : Y_i \to Z_i$ over $P_i$. 
Then $Y_i$ belongs to one of Table \ref{table-DE1} (Subsection \ref{ss-table-pic2}).  
The assertion (1) follows from Lemma \ref{l-exactly3 2P^1}.



Let us show (2). 
First of all, we prove that $V \simeq \P^3$ or $V \simeq Q$. 
Otherwise, it follows from  the classification list (Table \ref{table-DE1}) that 
$r_V = 2$ and $h_i : Y_i \to Z_i =\P^1$ is of type $D_1$ for each $i \in \{1, 2\}$. 
For the length $\mu_i$ of the extremal ray corresponding to $h_i : Y_i \to Z_i =\P^1$, 
we get $\mu_1 + \mu_2 = 1+1=2 =r_V$, which contradicts Lemma \ref{l-K-disjoint-blowup}. 
Thus $V \simeq \P^3$ or $V \simeq Q$. 
In what follows, we shall use the following:
\begin{equation}\label{e1-V-2dP}
0 < (-K_X)^3 = (-K_V)^3 - ( (-K_V)^3 -(-K_{Y_1})^3) - ( (-K_V)^3 -(-K_{Y_2})^3). 
\end{equation}

\medskip

Assume $V=Q$. 
Then $Y_i$ is 2-7 or 2-29 for each $i \in \{1, 2\}$ (Table \ref{table-DE1}). 
If $Y_i$ is 2-7 (resp. 2-29), then 
$(-K_V)^3 - (-K_{Y_i})^3 = 40$ (resp. $=14$). 
If $Y_1$ is 2-7, then (\ref{e1-V-2dP}) leads to the following contradiction: 
\[
0< (-K_X)^3 =  (-K_V)^3 - ( (-K_V)^3 -(-K_{Y_1})^3) - ( (-K_V)^3 -(-K_{Y_2})^3) 
\leq  54 -40 -14 =0. 
\]
By symmetry, both $Y_1$ and $Y_2$ must be 2-29. 
In this case, 
(c) holds except for the assertion on $\Delta_h$. 

\medskip

Assume $V= \P^3$. 
Then $Y_i$ is 2-4, 2-25, or 2-33 (Table \ref{table-DE1}). 
If $Y_1$ and $Y_2$ are 2-25, then we would get the following contradiction by  (\ref{e1-V-2dP}): 
\[
0< (-K_X)^3 =  (-K_V)^3 - ( (-K_V)^3 -(-K_{Y_1})^3) - ( (-K_V)^3 -(-K_{Y_2})^3) 
=64 -32 -32  =0. 
\]
Similarly, none of $Y_1$ nor $Y_2$ is 2-4. 
Therefore, $(Y_1, Y_2)$ is either (2-33, 2-33) or (2-33, 2-25), 
i.e., (a) or (b) holds except for the assertions on $h$ and $\Delta_h$. 

\medskip

It is enough to compute the bidegree $(a_1, a_2)$ of $\Delta_h$ for the cases (a)-(c). 
Recall that the following holds for any divisor $D$ on $S := Z_1 \times Z_2$ satisfying $(-K_S) \cdot D=2$ \cite[Proposition 3.16]{ATIII}:
\[
\Delta_h \cdot D = 4 (-K_S) \cdot D - (-K_X)^2 \cdot h^*D = 8 -  (-K_X)^2 \cdot h^*D. 
\]
For each $i \in \{1, 2\}$, let $P_i$ be a closed point on $Z_i = \P^1$. 
Set $D_i := {\rm pr}_i^*P_i$. 
We have  $(-K_S) \cdot D_i = 2$ and $\Delta_h \cdot D_i = a_{3-i}$ 
(note that $(i, 3-i) \in \{ (1, 2), (2, 1)\}$). 
Then 
\[
(-K_X)^2 \cdot h^*D_i = (-f_i^*K_{Y_i}- \Ex(f_i))^2 \cdot f_i^*h_i^*P_i 
= (-K_{Y_i})^2 \cdot (Y_i)_{P_i}- f_i(\Ex(f_i)) \cdot h_i^*P_i, 
\]
where 
$(Y_i)_{P_i} = h_i^*P_i$ and the last equality holds by 
$f_i^*K_Y \cdot \Ex(f_i)  \cdot f_i^*h_i^*P_i=0$ and 
$\Ex(f_i)^2 \cdot f_i^*h_i^*P_i = -f_i(\Ex(f_i)) \cdot h_i^*P_i$.

(a) 
By $-K_{Y_1} \sim h_1^*P_1 + 3g_1^*\MO_{\P^3}(1)$, 
we obtain 
$4 =-K_{\P^3} \cdot B_2 = -K_{Y_1} \cdot g_1^{-1}(B_2) = (h_1^*P_1 + 3g_1^*\MO_{\P^3}(1))\cdot g_1^{-1}(B_2) =h_1^*P_1 \cdot f_1(\Ex(f_1)) +3$. 
Thus $ h_1^*P_1  \cdot f_1(\Ex(f_1))=1$, 
which implies  $(-K_X)^2 \cdot h^*D_1 = 9 - 1 =8$ 
(indeed, for the generic point $\xi$ of $Z_1 = \P^1$, 
the equality $ h_1^*P_1  \cdot f_1(\Ex(f_1))=1$ implies that 
the birational morphism $(f_1)_{\xi} : X_{\xi} \to (Y_1)_{\xi}$ is a blowup at a rational point, and hence 
$(-K_X)^2 \cdot h^*D_1 = (-K_{X_{\xi}})^2 = (-K_{(Y_1)_{\xi}})^2 -1 = 9-1$). 
Hence 
\[
a_2 = \Delta_h \cdot D_1 = 
8 -(-K_X)^2 \cdot h^*D_1= 8 -8=0.
\]
By symmetry, we get $a_1=0$. Thus $\Delta_h =0$, i.e.,  $h$ is of type $C_2$.

(b) 
By $-K_{Y_1} \sim h_1^*P_1 + 3g_1^*\MO_{\P^3}(1)$, 
we obtain 
$16 =-K_{\P^3} \cdot B_2 = -K_{Y_1} \cdot g_1^{-1}(B_2) = (h_1^*P_1 + 3g_1^*\MO_{\P^3}(1))\cdot g_1^{-1}(B_2) =h_1^*P_1 \cdot f_1(\Ex(f_1)) +12$. 
Thus $ h_1^*P_1  \cdot f_1(\Ex(f_1))=4$, 
which implies  $(-K_X)^2 \cdot  h^*D_1 = 9 - 4 =5$. 
Hence 
\[
a_2 = \Delta_h \cdot D_1 = 8 -(-K_X)^2 \cdot h^*D_1= 8 -5=3.
\]

By $-K_{Y_2} \sim h_2^*P_2 + 2g_2^*\MO_{\P^3}(1)$, 
we obtain 
$4 =-K_{\P^3} \cdot B_1 = -K_{Y_2} \cdot g_2^{-1}(B_1) = 
(h_2^*P_2 + 2g_2^*\MO_{\P^3}(1))\cdot g_2^{-1}(B_1) 
=h_2^*P_2 \cdot f_2(\Ex(f_2)) +2$. 
Thus $ h_2^*P_2  \cdot f_2(\Ex(f_2))=2$, 
which implies  $(-K_X)^2 \cdot D_2 = 8 - 2 =6$. 
Hence 
\[
a_1 = \Delta_h \cdot D_2 = 8 -(-K_X)^2 \cdot h^*D_1= 8 -6=2.
\]
 Thus $\Delta_h$  is of bidegree $(2, 3)$ and $h$ is of type $C_1$.


(c) 
By $-K_{Y_1} \sim h_1^*P_1 + 2g_1^*\MO_{Q}(1)$, 
we obtain 
$6 =-K_{Q} \cdot B_2 = -K_{Y_1} \cdot g_1^{-1}(B_2) 
= (h_1^*P_1 + 2g_1^*\MO_{Q}(1))\cdot g_1^{-1}(B_2) =h_1^*P_1 \cdot f_1(\Ex(f_1)) +4$. 
which implies  $(-K_X)^2 \cdot  h^*D_1 = 8 - 2 =6$.  
Hence 
\[
a_2 = \Delta_h \cdot D_1 = 8 -(-K_X)^2 \cdot h^*D_1 =8-6 =2.
\]
By symmetry, we get $a_1=2$. Thus $\Delta_h$  is of bidegree $(2, 2)$ and $h$ is of type $C_1$. 
\qedhere



\end{proof}

  \begin{center}
\begin{longtable}{ccp{10cm}c}
No. & $(-K_{Y_i})^3$ & Description & types \\ \hline
2-1 & $4$ & $(-K_{Y_i})^2 \cdot (Y_i)_{P_i}=1$ & $D_1$\\ 
&  & blowup of $V_1$ along an elliptic curve of degree $1$ & $E_1$\\ \hline
2-3 & $8$ & $(-K_{Y_i})^2 \cdot (Y_i)_{P_i}=2$ & $D_1$\\ 
 &  & blowup of $V_2$ along an elliptic curve of degree $2$ & $E_1$\\ \hline
2-4 & $10$ & $(-K_{Y_i})^2 \cdot (Y_i)_{P_i}=3$ & $D_1$\\ 
 &  & blowup of $\P^3$ along a curve of genus $10$ degree $9$ & $E_1$\\ \hline
2-5 & $12$ & $(-K_{Y_i})^2 \cdot (Y_i)_{P_i} =3$ & $D_1$\\ 
&  & blowup of $V_3$ along an elliptic curve of degree $3$ & $E_1$\\ \hline
2-7 & $14$ & $(-K_{Y_i})^2 \cdot (Y_i)_{P_i}=4$ & $D_1$\\ 
 &  & blowup of $Q$ along a curve of genus $5$ degree $8$ & $E_1$\\ \hline
 2-10 & $16$ & $(-K_{Y_i})^2 \cdot (Y_i)_{P_i}=4$ & $D_1$\\ 
 &  & blowup of $V_4$ along an elliptic curve of degree $4$ & $E_1$\\ \hline
2-14 & $20$ & $(-K_{Y_i})^2 \cdot (Y_i)_{P_i} = 5$ & $D_1$\\ 
 &  & blowup of $V_5$ along an elliptic curve of degree $5$ & $E_1$\\ \hline
2-25 & $32$ & $(-K_{Y_i})^2 \cdot (Y_i)_{P_i} =8$ & $D_2$\\ 
 &  & blowup of $\P^3$ along an elliptic curve of degree $4$ & $E_1$\\ \hline
2-29 & $40$ & $(-K_{Y_i})^2 \cdot (Y_i)_{P_i} =8$ & $D_2$\\ 
 &  & blowup of $Q$ along a conic & $E_1$\\ \hline
2-33 & $54$ & $(-K_{Y_i})^2 \cdot (Y_i)_{P_i} =9$ & $D_3$\\ 
 &  & blowup of $\P^3$ along a line & $E_1$\\ \hline\\
\\
      \caption{Fano threefolds $Y_i$ with $\rho(Y_i)=2$ and 
      whose extremal rays are of type 
      $D+E_1$}\label{table-DE1}
      \end{longtable}
  \end{center}

\begin{lem}\label{l-Ex-differ}
Let $Y$ be a Fano threefold with $\rho(Y)=2$. 
For the extremal rays $R_1$ and $R_2$ of $\NE(Y)$, 
let $f_1 : Y \to Z_1$ and $f_2 : Y \to Z_2$ be the contractions of $R_1$ and $R_2$, respectively. 
Assume that both $R_1$ and $R_2$ are of type $E$. 
Then $\Ex(f_1) \neq \Ex(f_2)$. 
\end{lem}

\begin{proof}
Suppose $\Ex(f_1) = \Ex(f_2) =:D$. 
Let us derive a contradiction. 
Then both $f_1$ and $f_2$ are of type $E_1$, as otherwise 
one of $R_1$ and $R_2$ would contain the numerical class of every curve on $D$. 
In particular, each $Z_i$ is a Fano threefold. 
Moreover, we obtain $-f_1^*K_{Z_1} \sim -K_Y+D \sim -f_2^*K_{Z_2}$, which leads to the following contradiction: 
\[
R_1 = \NE(X) \cap  (-f_1^*K_{Z_1})^{\perp} = \NE(X) \cap (-f_2^*K_{Z_2})^{\perp} = R_2. 
\]
\end{proof}

\begin{lem}\label{l-P3-disjoint-E}
We use the same notation as in the statement of Lemma \ref{l-K-disjoint-blowup}. 
Assume that $h_1$ is birational. 
Then $V =\P^3$, $B_1$ is a conic, and $B_2$ is a line. 
\end{lem}

\begin{proof}
By $\Ex(g_1) \neq \Ex(h_1)$ (Lemma \ref{l-Ex-differ}), 
$g_1(\Ex(h_1))$ is an ample divisor on $V$. 
Then $\Ex(h_1)$ must intersect the curve $B_{2, Y_1} := g_1^{-1}(B_2)$, i.e., $\Ex(h_1) \cap B_{2, Y_1} \neq \emptyset$. 
By $B_1 \cap B_2 = \emptyset$, we have 
$\Ex(h_1) \cdot B_{2, Y_1} = g_1(\Ex(h_1)) \cdot B_{2}>0$. 
Hence $B_{2, Y_1}$ is not contracted by $h_1$. 

\begin{claim*}
$h_1 : Y_1 \to Z_1$ is of type $E_2$. 
\end{claim*}

\begin{proof}[Proof of Claim]
Suppose that $h_1$ is of type $E_3, E_4$, or $E_5$. 
Since $\Ex(h_1)$ is covered by the curves $C$ on $\Ex(h_1)$ satisfying $-K_{Y_1} \cdot C =1$, 
we can find a curve $C$ on $\Ex(h_1)$ 
such that $-K_{Y_1} \cdot C=1$ and $C$ intersects $B_{2, Y_1}$ properly. 
This contradicts Lemma \ref{l-line-meeting}. 
Suppose that $h_1$ is of type $E_1$. 
Since $B_{2, Y_1}$ is not contracted by $h_1 : Y_1 \to Z_1$, 
we can find a a one-dimensional fibre $C$ of 
$h_1 : Y_1 \to Z_1$ which intersects $B_{2, Y_1}$ properly. 
Again by $-K_{Y_1} \cdot C=1$ and Lemma \ref{l-line-meeting}, we get a contradiction. 
This completes the proof of Claim. 
\end{proof}

Hence the extremal rays of $Y_1$ are of type $E_1$ and $E_2$. 
By the classification list (Subsection \ref{ss-table-pic2}), 
$Y_1$ is of No.\ 2-30, $V = \P^3$, and $B_1$ is a conic. 


It is enough to show that $B_2$ is a line. 
Suppose $\deg B_2 \geq 2$. 
By $B_1 \cap B_2 = \emptyset$, $\langle B_1 \rangle \cap B_2$ is zero-dimensional with 
$\dim \MO_{\langle B_1 \rangle \cap B_2} \geq 2$. 
We can find a line $L$ on $\langle B_1 \rangle =\P^2$ 
satisfying $\dim \MO_{L \cap B_2} \geq 2$. 
We then get  
$\dim (\MO_{L \cap B_1} \oplus \MO_{L \cap B_2}) \geq 2 +2 =4$. 
Then the proper transform $L_X$ of $L$ on $X$ satisfies 
$K_X \cdot L_X = (f^*K_{\P^3}+E_1 +E_2) \cdot L_X \geq -4+4 =0$, 
which contradicts the ampleness of $-K_X$.
\qedhere

\end{proof}

\begin{prop}\label{p-pic2-structure}
Let $Y$ be a Fano threefold with $\rho(Y)=2$. 
Then one of the following holds. 
\begin{enumerate}
\item[(i)] $Y$ has a conic bundle structure over $\P^2$. 
\item[(ii)] $Y$ is isomorphic to a blowup of $\P^3$ along a smooth curve. 
\item[(iii)] $Y$ is isomorphic to a blowup of $Q$ along a smooth curve. 
\item[(iv)] 
The extremal rays of $Y$ are of type $E_1$ and $D$. Furthermore, the following hold. 
\begin{itemize}
\item 
The contraction of the extremal ray of type $E_1$ is a blowup $Y \to V_d$ along an elliptic curve of degree $d$. 
\item 
For the contraction $Y \to \P^1$ of the extremal ray of type $D$, 
it holds that $(-K_Y)^2 \cdot Y_t =d$ for a closed point $t$ on $\P^1$ and its fibre $Y_t$. 
\end{itemize}
\end{enumerate}
\end{prop}

\begin{proof}
This follows from Table \ref{table-pic2} in Subsetion \ref{ss-table-pic2} (note that the case (iv) corresponds to No.\ 2-1, 2-3, 2-5, 2-10, 2-14). 
\end{proof}

\begin{lem}\label{l-pic3-structure}
Let $X$ be a Fano threefold with $\rho(X)=3$. 
Then one of (1)--(3) holds. 
\begin{enumerate}
\item 
\begin{enumerate}
    \item $X$ is primitive, 
    \item $X$ has a conic bundle structure over $\P^2$, or 
   \item $X$ has a conic bundle structure over $\F_1$. 
\end{enumerate}
\item $X$ is isomorphic to a blowup of $\P^3$ along a disjoint union of 
two smooth curves. 
\item $X$ is isomorphic to a blowup of $Q$ along a disjoint union of 
two smooth curves. 
\end{enumerate}
\end{lem}
\begin{proof}
If $X$ is primitive, then (1) holds. 
In what follows, we assume that $X$ is imprimitive. 
Then there is a blowup $f: X \to Y$ of a Fano threefold $Y$ with $\rho(Y)=2$ 
along a smooth curve $\Gamma$ on $Y$. 
By Proposition \ref{p-pic2-structure}, 
one of (i)-(iv) in Proposition \ref{p-pic2-structure} holds. 

(i) If Proposition \ref{p-pic2-structure}(i) holds, then 
we have a Fano conic bundle $Y \to \P^2$, and hence (b) or (c) in (1) holds (Proposition \ref{p-FCB-centre}).

(ii) Assume Proposition \ref{p-pic2-structure}(ii). 
We have a blowup $g: Y \to \P^3$ along a smooth curve $B$ on $\P^3$. 
Set $E := \Ex(g)$. 
If $\Gamma \cap E =\emptyset$, then (2) holds. 
We may assume that $\Gamma \cap E \neq \emptyset$. 
If $\Gamma$ is not a fibre of the $\P^1$-bundle structure $E \to B$, 
then we would get a contradiction (Lemma \ref{l-line-meeting}). 
Hence $\Gamma$ is a fibre of $E \to B$. 
In this case, 
we have a blowup $X \to V_7$ along a smooth curve (Lemma \ref{l-P3-V7}). 
Then (1) holds  by Proposition \ref{p-FCB-centre}, as $V_7$ has a conic bundle structure over $\P^2$. 

(iii) 
Assume Proposition \ref{p-pic2-structure}(iii). 
We have a blowup $g: Y \to Q$ along a smooth curve $B$. 
Set $E := \Ex(g)$. 
If $\Gamma \cap E =\emptyset$, then (3) holds. 
We may assume that $\Gamma \cap E \neq \emptyset$. 
If $\Gamma$ is not a fibre of the $\P^1$-bundle structure $E \to B$, 
then we would get a contradiction (Lemma \ref{l-line-meeting}). 
Hence $\Gamma$ is a fibre of $E \to B$. 
In this case, 
we have a blowup $X \to \Bl_P\,Q$ along a smooth curve (Lemma \ref{l-P3-V7}), 
where $\Bl_P\,Q$ is the blowup of $Q$ at $P:=g(\Gamma)$, 
i.e., $\Bl_P Q$ is a Fano threefold of No.\ 2-30. 
Since  we have $\Bl_P\,Q \simeq \Bl_C\,\P^3$ for a smooth conic $C$ on $\P^3$, 
we may apply the case when (ii) holds. 

(iv) 
Assume Proposition \ref{p-pic2-structure}(iv). 
We have a blowup $g: Y \to V_d$ along an elliptic curve $B$. 
Set $E := \Ex(g)$. 
As above, either $\Gamma$ is a one-dimensional fibre of $g$ or $E \cap \Gamma = \emptyset$. 
Suppose that $\Gamma$ is a one-dimensional fibre of $g$. 
By Lemma \ref{l-P3-V7}, we obtain $X \simeq \Bl_{\Gamma_U}\,U$, 
where $U := \Bl_P\,V_d$ for $P :=g(\Gamma)$ and $\Gamma_U$ denotes the proper transform of $\Gamma$ on $U$. 
By $\Gamma_U \not\simeq \P^1$ and Lemma \ref{l-nonFano-blowdown}, $U$ is Fano. 
However, this is a contradiction, because $U = \Bl_P\,V_d$ can not be Fano 
by the classification list (Subsection \ref{ss-table-pic2}). 
Thus we get $E \cap \Gamma = \emptyset$. 
Hence $X \simeq \Bl_{B_1 \amalg B_2}\,V_d$ for a disjoint union of smooth curves $B_1$ and $B_2$ on $V_d$. 
In this case, each of  $\Bl_{B_1}\,V_d$ and $\Bl_{B_2}\,V_d$ is Fano (Corollary \ref{c-disjoint-blowup}). 
We may assume that both Fano threefolds $\Bl_{B_1}\,V_d$ and $\Bl_{B_2}\,V_d$ satisfy Proposition \ref{p-pic2-structure}(iv), as otherwise one of (1)--(3) holds. 
Then  
Proposition \ref{p-V-2dP} is applicable, which is absurd. 
\qedhere




\end{proof}

\begin{lem}\label{l-P3-V7}
Let $V$ be a smooth projective threefold. 
Fix a smooth curve $\Gamma$ on $V$ and a closed point $P \in \Gamma$. 
Then there exists an isomorphism 
$\Bl_F(\Bl_{\Gamma}V) \xrightarrow{\iota, \simeq} \Bl_{\Gamma'}(\Bl_P V)$ which satisfies the following commutative diagram 
\begin{equation}\label{e1-P3-V7}
\begin{tikzcd}
    & \Bl_F(\Bl_{\Gamma}\,V)\xrightarrow{\iota, \simeq} \Bl_{\Gamma'}(\Bl_P\,V)
    \arrow[ld, "\sigma"'] \arrow[rd, "\sigma'"]\\
    \Bl_{\Gamma}\,V \arrow[rd, "\tau"'] & & \Bl_P\,V \arrow[ld, "\tau'"]\\
    & V
\end{tikzcd}
\end{equation}
where each arrow except for $\iota$ is the induced blowup, 
$F$ is the fibre of $\Bl_{\Gamma}\,V \to V$ over $P$, and 
$\Gamma'$ denotes the proper transform of $\Gamma$ on $\Bl_P\,V$. 
\end{lem}

\begin{proof}
Set $X := \Bl_F(\Bl_{\Gamma}\,V)$. 
Then the scheme-theoretic inverse image of $P$ on $X$ is defined by an invertible ideal sheaf, because 
the scheme-theoretic inverse image of $P$ on $\Bl_{\Gamma}\,V$ is equal to $F$. 
By the universal property of blowups \cite[Ch. II, Proposition 7.14]{Har77}, we get a factorisation 
\[
\tau \circ \sigma: X \xrightarrow{f}  Y:= \Bl_P\,V \xrightarrow{\tau'} V. 
\]

It suffices to show that $f$ coincides with the blowup along $\Gamma'$. 
Note that $f : X \to Y$ is a birational morphism of smooth projective threefolds such that 
$\rho(X) = \rho(Y)+1$. 
In particular, $\rho(X/Y)=1$ (cf.\ Remark \ref{r-rel-pic}). 
For a curve $C$ on $X$ such that $(\tau \circ \sigma)(C)$ is a point satisfying  $(\tau \circ \sigma)(C) \in \Gamma \setminus P$, 
we get $K_X \cdot C = K_{\Bl_{\Gamma}\,V} \cdot \sigma(C) = -1$. 
Therefore, $f$ is a contraction of a $K_X$-negative extremal ray of $\overline{\NE}(X)$. 
Note that $Y$ is smooth and $f: X \to Y$ coincides with the blowup along $\Gamma'$ outside $\Ex(\tau')$. 
By the classification of extremal ray, 
$f$ coincides with the blowup along $\Gamma'$. 
\end{proof}

\begin{lem}\label{l-3-18}
Let $B_1$ and $B_2$ be a line and a conic on $\P^3$ such that $B_1 \cap B_2 = \emptyset$. 
For each $i \in \{1, 2\}$, let $g_i : Y_i :=\Bl_{B_i}\,\P^3 \to \P^3$ be the blowup along $B_i$ and let $\varphi : X :=\Bl_{B_1 \amalg B_2}\,\P^3 \to \P^3$ be the blowup along $B_1 \amalg B_2$. 
Assume that $X$ is a Fano threefold. 
Then the following hold. 
\begin{enumerate}
\item $(-K_X)^3 = 36$ ($X$ is 3-18). 
    \item $X$ has exactly three extremal rays. 
    \item 
    For each $i \in \{1, 2\}$, let $f_i :  X \to Y_i$ be the induced blowup. 
    Let $f_3 : X \to Y_3$ be the contraction of the extremal ray not corresponding to $f_1$ nor $f_2$. Then $f_3$ is of type $E_1$, $\Ex(f_3) \simeq \F_1$, and 
    $Y_3$ is a Fano threefold of No.\ 2-29. 
\end{enumerate}
\end{lem}

\begin{proof}
Recall that  $Y_1$ and $Y_2$ are Fano threefolds of No.\ 2-33 and No.\ 2-30, respectively (Subsection \ref{ss-table-pic2}). 
Then we have the following commutative diagram: 
\[
\begin{tikzcd}
& & X \arrow[ld, "f_1"'] \arrow[rd, "f_2"] \arrow[dd, "\varphi"] 
 \arrow[lldd, "\varphi_1"', bend right=40] \arrow[rrdd, "\varphi_2", bend left=40]\\
& Y_1  \arrow[ld, "h_1"'] \arrow[rd, "g_1"]& & Y_2 \arrow[ld, "g_2"'] \arrow[rd, "h_2"]\\
Z_1 :=\P^1 & & \P^3 & & Z_2 :=Q
\end{tikzcd}
\]
The assertion (1) follows from $(-K_X)^3 = (-K_{\P^3})^3 - 10 - 18=36$ 
(Lemma \ref{l-blowup-formula}(2)).

Let us show (2). 
It suffices to find a curve $C$ on $X$ such that $h_1(f_1(C))$ and $h_2(f_2(C))$ are points (Lemma \ref{l-exactly3}(3)). 
Recall that $h_2 : Y_2 \to Q$ is a blowup at a point $P$. 
Since $f_2^{-1}(h_2^{-1}(Q))$ is two-dimensional, we can find a curve $C$ on $X$ such that 
$C \subset f_2^{-1}(h_2^{-1}(Q))$ and $h_1(f_1(C))$ is a point. 
Thus (2) holds.

\[
\begin{tikzcd}
& X \arrow[ld, "f_2"'] \arrow[rd, "f_3"] \arrow[dd, "\varphi_2"]\\
Y_2\arrow[rd, "h_2"]& & Y_3 \arrow[ld]\\
& Z_2 =Q. 
\end{tikzcd}
\]

Let us show (3). By the factorisation $\varphi_2 : X \xrightarrow{f_3} Y_3 \to  Q$, 
the contraction $f_3:X \to Y_3$ is of type $E$. 
It is easy to see that $\Ex(h_2)$ is 
the proper transform $(g_2)_*^{-1}\langle B_2 \rangle$ of 
the plane $\langle B_2 \rangle$ containing the conic $B_2$. 
In particular, $\Ex(h_2) \cdot g_2^{-1}(B_1)  = \langle B_2 \rangle \cdot B_1=1$. 
Thus $f_2^{-1}(\Ex(h_2)) \simeq \F_1$. 
Then the image  $\varphi_2(\Ex(f_2))$ of $\Ex(f_2)$ to $Z_2 = Q$ is a curve. 
Hence we obtain $\Ex(f_2) \neq \Ex(f_3)$ (as otherwise, $\Ex(f_2) \simeq \P^1 \times \P^1$ and $\varphi_2(\Ex(f_2))$ would be a point). 
Therefore, $\Ex(f_3) = f_2^{-1}(\Ex(h_2)) \simeq \F_1$, and hence $f_3$ is of type $E_1$ and  $Y_3$ is Fano (Lemma \ref{l-nonFano-blowdown}). 
By (2),  we have  a contraction $Y_3 \to Z_1 =\P^1$ and a birational morphism 
$Y_3 \to Z_2 =Q$, 
In particular, the extremal rays of $Y_3$ are type $D$ and $E_1$, 
and the type-E contraction is a birational morphism to $Q$. 
By the classification table Subsection \ref{ss-table-pic2} 
(or the table immediately after the proof of Proposition \ref{p-V-2dP}), 
$Y_3$ is either 2-7 or 2-29. 
By $(-K_X)^3 = 36 < (-K_{Y^3})^3$ (Lemma \ref{l-blowup-formula2}),  
the case 2-7 is impossible. 
Hence $Y_3$ is 2-29. Thus (3) holds. 
\end{proof}


We are ready to prove the main result of this subsection.

\begin{thm}\label{t-pic3-structure}
Let $X$ be a Fano threefold with $\rho(X)=3$. 
Then one of the following holds. 
\begin{enumerate}
\item[(I)] $X$ has a conic bundle structure over $\P^2$. 
\item[(II)] $X$ has a conic bundle structure over $\F_1$. 
\item[(III)] $X$ is primitive. In particular, $X$ has a conic bundle structure over $\P^1 \times \P^1$. 
\item[(IV)] 
There exist a line $L$ and a conic $C$ on $\P^3$ such that $L \cap C = \emptyset$ 
and $X \simeq \Bl_{L \amalg C} \P^3$ ($X$ is 3-18). 
\item[(V)] $X$ satisfies the assumption and the conclusion of Proposition \ref{p-V-2dP}  
($X$ is one of 3-6, 3-10, 3-25). 
In particular, $X$ has a conic bundle structure over $\P^1 \times \P^1$. 
\end{enumerate}
\end{thm}


\begin{proof}
By Lemma \ref{l-pic3-structure}, 
we may assume that $X$ is  a blowup of $V \in \{\P^3, Q\}$ along 
a disjoint union $C_1 \amalg C_2$ of smooth curves $C_1$ and $C_2$ on $V$. 
Let $g_i : Y_i \to V$ be the blowup along $C_i$. 
Since each $Y_i$ is a Fano threefold with $\rho(Y_i)=2$ (Corollary \ref{c-disjoint-blowup}), 
we have the extremal ray of $\NE(Y_i)$ not corresponding to $g_i$. 
Let $h_i : Y_i \to Z_i$ be its contraction. 
To summarise, we get the following commutative diagram. 
\[
\begin{tikzcd}
& & X \arrow[ld, "f_1"'] \arrow[rd, "f_2"] \arrow[dd, "f"]\\
& Y_1  \arrow[ld, "h_1"'] \arrow[rd, "g_1"]& & Y_2 \arrow[ld, "g_2"'] \arrow[rd, "h_2"]\\
Z_1 & & V & & Z_2
\end{tikzcd}
\]
We may assume that none of $h_1$ and $h_2$ is of type C, as otherwise (I) or (II) holds. 
If both $h_1$ and $h_2$ are of type D, then (V) holds 
(Proposition \ref{p-V-2dP}). 
After permuting $h_1$ and $h_2$ if necessary, the problem is reduced to the case when  $h_1$ is of type $E$. 
Then (IV) holds by 
Lemma \ref{l-P3-disjoint-E}. 
\qedhere

\end{proof}

\subsection{Fano conic bundles over $\P^2$ ($\rho=3$)}\label{ss-pic3-P2}



The purpose of this subsection is to classify Fano conic bundles $X \to \P^2$  with $\rho(X)=3$ (Theorem \ref{t-ele-tr-P2}). 
In this case, we have an elementary transform as in Notation \ref{n-ele-tr-P2}. 
Since a Fano $\P^1$-bundle $Y \to \P^2$ is classified (Remark \ref{r-ele-tr-P2}), we will apply case study depending on No.\ of $Y$. 

\begin{lem}\label{l-FCB-P^2-1}
Let $f: X \to S:=\P^2$ be a Fano conic bundle. 
Then the following hold. 
\begin{enumerate}
    \item $\rho(X)=2$ or $\rho(X)=3$. 
    \item If $\rho(X)=3$, then there exist a Fano $\P^1$-bundle $g: Y \to S$ with $\rho(Y)=2$ and 
    a blowup $\sigma: X \to Y$ along a regular subsection $B_Y$ of $g$: 
    \[
    f: X \xrightarrow{\sigma} Y \xrightarrow{g} S= \P^2. 
    \]
\end{enumerate}
\end{lem}

\begin{proof}
Since the intersection of two curves on $\P^2$ is not empty, 
the assertion follows from Proposition \ref{p-ele-tf}, 
 Proposition \ref{p-ele-tf-1Fano}, Proposition \ref{p-FCB-centre}, and Lemma \ref{l-FCB-pic-irre}. 
\end{proof}


\begin{nota}\label{n-ele-tr-P2}
Let $g: Y \to S:=\P^2$ be a Fano $\P^1$-bundle. 
Let $B_Y$ be a regular subsection of $g$ and let $\sigma : X \to Y$ be the blowup along $B_Y$. 
Assume that $X$ is Fano. 
Let $Y'$ be the elementary transfom of $f: X \xrightarrow{\sigma} Y \xrightarrow{g} S$. 
Set $B:= g(B_Y)$ and $B_{Y'} := \sigma'(\Ex(\sigma'))$, which implies $B_Y \simeq B \simeq B_{Y'}$. Set $d := \deg B$. 
Note that $Y$ has exactly two extremal rays. 
Let $h: Y \to Z$ be the contraction of the extremral ray $R_h$ not corresponding to $g$. 
\begin{itemize}
    \item Let $\ell_h$ be an extremal rational curve on $Y$ with $R_h = \R_{\geq 0}[\ell_h]$. 
    \item Set $\mu_h := -K_Y \cdot \ell_h$, which is the length of $R_h$. 
    \item Let $H_Z$ an ample generator of the Picard group $\Pic\,Z \simeq \Z$. 
    \item It holds that $-K_Y \sim 2h^*H_Z +  g^*\MO_{\P^2}(\mu_h)$ (Proposition \ref{p-pic2-Pic}(2)). 
\end{itemize}
Similarly, if also $Y'$ is Fano, then 
we have the contraction $h': Y' \to Z'$ whose extremal ray does not correspond to $g'$. 
In this case, we define $\ell_{h'}, \mu_{h'}, H_{Z'}$ in a similar way. 
\[
\begin{tikzcd}
& & X \arrow[ld, "\sigma"'] \arrow[rd, "\sigma'"] \arrow[dd, "f"] 
 \arrow[lldd, "\psi"', bend right=40] \arrow[rrdd, "\psi'", bend left=40]\\
& Y  \arrow[ld, "h"'] \arrow[rd, "g"]& & Y' \arrow[ld, "g'"'] \arrow[rd, "h'"]\\
Z & & S= \P^2 & & Z'
\end{tikzcd}
\]


\end{nota}

\begin{rem}\label{r-ele-tr-P2} 
We use Notation \ref{n-ele-tr-P2}. 
Since $g$ is of type $C_2$, 
the possibilities for $Y$ are as in Table \ref{table-P1bdl/P2} (Subsection \ref{ss-table-pic2}). 
In particular, the No.\ of $Y$ is determined only by $(-K_Y)^3$. 
If $Y'$ is Fano, then also the possibilities for $Y'$ are as in Table \ref{table-P1bdl/P2}. 

  \begin{center}
\begin{longtable}{ccp{12cm}c}
No. & $(-K_Y)^3$ & descriptions and extremal rays  & \\ \hline
2-24 & $30$ & $Y$ is a  divisor on $\mathbb{P}^2\times \mathbb{P}^2$ of bidegree $(1,2)$ & \\
 &  & $C_1: \deg \Delta = 3$ & \\ 
&  & $C_2$ & \\ \hline
2-27 & $38$ &$C_2$ & \\ 
 &  & $E_1:$ blowup of $\P^3$ along a cubic rational curve & \\ \hline
2-31 & $46$ & $C_2$ & \\ 
&  & $E_1:$ blowup of $Q$ along a line & \\ \hline
2-32 & $48$ & $Y$ is  a divisor $W$ on $\P^2 \times \P^2$ of bidegree $(1, 1)$ &\\
 &  & $C_2: W \hookrightarrow \P^2 \times \P^2 \xrightarrow{\pr_1} \P^2$ & \\ 
 &  & $C_2: W \hookrightarrow \P^2 \times \P^2 \xrightarrow{\pr_1} \P^2$ & \\ \hline
2-34 & $54$ & $Y=\P^2 \times \P^1$ \\
 & &  $C_2:$ the projection $\P^2 \times \P^1 \to \P^2$ & \\ 
 &  & $D_3:$ the projection $\P^2 \times \P^1 \to \P^1$ & \\ \hline
2-35 & $56$ & $Y=V_7 =\P_{\P^2}(\MO_{\P^2} \oplus \MO_{\P^2}(1))$ \\
 & & $C_2:$ the projection $\mathbb{P}_{\P^2} (\MO_{\mathbb{P}^2}\oplus \MO_{\mathbb{P}^2}(1)) \to \P^2$ & \\ 
 &  & $E_2:$ blowup of $\P^3$ at a point &\\ \hline
2-36 & $62$  & 
$Y=\mathbb{P}_{\P^2}(\MO_{\mathbb{P}^2}\oplus \MO_{\mathbb{P}^2}(2))$\\
& & 
$C_2:$ the projection $\mathbb{P}_{\P^2}(\MO_{\mathbb{P}^2}\oplus \MO_{\mathbb{P}^2}(2)) \to \P^2$  & \\ 
&  &  $E_5:$ blowup at the singular point of the cone over the Veronese surface
& \\ \hline
\\
    \caption{Fano $\P^1$-bundles $Y$ over $\P^2$}\label{table-P1bdl/P2}
      \end{longtable}
  \end{center}


\end{rem}


\begin{lem}\label{l-P2-2-24}
We use Notation \ref{n-ele-tr-P2}. 
Assume that $Y$ is of No.\ 2-24. 
Then 
$Y'$ is a Fano threefold of No.\ 2-34, $(-K_X)^3=24$, $\deg B = 2$, $p_a(B)=0$, 
$-K_Y \cdot B_Y = 2$, and $-K_{Y'} \cdot B'=14$ ($X$ is 3-8). 
\end{lem}

\begin{proof}
Note that $g: Y  \to S = \P^2$ is of type $C_2$ and $h : Y \to Z \simeq \P^2$ is of type $C_1$. 
Hence $B_Y$ is a smooth fibre of $h$ (Proposition \ref{p-FCB-centre}). 
In particular, $B \simeq B_Y \simeq \P^1$, $-K_Y \cdot B_Y = 2$, and 
$(-K_X)^3 =(-K_Y)^3 - 6 = 30 -6 =24$ (Lemma \ref{l-blowup-formula}). 
Recall that $-K_Y \sim 2h^*H_Z + g^*\MO_{\P^2}(1)$ and  $B_Y \equiv (h^*H_Z)^2 \equiv 2\ell_h$. 
We then get $\deg B = B \cdot \MO_{\P^2}(1) =  B_Y \cdot g^*\MO_{\P^2}(1)  
= 2\ell_{h}\cdot g^*\MO_{\P^2}(1) = 2$ (Proposition \ref{p-pic2-Pic}). 
We obtain the following (Proposition \ref{p-ele-tf-numbers}): 
\[
-K_{Y/S} \cdot B_Y = -K_Y \cdot B_Y + K_S \cdot B= 2 -6 =-4. 
\]
\[
(-K_{Y'})^3 = (-K_Y)^3 -4 (-K_{Y/S}) \cdot B_Y + 2B^2
=30 +16 + 8= 54, 
\]
\[
-K_{Y'} \cdot B_{Y'} = B^2 +2 (- K_S \cdot B) -(-K_{Y} \cdot B_Y) 
= 4 +12 - 2 =14. 
\]
By $-K_{Y/S} \cdot B_Y =-4 \neq 10 = 2(B^2+1)$, 
$Y'$ is a Fano threefold (Proposition \ref{p-ele-tf-1Fano}), which is of No.\ 2-34 (Remark \ref{r-ele-tr-P2}). 
\qedhere

\end{proof}

\begin{lem}\label{l-P2-nonFano}
We use Notation \ref{n-ele-tr-P2}. 
Assume that $Y'$ is not Fano. 
Then $Y$ is of No.\ 2-34 and one of the following holds. 
\begin{enumerate}
\item $(-K_X)^3=38$, $(-K_{Y'})^3 =40$, $\deg B = 1$, $p_a(B)=0$, 
$-K_Y \cdot B_Y = 7$, and $-K_{Y'} \cdot B_{Y'}=0$ ($X$ is  3-21). 
\item $(-K_X)^3=20$, $(-K_{Y'})^3 =22$, $\deg B = 2$, $p_a(B)=0$, 
$-K_Y \cdot B_Y = 16$, and $-K_{Y'} \cdot B_{Y'}=0$ ($X$ is   3-5). 
\end{enumerate}
\end{lem}

\begin{proof}
Proposition \ref{p-ele-tf-1Fano} implies $B \simeq \P^1, -K_{Y'} \cdot B_{Y'} =0, g^{-1}(B)=\P^1 \times \P^1$, and $-K_{Y/S}\cdot B_Y =2 (B^2 +1 )$. 
Since $B$ is a smooth rational curve on $\P^2$, $B$ is either a line or a conic on $S = \P^2$. 
If $B$ is a line (resp. a conic), then 
$-K_Y\cdot B_Y =-K_S \cdot B + 2 (B^2 +1 )=7$ (resp. $=16$). 
In particular, $B_Y$ is disjoint from any curve $L$ on $Y$ satisfying $-K_Y \cdot L =1$ (Lemma \ref{l-line-meeting}).

\begin{claim*}
$Y$ is of No.\ 2-34.
\end{claim*}

\begin{proof}[Proof of Claim]
Suppose that $Y$ is of No.\ 2-35 or 2-36. 
Then we have $Y = \P_{\P^2}(\MO_{\P^2} \oplus \MO_{\P^2}(n))$ with $n \in \{1, 2\}$ and $g: Y \to S =\P^2$ is the projection. 
Since $B$ is either a line or a conic on $S = \P^2$, 
we get $(\MO_{\P^2} \oplus \MO_{\P^2}(n))|_B \simeq \MO_{\P^1} \oplus \MO_{\P^1}(m)$ for some $m \in \{ 1, 2, 4\}$, i.e., $g^{-1}(B) \simeq \P_{\P^1}(\MO_{\P^1} \oplus \MO_{\P^1}(m))$. 
In any case, we get $g^{-1}(B) \not\simeq \P^1 \times \P^1$, which is absurd. 
Thus $Y$ is not of No.\ 2-35 nor 2-36. 

Suppose that $Y$ is of type 2-32. 
Recall that $Y$ is a prime divisor on $\P^2 \times \P^2$ of bidegree $(1, 1)$ and 
$g: Y \to S=\P^2$ is the composition 
$Y \hookrightarrow \P^2 \times \P^2 \xrightarrow{ {\rm pr}_1} \P^2$. 
Set $D := g^{-1}(B)$. 
By $-K_Y \simeq \MO_{\P^2 \times \P^2}(2, 2)|_Y$, we have that 
$\MO_D(-K_D) \sim \MO_Y(-K_Y-D)|_D \sim 
\MO_{\P^2 \times \P^2}(2-d, 2)|_D$. 
By $D \simeq \P^1 \times \P^1$, $-K_D$ is ample, and hence $d =1$, i.e., $B$ is a line. 
By the Euler sequence $0 \to \MO_{\P^2} \to \MO_{\P^2}(1)^{\oplus 3} \to T_{\P^2} \to 0$, 
we get an exact sequence 
$0 \to \MO_{\P^1} \to \MO_{\P^1}(1)^{\oplus 3} \to T_{\P^2}|_B \to 0$, which implies $\deg (T_{\P^2}|_B)=3 \not\in 2\Z$. 
By $Y \simeq \P_{\P^2}(T_{\P^2})$, 
we obtain $D \simeq \P_B(T_{\P^2}|_B)$. 
Then  $D \not\simeq \P^1 \times \P^1$, which is absurd. 
Thus $Y$ is not of No.\ 2-32.

Suppose that $Y$ is of No.\ 2-31. 
Recall that $h : Y \to Z = Q$ is of type $E_1$. 
As $B_Y$ is disjoint from any curve $L$ satisfying $-K_Y \cdot L=1$, $B_Y$ is disjoint from $\Ex(h)$. 
Then $-K_Y \cdot B_Y = -h^*K_{Q} \cdot B_Y \in 3\Z$. 
By $-K_Y \cdot B_Y \in \{7, 16\}$, this is a contradiction.

Suppose that $Y$ is of No.\ 2-27. 
Recall that $h : Y \to Z=\P^3$ is a blowup along  a smooth cubic rational curve $\Gamma$. 
It follows from $\Ex(h) \cap B_Y = \emptyset$ that 
$-K_Y \cdot B_Y = -h^*K_{\P^3} \cdot B_Y \in 4\Z$. 
By $-K_Y \cdot B_Y \in \{7, 16\}$, we obtain $-K_Y \cdot B_Y = 16$, i.e., 
$B_{\P^3} := h(B_Y)$ is a smooth rational curve of degree $4$. 
Hence $\psi : X \to Z = \P^3$ is the blowup along the disjoint union $\Gamma \amalg B_{\P^3}$ 
of a cubic rational curve $\Gamma$ and a quartic rational curve $B_{\P^3}$, 
which contradicts Lemma \ref{l-P3-disjoint-E} 
(this is applicable because $\Bl_{B_{\P^3}}\P^3$ is a Fano threefold of No.\ 2-22 (Subsection \ref{ss-table-pic2})).

Finally, $Y$ is not of No.\ 2-24 by Lemma \ref{l-P2-2-24}. This completes the proof of Claim. 
\end{proof}

Since $Y$ is of No.\ 2-34, we have $Y=\P^2 \times \P^1$. 
Then the following hold  (Proposition \ref{p-ele-tf-numbers}, 
Proposition \ref{p-ele-tf-1Fano}, Lemma \ref{l-blowup-formula}): 
\[
(-K_{Y'})^3 = (-K_Y)^3 -4 (-K_{Y/S}) \cdot B_Y + 2B^2
= 54 -8 (B^2+1) +2B^2 =  46 - 6d^2, 
\]
\[
(-K_X)^3 = (-K_Y)^3 - 2  (-K_Y) \cdot B_Y -2 = 52 -2 (-K_Y) \cdot B_Y. 
\]

(1) Assume $\deg B =1$, i.e., $d=1$. 
Then $-K_Y \cdot B_Y = 7, (-K_{Y'})^3 = 40, (-K_X)^3 = 52-14 = 38$. 

(2) Assume $\deg B =2$, i.e., $d=2$. 
Then $-K_Y \cdot B_Y = 16, (-K_{Y'})^3 = 22, (-K_X)^3 = 52-32 = 20$.  
\end{proof}

\begin{lem}\label{l-P2-2-36}
We use Notation \ref{n-ele-tr-P2}. 
Assume that $Y$ is of No.\ 2-36. 
Set $D:=\Ex(h)$.
Then $D \cap B_Y=\emptyset$ and one of the following holds. 
\begin{enumerate}
\item 
$Y'$ is a Fano threefold of No.\ 2-35, $(-K_X)^3=50$, $\deg B = 1$, $p_a(B)=0$, 
$-K_Y \cdot B_Y = 5$, and $-K_{Y'} \cdot B'=2$ ($X$ is 3-29). 
\item 
$Y'$ is a Fano threefold of No.\ 2-34, $(-K_X)^3=40$, $\deg B = 2$, $p_a(B)=0$, 
$-K_Y \cdot B_Y = 10$, and $-K_{Y'} \cdot B'=6$ ($X$ is 3-22). 
\item 
$Y'$ is a Fano threefold of No.\ 2-35, $(-K_X)^3=32$, $\deg B = 3$, $p_a(B)=1$, 
$-K_Y \cdot B_Y = 15$, and $-K_{Y'} \cdot B'=12$ ($X$ is 3-14). 
\item 
$Y'$ is a Fano threefold of No.\ 2-36, $(-K_X)^3=26$, $\deg B = 4$, $p_a(B)=3$, 
$-K_Y \cdot B_Y = 20$, and $-K_{Y'} \cdot B'=20$  ($X$ is 3-9). 
\end{enumerate}
\end{lem}


\begin{proof}
Since $h$ is of type $E_5$, we have $-K_Y \cdot L=1$ for any line $L$ on $D=\P^2$ (cf.\ \cite[Proposition 3.22]{TanII}). 
Therefore, we get $D \cap B_Y = \emptyset$ (Lemma \ref{l-line-meeting}). 

We now show that  $-K_{Y/S} \sim 2D + g^*\MO_{\P^2}(2)$. 
Since $D$ is a section of $g: Y \to \P^2$, 
we can write $-K_{Y/S} \sim 2D + g^*\MO_{\P^2}(n)$ for some $n \in \Z$. 
As $D \hookrightarrow Y \xrightarrow{g} S$ is an isomorphism, we get 
\[
0 \sim -K_D +K_D \sim -(K_Y+D)|_D +g^*K_S|_D \sim (D + g^*\MO_{\P^2}(n))|_D \sim \MO_{\P^2}(-2+n). 
\]
Thus $n=2$, which completes the proof of  $-K_{Y/S} \sim 2D + g^*\MO_{\P^2}(2)$.

Then $-K_{Y/S} \cdot B_Y = 2D \cdot B_Y + g^*\MO_{\P^2}(2) \cdot B_Y = 0 + 2 \deg B = 2d$. 
It follows from Proposition \ref{p-ele-tf-numbers} that 
\[
(-K_{Y'})^3 =(-K_Y)^3 -4 (-K_{Y/S} \cdot B_Y) +2B^2 = 62 - 8d +2d^2. 
\]
By Lemma \ref{l-P2-nonFano}, $Y'$ is Fano. 
We then have $(-K_{Y'})^3 \in \{ 30, 38, 46, 48, 54, 56, 62\}$ 
(Remark \ref{r-ele-tr-P2}), and hence $(-K_{Y'})^3 \leq 62$. 
Since $d \geq 5$ would imply 
$(-K_{Y'})^3 =(-K_Y)^3 -4 (-K_{Y/S} \cdot B_Y) +2B^2 = 62 - 8d +2d^2 >62$, 
we get $d \leq 4$. 
Moreover, the following hold 
(Lemma \ref{l-blowup-formula}, Proposition \ref{p-ele-tf-numbers}):  
\[
-K_{Y} \cdot B_Y = -K_{Y/S} \cdot B_Y -K_S \cdot B = 2d+3d =5d,  
\]
\[
(-K_X)^3 = (-K_Y)^3 - 2 (-K_Y) \cdot B_Y +2p_a(B) -2 
= 60 -2 (-K_Y) \cdot B_Y +2p_a(B), 
\]
\[
-K_{Y'} \cdot B_{Y'} = B^2 +2 (- K_S \cdot B) -(-K_{Y} \cdot B_Y) 
= d^2 +6d -(-K_{Y} \cdot B_Y). 
\]

(1) Assume $d =1$. Then $(-K_{Y'})^3= 62 -8+2  =56$, $p_a(B)=0$, $-K_{Y} \cdot B_Y =5$, $(-K_X)^3 = 50$, 
$-K_{Y'} \cdot B_{Y'} =1+6 - 5=2$.

(2) Assume $d =2$. Then $(-K_{Y'})^3= 62 -16+8  =54$, $p_a(B)=0$, $-K_{Y} \cdot B_Y =10$, $(-K_X)^3 = 40$, 
$-K_{Y'} \cdot B_{Y'} =4+12 - 10=6$. 

(3) Assume $d =3$. Then $(-K_{Y'})^3= 62 -24+18  =56$, $p_a(B)=1$, $-K_{Y} \cdot B_Y =15$, $(-K_X)^3 = 32$, 
$-K_{Y'} \cdot B_{Y'} =9 +18 - 15=12$.

(4) Assume $d =4$. Then $(-K_{Y'})^3= 62 -32+32  =62$,  $p_a(B)=3$, $-K_{Y} \cdot B_Y =20$, $(-K_X)^3 = 26$, 
$-K_{Y'} \cdot B_{Y'} =16 +24 - 20=20$. 
\end{proof}

\begin{lem}\label{l-P2-2-27}
We use Notation \ref{n-ele-tr-P2}. 
Assume that $Y$ is of No.\ 2-27. 
Set $D := \Ex(h)$. 
Then one of the following holds. 
\begin{enumerate}
\item 
$B_Y$ is a fibre of the $\P^1$-bundle $D \to h(D)$,  $D \cdot B_Y = -1$, 
$Y'$ is a Fano threefold of No.\ 2-32, $(-K_X)^3=34$, $\deg B = 1$, $p_a(B)=0$, 
$-K_Y \cdot B_Y = 1$, $-K_{Y'} \cdot B'=6$ ($X$ is 3-16). 
\item  
$D \cap B_Y = \emptyset$, 
$Y'$ is a Fano threefold of No.\ 2-34, $(-K_X)^3=28$, $\deg B = 2$, $p_a(B)=0$, 
$-K_Y \cdot B_Y = 4$, $-K_{Y'} \cdot B'=12$ ($X$ is 3-12). 
\end{enumerate}
\end{lem}

\begin{proof}
Recall that $g : Y \to S = \P^2$ is of type $C_2$ and 
$h : Y \to Z = \P^3$ is the blowup along a smooth cubic rational curve. 
We have $-K_Y \sim h^*\MO_{\P^3}(2) + g^*\MO_{\P^2}(1)$ (Proposition \ref{p-pic2-Pic}(2)). 
We have 
\[
K_Y = h^* K_{\P^3} +D = h^*\MO_{\P^3}(-4)  +D. 
\]
By 
$-K_Y \sim h^*\MO_{\P^3}(2) + g^*\MO_{\P^2}(1) \equiv g^*\MO_{\P^2}(1) + \frac{-K_Y+D}{2}$, we get 
\[
-K_Y \sim D +g^*\MO_{\P^2}(2), \qquad -K_{Y/S} \sim D +g^*\MO_{\P^2}(-1). 
\]

(1) Assume $B_Y \subset D$. 
Then $B_Y$ is a fibre of the $\P^1$-bundle $D \to h(D)$ (Lemma \ref{l-line-meeting}). 
In this case, $K_Y \cdot B_Y =D \cdot B_Y = -1$. 
Hence $\deg B = \MO_{\P^2}(1) \cdot B = g^*\MO_{\P^2}(1) \cdot B_Y =  \frac{1}{2}(-K_Y \cdot B_Y -D \cdot B_Y)=1$, $p_a(B)=0$, 
$(-K_X)^3 =(-K_Y)^3 - 4 = 38-4=34$, 
$(-K_{Y/S}) \cdot B_Y = -1 -\deg B = -2$, 
$(-K_{Y'})^3 =(-K_Y)^3 -4 (-K_{Y/S} \cdot B_Y) +2B^2 = 38 -4 \cdot (-2) + 2 = 48$ (Proposition \ref{p-ele-tf-numbers}). 
In particular, $Y'$ is a Fano threefold of No.\ 2-32 
(Remark \ref{r-ele-tr-P2}, Lemma \ref{l-P2-nonFano}). 
It follows from Proposition \ref{p-ele-tf-numbers} that 
$-K_{Y'} \cdot B_{Y'} = B^2 +2 (- K_S \cdot B) -(-K_{Y} \cdot B_Y) 
= 1 + 6 -1 =6.$

(2) Assume $B_Y \not\subset D$. 
Since $D$ is covered by curves $L$ satisfying $-K_Y \cdot L=1$, 
we obtain $B_Y \cap D= \emptyset$  (Lemma \ref{l-line-meeting}). 
Then $-K_{Y/S} \cdot B_Y =( D+ g^*\MO_{\P^2}(-1)) \cdot B_Y = -d$, 
$-K_Y \cdot B_Y = (D+ g^*\MO_{\P^2}(2)) = 2d$, and 
\[
(-K_{Y'})^3 = (-K_Y)^3 -4 (-K_{Y/S}) \cdot B_Y +2B^2 = 38
+4d +2d^2. 
\]
For $d=1,2,3, ...$, we get $(-K_{Y'})^3 = 44, 54, 68, ...$ 
By 
$(-K_{Y'})^3 \in \{ 30, 38, 46, 48, 54, 56, 62 \}$ (Remark \ref{r-ele-tr-P2}, Lemma \ref{l-P2-nonFano}), 
we obtain $(d, (-K_{Y'})^3) =(2, 54)$, and hence $Y'$ is a Fano threefold of No.\ 2-34. 
We have $p_a(B)=0$, $(-K_X)^3 =(-K_Y)^3 - 10 = 38-10=28$,  and 
$-K_{Y'} \cdot B_{Y'} = B^2 +2 (- K_S \cdot B) -(-K_{Y} \cdot B_Y) 
= 4 +12 - 4 =12$ (Proposition \ref{p-ele-tf-numbers}).  
\end{proof}

\begin{lem}\label{l-P2-2-31}
We use Notation \ref{n-ele-tr-P2}. 
Assume that $Y$ is of No.\ 2-31. 
Set $D := \Ex(h)$. 
Then one of the following holds. 
\begin{enumerate}
\item 
$B_Y$ is a fibre of the $\P^1$-bundle $D \to h(D)$, $D \cdot B_Y  =-1$, 
$Y'$ is a Fano threefold of No.\ 2-35, $(-K_X)^3=42$, $\deg B = 1$, $p_a(B)=0$, 
$-K_Y \cdot B_Y = 1$, and $-K_{Y'} \cdot B_{Y'}=6$ ($X$ is  3-23). 
\item 
$D \cap B_Y = \emptyset$, 
$Y'$ is a Fano threefold of No.\ 2-32, $(-K_X)^3=38$, $\deg B = 1$, $p_a(B)=0$, 
$-K_Y \cdot B_Y = 3$, and $-K_{Y'} \cdot B_{Y'}=4$ ($X$ is  3-20). 
\item 
$D \cap B_Y = \emptyset$, 
$Y'$ is a Fano threefold of No.\ 2-34, $(-K_X)^3=32$, $\deg B = 2$, $p_a(B)=0$, 
$-K_Y \cdot B_Y = 6$, and $-K_{Y'} \cdot B_{Y'}=10$ ($X$ is  3-15). 
\end{enumerate}
\end{lem}

\begin{proof}
Recall that $g : Y \to S=\P^2$ is of type $C_2$ and 
$h : Y \to Z= Q$ is the blowup along a line. 
Then $-K_Y \sim g^*\MO_{\P^2}(1) +h^*\MO_Q(2)$ (Proposition \ref{p-pic2-Pic}(2)).  
We have 
\[
K_Y = h^* K_{Q} +D = h^*\MO_Q(-3) +D. 
\]
By 
$-K_Y \sim g^*\MO_{\P^2}(1) +h^*\MO_Q(2) \equiv g^*\MO_{\P^2}(1) + \frac{2}{3}(-K_Y+D)$, we get 
\[
-K_Y \sim 2D +g^*\MO_{\P^2}(3), \qquad -K_{Y/S} \sim 2D. 
\]

(1) Assume $B_Y \subset D$. 
Then $B_Y$ is a fibre of the $\P^1$-bundle $D \to h(D)$ (Lemma \ref{l-line-meeting}). 
In this case, $K_Y \cdot B_Y =D \cdot B_Y = -1$. 
Hence $\deg B = \MO_{\P^2}(1) \cdot B = g^*\MO_{\P^2}(1) \cdot B_Y = 
\frac{1}{3}(-K_Y \cdot B_Y -2D \cdot B_Y)=1$. 
Thus $p_a(B)=0$, 
$(-K_X)^3 =(-K_Y)^3 - 4 = 46-4=42$, 
$(-K_{Y/S}) \cdot B_Y = -2$, 
$(-K_{Y'})^3 =(-K_Y)^3 -4 (-K_{Y/S} \cdot B_Y) +2B^2 = 46 -4 \cdot (-2) + 2 = 56$ (Proposition \ref{p-ele-tf-numbers}). 
In particular, $Y'$ is a Fano threefold of No.\ 2-35 (Remark \ref{r-ele-tr-P2}, Lemma \ref{l-P2-nonFano}). 
It follows from Proposition \ref{p-ele-tf-numbers} that  $-K_{Y'} \cdot B_{Y'} = B^2 +2 (- K_S \cdot B) -(-K_{Y} \cdot B_Y) 
=1 +6-1=6.$  

(2), (3) Assume $B_Y \not\subset D$. 
Since $D$ is covered by curves $L$ satisfying $-K_Y \cdot L=1$, 
we obtain $B_Y \cap D= \emptyset$  (Lemma \ref{l-line-meeting}). 
Then $-K_{Y/S} \cdot B_Y = 0$, $-K_Y \cdot B_Y = 3d$, and 
\[
(-K_{Y'})^3 = (-K_Y)^3 -4 (-K_{Y/S}) \cdot B_Y +2B^2 = 46 +2d^2. 
\]
For $d=1,2,3, ...$, we get $(-K_{Y'})^3 = 48, 54, 64...$ By 
$(-K_{Y'})^3 \in \{ 30, 38, 46, 48, 54, 46, 62 \}$ (Remark \ref{r-ele-tr-P2}, Lemma \ref{l-P2-nonFano}), 
we obtain $(d, (-K_{Y'})^3) \in \{(1, 48), (2, 54)\}$. 
\begin{enumerate}
\item[(2)] If $(d, (-K_{Y'})^3) \in (1, 48)$, then $p_a(B)=0, 
-K_Y \cdot B_Y = 3, (-K_X)^3 =(-K_Y)^3 -8 =38, 
$ and $-K_{Y'} \cdot B_{Y'} = B^2 +2 (- K_S \cdot B) -(-K_{Y} \cdot B_Y) 
=1 +6-3=4.$  
In this case, $Y'$ is a Fano threefold of No.\ 2-32. 
\item[(3)] If $(d, (-K_{Y'})^3) \in (2, 54)$, then $p_a(B)=0, -K_Y \cdot B_Y = 6, (-K_X)^3 =(-K_Y)^3 -14 =32,$ and 
$-K_{Y'} \cdot B_{Y'} = B^2 +2 (- K_S \cdot B) -(-K_{Y} \cdot B_Y) 
=4 +12-6=10.$ 
In this case, $Y'$ is a Fano threefold of No.\ 2-34. 
\end{enumerate}
\end{proof}

\begin{lem}\label{l-P2-2-35}
We use Notation \ref{n-ele-tr-P2}. Set $D := \Ex(h)$. 
Assume that $Y$ is of No.\ 2-35 and $(-K_{Y'})^3 \in \{ 48, 54, 56\}$. 
Then one of the following holds. 
\begin{enumerate}
\item 
$D \cap B_Y = \emptyset$, 
$Y'$ is a Fano threefold of No.\ 2-34, $(-K_X)^3=46$, $\deg B = 1$, $p_a(B)=0$, 
$-K_Y \cdot B_Y = 4$, and $-K_{Y'} \cdot B_{Y'}=3$ ($X$ is 3-26). 
\item 
$D \cap B_Y = \emptyset$, 
$Y'$ is a Fano threefold of No.\ 2-35, $(-K_X)^3=38$, $\deg B = 2$, $p_a(B)=0$, 
$-K_Y \cdot B_Y = 8$, and $-K_{Y'} \cdot B_{Y'}=8$ ($X$ is 3-19). 
\item 
$B_Y \not\subset D$, $D \cdot B_Y =1$, 
$Y'$ is a Fano threefold of No.\ 2-32, $(-K_X)^3=34$, $\deg B = 2$, $p_a(B)=0$, 
$-K_Y \cdot B_Y = 10$, and $-K_{Y'} \cdot B_{Y'}=6$ ($X$ is 3-16). 
\item 
$B_Y \not\subset D$, $D \cdot B_Y =1$, 
$Y'$ is a Fano threefold of No.\ 2-34, $(-K_X)^3=28$, $\deg B = 3$, $p_a(B)=1$, 
$-K_Y \cdot B_Y = 14$, and $-K_{Y'} \cdot B_{Y'}=13$ ($X$ is 3-11). 
\end{enumerate}
\end{lem}

\begin{proof}
Note that $g : Y = \P_{\P^2}(\MO_{\P^2} \oplus \MO_{\P^2}(1)) \to S = \P^2$ is the projection and $h : Y \to \P^3$ is a blowup at a point of $Z=\P^3$. 
Since $D$ is a section of $g$, we can write  
\[
-K_{Y/S} \sim 2D +g^*\MO_{\P^2}(n)
\]
for some $n \in \Z$. 
Let $L$ be a line on $D = \P^2$. 
Then $D \cdot L =-1$ implies 
\[
3 = -K_D \cdot L = -(K_{Y}+D) \cdot L = D \cdot L +g^*\MO_{\P^2}(n+3) \cdot L = -1 +n+3. 
\]
Hence $n=1$ and $-K_{Y/S} \sim 2D +g^*\MO_{\P^2}(1)$. 
We get 
\[
-K_{Y/S} \cdot B_{Y} = 2D \cdot B_{Y} + d = d+2u, 
\]
\[
-K_Y \cdot B_Y = 4d+2u, \qquad u :=  D \cdot B_Y. 
\]

\begin{claim*}
    The following hold. 
    \begin{enumerate}
    \renewcommand{\labelenumi}{(\roman{enumi})}
    \item $u =-1$ $\Leftrightarrow$ $B_Y$ is a line on $D=\P^2$ $\Leftrightarrow$ $B_Y \subset D$. 
    \item $u=0$ $\Leftrightarrow$ $D \cap B_Y = \emptyset$. 
    \item $u=1$ $\Leftrightarrow$ $D \cap B_Y \neq \emptyset$ and $B_Y \not\subset D$. 
    \end{enumerate}
    In particular, $u =D \cdot B_Y \in \{ -1, 0, 1\}$.  
\end{claim*}

\begin{proof}[Proof of Claim] 
Let us show (i). 
By $D|_D \simeq \MO_{\P^2}(-1)$, 
it holds that $u =-1$ $\Leftrightarrow$ $B_Y$ is a line on $D=\P^2$. 
If $B_Y$ is a line on $D$, then we get $B_Y \subset D$.  
To prove the converse, assume $B_Y \subset D$. 
Then $\dim \MO_{B_Y \cap L} < -K_Y \cdot L =2$ for a general line $L$ on $D=\P^2$ (Lemma \ref{l-line-meeting}), which implies that $B_Y$ is a line, i.e., $D \cdot B_Y =-1$. 
Thus (i) holds.

Let us show (ii). 
The implication $\Leftarrow$ is obvious. 
To prove the converse, assume $u= D \cdot B_Y = 0$. 
By (i), we get $B_Y \not\subset D$. Then $D \cap B_Y = \emptyset$. 
Thus (ii) holds. 

Let us show (iii). 
If $u=1$, then it follows from (i) and (ii) that 
$D \cap B_Y \neq \emptyset$ and $B_Y \not\subset D$. 
Conversely, assume that 
$D \cap B_Y \neq \emptyset$ and $B_Y \not\subset D$. 
If $D \cdot B_Y \geq 2$, then we can find a line $L$ on $D=\P^2$ 
such that $B_Y \cap L$ is zero-dimensional and 
$\dim_k \MO_{B_Y \cap L} \geq 2$, which is absurd (Lemma \ref{l-line-meeting}). 
Hence $D \cdot B_Y \leq 1$. 
By $D \cap B_Y \neq \emptyset$ and $B_Y \not\subset D$, we get $D \cdot B_Y = 1$. 
Thus (iii) holds. 
This completes the proof of Claim. 
\end{proof}

It follows from Proposition \ref{p-ele-tf-numbers} that  
\[
\{ 48, 54, 56\} \ni (-K_{Y'})^3 = (-K_{Y})^3 -4 (-K_{Y/S}) \cdot B_Y +2B^2 
\]
\[
= 56 -4 (d+2u) + 2d^2 = 56 +2 ( d^2-2d -4u). 
\]
If $d \geq 4$, then 
we would get $d^2-2d -4u \geq 16-8-4u>0$, which contradicts 
$(-K_{Y'})^3 \in \{ 48, 54, 56\}$. 
Hence $d \in \{1, 2, 3\}$.

(1) Assume $d=1$. Then 
\[
\{ 48, 54, 56\} \ni (-K_{Y'})^3 = 54  -8u. 
\]
Hence $(d, u, (-K_{Y'})^3) = (1, 0, 54)$. 
Thus $p_a(B_Y) = 0, -K_Y \cdot B_Y = 4,$ 
$(-K_X)^3 = (-K_Y)^3 -2 (-K_Y) \cdot B_Y + 2p_a(B) -2 = 56 -8+0-2 = 46$, 
$-K_{Y'} \cdot B_{Y'} = B^2 +2 (- K_S \cdot B) -(-K_{Y} \cdot B_Y) = 1 + 6 - 4 =3$ (Proposition \ref{p-ele-tf-numbers}). 
Then $Y'$ is a Fano threefold of No.\ 2-34 (Remark \ref{r-ele-tr-P2}, Lemma \ref{l-P2-nonFano}). 

(2), (3) 
Assume $d=2$. Then 
\[
\{ 48, 54, 56\} \ni (-K_{Y'})^3 = 56  -8u. 
\]
Hence $(d, u, (-K_{Y'})^3) = (2, 0, 56), (2, 1, 48)$. 
\begin{enumerate}
\item[(2)] 
Assume $(d, u, (-K_{Y'})^3) = (2, 0, 56)$. 
Then $p_a(B_Y) = 0, -K_Y \cdot B_Y = 8,$ 
$(-K_X)^3 = (-K_Y)^3 -2 (-K_Y) \cdot B_Y + 2p_a(B) -2 = 56 -16+0-2 = 38$, 
$-K_{Y'} \cdot B_{Y'} = B^2 +2 (- K_S \cdot B) -(-K_{Y} \cdot B_Y) = 4 + 12 - 8 =8$ (Proposition \ref{p-ele-tf-numbers}). 
Hence $Y'$ is a Fano threefold  of No.\ 2-35  (Remark \ref{r-ele-tr-P2}, Lemma \ref{l-P2-nonFano}). 
\item[(3)] 
Assume $(d, u, (-K_{Y'})^3) = (2, 1, 48)$. 
Then $p_a(B_Y) = 0, -K_Y \cdot B_Y = 10,$ 
$(-K_X)^3 = (-K_Y)^3 -2 (-K_Y) \cdot B_Y + 2p_a(B) -2 = 56 -20+0-2 = 34$, 
$-K_{Y'} \cdot B_{Y'} = B^2 +2 (- K_S \cdot B) -(-K_{Y} \cdot B_Y) = 4 + 12 - 10=6$ (Proposition \ref{p-ele-tf-numbers}). 
Hence $Y'$ is a Fano threefold  of No.\ 2-32  (Remark \ref{r-ele-tr-P2}, Lemma \ref{l-P2-nonFano}).  
\end{enumerate}

(4) Assume $d=3$. Then 
\[
\{ 48, 54, 56\} \ni (-K_{Y'})^3 = 56 +2 ( d^2-2d -4u) = 62-8u. 
\]
Hence $(d, u, (-K_{Y'})^3) =(3, 1, 54)$. 
Then $p_a(B_Y) = 1, -K_Y \cdot B_Y = 14,$ $(-K_X)^3 = (-K_Y)^3 -2 (-K_Y) \cdot B_Y + 2p_a(B) -2 = 56 -28+2-2 = 28$, 
$-K_{Y'} \cdot B_{Y'} = B^2 +2 (- K_S \cdot B) -(-K_{Y} \cdot B_Y) = 9 + 18 - 14=13$ (Proposition \ref{p-ele-tf-numbers}). 
Hence $Y'$ is a Fano threefold  of No.\ 2-34  (Remark \ref{r-ele-tr-P2}, Lemma \ref{l-P2-nonFano}).  
\qedhere


\end{proof}

\begin{lem}\label{l-P2-2-32}
We use Notation \ref{n-ele-tr-P2}. 
Assume that both $Y$ and $Y'$ are Fano threefolds of No.\ 2-32. 
Then 
$(-K_X)^3=30$, $\deg B = 2$, $p_a(B)=0$, 
$-K_Y \cdot B_Y =8$, and $-K_{Y'} \cdot B_{Y'}=8$ ($X$ is 3-13). 
\end{lem}

\begin{proof}
We have $Y  \in |\MO_{\P^2 \times \P^2}(1, 1)|$. 
Set $\alpha_1 := g : Y \to S = \P^2$ and $\alpha_2 := h : Y \to Z = \P^2$. 
For each $i \in \{1, 2\}$, let $\ell_i$ be a fibre of $\alpha_i$ and 
set $H_i := \alpha_i^*\MO_{\P^2}(1)$. 
Recall that we have $H_1 \cdot \ell_2 =H_2 \cdot \ell_1 = 1$ and $-K_Y \sim 2H_1 + 2H_2$ (Proposition \ref{p-pic2-Pic}). 
Hence $-K_{Y/S} \sim 2H_1 -H_2$. 
Set $d_1 :=d = H_1 \cdot B_Y$ and $d_2 := H_2 \cdot B_Y$, i.e., $(d_1, d_2)$ is the bidegree of the curve $B_Y$ in $\P^2 \times \P^2$. 
It follows from Proposition \ref{p-ele-tf-numbers} that  
\[
-K_{Y/S} \cdot B_Y = (2H_1 -H_2) \cdot B_Y = 2d_1 -d_2
\]
\[
(-K_Y)^3 = (-K_{Y'})^3 = (-K_{Y})^3 -4 (-K_{Y/S}) \cdot B_Y +2B^2 
=(-K_{Y})^3 -8d_1 +4d_2 +2d_1^2.
\]
Thus $d_1^2 -4d_1+2d_2=0$. 

\begin{claim*}
One of the following holds. 
\begin{enumerate}
\renewcommand{\labelenumi}{(\roman{enumi})}
\item $(d_1, d_2) = (1, 0)$, i.e., $B_Y$ is a fibre of $\alpha_2 : Y \to \P^2$. 
\item $(d_1, d_2) = (1, 2)$.  
\item $(d_1, d_2)=(2, 1)$. 
\item $(d_1, d_2) = (d, d)$. 
\end{enumerate}
\end{claim*}

\begin{proof}[Proof of Claim]
Recall that $B_Y$ is a regular subsection of $\alpha_1=g$. 
By Proposition \ref{p-FCB-centre}, 
$B_Y$ is either a fibre of $\alpha_2$ or a regular subsection of $\alpha_2$. 
If $B_Y$ is a fibre of $\alpha_2$, then (i) holds, 
because $H_2 \cdot B_Y = 0$ and $H_1 \cdot B_Y = H_1 \cdot \ell_2 = 1$. 
Assume that $B_Y$ is  a regular subsection of $\alpha_2$. 
Then we get $B_Y \simeq \alpha_1(B_Y) \simeq \alpha_2(B_Y)$, 
which implies $p_a(B_Y) = \frac{(d_1-1)(d_1-2)}{2} = \frac{(d_2-1)(d_2-2)}{2}$, 
because $\alpha_i(B_Y)$ is a smooth plane curve of degree $d_i$. 
Then one of (ii)-(iv) holds. This completes the proof of  Claim. 
\end{proof}

It is easy to see that each of (i)--(iii) is not a solution of 
$d_1^2 -4d_1 +2d_2=0$. 
Hence (iv) holds. 
We then get $d^2 -2d=0$, which implies $d=2$, i.e., $(d_1, d_2)=(2, 2)$. 
Hence $-K_Y \cdot B_Y = 2(d_1+d_2) = 8$, $-K_{Y/S}\cdot B_Y = 2d_1 - d_2 =2$, 
$p_a(B)=0$, $(-K_X)^3 = (-K_{Y})^3 -2(-K_{Y}) \cdot B_Y+2p_a(B)-2 =48 -16 +0-2 =30$,  
$-K_{Y'} \cdot B_{Y'} = B^2 +2 (- K_S \cdot B) -(-K_{Y} \cdot B_Y) = 4 +12 - 8 =8$ (Proposition \ref{p-ele-tf-numbers}). 
\qedhere

\end{proof}

\begin{lem}\label{l-P2-2-32-same}
We use Notation \ref{n-ele-tr-P2}. 
Assume that $Y$ is of No.\ 2-32 and $(-K_{Y'})^3 = 54$. 
Then $\deg B =1$ or $\deg B =3$. 
\end{lem}

\begin{proof}
The following proof is identical to that of Lemma \ref{l-P2-2-32}. 
Set $d_1 :=\deg B$ and $d_2 := h^*\MO_{\P^2}(1) \cdot B_Y$. 
The same argument as in the proof of Lemma \ref{l-P2-2-32} (the argument before Claim) implies $-K_{Y/S} \cdot B_Y  = 2d_1 -d_2$ and 
\[
54 = (-K_{Y'})^3 
=(-K_{Y})^3 -8d_1 +4d_2 +2d_1^2 = 48-8d_1 +4d_2 +2d_1^2, 
\]
i.e., $d_1^2 -4d_1 + 2d_2 = 3$. 
In particular, $d_1 \not\in 2\Z$ and $d_1 <5$. Hence we get $\deg B = d_1 \in \{1, 3\}$, as required. 
\end{proof}

\begin{lem}\label{l-P2-2-34}
We use Notation \ref{n-ele-tr-P2}. 
Assume that $Y$ is of No.\ 2-34 and $(-K_{Y'})^3 \in \{ 48, 54\}$. 
Let $D$ be a fibre of the second projection ${\rm pr}_2 : Y =\P^2 \times \P^1 \to \P^1$. 
Then one of the following holds. 
\begin{enumerate}
\item 
$D \cdot B_Y = 1$, 
$Y'$ is a Fano threefold of No.\ 2-34, $(-K_X)^3=36$, $\deg B = 2$, $p_a(B)=0$, 
$-K_Y \cdot B_Y = 8$, and $-K_{Y'} \cdot B'=8$ ($X$ is 3-17). 
\item 
$D \cdot B_Y =4$, 
$Y'$ is a Fano threefold of No.\ 2-34, $(-K_X)^3=18$, $\deg B = 4$, $p_a(B)=3$, 
$-K_Y \cdot B_Y = 20$, and $-K_{Y'} \cdot B'=20$ ($X$ is 3-3). 
\item 
$D \cdot B_Y =1$, 
$Y'$ is a Fano threefold of No.\ 2-32, $(-K_X)^3=42$, $\deg B = 1$, $p_a(B)=0$, 
$-K_Y \cdot B_Y = 5$, and $-K_{Y'} \cdot B'=2$ ($X$ is 3-24). 
\item 
$D \cdot B_Y =3$, 
$Y'$ is a Fano threefold of No.\ 2-32, $(-K_X)^3=24$, $\deg B = 3$, $p_a(B)=1$, 
$-K_Y \cdot B_Y = 15$, and $-K_{Y'} \cdot B'=12$ ($X$ is 3-7). 
\end{enumerate}
\end{lem}

\begin{proof}
Note that $g={\rm pr}_1 : Y=\P^2 \times \P^1 \to \P^2=S$ and 
$h = {\rm pr}_2 : \P^2 \times \P^1 \to \P^1=Z$.  
We obtain $K_{Y} = g^*K_S + h^*K_{\P^1}$, i.e., 
$-K_{Y/S} \sim -h^*K_{\P^1} \sim 2D$. 
We then get $-K_{Y/S} \cdot B_{Y} = 2D \cdot B_{Y}$ and $-K_Y \cdot B_Y = 2D \cdot B_Y +3d$. 
It follows from Proposition \ref{p-ele-tf-numbers} that  
\begin{equation}\label{e1-P2-2-34}
(-K_{Y'})^3 = (-K_{Y})^3 -4 (-K_{Y/S}) \cdot B_Y +2B^2 
=54 - 8 D \cdot B_{Y} +2d^2. 
\end{equation}
Since $D$ is nef, we have $D \cdot B_Y \geq 0$. 
By 
\[
(-K_{Y})^2 \cdot D = \MO_{\P^2 \times \P^1}(3, 2)^2 \cdot 
\MO_{\P^2 \times \P^1}(0, 1) = 9, 
\]
we obtain $0 \leq D \cdot B_Y < (-K_Y)^2 \cdot D=9$ (Lemma \ref{l-DB-bound}).

(1), (2) Assume $(-K_{Y'})^3 =54$, i.e., $Y'$ is a Fano threefold of No.\ 2-34 (Remark \ref{r-ele-tr-P2}, Lemma \ref{l-P2-nonFano}).  
Then (\ref{e1-P2-2-34}) implies $d^2  = 4 D \cdot B_Y$. 
By  $0 \leq D \cdot B_Y < 9$, we get $(D \cdot B_Y, d)  \in \{(1, 2), (4, 4)\}$. 
\begin{enumerate}
\item 
If $(D \cdot B_Y, d) =(1, 2)$, then 
$p_a(B_Y) = 0, -K_Y \cdot B_Y = 2 + 6 =8,$ 
$(-K_X)^3 = (-K_Y)^3 -2 (-K_Y) \cdot B_Y + 2p_a(B) -2 = 54 -16+0-2 = 36$, 
$-K_{Y'} \cdot B_{Y'} = B^2 +2 (- K_S \cdot B) -(-K_{Y} \cdot B_Y) = 4 + 12 - 8=8$  (Proposition \ref{p-ele-tf-numbers}). 
\item 
If $(D \cdot B_Y, d) =(4, 4)$, then 
$p_a(B_Y) = 3, -K_Y \cdot B_Y = 8 + 12 =20,$ 
$(-K_X)^3 = (-K_Y)^3 -2 (-K_Y) \cdot B_Y + 2p_a(B) -2 = 54 -40+6-2 = 18$, 
$-K_{Y'} \cdot B_{Y'} = B^2 +2 (- K_S \cdot B) -(-K_{Y} \cdot B_Y) = 16 + 24 - 20=20$  (Proposition \ref{p-ele-tf-numbers}). 
\end{enumerate}

(3), (4) 
Assume $(-K_{Y'})^3 =48$, i.e., $Y'$ is of No.\ 2-32  (Remark \ref{r-ele-tr-P2}, Lemma \ref{l-P2-nonFano}). 
Applying Lemma \ref{l-P2-2-32-same} after switching $Y$ and $Y'$, 
it holds that  $d = 1$ or $d=3$. 
By  $4D \cdot B_Y = 3+ d^2$ (\ref{e1-P2-2-34}), 
we get $(D \cdot B_Y, d) \in \{(1, 1), (3, 3)\}$. 
\begin{enumerate}
\item[(3)] Assume $(D \cdot B_Y, d) =(1, 1)$. 
Then $p_a(B_Y) = 0, -K_Y \cdot B_Y = 2 + 3 =5,$ 
$(-K_X)^3 = (-K_Y)^3 -2 (-K_Y) \cdot B_Y + 2p_a(B) -2 = 54 -10+0-2 = 42$, 
$-K_{Y'} \cdot B_{Y'} = B^2 +2 (- K_S \cdot B) -(-K_{Y} \cdot B_Y) = 1 + 6 - 5=2$ (Proposition \ref{p-ele-tf-numbers}).  
\item[(4)] Assume $(D \cdot B_Y, d) =(3, 3)$. 
Then $p_a(B_Y) = 1, -K_Y \cdot B_Y = 6+9  =15,$ 
$(-K_X)^3 = (-K_Y)^3 -2 (-K_Y) \cdot B_Y + 2p_a(B) -2 = 54 -30+2-2 = 24$, 
$-K_{Y'} \cdot B_{Y'} = B^2 +2 (- K_S \cdot B) -(-K_{Y} \cdot B_Y) = 9 + 18- 15=12$ (Proposition \ref{p-ele-tf-numbers}). 
\end{enumerate}
\end{proof}







\begin{thm}\label{t-ele-tr-P2}
Let 
\[
\begin{tikzcd}
& X \arrow[ld, "\sigma"'] \arrow[rd, "\sigma'"] \arrow[dd, "f"]\\
Y \arrow[rd, "g"']& & Y' \arrow[ld, "g'"]\\
& S:=\P^2. 
\end{tikzcd}
\]
be an elementary transform (cf.\ Definition \ref{d-ele-tf}), 
where $f: X \to S=\P^2$ and $g : Y \to S=\P^2$ are Fano conic bundles, 
and $\sigma: X \to Y$ is a blowup along a smooth curve $B_Y$. 
Set $B:=g(B_Y)$. 
If also $Y'$ is Fano, then we assume that $(-K_Y)^3 \leq (-K_{Y'})^3$. 
Then one of the following holds. 
    \begin{center}
\begin{longtable}{ccccccccc}
$X$ & $Y$ & $Y'$ & $(-K_X)^3$& $\deg B$ & $p_a(B)$ & $-K_Y \cdot B_Y$ & $-K_{Y'} \cdot B_{Y'}$ &\\ \hline
3-3 & 2-34 & 2-34 & $18$ & $4$ & $3$ & $20$ & $20$\\ \hline
3-5 & 2-34 & non-Fano & $20$ & $2$ & $0$ & $16$ & $0$\\ \hline
 3-7 &  2-32 & 2-34 & $24$ & $3$ & $1$ & $12$ & $15$\\ \hline
3-8 &  2-24 & 2-34 & $24$ & $2$ & $0$ & $2$ & $14$\\ \hline
 3-9 &  2-36 & 2-36 & $26$ & $4$ & $3$ & $20$ & $20$\\ \hline
   3-11 &  2-34 & 2-35 & $28$ & $3$ & $1$ & $13$ & $14$\\ \hline
3-12 &  2-27 & 2-34 & $28$ & $2$ & $0$ & $4$ & $12$\\ \hline
 3-13 &  2-32 & 2-32 & $30$ & $2$ & $0$ & $8$ & $8$\\ \hline
 3-14 &  2-35 & 2-36 & $32$ & $3$ & $1$ & $12$ & $15$\\ \hline
 3-15 &  2-31 & 2-34 & $32$ & $2$ & $0$ & $6$ & $10$\\ \hline
 3-16 &  2-27 & 2-32 & $34$ & $1$ & $0$ & $1$ & $6$\\ \hline
  3-16 &  2-32 & 2-35 & $34$ & $2$ & $0$ & $6$ & $10$\\ \hline
   3-17 &  2-34 & 2-34 & $36$ & $2$ & $0$ & $8$ & $8$\\ \hline
  3-19 &  2-35 & 2-35 & $38$ & $2$ & $0$ & $8$ & $8$\\ \hline
 3-20 &  2-31 & 2-32 & $38$ & $1$ & $0$ & $3$ & $4$\\ \hline
3-21 & 2-34 & non-Fano & $38$ & $1$ & $0$ & $7$ & $0$\\ \hline
 3-22 &  2-34 & 2-36 & $40$ & $2$ & $0$ & $6$ & $10$\\ \hline
 3-23 &  2-31 & 2-35 & $42$ & $1$ & $0$ & $1$ & $6$\\ \hline
 3-24 &  2-32 & 2-34 & $42$ & $1$ & $0$ & $2$ & $5$\\ \hline
 3-26 &  2-34 & 2-35 & $46$ & $1$ & $0$ & $3$ & $4$\\ \hline
 3-29 &  2-35 & 2-36 & $50$ & $1$ & $0$ & $2$ & $5$\\ \hline\\
 \\
 \caption{Elementary transforms over $\P^2$}\label{table-ele-tr-P2}
      \end{longtable}
  \end{center}

\end{thm}

We say that the above diagram is called an elemental transform over $\P^2$ 
{\em of type 2-xx-vs-2-yy} if $Y$ is 2-xx and $Y'$ is 2-yy. 
In this case, we say that $X$ has a conic bundle structure over $\P^2$ {\em of type 2-xx-vs-2-yy}. 

\begin{proof}
If $Y'$ is non-Fano, then 
the assertion follows from Lemma \ref{l-P2-nonFano}. 
Hence we may assume that $Y'$ is Fano. 
By Remark \ref{r-ele-tr-P2}, each of $Y$ and $Y'$ is one of 
\[
\text{2-24,\quad 2-27,\quad  2-31,\quad  2-32,\quad  2-34,\quad  2-35,\quad  2-36.} 
\]
If $Y$ is 2-24 (resp. 2-27, resp. 2-31), then apply Lemma \ref{l-P2-2-24} (resp. Lemma \ref{l-P2-2-27}, resp. Lemma \ref{l-P2-2-31}). 
By $(-K_Y)^3 \leq (-K_{Y'})^3$, 
we may assume that each of $Y$ and $Y'$ is one of 2-32,  2-34,  2-35,  2-36. 
Depending on No.\ of $Y'$,  the assertion holds by applying lemmas as follows. 
\begin{itemize}
\item $Y'$ is 2-32:  Lemma \ref{l-P2-2-32}.
\item $Y'$ is 2-34:  Lemma \ref{l-P2-2-34}.
\item $Y'$ is 2-35:  Lemma \ref{l-P2-2-35}.
\item $Y'$ is 2-36:  Lemma \ref{l-P2-2-36}.
\end{itemize}
\end{proof}


\subsection{Fano conic bundles over $\F_1$ ($\rho=3$)}\label{ss-pic3-F1}

The purpose of this subsection is to classify Fano conic bundles $X \to \F_1$ with $\rho(X)=3$. 
Note that such a Fano conic bundle is obtained 
from a Fano conic bundle $X' \to S'=\P^2$ as follows. 

\begin{lem}\label{l-F1-pic3}
Let $f:  X \to S:=\F_1$ be a Fano conic bundle with $\rho(X)=3$. 
Then there exists a cartesian diagram 
\[
\begin{tikzcd}
X \arrow[r, "\sigma"] \arrow[d, "f"] & X' \arrow[d, "f'"] \\
S=\F_1 \arrow[r, "\tau"] &  S':=\P^2
\end{tikzcd}
\]
such that 
\begin{enumerate}
\item $f' : X' \to S' =\P^2$ is a Fano conic bundle, 
\item $\tau : \F_1 \to \P^2$ is the blowdown of the $(-1)$-curve $\Gamma$ on $S=\F_1$, 
\item $P := \tau(\Gamma) \not\in \Delta_{f'}$, and 
\item $\sigma$ is the blowup along the smooth fibre $f'^{-1}(P)$. 
\end{enumerate}
\end{lem}

\begin{proof}
By $\rho(X) = 3 =2+1 =\rho(S)+1$, $f^{-1}(\Gamma)$ is irreducible 
(Lemma \ref{l-FCB-pic-irre}). 
Then the assertion follows from Proposition \ref{p-FCB-(-1)}. 
\end{proof}

The main result of this subsection is Theorem \ref{t-F1-pic3}. 
To this end, we now establish the following auxiliary result. 

\begin{lem}\label{l-CE-avoiding-Ex}
Let $Y$ be a Fano threefold with $\rho(Y)=2$ and 
let $R_1$ and $R_2$ be the extremal rays of $\NE(Y)$. 
Assume that 
$R_1$ is of type $C$ and $R_2$ is of type $E$. 
Let $f_1 : Y \to T$ and $f_2: Y \to Z$ be the contractions of $R_1$ and $R_2$, respectively. 
Then one of the following holds. 
\begin{enumerate}
    \item $\Ex(f_2)$ dominates $T$. 
    \item $f_1$ is of type $C_2$, $f_2$ is of type $E_1$, and 
    $Z$ is a Fano threefold of index two. 
\end{enumerate}
\end{lem}

\begin{proof}
Set $D := \Ex(f_2)$. 
Assume that (1) does not hold, i.e., $f_1(D) \subsetneq T$. 
It is enough to show (2). 
By $f_1(D) \subsetneq T$, there exists a curve $C$ on $D$ contracted by $f_1$. 
Then $f_2$ is of type $E_1$, as otherwise $f_2(D)$ would be a point. 
In particular, $Z$ is a Fano threefold of index $r_Z \geq 2$ \cite[Proposition 5.8]{ATIII}. 
We have 
\[
K_Y = f_2^* K_Z + D. 
\]
Let $\ell_Y$ be an extremal rational curve of $f_1 : Y \to T=\P^2$. 
By $D \cdot \ell_Y=0$, the following holds: 
\[
\{1, 2\} \ni -K_Y \cdot \ell_Y = -f_2^*K_Z \cdot \ell_Y = 
-K_Z \cdot (f_2)_*(\ell_Y) \in r_Z \Z. 
\]
By $r_Z \geq 2$, it holds that $-K_Y \cdot \ell_Y = 2$ and $r_Z =2$. 
Hence $f_1$ is of type $C_2$. 
\end{proof}

We are ready to prove the main theorem of this subsection.

\begin{thm}\label{t-F1-pic3}
Set $S := \F_1$ and $S' :=\P^2$. 
Let 
\[
\begin{tikzcd}
X \arrow[r, "\sigma"] \arrow[d, "f"] & X' \arrow[d, "f'"] \\
S :=\F_1 \arrow[r, "\tau"] & S' := \P^2, 
\end{tikzcd}
\]
be a cartesian diagram, where $f$ and $f'$ are Fano conic bundles and 
$\tau : S =\F_1 \to \P^2 =S'$ is the blowdown of the $(-1)$-curve on $S=\F_1$. 
Assume that $\rho(X)=3$. 
Then one of the following holds. 
  \begin{center}
\begin{longtable}{cccccc}
$X$ & $(-K_X)^3$ & $X'$ & $(-K_{X'})^3$ & $X'\to S'$ & $\deg \Delta_{f'}$\\ \hline
3-4 & $18$ & 2-18 & $24$ & $X' \xrightarrow{2:1} \P^2\times \P^1 \xrightarrow{\pr_1} \P^2$ & $4$\\ \hline
3-8 & $24$ & 2-24 & $30$ & $X' \hookrightarrow \P^2\times \P^2 \xrightarrow{\pr_1} \P^2, X' \in |\MO(1, 2)|$ & $3$\\ \hline
3-24 & $42$ & 2-32 & $48$ & $X' =W \hookrightarrow \P^2 \times \P^2 \xrightarrow{{\rm pr}_i} \P^2$ & $0$\\ \hline
3-28 & $48$ & 2-34 & $54$ & ${\rm pr}_1 : X'=\P^2 \times \P^1 \to \P^2$ & $0$ \\ \hline
3-30 & $50$ & 2-35 & $56$ & $X'=\P_{\P^2}(\MO_{\P^2} \oplus \MO_{\P^2}(1)) \to \P^2$ & $0$\\ \hline
      \end{longtable}
  \end{center} 
\end{thm}

\begin{proof}
Note that the centre $B_{X'}$ of the blowup $\sigma : X \to X'$ is 
a smooth fibre of $f' : X' \to S'=\P^2$ (Proposition \ref{p-FCB-centre}), i.e., a fibre over a closed point of $S' \setminus \Delta_{f'}$. 
In particular, 
\[
(-K_{X'}) \cdot B_{X'} =2,  \qquad 
(-K_X)^3 = (-K_{X'})^3 -6,
\]
and 
\begin{enumerate}
\item $B_{X'}$ is disjoint from any curve $L$ on $X'$ satisfying 
$-K_{X'} \cdot L =1$ (Lemma \ref{l-line-meeting}). 
\end{enumerate}
Recall that the contraction $f' :X' \to S'=\P^2$ is of type $C_1$ or $C_2$. 
Let $h : X' \to Z$ be the contraction of the extremal ray not corresponding to $f'$. 
We now prove (2)-(4) below. 
\begin{enumerate}
\setcounter{enumi}{1}
    \item $h$ is not of type $E_3, E_4$, nor $E_5$. 
    \item If $h$ is of type $E_1$, then $f'$ is of type $C_2$ and $r_Z =2$, 
    where $r_Z$ denotes the index of the Fano threefold $Z$. 
    \item $h$ is not of type $C_1$.  
\end{enumerate}
Let us show (2) and (3). 
Assume that the type of $h$ is one of $E_1, E_3, E_4, E_5$. 
Set $D :=\Ex(h)$. Since the extremal ray corresponding to $h : X' \to Z$ is of length $1$, 
$D$ is covered by curves $L$ satisfying $-K_{X'} \cdot L=1$. 
Hence $D \cap B_{X'}=\emptyset$ by (1). 
In particular, $D$ does not dominate $S'$. 
Then Lemma \ref{l-CE-avoiding-Ex} implies that (2) and (3) hold. 
Let us show (4). 
Suppose that $h$ is of type $C_1$, i.e., 
$h : X' \to \P^2$ is a Fano conic bundle with $\Delta_h \neq \emptyset$. 
Then the blowup centre $B_{X'}$ of $\sigma : X \to X'$ 
must be a smooth fibre of $h$ (Proposition \ref{p-FCB-centre}). 
However, this would imply $[B_{X'}] \in R_{f'} \cap R_h = \{0\}$ for the extremal rays $R_{f'}$ and $R_h$ of $f'$ and $h$, respectively. 
This is absurd. 
This completes the proofs of (2)-(4).


\medskip

Assume that  $f' :X' \to S'=\P^2$ is of type $C_1$, i.e., 
$\deg \Delta_{f'} \neq 0$. 
The possibilities for $X' \to S'$ are as follows (Subsection \ref{ss-table-pic2}): 

  \begin{center}
\begin{longtable}{ccp{10cm}c}
No.\ & $(-K_{X'})^3$ & Description  & Extremal rays\\ \hline
2-2 & $6$ & a split double cover of $\mathbb{P}^2\times \mathbb{P}^1$ with $\mathcal L \simeq \MO(2, 1)$ & $C_1+D_1$\\ \hline
2-6 & $12$ & 
a smooth divisor on $\mathbb{P}^2\times \mathbb{P}^2$ of bidegree $(2,2)$, or 
a split double cover of $W$ with $\mathcal L^{\otimes 2} \simeq \omega_W^{-1}$ & $C_1+C_1$\\ \hline
2-8 & $14$ &  a split double cover of $V_7=\mathbb{P}(\MO_{\mathbb{P}^2}\oplus \MO_{\mathbb{P}^2}(1))$ & $C_1+E_3\,{\rm or}\,E_4$\\ \hline
2-9 & $16$ & blowup of $\P^3$ along a curve of genus $5$ and degree $7$ & $C_1+E_1$\\ \hline
2-11 & $18$ & blowup of $V_3$ along a line & $C_1+E_1$\\ \hline
2-13 & $20$ & blowup of $Q$ along a curve of genus $2$ and degree $6$ & $C_1+E_1$\\ \hline
2-16 & $22$ & blowup of $V_4$ along a conic & $C_1+E_1$\\ \hline
2-18 & $24$ & a split double cover of $V_7=\mathbb{P}(\MO_{\mathbb{P}^2}\oplus \MO_{\mathbb{P}^2}(1))$ & $C_1+D_2$\\ \hline
2-20 & $26$ & blowup of $V_5$ along a cubic rational curve & $C_1+E_1$\\ \hline
2-24 & $30$ & a smooth divisor on $\mathbb{P}^2\times \mathbb{P}^2$ of bidegree $(1,2)$ & $C_1+C_2$\\ \hline
      \end{longtable}
  \end{center} 
By (2)-(4),  
the case when the type of $h:X' \to Z$ is one of $C_1, E_1, E_3, E_4$ 
does not occur. 
It follows from $0<(-K_X)^3 = (-K_{X'})^3 -6$ that also 2-2 is impossible. 
The remaining possibility is when $X'$ is 2-18 or 2-24. 
This completes the proof for the case when $f'$ is of type $C_1$. 

\medskip

Assume that $f' :X' \to S'=\P^2$ is of type $C_2$, i.e., 
$\deg \Delta_{f'} = 0$. 
The possibilities for $f': X' \to S'$ are as follows (Subsection \ref{ss-table-pic2}):

  \begin{center}
\begin{longtable}{ccp{10cm}c}
No. & $(-K_{X'})^3$ & Description  & Extremal rays\\ \hline
2-24 & $30$ & a smooth divisor on $\mathbb{P}^2\times \mathbb{P}^2$ of bidegree $(1,2)$ & $C_2+C_1$\\ \hline
2-27 & $38$ & blowup of $\P^3$ along a cubic rational curve & $C_2+E_1$\\ \hline
2-31 & $46$ & blowup of $Q$ along a line & $C_2+E_1$\\ \hline
2-32 & $48$ & $W$ & $C_2+C_2$\\ \hline
2-34 & $54$ & $\P^2 \times \P^1$ & $C_2+D_3$\\ \hline
2-35 & $56$ & $V_7=\mathbb{P}(\MO_{\mathbb{P}^2}\oplus \MO_{\mathbb{P}^2}(1))$ & $C_2+E_2$\\ \hline
2-36 & $62$ & $\mathbb{P}(\MO_{\mathbb{P}^2}\oplus \MO_{\mathbb{P}^2}(2))$ & $C_2+E_5$\\ \hline
      \end{longtable}
  \end{center} 
By (2)-(4),  
the case when the type of $h:X' \to Z$ is one of $C_1, E_1, E_5$ is  impossible. 
The remaining possibility is when $X'$ is 2-32, 2-34, 2-35. 
We are done. 
\end{proof}

\subsection{Classification ($\rho=3$)}\label{ss-pic3-classify}

\begin{nasi}\label{n-pic3-vol}
Let $X$ be a Fano threefold with $\rho(X)=3$. 
By Theorem \ref{t-pic3-structure}, one of the following holds. 
\begin{enumerate}
\item[(I)] $X$ has a conic bundle structure over $\P^2$. 
In this case, the following holds (Theorem \ref{t-ele-tr-P2}):   
\[
(-K_X)^3 \in \{18, 20, 24, 26, 28, 30, 32, 34, 36, 38, 40, 42, 46, 50\}. 
\]
\item[(II)] $X$ has a conic bundle structure over $\F_1$. 
We have $(-K_X)^3 \in \{18, 24, 42, 48, 50\}$ (Theorem \ref{t-F1-pic3}). 
\item[(III)] $X$ is primitive. 
Then $(-K_X)^3 \in \{ 12, 14, 48, 52\}$ and $X$ has a conic bundle structure over $\P^1 \times \P^1$ 
\cite[Theorem 1.1, Theorem 4.17]{ATIII}. 
\item[(IV)] 
There exist a line $L$ and a conic $C$ on $\P^3$ such that $L \cap C = \emptyset$ 
and $X \simeq \Bl_{L \amalg C} \P^3$. In particular, $(-K_X)^3 = 36$ 
(Lemma \ref{l-3-18}). 
\item[(V)] $X$ satisfies the assumption and the conclusion of Proposition \ref{p-V-2dP}.  
In particular, $(-K_X)^3 \in \{ 22, 26, 44\}$ and 
$X$ has a conic bundle structure over $\P^1 \times \P^1$. 
\end{enumerate}
\end{nasi}



\begin{nota}\label{n-exactly3}
Let $X$ be a Fano threefold with $\rho(X)=3$. 
Assume that there exist exactly three (resp. four) extremal rays $R_1, R_2, R_3$ (resp. $R_1. R_2, R_3, R_4$). 
We have exactly three (resp. four) two-dimensional extremal faces $F_1 := R_3+R_1, F_2 := R_1 +R_2, F_3 := R_2+R_3$ 
(resp.  $F_1 := R_4+R_1, F_2 := R_1 +R_2, F_3 := R_2+R_3, F_4 :=R_3+R_4$). 
\begin{itemize}
\item For each $i \in \{1, 2, 3\}$ (resp. $i \in \{ 1, 2, 3, 4\}$), let $f_i : X \to Y_i$ be the contraction of $f_i$. 
\item For each $i \in \{1, 2, 3\}$ (resp. $i \in \{ 1, 2, 3, 4\}$), let $\varphi_i : X \to Z_i$ be the contraction of $\varphi_i$. Although the existence of $\varphi_i$ is not clear, we shall prove it. 
\item Let $g_{ij} : Y_i \to Z_j$ be the induced morphism whenever it exists. 
\end{itemize}
Let $H_{Z_i}$ be an ample Cartier divisor which generates $\Pic\,Z_i \simeq \Z$. 
Set $H_i := \varphi_i^*H_{Z_i}$.  
\begin{itemize}
\item When $f_i : X \to Y_i$ is of type $C_1$, $\Delta_{f_i}$ denotes its discriminant divisor. 
\item When $f_i : X \to Y_i$ is of type $E_1$, 
$B_i$ denotes its blowup centre. 
\item 
When $Z_i =\P^2$ and $\varphi_i : X \to Z_i$ is a conic bundle, 
the square diagram consisting of $X, Z_i, Y_{i-1}, Y_i$ (where $Y_0:=Y_3$) 
is a elementrary transform. 
In this case, $B_{Z_i} := g_{ii}(B_i)$. 
\end{itemize}
\end{nota}
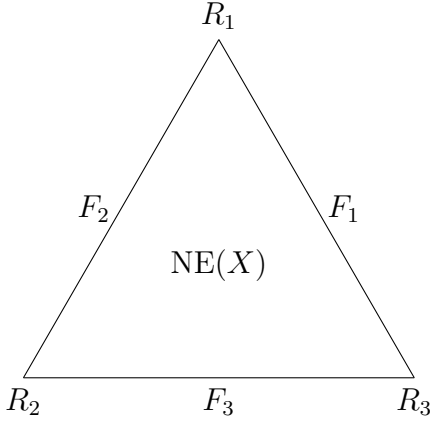
\[
\begin{tikzpicture}[commutative diagrams/every diagram]
\coordinate (X)  at (0,0);
\coordinate (Y1) at (0,3);
\coordinate (Y2) at (-{(3/2)*sqrt(3)}, -3/2);
\coordinate (Y3) at ({(3/2)*sqrt(3)}, -3/2);

\draw (Y1) -- node[midway, left] {$F_2$} (Y2);
\draw (Y2) -- node[midway, below] {$F_3$} (Y3);
\draw (Y3) -- node[midway, right] {$F_1$} (Y1);

\node at (X) {$\NE(X)$};
\node[above] at (Y1) {$R_1$};
\node[below] at (Y2) {$R_2$};
\node[below] at (Y3) {$R_3$};
\end{tikzpicture}
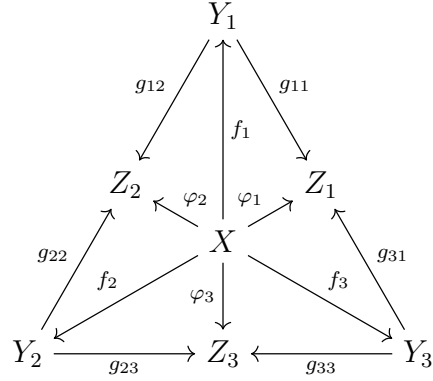
\hspace{20mm}
\begin{tikzpicture}[commutative diagrams/every diagram,
    declare function={R=3;Rs=R*cos(60);}]
 \path 
  (0,0)  node(X) {$X$} 
  (90:R) node (Y1) {$Y_1$}
  (210:R) node (Y2) {$Y_2$}
  (-30:R) node (Y3) {$Y_3$}  
  (30:Rs) node(Z1) {$Z_1$} 
  (150:Rs) node(Z2) {$Z_2$} 
  (270:Rs) node(Z3) {$Z_3$};
 \path[commutative diagrams/.cd, every arrow, every label]
 (X) edge[swap] node {$f_1$} (Y1)
 (X) edge[swap] node {$f_2$} (Y2)
 (X) edge node {$f_3$} (Y3)
 (X) edge node {$\varphi_1$} (Z1)
 (X) edge[swap] node {$\varphi_2$} (Z2)
 (X) edge[swap] node {$\varphi_3$} (Z3)
 (Y1) edge node {$g_{11}$} (Z1)
 (Y1) edge[swap] node {$g_{12}$} (Z2)
 (Y2) edge node {$g_{22}$} (Z2)
 (Y2) edge[swap] node {$g_{23}$} (Z3)
 (Y3) edge node {$g_{33}$} (Z3)
 (Y3) edge[swap] node {$g_{31}$} (Z1);
\end{tikzpicture}
\]
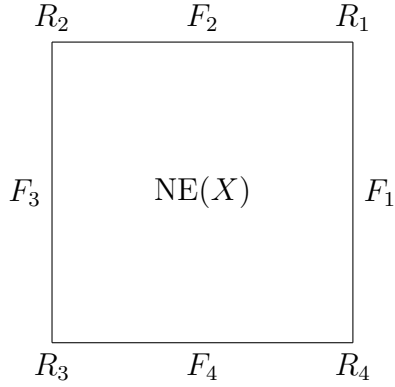
\[
\begin{tikzpicture}[commutative diagrams/every diagram]
\coordinate (X)  at (0,0);
\coordinate (Y1) at (2, 2);
\coordinate (Y2) at (-2, 2);
\coordinate (Y3) at (-2, -2);
\coordinate (Y4) at (2, -2);

\draw (Y1) -- node[midway, above] {$F_2$} (Y2);
\draw (Y2) -- node[midway, left] {$F_3$} (Y3);
\draw (Y3) -- node[midway, below] {$F_4$} (Y4);
\draw (Y4) -- node[midway, right] {$F_1$} (Y1);

\node at (X) {$\NE(X)$};
\node[above] at (Y1) {$R_1$};
\node[above] at (Y2) {$R_2$};
\node[below] at (Y3) {$R_3$};
\node[below] at (Y4) {$R_4$};
\end{tikzpicture}
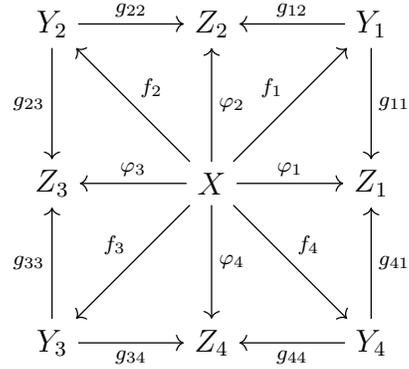
\hspace{20mm}
\begin{tikzpicture}[commutative diagrams/every diagram,
    declare function={R=3;Rs=R*cos(45);}]
 \path 
  (0,0)  node(X) {$X$} 
  (45:R) node (Y1) {$Y_1$}
  (135:R) node (Y2) {$Y_2$}
  (-135:R) node (Y3) {$Y_3$}  
  (-45:R) node (Y4) {$Y_4$}  
  (0:Rs) node(Z1) {$Z_1$} 
  (90:Rs) node(Z2) {$Z_2$} 
  (180:Rs) node(Z3) {$Z_3$}
  (270:Rs) node(Z4) {$Z_4$};
 \path[commutative diagrams/.cd, every arrow, every label]
 (X) edge node {$f_1$} (Y1)
 (X) edge[swap] node {$f_2$} (Y2)
 (X) edge[swap] node {$f_3$} (Y3)
 (X) edge node {$f_4$} (Y4)
 (X) edge node {$\varphi_1$} (Z1)
 (X) edge[swap] node {$\varphi_2$} (Z2)
 (X) edge[swap] node {$\varphi_3$} (Z3)
 (X) edge node {$\varphi_4$} (Z4)
 (Y1) edge node {$g_{11}$} (Z1)
 (Y1) edge[swap] node {$g_{12}$} (Z2)
 (Y2) edge node {$g_{22}$} (Z2)
 (Y2) edge[swap] node {$g_{23}$} (Z3)
 (Y3) edge node {$g_{33}$} (Z3)
 (Y3) edge[swap] node {$g_{34}$} (Z4)
 (Y4) edge node {$g_{44}$} (Z4)
 (Y4) edge[swap] node {$g_{41}$} (Z1);
\end{tikzpicture}
\]


\begin{lem}\label{l-MM-formula}
Let $f: X \to S$ be a threefold conic bundle. 
The following holds for a Cartier divisor $D$ on $S$: 
\[
\Delta_f \cdot D = -4 K_S \cdot D - (-K_X)^2 \cdot f^*D. 
\]
\end{lem}
\begin{proof}
See \cite[Proposition 3.16]{ATIII}. 
\end{proof}

\begin{lem}\label{l-P1P1-Delta-bideg}
Let $f: X \to \P^1 \times \P^1$ be a threefold conic bundle. 
For each $i \in \{1, 2\}$ and the contraction $\pr_i \circ f : X \xrightarrow{f} \P^1 \times \P^1 \xrightarrow{\pr_i} \P^1$, 
we set $H_i := (\pr_i \circ f)^*\MO_{\P^1}(1)$. 
Assume $-K_X \equiv a_1 H_1 + a_2 H_2 + D$ for some $a_1, a_2 \in \Q$ and $\Q$-divisor $D$. 
Let $\deg \Delta_f =(d_1, d_2)$. 
Then the following hold. 
\begin{enumerate}
\item $(-K_X)^2 \cdot H_1 = 2a_2 H_1 \cdot H_2 \cdot D +  H_1 \cdot D^2$. 
\item $(-K_X)^2 \cdot H_2 = 2a_1 H_1 \cdot H_2 \cdot D +  H_2 \cdot D^2$. 
\item $d_1 = 8 - (-K_X)^2 \cdot H_2$ and $d_2 = 8 - (-K_X)^2 \cdot H_1$. 
\end{enumerate}
\end{lem}

\begin{proof}
The assertion (1) follows from 
\[
(-K_X)^2 \cdot H_1 = (a_1 H_1 + a_2 H_2 + D)^2 \cdot H_1 = (a_2H_2+D)^2 \cdot H_1 
\]
\[
= ( 2 a_2 H_2 \cdot D +D^2 ) \cdot H_1 = 2a_2 H_1 \cdot H_2 \cdot D +  H_1 \cdot D^2. 
\]
By symmetry, (2) holds. 
The assertion (3) follows from Lemma \ref{l-MM-formula}.  
\end{proof}


\begin{lem}\label{l rho3 primitive type}
Let $X$ be a primitive Fano threefold with $\rho(X)=3$. 
Then the following hold. 
\begin{enumerate}
\item $(-K_X)^3 =12 \Leftrightarrow$ all the extremal rays are of type $C_1$. 
\item $(-K_X)^3 =14 \Leftrightarrow$ 
$X$ has extremal rays $R_1$ and $R_2$ such that 
$R_1$ is of type $C_1$ and $R_2$ is of type $E_1$. 
\item 
$(-K_X)^3 =48 \Leftrightarrow$ all the extremal rays are of type $C_2$. 
\item 
$(-K_X)^3 =52 \Leftrightarrow$ 
$X$ has extremal rays $R_1$ and $R_2$ such that 
$R_1$ is of type $C_2$ and $R_2$ is of type $E_1$. 
\end{enumerate}
\end{lem}

\begin{proof}
All the implications in direction $\lq\lq \Leftarrow"$ follow from \cite[Theorem 6.7 and Theorem 6.17]{ATIII}. 
The opposite implications  hold by $(-K_X)^3 \in \{12, 14, 48, 52\}$ 
 \cite[Theorem 6.1]{ATIII}. 
\end{proof}

\begin{prop}[No.\ {\hyperref[table-3-1]{3-1}}]\label{p-pic3-1}
Let $X$ be a Fano threefold with $\rho(X)=3$ and $(-K_X)^3=12$. 
Then the following hold. 
\begin{enumerate}
\item $X$ has exactly three extremal rays. 
In what follows, we use Notation \ref{n-exactly3}. 
\item 
The contractions of the extremal faces are as in the following diagram. 
\begin{enumerate}
\item $f_1$ is of type $C_1$, $\deg \Delta_{f_1} =(4, 4)$.
\item $f_2$ is of type $C_1$, $\deg \Delta_{f_2} =(4, 4)$. 
\item $f_3$ is of type $C_1$, $\deg \Delta_{f_3} =(4, 4)$. 
\end{enumerate}
\item $-K_X \sim H_1 + H_2 + H_3$. 
\item 
$\varphi_1 \times \varphi_2 \times \varphi_3 : X \to \P^1_1 \times \P^1_2 \times \P^1_3$ is a split double cover satisfying $ (\varphi_1 \times \varphi_2 \times \varphi_3)_*\MO_X/\MO_{\P^1 \times \P^1 \times \P^1}  \simeq \MO_{\P^1 \times \P^1 \times \P^1}(-1, -1, -1)$. 
\end{enumerate}
\end{prop}

\[
\begin{tikzpicture}[commutative diagrams/every diagram,
    declare function={R=3;Rs=R*cos(60);}]
 \path 
  (0,0)  node(X) {$X$} 
  (90:R) node (Y1) {$\P^1_1\times \P^1_2$}
  (210:R) node (Y2) {$\P^1_2\times \P^1_3$}
  (-30:R) node (Y3) {$\P^1_3\times \P^1_1$}  
  (30:Rs) node(Z1) {$\P^1_1$} 
  (150:Rs) node(Z2) {$\P^1_2$} 
  (270:Rs) node(Z3) {$\P^1_3$};
 \path[commutative diagrams/.cd, every arrow, every label]
 (X) edge[swap] node {$f_1$} (Y1) 
 (X) edge[swap] node {$f_2$} (Y2)
 (X) edge node {$f_3$} (Y3)
 (X) edge node {$\varphi_1$} (Z1)
 (X) edge[swap] node {$\varphi_2$} (Z2)
 (X) edge[swap] node {$\varphi_3$} (Z3)
 (Y1) edge node {$g_{11}=\pr_1$} (Z1)
 (Y1) edge[swap] node {$g_{12}=\pr_2$} (Z2)
 (Y2) edge node {$g_{22}=\pr_1$} (Z2)
 (Y2) edge[swap] node {$g_{23}=\pr_2$} (Z3)
 (Y3) edge node {$g_{33}=\pr_1$} (Z3)
 (Y3) edge[swap] node {$g_{31}=\pr_2$} (Z1);
\end{tikzpicture}
\]

\begin{proof}
Only (III) of (\ref{n-pic3-vol}) holds. 
Then all the extremal rays are of type $C_1$ 
(Lemma \ref{l rho3 primitive type}). 
Hence all the assertions follow from \cite[Theorem 6.7 and Remark 6.8]{ATIII}. 
\end{proof}

\begin{prop}[No.\  {\hyperref[table-3-2]{3-2}}]\label{p-pic3-2}
Let $X$ be a Fano threefold with $\rho(X)=3$ and $(-K_X)^3=14$. 
Then the following hold. 
\begin{enumerate}
\item $X$ has exactly three extremal rays. 
In what follows, we use Notation \ref{n-exactly3}. 
\item The contractions of the extremal faces are as in the following diagram.
\begin{enumerate}
\item $f_1$ is of type $C_1$, $\deg \Delta_{f_1} = (2, 5)$. 
\item $f_2$ is of type $E_1$, $p_a(B_2)= 0$, $-K_{Y_2} \cdot B_2 = 0$. 
\item $f_3$ is of type $E_1$, $p_a(B_3)= 0$, $-K_{Y_3} \cdot B_3 = 0$. 
\item 
$D := \Ex (\varphi_3) = \Ex(f_2) = \Ex(f_3) \simeq \P^1 \times \P^1$, $\varphi_3(D)$ is a point, 
and $f_1^*\MO_{\P^1 \times \P^1}(a, b)|_D \simeq \MO_{\P^1 \times \P^1}(a, 2b)$ for all $a, b \in \Z$. 
\end{enumerate}
\item $-K_X \sim 2H_1 + H_2 + D$. 
\item 
$X$ is a member of the complete linear system $|\MO_P(2)\otimes\pi^*\MO_{\mathbb{P}^1\times \mathbb{P}^1}(2,3)|$ on the $\mathbb{P}^2$-bundle $\pi\colon P=\mathbb{P}(\MO_{\mathbb{P}^1\times \mathbb{P}^1}\oplus \MO_{\mathbb{P}^1\times \mathbb{P}^1}(-1,-1)^{\oplus 2})\to \mathbb{P}^1\times \mathbb{P}^1$ such that 
$f_1$ coincides with the composition $X \hookrightarrow P \xrightarrow{\pi} \P^1 \times \P^1$. 
\end{enumerate}
\end{prop}
\[
\begin{tikzpicture}[commutative diagrams/every diagram,
    declare function={R=3;Rs=R*cos(60);}]
 \path 
  (0,0)  node(X) {$X$} 
  (90:R) node (Y1) {$\P^1\times \P^1$}
  (210:R) node (Y2) {$Y_{\text{non-Fano}}$}
  (-30:R) node (Y3) {$Y'_{\text{non-Fano}}$}  
  (30:Rs) node(Z1) {$\P^1$} 
  (150:Rs) node(Z2) {$\P^1$} 
  (270:Rs) node(Z3) {$Z_3$};
 \path[commutative diagrams/.cd, every arrow, every label]
 (X) edge[swap] node {$f_1$} (Y1)
 (X) edge[swap] node {$f_2$} (Y2)
 (X) edge node {$f_3$} (Y3)
 (X) edge node {$\varphi_1$} (Z1)
 (X) edge[swap] node {$\varphi_2$} (Z2)
 (X) edge[swap] node {$\varphi_3$} (Z3)
 (Y1) edge node {$g_{11}=\pr_1$} (Z1)
 (Y1) edge[swap] node {$g_{12}=\pr_2$} (Z2)
 (Y2) edge node {$g_{22}$} (Z2)
 (Y2) edge[swap] node {$g_{23}$} (Z3)
 (Y3) edge node {$g_{33}$} (Z3)
 (Y3) edge[swap] node {$g_{31}$} (Z1);
\end{tikzpicture}
\]
\begin{proof}
Only (III) of (\ref{n-pic3-vol}) holds. 
Then $X$ has an extremal rays $R_1$ and $R_2$ such that 
$R_1$ is of type $C_1$ and $R_2$ is of type $E_1$ 
(Lemma \ref{l rho3 primitive type}). 
It follows from \cite[Proposition 6.16]{ATIII} that (4) holds. 
Moreover, we get the above commutative diagram except for $f_3, g_{31}, g_{33}$ (note that the target $Y_2$ of $f_2$ is non-Fano, because 
$X$ is primitive and $f_2$ is of type $E_1$). 
By \cite[Lemma 6.13]{ATIII}, Proposition \ref{p-nonFano-flop}, and Lemma \ref{l-exactly3}, 
we get the above commutative diagram 
and 
all the assertions hold except for (3) and $\deg \Delta_{f_1} =(2, 5)$. 

Let us show (3). 
Since $H_1, H_2, D$ are linearly independent in $\Pic(X) (\simeq \Z^3)$, 
we can write $-K_X \equiv a_1 H_1 + a_2 H_2 + b D$ for some $a_1, a_2, b \in \Q$. 
Since $f_1|_D : D \to \P^1 \times \P^1$ is a double cover \cite[Lemma 6.13]{ATIII}, 
we obtain $b=1$ (consider the intersection with a fibre of $f_1 : X \to \P^1 \times \P^1$). 
Recall that $-K_X|_D =-D|_D =\MO_{\P^1 \times \P^1}(1, 1)$ (Lemma \ref{l-nonFano-blowdown}). 
Hence 
\[
\MO_{\P^1 \times \P^1}(1, 1) = -K_X|_D \equiv 
(a_1 H_1 + a_2 H_2 + D)|_D \equiv \MO_{\P^1 \times \P^1}(a_1 -1, 2a_2 -1), 
\]
which implies $a_1  =2$ and $a_2  =1$. 
Thus (3) holds.

It suffices to compute $\deg \Delta_{f_1}=(d_1, d_2)$. 
The following holds for any divisor $D$ on $\P^1 \times \P^1$ (Lemma \ref{l-MM-formula}):
\[
\Delta_{f_1} \cdot D = -4K_{\P^1 \times \P^1} \cdot D -
(-K_X)^2 \cdot f_1^*D.
\]
We have 
$H_1^2 \equiv H_2^2 \equiv 0$, $H_1 \cdot H_2 \cdot D = \deg (f_1|_D :D \to \P^1 \times \P^1) =2$, $D^2 \cdot H_1 = \MO_{\P^1 \times \P^1}(-1, -1) \cdot \MO_{\P^1 \times \P^1}(1, 0)=-1,$ and 
$D^2 \cdot H_2 = \MO_{\P^1 \times \P^1}(-1, -1) \cdot \MO_{\P^1 \times \P^1}(0, 2)=-2$. 
Hence 
\[
(-K_X)^2 \cdot f_1^*\MO_{\P^1 \times \P^1}(1, 0) = (2H_1+H_2+D)^2 \cdot H_1 
\]
\[
=(H_2+D)^2 \cdot H_1 = (2H_2 \cdot D + D^2) \cdot H_1 = 4 -1 =3, 
\]
\[
d_2 = \Delta_{f_1} \cdot \MO_{\P^1 \times \P^1}(1, 0) = 8-3 =5, 
\]
\[
(-K_X)^2 \cdot f_1^*\MO_{\P^1 \times \P^1}(0, 1) = (2H_1+H_2+D)^2 \cdot H_2 
\]
\[
=(2H_1+D)^2 \cdot H_2 = (4H_1 \cdot D + D^2) \cdot H_2 = 8 -2 =6, 
\]
\[
d_1 = \Delta_{f_1} \cdot \MO_{\P^1 \times \P^1}(0, 1) = 8-6 =2. 
\]
Thus $\deg \Delta_{f_1} = (2, 5)$. 
\end{proof}

\begin{prop}[No.\  {\hyperref[table-3-3]{3-3}}]\label{p-pic3-3}
Let $X$ be a Fano threefold with $\rho(X)=3$ and $(-K_X)^3=18$. 
Assume that $X$ has a conic bundle structure over $\P^2$. 
Then the following hold. 
\begin{enumerate}
\item $X$ has exactly three extremal rays. 
In what follows, we use Notation \ref{n-exactly3}. 
\item The contractions of the extremal faces are as in the following diagram. 
\begin{enumerate}
\item $f_1$ is of type $C_1$, $\deg \Delta_{f_1} = (3, 3)$. 
\item $f_2$ is of type $E_1$, $p_a(B_2)= 3$, $-K_{Y_2} \cdot B_2 = 20$. 
\item $f_3$ is of type $E_1$, $p_a(B_3)= 3$, $-K_{Y_3} \cdot B_3 = 20$. 
\item $\deg B_{Z_3} = 4$. 
\end{enumerate}
\item $-K_X \sim H_1 + H_2 + H_3$. 
\item 
$\varphi_1 \times \varphi_2 \times \varphi_3 : X \to \P^1 \times \P^1 \times \P^2$ is a closed immersion whose image $(\varphi_1 \times \varphi_2 \times \varphi_3)(X)$ 
is of tridegree $(1, 1, 2)$, i.e., linearly equivalent to $\MO_{\P^1 \times \P^1 \times \P^2}(1, 1, 2)$. 
\end{enumerate}
\end{prop}
\[
\begin{tikzpicture}[commutative diagrams/every diagram,
    declare function={R=3;Rs=R*cos(60);}]
 \path 
  (0,0)  node(X) {$X$} 
  (90:R) node (Y1) {$\P^1\times \P^1$}
  (210:R) node (Y2) {$\P^2\times \P^1$}
  (-30:R) node (Y3) {$\P^2\times \P^1$}  
  (30:Rs) node(Z1) {$\P^1$} 
  (150:Rs) node(Z2) {$\P^1$} 
  (270:Rs) node(Z3) {$\P^2$};
 \path[commutative diagrams/.cd, every arrow, every label]
 (X) edge[swap] node {$f_1$} (Y1)
 (X) edge[swap] node {$f_2$} (Y2)
 (X) edge node {$f_3$} (Y3)
 (X) edge node {$\varphi_1$} (Z1)
 (X) edge[swap] node {$\varphi_2$} (Z2)
 (X) edge[swap] node {$\varphi_3$} (Z3)
 (Y1) edge node {$g_{11} =\pr_1$} (Z1)
 (Y1) edge[swap] node {$g_{12}=\pr_2$} (Z2)
 (Y2) edge node {$g_{22}=\pr_2$} (Z2)
 (Y2) edge[swap] node {$g_{23}=\pr_1$} (Z3)
 (Y3) edge node {$g_{33}=\pr_1$} (Z3)
 (Y3) edge[swap] node {$g_{31}=\pr_2$} (Z1);
\end{tikzpicture}
\]
\begin{proof}
By our assumption, we obtain the elementary transform over $\P^2$ consisting of 
$f_2, g_{23}, f_3, g_{33}$ 
such that (b), (c), (d) hold (Theorem \ref{t-ele-tr-P2}). 
For the second projections 
$g_{31} :\P^2 \times \P^1 \to \P^1 =:Z_1$ and 
$g_{22} : \P^2 \times \P^1 \to \P^1 =:Z_2$, 
we get the above diagram except for $f_1, g_{11}, g_{12}$. 
By Lemma \ref{l-exactly3 2P^1}, 
(1) holds, 
$f_1 :=\varphi_1 \times \varphi_3 : X \to Z_1 \times Z_3$ is a contraction,  and we obtain the above commutative diagram consisting of contractions. 
The assertion  (3) holds by the following (cf. Proposition \ref{p-pic2-Pic} and Lemma \ref{l-ele-tf-K-relation}): 
\[
-2K_X \sim -f_2^*K_{\P^2 \times \P^1} -f_3^*K_{\P^2 \times \P^1} -\varphi_3^*B_{\P^2} 
\]
\[
= (2H_2 + 3H_3) +(3H_3+2H_1)  - 4H_3 = 2(H_1+H_2+H_2). 
\]

Let us compute $\deg \Delta_{f_1} = (d_1, d_2)$. 
We have  $H_1 \cdot H^2_3=1$ and $H_2 \cdot H_3^2 =1$, because 
$f_2 : X \to \P^2 \times \P^1$ and $f_3 : X \to \P^2 \times \P^1$ are birational. 
It holds that 
$H_1 \cdot H_2 \cdot H_3 = H_1 \cdot H_2 \cdot (H_1 + H_2 + H_3) = -K_X \cdot \zeta =2$, 
where $\zeta$ denotes a fibre of $X \to \P^1 \times \P^1$. 
Hence 
\[
(-K_X)^2 \cdot H_1 = (H_1 + H_2 +H_3)^2 \cdot H_1 = (H_2+H_3)^2 \cdot H_1 
= 2 H_1 \cdot H_2 \cdot H_3 + H_3^2 \cdot H_1= 4 +1  = 5. 
\]
Lemma \ref{l-MM-formula} implies 
\[
d_2 = \Delta_{f_1} \cdot \MO_{\P^1 \times \P^1}(1, 0) = -4K_S \cdot \MO_{\P^1 \times \P^1}(1, 0) -(-K_X)^2 \cdot H_1=8-5=3. 
\]
By symmetry, we get $\deg \Delta_{f_1} = (3, 3)$. 
Thus (2) holds.

Let us show (4). 
Note that $\varphi := \varphi_1 \times \varphi_2 \times \varphi_3 : X \to \P^1 \times \P^1 \times \P^2$ is a finite morphism, because 
the intersection of the corresponding three extremal faces is equal to $\{0\}$. 
For $X' := \varphi(X)$, let $\psi: X \to X'$ be the induced morphism. 
Then $\psi$ is birational, because $f_2 : X \to \P^2 \times \P^1$ factors through $\psi$. 
Set $H'_1 := \MO_{\P^1 \times \P^1 \times \P^2}(1, 0, 0)|_{X'},$ 
$H'_2 := \MO_{\P^1 \times \P^1 \times \P^2}(0, 1, 0)|_{X'},$ and 
$H'_3 := \MO_{\P^1 \times \P^1 \times \P^2}(0, 0, 1)|_{X'}.$ 
By 
\[
H_2' \cdot H_3'^2 = H_2 \cdot H_3^2 = 1,\quad 
H_1' \cdot H_3'^2 = H_1 \cdot H_3^2 = 1,\quad  
H'_1 \cdot H'_2 \cdot H'_3  =  H_1 \cdot H_2 \cdot H_3  =2, 
\]
$X'$ is of tridegree $(1, 1, 2)$. 
For the conductor $C$ of the normalisation $\psi : X \to X'$, 
we obtain $\MO_X(K_X+C) \simeq \psi^*\omega_{X'}$. 
By (3) and the adjunction formula $(\omega_{\P^1 \times \P^1 \times \P^2} \otimes \MO(X')|_{X'} \simeq \omega_{X'}$, we get $C \equiv 0$, which implies $C=0$. 
Hence $\varphi$ is a closed immersion. 
Thus (4) holds. 
\end{proof}

\begin{lem}\label{l-F1cart-KX}
Let $g: Y \to \P^2$ be a Fano conic bundle. 
Take a  blowup $\tau : \F_1 \to \P^2$ at a point $P \in \P^2 \setminus \Delta_g$. 
Consider the following commutative diagram: 
\[
\begin{tikzcd}
& X \arrow[d, "f"] \arrow[r, "\sigma"] \arrow[rd, "\varphi"]& Y \arrow[d, "g"]\\
\P^1 & \F_1 \arrow[l, "\pi"'] \arrow[r, "\tau"] & \P^2
\end{tikzcd}
\]
where the square diagram is cartesian and $\pi : \F_1 \to \P^1$ denotes the $\P^1$-bundle. 
Then it holds that 
\[
-K_X \sim -\sigma^*K_Y -H_{\P^2} +H_{\P^1}, 
\]
where 
$H_{\P^2}$ and $H_{\P^1}$ are the pullbacks of $\MO_{\P^2}(1)$ and $\MO_{\P^1}(1)$ 
by $\varphi: X \to \P^2$ and $ X \xrightarrow{f} \F_1 \xrightarrow{\pi} \P^1$, 
respectively. 
\end{lem}

\begin{proof}
Let $E_{X/Y}$ and $E_{\F_1/\P^2}$ be the exceptional prime divisors 
of $\sigma : X \to Y$ and $\tau : \F_1 \to \P^2$, respectively. 
We have $K_X \sim \sigma^*K_Y +E_{X/Y}$,  $K_{\F_1} = \tau^*K_{\P^2}+E_{\F_1/\P^2}$, and $E_{X/Y} = f^*E_{\F_1/\P^2}$. 
Recall that $K_{\F_1} +2E_{\F_1/\P^2} \sim \pi^*\MO_{\P^1}(-3)$. 
Hence 
\[
\pi^*\MO_{\P^1}(-3) \sim K_{\F_1} +2E_{\F_1/\P^2} \sim \tau^*K_{\P^2}+3E_{\F_1/\P^2} \sim \tau^*\MO_{\P^2}(-3) + 3E_{\F_1/\P^2}, 
\]
which implies $E_{\F_1/\P^2} \sim \tau^*\MO_{\P^2}(1) -\pi^*\MO_{\P^1}(1)$. 
To summarise, 
\[
-K_X \sim -\sigma^*K_Y -E_{X/Y} 
 = -\sigma^*K_Y -f^*E_{\F_1/\P^2}
 \sim -\sigma^*K_Y -H_{\P^2} +H_{\P^1}. 
\]
\end{proof}

\begin{prop}[No.\  {\hyperref[table-3-4]{3-4}}]\label{p-pic3-4}
Let $X$ be a Fano threefold with $\rho(X)=3$ and $(-K_X)^3=18$. 
Assume that $X$ has no conic bundle structure over $\P^2$. 
Then the following hold. 
\begin{enumerate}
\item $X$ has exactly three extremal rays. 
In what follows, we use Notation \ref{n-exactly3}. 
\item The contractions of the extremal faces are as in the following diagram. 
\begin{enumerate}
\item $f_1$ is of type $C_1$, $\deg \Delta_{f_1} = (4, 2)$.  
\item $f_2$ is of type $C_1$, $\Delta_{f_2} \in |\tau^* \MO_{\P^2}(4)|$. 
\item $f_3$ is of type $E_1$ and the blowup centre of $f_3$ is a smooth fibre of  $g_{33} : Y_{\text{2-18}} \to \P^2$. 
\end{enumerate}
\item $-K_X \sim H_1 +H_2 + H_3$. 
\end{enumerate}
\end{prop}
\[
\begin{tikzpicture}[commutative diagrams/every diagram,
    declare function={R=3;Rs=R*cos(60);}]
 \path 
  (0,0)  node(X) {$X$} 
  (90:R) node (Y1) {$\P^1\times \P^1$}
  (210:R) node (Y2) {$\F_1$}
  (-30:R) node (Y3) {$Y_{\text{2-18}}$}  
  (30:Rs) node(Z1) {$\P^1$} 
  (150:Rs) node(Z2) {$\P^1$} 
  (270:Rs) node(Z3) {$\P^2$};
 \path[commutative diagrams/.cd, every arrow, every label]
 (X) edge[swap] node {$f_1$} (Y1)
 (X) edge[swap] node {$f_2$} (Y2)
 (X) edge node {$f_3$} (Y3)
 (X) edge node {$\varphi_1$} (Z1)
 (X) edge[swap] node {$\varphi_2$} (Z2)
 (X) edge[swap] node {$\varphi_3$} (Z3)
 (Y1) edge node {$g_{11}=\pr_1$} (Z1)
 (Y1) edge[swap] node {$g_{12}=\pr_2$} (Z2)
 (Y2) edge node {$g_{22}$} (Z2)
 (Y2) edge[swap] node {$g_{23} =\tau$} (Z3)
 (Y3) edge node {$g_{33}$} (Z3)
 (Y3) edge[swap] node {$g_{31}$} (Z1);
\end{tikzpicture}
\]

\begin{proof}
Only (II) of  (\ref{n-pic3-vol}) holds. 
Hence $X$ has a conic bundle structure over $\F_1$. 
By Theorem \ref{t-F1-pic3}, we obtain $X \simeq Y_{\text{2-18}} \times_{\P^2} \F_1$. 
We then get the above commutative diagram except for $f_1, g_{11}, g_{12}$. 
Moreover, (b) and (c) hold. 
By Lemma \ref{l-exactly3 2P^1}, 
(1) holds, 
$f_1 :=\varphi_1 \times \varphi_2 : X \to Z_1 \times Z_2$ is a contraction,  and we obtain the above commutative diagram consisting of contractions.


Let us show (3). 
By Lemma \ref{l-F1cart-KX} and Proposition \ref{p-pic2-Pic}(2), we obtain 
\[
-K_X \sim -f_3^*K_{Y_{\text{2-18}}} - H_3 + H_2 \sim (H_1 + 2H_3)-H_3+H_2 = H_1 + H_2 + H_3.
\]
Thus (3) holds. 

It suffices to show $\deg \Delta_{f_1} = (4, 2)$. 
We have $H_1 \cdot H_3^2 = g_{31}^*\MO_{\P^1}(1) \cdot (2\ell_{g_{33}}) = 2$ for an extremal rational curve $\ell_{g_{33}}$ of $g_{33}$ (Proposition \ref{p-pic2-Pic}). 
It holds that $H_2 \cdot H_3^2 = 0$, because the divisors 
$H_2$ and $H_3$ come from $\F_1$. 
We obtain $H_1 \cdot H_2 \cdot H_3 =2$ by 
\[
18 = (-K_X)^3 = (H_1 + H_2 + H_3)^3 = 6H_1 \cdot H_2 \cdot H_3 + 3H_1 \cdot H_3^2 + 3H_2 \cdot H_3^2 = 6H_1 \cdot H_2 \cdot H_3 +6. 
\]
Then the following hold (Lemma \ref{l-P1P1-Delta-bideg}):
\[
(-K_X)^2 \cdot H_i = 2 H_1 \cdot H_2 \cdot H_3 + H_i \cdot H_3^2 =
\begin{cases}
6 \qquad (i=1)\\
4 \qquad (i=2)
\end{cases}
\]
\[
d_2 = 8 - 6 =2, \qquad d_1 = 8-4=4.
\]

\end{proof}

\begin{prop}[No.\  {\hyperref[table-3-5]{3-5}}]\label{p-pic3-5}
Let $X$ be a Fano threefold with $\rho(X)=3$ and $(-K_X)^3=20$. 
Then the following hold. 
\begin{enumerate}
\item $X$ has exactly three extremal rays. 
In what follows, we use Notation \ref{n-exactly3}. 
\item The contractions of the extremal faces are as in the following diagram. 
\begin{enumerate}
\item $f_1$ is of type $E_1$, $p_a(B_1)= 0$, $-K_{Y_1} \cdot B_1 = 16$, and $B_1$ is a regular subsection of $g_{12}$ which is of bidegree $(2, 5)$. 
\item $f_2$ is of type $E_1$, $p_a(B_2)= 0$, $-K_{Y_2} \cdot B_2 = 0$. 
\item $f_3$ is of type $E_1$, $p_a(B_3)= 0$, $-K_{Y_3} \cdot B_3 = 0$. 
\item $\deg B_{Z_2} = 2$. 
\item $\varphi_3$ is birational, 
$D:= \Ex(\varphi_3) = \Ex(f_2) = \Ex(f_3) \simeq \P^1 \times \P^1$, and $\varphi_3(D)$ is a point. 
\end{enumerate}
\item $-K_X \sim 2H_1 + H_2 +D$. 
\end{enumerate}
\end{prop}
\[
\begin{tikzpicture}[commutative diagrams/every diagram,
    declare function={R=3;Rs=R*cos(60);}]
 \path 
  (0,0)  node(X) {$X$} 
  (90:R) node (Y1) {$\P^2\times \P^1$}
  (210:R) node (Y2) {$Y_{\text{non-Fano}}$}
  (-30:R) node (Y3) {$Y'_{\text{non-Fano}}$}  
  (30:Rs) node(Z1) {$\P^1$} 
  (150:Rs) node(Z2) {$\P^2$} 
  (270:Rs) node(Z3) {$Z_3$};
 \path[commutative diagrams/.cd, every arrow, every label]
 (X) edge[swap] node {$f_1$} (Y1)
 (X) edge[swap] node {$f_2$} (Y2)
 (X) edge node {$f_3$} (Y3)
 (X) edge node {$\varphi_1$} (Z1)
 (X) edge[swap] node {$\varphi_2$} (Z2)
 (X) edge[swap] node {$\varphi_3$} (Z3)
 (Y1) edge node {$g_{11}=\pr_2$} (Z1)
 (Y1) edge[swap] node {$g_{12}=\pr_1$} (Z2)
 (Y2) edge node {$g_{22}$} (Z2)
 (Y2) edge[swap] node {$g_{23}$} (Z3)
 (Y3) edge node {$g_{33}$} (Z3)
 (Y3) edge[swap] node {$g_{31}$} (Z1);
\end{tikzpicture}
\]
\begin{proof}
Only (I) of  (\ref{n-pic3-vol}) holds. 
By Theorem \ref{t-ele-tr-P2}, $X$ has a conic bundle structure over $\P^2$ of type 2-34-vs-non-Fano, 
and hence we get the above square diagram consisting of $f_1, g_{12}, f_2, g_{22}$. 
Since $f_2 : X \to Y_{\text{non-Fano}}$ is a contraction of type $E_1$ to a non-Fano threefold $Y_{\text{non-Fano}}$, 
we obtain the square diagram consisting of $f_2, g_{23}, f_3, g_{33}$ (Proposition \ref{p-nonFano-flop}). 
Then (b)-(e) hold. 
Recall that $D = \Ex(f_3) = \Ex(f_2)$ coincides with the proper transform of $g_{12}^{-1}(B_{Z_2})$. 
The restriction of the second projection $g_{11} : \P^2 \times \P^1 \to \P^1$ to 
$g_{12}^{-1}(B_{Z_2}) \simeq B_{Z_2} \times \P^1$ is again the second projection 
$g_{11}|_{g_{12}^{-1}(B_{Z_2})} : g_{12}^{-1}(B_{Z_2}) \simeq B_{Z_2} \times \P^1 \to \P^1$. 
Therefore, $X$ has a curve $C$ contracted by each of $X \xrightarrow{f_1} \P^2 \times \P^1 \xrightarrow{g_{11}} \P^1$ and $f_3 : X \to Y'_{\text{non-Fano}}$ (as otherwise, there is a curve $C$ contracted by each of $X \xrightarrow{f_1} \P^2 \times \P^1 \xrightarrow{g_{11}} \P^1$ and $f_2 : X \to Y_{\text{non-Fano}}$, which is absurd because $Y_{\text{non-Fano}}$ has exactly two non-trivial contractions and $g_{22}$ and $g_{23}$). 
Hence (1) holds (Lemma \ref{l-exactly3}) and we obtain the above commutative diagram. 
Note that the bidegree $(d_1, d_2)$ of $B_1 \subset \P^2 \times \P^1$ can be computed by $d_1 = \pr_1^*\MO_{\P^2}(1) \cdot B_1 =\deg B_{Z_2} =2$ and 
$16 = -K_{\P^2 \times \P^1} \cdot B_1 = 3d_1 + 2d_2 = 6+2d_2$. 
Thus (2) holds. 

It suffices to show (3). 
The following holds (Lemma \ref{l-ele-tf-K-relation}):
\[
-2K_X \sim -f_1^*K_{\P^2 \times \P^1} -f_2^*K_{Y_{\text{non-Fano}}}-\varphi_2^*B_{Z_2}. 
\]
We have $-f_1^*K_{\P^2 \times \P^1} \sim 2H_1  +3H_2$, $\varphi_2^*B_{Z_2} \sim 2H_2$, 
and $K_X =f_2^*K_{Y_{\text{non-Fano}}} +D$. 
Hence 
\[
-K_X \sim 
-f_1^*K_{\P^2 \times \P^1} + (K_X -f_2^*K_{Y_{\text{non-Fano}}})-\varphi_2^*B_{Z_2}
\]
\[
\sim (2H_1+3H_2) + D - 2H_2 = 2H_1 + H_2 +D.
\]
Thus (3) holds. 
\qedhere


\end{proof}

\begin{prop}[No.\  {\hyperref[table-3-6]{3-6}}]\label{p-pic3-6}
Let $X$ be a Fano threefold with $\rho(X)=3$ and $(-K_X)^3=22$. 
Then the following hold. 
\begin{enumerate}
\item $X$ has exactly three extremal rays. 
In what follows, we use Notation \ref{n-exactly3}. 
\item The contractions of the extremal faces are as in the following diagram. 
\begin{enumerate}
\item $f_1$ is of type $C_1$, $\deg \Delta_{f_1} = (3, 2)$. 
\item $f_2$ is of type $E_1$, $p_a(B_2)= 1$, $-K_{Y_2} \cdot B_2 = 16$. 
\item $f_3$ is of type $E_1$, $p_a(B_3)= 0$, $-K_{Y_3} \cdot B_3 = 4$. 
\item $\varphi_3$ is a blowup along a disjoint union of a line and an elliptic curve of degree four. 
\end{enumerate}
\item $-K_X \sim H_1 + H_2 + H_3$ under Notation \ref{n-exactly3}. 
\end{enumerate}
\end{prop}
\[
\begin{tikzpicture}[commutative diagrams/every diagram,
    declare function={R=3;Rs=R*cos(60);}]
 \path 
  (0,0)  node(X) {$X$} 
  (90:R) node (Y1) {$\P^1\times \P^1$}
  (210:R) node (Y2) {$Y_{\text{2-25}}$}
  (-30:R) node (Y3) {$Y_{\text{2-33}}$}  
  (30:Rs) node(Z1) {$\P^1$} 
  (150:Rs) node(Z2) {$\P^1$} 
  (270:Rs) node(Z3) {$\P^3$};
 \path[commutative diagrams/.cd, every arrow, every label]
 (X) edge[swap] node {$f_1$} (Y1)
 (X) edge[swap] node {$f_2$} (Y2)
 (X) edge node {$f_3$} (Y3)
 (X) edge node {$\varphi_1$} (Z1)
 (X) edge[swap] node {$\varphi_2$} (Z2)
 (X) edge[swap] node {$\varphi_3$} (Z3)
 (Y1) edge node {$g_{11} =\pr_1$} (Z1)
 (Y1) edge[swap] node {$g_{12}=\pr_2$} (Z2)
 (Y2) edge node {$g_{22}$} (Z2)
 (Y2) edge[swap] node {$g_{23}$} (Z3)
 (Y3) edge node {$g_{33}$} (Z3)
 (Y3) edge[swap] node {$g_{31}$} (Z1);
\end{tikzpicture}
\]
\begin{proof}
Only (V) of  (\ref{n-pic3-vol}) holds.  
All the assertions except for (3) follow from Proposition \ref{p-V-2dP}. 
The assertion (3) holds by the following (Lemma \ref{l-K-disjoint-blowup}, Proposition \ref{p-pic2-Pic}): 
\[
-K_X \sim -f_2^*K_{Y_{\text{2-25}}}-f_3^*K_{Y_{\text{2-33}}} + \varphi_3^*K_{\P^3} 
\sim (H_2 +2H_3) + (H_1+ 3H_3)-4H_3 = H_1 + H_2 +H_3. 
\]
\end{proof}

\begin{lem}\label{l-dP-fib-CI}
Let $\sigma: X \to Y$ be a blowup along a smooth curve $\Gamma$ on a Fano threefold $Y$. 
Let $\pi:X \to \P^1$ be a contraction. 
Take a Cartier divisor $D$ on $Y$. 
Assume that {\rm (i)} and {\rm (ii)} hold. 
\begin{enumerate}
\renewcommand{\labelenumi}{(\roman{enumi})}
\item $(-K_Y) \cdot D^2= (-K_Y) \cdot \Gamma$. 
\item 
$D \sim \sigma_*F$ for a fibre $F$ of $\pi$. 
\end{enumerate}
Then $\Gamma$ is a complete intersection of two members of $|D|$. 
\end{lem}

\begin{proof}
See \cite[Lemma 5.42]{ATIII}. 
\end{proof}

\begin{prop}[No.\  {\hyperref[table-3-7]{3-7}}]\label{p-pic3-7}
Let $X$ be a Fano threefold with $\rho(X)=3$ and $(-K_X)^3=24$. 
Assume that $X$ has a conic bundle structure over $\P^2$ 
of type 2-32-vs-2-34.
Then the following hold. 
\begin{enumerate}
\item $X$ has exactly three extremal rays. 
In what follows, we use Notation \ref{n-exactly3}. 
\item The contractions of the extremal faces are as in the following diagram.
\begin{enumerate}
\item $f_1$ is of type $E_1$, $p_a(B_1)= 1$, $-K_{Y_1} \cdot B_1 = 12$. 
\item $f_2$ is of type $E_1$, $p_a(B_2)= 1$, $-K_{Y_2} \cdot B_2 = 15$. 
\item $f_3$ is of type $E_1$, $p_a(B_3)= 1$, $-K_{Y_3} \cdot B_3 = 15$. 
\item $\deg B_{Z_1} = 3$. 
\item $\deg B_{Z_2} = 3$. 
\end{enumerate}
\item $-K_X \sim H_1 + H_2 + H_3$. 
\item 
The blowup centre $B_1$ of $f_1 : X \to W$ is an elliptic curve which is a complete intersection of two members of $|-\frac{1}{2}K_W|$. 
\end{enumerate}
\end{prop}
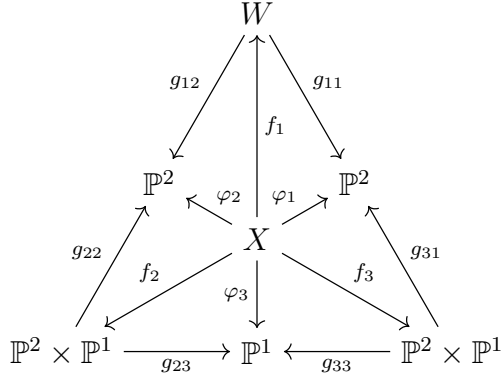
\[
\begin{tikzpicture}[commutative diagrams/every diagram,
    declare function={R=3;Rs=R*cos(60);}]
 \path 
  (0,0)  node(X) {$X$} 
  (90:R) node (Y1) {$W$}
  (210:R) node (Y2) {$\P^2 \times \P^1$}
  (-30:R) node (Y3) {$\P^2 \times \P^1$}  
  (30:Rs) node(Z1) {$\P^2$} 
  (150:Rs) node(Z2) {$\P^2$} 
  (270:Rs) node(Z3) {$\P^1$};
 \path[commutative diagrams/.cd, every arrow, every label]
 (X) edge[swap] node {$f_1$} (Y1)
 (X) edge[swap] node {$f_2$} (Y2)
 (X) edge node {$f_3$} (Y3)
 (X) edge node {$\varphi_1$} (Z1)
 (X) edge[swap] node {$\varphi_2$} (Z2)
 (X) edge[swap] node {$\varphi_3$} (Z3)
 (Y1) edge node {$g_{11}$} (Z1)
 (Y1) edge[swap] node {$g_{12}$} (Z2)
 (Y2) edge node {$g_{22}=\pr_1$} (Z2)
 (Y2) edge[swap] node {$g_{23}=\pr_2$} (Z3)
 (Y3) edge node {$g_{33}=\pr_2$} (Z3)
 (Y3) edge[swap] node {$g_{31}=\pr_1$} (Z1);
\end{tikzpicture}
\]
\begin{proof}
For $Y_1 := W$ (No.\ 2-32) and $Y_2 := \P^2 \times \P^1$ (No.\ 2-34), 
we get the above commutative diagram except for $f_3, g_{31}, g_{33}$. 
It follows from 
Lemma \ref{l-ele-tf-K-relation}  and Proposition \ref{p-pic2-Pic} that 
\[
-2K_X \sim -f_1^*K_W -f_2^*K_{\P^2 \times \P^1} - \varphi_2^*K_{B_{Z_2}}
\]
\[
= (2H_1 + 2H_2) + (3H_2 + 2H_3) - 3H_2 = 2(H_1+H_2+H_3). 
\]
Thus (3) holds. 
Recall that $H_1^3 =0$ and  $H_1^2 \cdot H_2 = \ell_{g_{12}} \cdot H_2 = 1$ for 
an extremal rational curve $\ell_{g_{12}}$ of $g_{12}$ 
(Proposition \ref{p-pic2-Pic}). 
It holds that 
\[
H_1^2 \cdot H_3 \overset{{\rm (3)}}{=} H_1^2 \cdot (-K_X -H_1 -H_2) = 2 + 0 -1=1. 
\]
For $\varphi_1 : X \to \P^2 =:Z_1$ and $\varphi_3 : X \to \P^1 =:Z_3$, 
the induced morphism $\varphi_1 \times \varphi_3 : X \to Z_1 \times Z_3 = \P^2 \times \P^1$ is birational 
by 
$1 = H_1^2 \cdot H_3 = \deg (\varphi_1 \times \varphi_3) \times (\MO_{\P^2 \times \P^1}(1, 0)^2 \cdot \MO_{\P^2 \times \P^1}(0, 1))$. 
Hence $\varphi_1 \times \varphi_3$ is a contraction. 
By $\rho(X)=3 > 2 =\rho(\P^2 \times \P^1)= \rho(Z_1 \times Z_3)$, 
there exists a curve $C$ on $X$ contracted by $\varphi_1 \times \varphi_3 : X \to Z_1 \times Z_3$. 
Hence (1) holds (Lemma \ref{l-exactly3}) and we obtain the above commutative diagram. 
Moreover, (2) holds by Theorem \ref{t-ele-tr-P2}. 

Let us show (4). 
Note that 
$W$ is a Fano threefold of index $2$ and 
the blowup centre $B_1$ of $f_1$ is an elliptic curve with $-K_W \cdot B_1=12$ by (2). 
Let  $D$ be a Cartier divisor satisfying  $-K_W \sim 2D$. 
It is enough to check that (i) and (ii) of  Lemma \ref{l-dP-fib-CI} holds. 
Lemma \ref{l-dP-fib-CI}(i) holds by 
$(-K_W) \cdot D^2 = (-K_W)^3/4 =12 = (-K_W) \cdot B_1$. 
By (3) and $f_1^*(2D) \sim -f_1^*K_W \sim f_1^*(g_{11}^*\MO_{\P^2}(2) +g_{12}^*\MO_{\P^2}(2)) = 2H_1 + 2H_2$, 
we have $-K_X \sim H_1 +H_2 +  F \sim f_1^*D +F$ for a fibre $F$ of $\varphi_3 : X \to \P^1$. 
Hence $(f_1)_*F \sim (f_1)_*(-K_X -f_1^*D) = -K_W -D \sim D$. 
Thus  Lemma \ref{l-dP-fib-CI}(ii) holds, which completes the proof of (4). 
\qedhere


\end{proof}

\begin{prop}[No.\  {\hyperref[table-3-8]{3-8}}]\label{p-pic3-8}
Let $X$ be a Fano threefold with $\rho(X)=3$ and $(-K_X)^3=24$. 
Assume that $X$ has no conic bundle structure over $\P^2$ 
of type 2-32-vs-2-34. 
Then the following hold. 
\begin{enumerate}
\item $X$ has exactly three extremal rays. 
In what follows, we use Notation \ref{n-exactly3}. 
\item The contractions of the extremal faces are as in the following diagram. 
\begin{enumerate}
\item $f_1$ is of type $C_1$, $\Delta_{f_1} \in |\tau^*\MO_{\P^2}(3)|$. 
\item $f_2$ is of type $E_1$, $p_a(B_2)= 0$, $-K_{Y_2} \cdot B_2 = 2$, and $B_2$ is a smooth fibre of $g_{22}$. 
\item $f_3$ is of type $E_1$, $p_a(B_3)= 0$, $-K_{Y_3} \cdot B_3 = 14$. 
\item $\deg B_{Z_3} = 2$. 
\end{enumerate}
\item $-K_X \sim H_1 + H_2 + H_3$. 
\item 
$f_1 \times \varphi_3 : X \to \F_1 \times \P^2$ is a closed immersion and its image is a divisor linearly equivalent to 
$\pr_1^*\tau^*\MO_{\P^2}(1) \otimes \pr_2^*\MO_{\P^2}(2)$. 
\end{enumerate}
\end{prop}
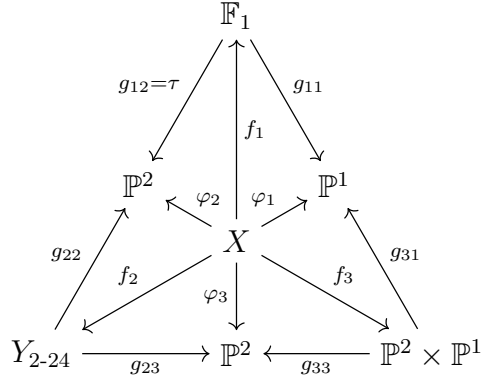
\[
\begin{tikzpicture}[commutative diagrams/every diagram,
    declare function={R=3;Rs=R*cos(60);}]
 \path 
  (0,0)  node(X) {$X$} 
  (90:R) node (Y1) {$\F_1$}
  (210:R) node (Y2) {$Y_{\text{2-24}}$}
  (-30:R) node (Y3) {$\P^2 \times \P^1$}  
  (30:Rs) node(Z1) {$\P^1$} 
  (150:Rs) node(Z2) {$\P^2$} 
  (270:Rs) node(Z3) {$\P^2$};
 \path[commutative diagrams/.cd, every arrow, every label]
 (X) edge[swap] node {$f_1$} (Y1)
 (X) edge[swap] node {$f_2$} (Y2)
 (X) edge node {$f_3$} (Y3)
 (X) edge node {$\varphi_1$} (Z1)
 (X) edge[swap] node {$\varphi_2$} (Z2)
 (X) edge[swap] node {$\varphi_3$} (Z3)
 (Y1) edge node {$g_{11}$} (Z1)
 (Y1) edge[swap] node {$g_{12} =\tau$} (Z2)
 (Y2) edge node {$g_{22}$} (Z2)
 (Y2) edge[swap] node {$g_{23}$} (Z3)
 (Y3) edge node {$g_{33} =\pr_1$} (Z3)
 (Y3) edge[swap] node {$g_{31}=\pr_2$} (Z1);
\end{tikzpicture}
\]

\begin{proof}
First of all, we prove that there exists a contraction $X \to Y_{\text{2-24}}$ to a Fano threefold 
$Y_{\text{2-24}}$ of No.\ 2-24. 
None of (III)-(V) of  (\ref{n-pic3-vol}) holds. 
Hence (I) or (II) of (\ref{n-pic3-vol}) holds. 
If (II) holds, then we get a contraction $X \to Y_{\text{2-24}}$ by Theorem \ref{t-F1-pic3}. 
If (I) holds, then Theorem \ref{t-ele-tr-P2} and our assumption imply that 
$X$ has a contraction $X \to Y_{\text{2-24}}$. 
This completes the proof of the existence of a contraction $X \to Y_{\text{2-24}}$. 

Fix a contraction $f_2 : X \to Y_{\text{2-24}}$. 
Then we get the commutative diagram consisting of $\varphi_2, \varphi_3, g_{22}, g_{23}$ as in the above diagram. 
Since the extremral rays of $Y_{\text{2-24}}$ are of type $C_1$ and $C_2$, 
we may assume that $g_{22}$ is of type $C_1$. 
Then the blowup centre $B_2$ of $f_2$ must be a smooth fibre of $g_{22}$. 
We then obtain the above commutative diagram except for $f_3, g_{31}, g_{33}$. 
Lemma \ref{l-F1cart-KX} and Proposition \ref{p-pic2-Pic} imply 
\[
-K_X \sim -f_2^*K_{Y_{\text{2-24}}} -H_2 +H_1 \sim (2H_2+H_3) -H_2 + H_1 = H_1 + H_2 +H_3. 
\]
Thus (3) holds. 
We have $H_3^3 =0$ and $H_2 \cdot H_3^2 = g_{22}^*\MO_{\P^2}(1) \cdot \ell_{g_{23}} =1$ 
for an extremal rational curve $\ell_{g_{23}}$ of $g_{23}$ (Proposition \ref{p-pic2-Pic}). 
It holds that 
\[
H_1 \cdot H_3^2 \overset{{\rm (3)}}{=} (-K_X -H_2 -H_3) \cdot H_3^2 = 2  -1+ 0 =1. 
\]
Hence $\varphi_3 \times \varphi_1 : X \to \P^2 \times \P^1$ is birational, 
i.e.,  a contraction. 
By $\rho(X) > \rho(\P^2 \times \P^1)$, there is a curve $C$ on $X$ contracted by $\varphi_3 \times \varphi_1$. 
Thus (1) holds (Lemma \ref{l-exactly3}) and we obtain the above commutative diagram. 
Since $B_2$ is a fibre of $g_{22}$, $B_2$ must be a subsection of $g_{23}$. 
Hence the square diagram consisting of $f_2, g_{23}, f_3, g_{33}$ 
is an elementary transform over $\P^2$ of type 2-24-vs-2-34. 
Hence $-K_{Y_3} \cdot B_3 = 14$ (Theorem \ref{t-ele-tr-P2}). It holds that $\Delta_{f_1} \in |\tau^*\MO_{\P^2}(3)|$ (Theorem \ref{t-F1-pic3}). 
Thus (2) holds.

Let us show (4). 
Since the extremal ray $R_{f_1}$ of $f_1$ is not contained in the extremal face $F_{\varphi_3}$ 
of $\varphi_3$, 
we obtain $R_{f_1} \cap F_{\varphi_3} = \{0\}$, which implies that 
$\varphi := f_1 \times \varphi_3 : X \to \F_1 \times \P^2$ is a finite morphism. 
For $X' := \varphi(X)$, we have the induced finite morphisms: 
\[
\varphi : X \xrightarrow{\psi} X' \hookrightarrow \F_1 \times \P^2. 
\]
Then $\psi : X \to X'$ is birational, because we have the following factorisation: 
\[
f_2 : X \xrightarrow{\psi} X' \to Y_{\text{2-24}} \subset \P^2 \times \P^2. 
\]
For $\tau \times {\rm id} :\F_1 \times \P^2 \to \P^2 \times \P^2$, 
we obtain $X' = (\tau\times {\rm id})^{-1}(Y_{\text{2-24}})$, 
because the blowup centre $\{ t\} \times \P^2$ of $\tau \times {\rm id}$ is not contained in 
$Y_{\text{2-24}}$ (recall that each contraction $Y_{\text{2-24}} \to \P^2$ is flat). 
Since $Y_{\text{2-24}}$ is a divisor on $\P^2 \times \P^2$ of bidegree $(1, 2)$, we obtain $X' \sim \pr_1^*\tau^*\MO_{\P^2}(1) \otimes \pr_2^*\MO_{\P^2}(2)$. 
For $H'_2 := \pr_1^*\tau^*\MO_{\P^2}(1)|_{X'}$ and $H'_3 := \pr_2^*\MO_{\P^2}(1)|_{X'}$, it follows from the adjunction formula that 
\[
\omega_{X'} \simeq (\omega_{\F_1 \times \P^2} \otimes \MO_{\F_1 \times \P^2}(X'))|_{X'} 
\]
\[
\simeq \MO_{X'}( (-3H'_2  +\Gamma_{X'} -3H'_3) + (H'_2 +2H'_3)) 
=\MO_{X'}(-2H'_2 -H'_3 +\Gamma_{X'}), 
\]
where $\Gamma_{X'} \subset X'$ denotes the pullback of the $(-1)$-curve $\Gamma$ on $\F_1$. 
Note that $\Gamma_X := \psi^*\Gamma_{X'}$ is the exceptional divisor of 
the blowup $f_2 : X \to Y_{\text{2-24}}$. 
Hence 
\[
-K_X \sim -f_2^*K_{Y_{\text{2-24}}} -\Gamma_{X} \sim 2H_2 +H_3  -\Gamma_X. 
\]
For the conductor $C$ of the normalisation $\psi :X \to X'$, 
we obtain $C\equiv K_X -\psi^*\omega_{X'} \sim 0$, and hence $C=0$. 
Thus $\psi$ is an isomorphism. 
Hence (4) holds. 
\end{proof}

\begin{lem}\label{l-disjoint-cont}
Let $X$ be a Fano threefold. 
Let $f_1 : X \to Y_1$ and $f_2: X \to Y_2$ be birational contractions of extremal rays. 
Assume that 
$\Ex(f_1) \cap \Ex(f_2) = \emptyset$ and 
one of $f_1(\Ex(f_1))$ and $f_2(\Ex(f_2))$ is a point. 
Then there exists a commutative diagram 
\[
\begin{tikzcd}
    & X \arrow[ld, "f_1"'] \arrow[rd, "f_2"] \arrow[dd, "\varphi"]\\
    Y_1 \arrow[rd, "g_1"'] & & Y_2\arrow[ld, "g_2"]\\
    & Z
\end{tikzcd}
\]
such that 
\begin{enumerate}
\item $Z$ is a projective normal threefold, 
\item $\varphi : X \to Z$ is a bitational contraction such that 
$\Ex(\varphi) = \Ex(f_1) \amalg \Ex(f_2)$, and 
\item for a curve $C$ on $X$, $\varphi(C)$ is a point if and only if $f_1(C)$ is a point or $f_2(C)$ is a point. 
\end{enumerate}
\end{lem}

\begin{proof}
Possibly after permuting $Y_1$ and $Y_2$, we may assume that $f_1(\Ex(f_1))$ is a point. 
Fix an ample Cartier divisor $A_{Y_2}$ on $Y_2$. 
Set $A_X := f_2^*A_{Y_2}$ and $E_1 := \Ex(f_1)$. 
Let $\lambda \in \Q_{>0}$ be the nef threshold: $A_X + \lambda E_1$, i.e., 
we define $\lambda$ as the largest rational number such that $A_X + \lambda E_1$ is nef. 
Since $f_1(E_1)$ is a point, we get $(A_X +\lambda E_1)|_{E_1} \equiv 0$, 
which implies  $(A_X +\lambda E_1)|_{E_1} \sim_{\Q} 0$. 
By Keel's theorem \cite[Theorem 0.2]{Kee99}, $A_X +\lambda E_1$ is semi-ample, which induces the contraction $\varphi : X \to Z$ satisfying (1)-(3). 
Then it is clear that we have the above commutative digaram. 
\end{proof}

\begin{prop}[No.\  {\hyperref[table-3-9]{3-9}}]\label{p-pic3-9}
Let $X$ be a Fano threefold with $\rho(X)=3$ and $(-K_X)^3=26$. 
Assume that $X$ has a conic bundle structure over $\P^2$.
Then the following hold. 
\begin{enumerate}
\item $X$ has exactly four extremal rays. 
In what follows, we use Notation \ref{n-exactly3}. 
\item The contractions of the extremal faces are as in the following diagram. 
\begin{enumerate}
\item $f_1$ is of type $E_1$, $p_a(B_1)= 3$, $-K_{Y_1} \cdot B_1 = 20$, and the blowup centre 
$B_1$ of $f_1$ is contained in the section $f_1(\Ex(f_3))$ of $g_{12}$. 
\item $f_2$ is of type $E_1$, $p_a(B_2)=3 $, $-K_{Y_2} \cdot B_2 = 20$, and the blowup centre 
$B_2$  of $f_2$ is contained in the section $f_2(\Ex(f_4))$ of $g_{22}$. 
\item $f_3$ is of type $E_5$. 
\item $f_4$ is of type $E_5$.
\item $Z_1$ is the cone over the Veronese surface $S \subset \P^5$, and 
$\varphi_1 : X \to Z_1$ is a blowup along a disjoint union of the singular point and 
a smooth curve $B_{Z_1}$ of genus $3$ satisfying $-K_{Z_1} \cdot B_{Z_1} =20$. 
\item $\deg B_{Z_2} = 4$. 
\item $Z_3$ is the cone over the Veronese surface $S \subset \P^5$, and 
$\varphi_3 : X \to Z_3$ is a blowup along a disjoint union of the singular point and 
a smooth curve $B_{Z_3}$ of genus $3$ satisfying $-K_{Z_3} \cdot B_{Z_3} =20$. 
\end{enumerate}
 \item 
$X \simeq \Bl_C\,Y_{\text{2-36}}$, 
where $C$ is a smooth curve on a section $T$ of the $\P^1$-bundle $\pi : Y_{\text{2-36}} = \P_{\P^2}(\MO \oplus \MO(2)) \to \P^2$ such that $\pi(C)$ is a smooth quartic curve and 
$T$ is disjoint from the section $S$ with $\MO_{Y_{\text{2-36}}}(-S)|_S$ ample. 
\item $-K_X \sim H_1 -H_2 +H_3$. 
\end{enumerate}
\end{prop}
\[
\begin{tikzpicture}[commutative diagrams/every diagram,
    declare function={R=3.5;Rs=R*cos(45);}]
 \path 
  (0,0)  node(X) {$X$} 
  (45:R) node (Y1) {$\P_{\P^2}(\MO \oplus \MO(2))$}
  (135:R) node (Y2) {$\P_{\P^2}(\MO \oplus \MO(2))$}
  (-135:R) node (Y3) {$Y_3$}  
  (-45:R) node (Y4) {$Y_4$}  
  (0:Rs) node(Z1) {$Z_1$} 
  (90:Rs) node(Z2) {$\P^2$} 
  (180:Rs) node(Z3) {$Z_3$}
  (270:Rs) node(Z4) {$Z_4$};
 \path[commutative diagrams/.cd, every arrow, every label]
 (X) edge node {$f_1$} (Y1)
 (X) edge[swap] node {$f_2$} (Y2)
 (X) edge[swap] node {$f_3$} (Y3)
 (X) edge node {$f_4$} (Y4)
 (X) edge node {$\varphi_1$} (Z1)
 (X) edge[swap] node {$\varphi_2$} (Z2)
 (X) edge[swap] node {$\varphi_3$} (Z3)
 (X) edge node {$\varphi_4$} (Z4)
 (Y1) edge node {$g_{11}$} (Z1)
 (Y1) edge[swap] node {$g_{12}$} (Z2)
 (Y2) edge node {$g_{22}$} (Z2)
 (Y2) edge[swap] node {$g_{23}$} (Z3)
 (Y3) edge node {$g_{33}$} (Z3)
 (Y3) edge[swap] node {$g_{34}$} (Z4)
 (Y4) edge node {$g_{44}$} (Z4)
 (Y4) edge[swap] node {$g_{41}$} (Z1);
\end{tikzpicture}
\]

\begin{proof}
By Theorem \ref{t-ele-tr-P2}, $X$ has a conic bundle structure over $\P^2$ of type 2-36-vs-2-36. 
We then obtain the commutative diagram consisting of $f_1, g_{12}, f_2, g_{22}$, 
which is an elementary transform over $\P^2$. 
Moreover, (a)', (b)', and (f) hold: 
\begin{enumerate}
\item[(a)'] $f_1$ is of type $E_1$, $p_a(B_2)= 3$, and $-K_{Y_1} \cdot B_1 = 20$.
\item[(b)'] $f_2$ is of type $E_1$, $p_a(B_2)=3 $, and $-K_{Y_2} \cdot B_2 = 20$. 
\end{enumerate}
Let $g_{11} : \P_{\P^2}(\MO \oplus \MO(2)) \to Z_1$ and $g_{23} :  \P_{\P^2}(\MO \oplus \MO(2)) \to Z_3$ be the birational contrations. 
In particular, each of $Z_1$ and $Z_3$ is the cone over the Veronese surface.  
Set $\varphi_1 := g_{11} \circ f_1$ and $\varphi_3 := g_{23} \circ f_2$. 
Since $D_{Y_1} := \Ex(g_{11}) \simeq \P^2$ is disjoint from the blowup centre $B_1$ of $f_1$ 
(Lemma \ref{l-P2-2-36}), 
$\varphi_1 : X \to Z_1$ is a blowup along a disjoint union of the singular point and a 
smooth curve $B_{Z_1}$ of genus $3$ satisfying $-K_{Z_1} \cdot B_{Z_1} = 20$. 
Therefore, we obtain the contraction $f_4 : X \to Y_4$ of type $E_5$ such that $D_{Y_1}^X =\Ex(f_4)$, where $D_{Y_1}^X$ denotes the proper transform of $D_{Y_1}$ on $X$. 
By symmetry, we obtain another contraction $f_3 : X \to Y_3$ of type $E_5$
such that $D_{Y_2}^X =\Ex(f_3)$, where $D_{Y_2}^X$ denotes the proper transform of $D_{Y_2} := \Ex(g_{23})$ on $X$. 
Hence  (c)-(g) hold and we obtain the above commutative diagram except for $\varphi_4, g_{34}, g_{44}$. 

Note that $D^X_{Y_1}$ and $D^X_{Y_2}$ are distinct, because 
$D^X_{Y_1}$ and $D^X_{Y_2}$ intersect the different prime divisors lying over 
the smooth curve $B_{Z_2} \subset \P^2$ (note that $\varphi_2^{-1}(B_{Z_2})$ consists of two prime divisors $F_1, F_2$ and 
each of $D^X_{Y_1}$ and $D^X_{Y_2}$ is disjoint from $F_1 \cap F_2$). 
Therefore, the extremal rays of $f_3$ and $f_4$ are different. 
Since $f_3(D^X_{Y_2})$ and $f_4(D^X_{Y_1})$ are points, 
we obtain $D^X_{Y_1} \cap D^X_{Y_2} = \emptyset$. 
Then there exists the birational contraction $\varphi_4 : X \to Z_4$ 
to a projective normal threefold $Z_4$ such that 
$\Ex(\varphi_4) = D^X_{Y_1} \amalg D^X_{Y_2}$, and 
both $\varphi_4(D^X_{Y_1})$ and $\varphi_4(D^X_{Y_2})$ are points (Lemma \ref{l-disjoint-cont}). 
Then we get the above commutative diagram, and hence (1) holds. 
By construction, $B_1$ (resp. $B_2$) is contained in the image of $D^X_{Y_2} =\Ex(f_3)$ (resp. $D^X_{Y_1} =\Ex(f_4)$). 
Hence (a) and (b) hold, and hence (2) holds. 
Moreover, (3) follows from (a) and 
$S \cap T= \emptyset$, 
where $S := \Ex(g_{11}) = D_{Y_1}$ and $T := f_1(D^X_{Y_2}) = f_1(\Ex(f_3))$. 
Lemma \ref{l-ele-tf-K-relation}  and Proposition \ref{p-pic2-Pic} imply 
\[
-2K_X \sim -f_1^*K_{Y_1} - f_2^*K_{Y_2}-\varphi_2^*B_{Z_2} 
\sim (2H_1+H_2) + (H_2+2H_3) - 4H_2 \sim 2H_1 -2H_2 +2H_3. 
\]
Thus (4) holds. 
\qedhere 
\end{proof}

\begin{prop}[No.\  {\hyperref[table-3-10]{3-10}}]\label{p-pic3-10} 
Let $X$ be a Fano threefold with $\rho(X)=3$ and $(-K_X)^3=26$. 
Assume that $X$ has no conic bundle structure over $\P^2$.
Then the following hold. 
\begin{enumerate}
\item $X$ has exactly three extremal rays. 
In what follows, we use Notation \ref{n-exactly3}. 
\item The contractions of the extremal faces are as in the following diagram. 
\begin{enumerate}
\item $f_1$ is of type $C_1$, $\deg \Delta_{f_1} = (2, 2)$. 
\item $f_2$ is of type $E_1$, $p_a(B_2)= 0$, $-K_{Y_2} \cdot B_2 = 6$. 
\item $f_3$ is of type $E_1$, $p_a(B_3)= 0$, $-K_{Y_3} \cdot B_3 = 6$. 
\item $\varphi_3 : X \to Q$ is a blowup of $Q$ along a disjoint union of two conics. 
\end{enumerate}
\item $-K_X \sim H_1 + H_2 + H_3$. 
\end{enumerate}
\end{prop}
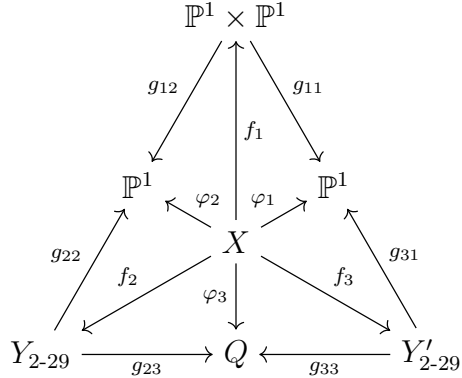
\[
\begin{tikzpicture}[commutative diagrams/every diagram,
    declare function={R=3;Rs=R*cos(60);}]
 \path 
  (0,0)  node(X) {$X$} 
  (90:R) node (Y1) {$\P^1 \times \P^1$}
  (210:R) node (Y2) {$Y_{\text{2-29}}$}
  (-30:R) node (Y3) {$Y'_{\text{2-29}}$}  
  (30:Rs) node(Z1) {$\P^1$} 
  (150:Rs) node(Z2) {$\P^1$} 
  (270:Rs) node(Z3) {$Q$};
 \path[commutative diagrams/.cd, every arrow, every label]
 (X) edge[swap] node {$f_1$} (Y1)
 (X) edge[swap] node {$f_2$} (Y2)
 (X) edge node {$f_3$} (Y3)
 (X) edge node {$\varphi_1$} (Z1)
 (X) edge[swap] node {$\varphi_2$} (Z2)
 (X) edge[swap] node {$\varphi_3$} (Z3)
 (Y1) edge node {$g_{11}=\pr_1$} (Z1)
 (Y1) edge[swap] node {$g_{12}=\pr_2$} (Z2)
 (Y2) edge node {$g_{22}$} (Z2)
 (Y2) edge[swap] node {$g_{23}$} (Z3)
 (Y3) edge node {$g_{33}$} (Z3)
 (Y3) edge[swap] node {$g_{31}$} (Z1);
\end{tikzpicture}
\]

\noindent 
The following proof works even if no fire of $f_1 : X \to \P^1 \times \P^1$ is smooth (cf. Theorem \ref{t-wild-cb}).

\begin{proof}
Only (V) of  (\ref{n-pic3-vol}) holds.  
All the assertions except for (3) follow from Proposition \ref{p-V-2dP}. 
Lemma \ref{l-K-disjoint-blowup} and Proposition \ref{p-pic2-Pic} imply 
\[
-K_X \sim -f_2^*K_{Y_{\text{2-29}}}-f_3^*K_{Y'_{\text{2-29}}} + \varphi_3^*K_{Q} 
\sim (H_2 +2H_3) + (H_1+ 2H_3)-3H_3 = H_1 + H_2 +H_3. 
\]
Thus (3) holds. 
\end{proof}

\begin{prop}[No.\ {\hyperref[table-3-11]{3-11}}]\label{p-pic3-11}
Let $X$ be a Fano threefold with $\rho(X)=3$ and $(-K_X)^3=28$. 
Assume $X$ has a conic bundle structure over $\P^2$ of type 2-34-vs-2-35.
Then the following hold. 
\begin{enumerate}
\item $X$ has exactly three extremal rays. 
In what follows, we use Notation \ref{n-exactly3}. 
\item The contractions of the extremal faces are as in the following diagram. 
\begin{enumerate}
\item $f_1$ is of type $E_1$, $p_a(B_1)= 0$, $-K_{Y_1} \cdot B_1 = 1$. 
\item $f_2$ is of type $E_1$, $p_a(B_2)= 1$, $-K_{Y_2} \cdot B_2 = 13$. 
\item $f_3$ is of type $E_1$, $p_a(B_3)= 1$, $-K_{Y_3} \cdot B_3 = 14$. 
\item $\deg B_{Z_3} =3$. 
\end{enumerate}
\item $-K_X \sim H_1 + H_2 + H_3$. 
\item 
The blowup centre $B_3$ of $f_3$ is an elliptic curve which is a complete intersection of two members of $|-\frac{1}{2}K_{V_7}|$. 
\end{enumerate}
\end{prop}
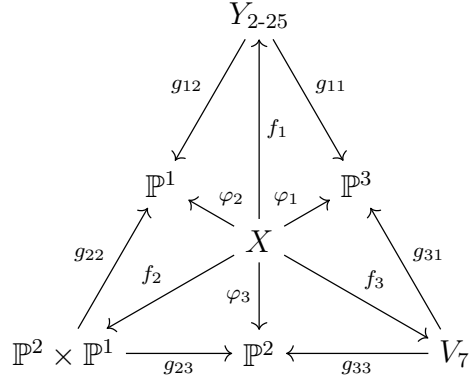
\[
\begin{tikzpicture}[commutative diagrams/every diagram,
    declare function={R=3;Rs=R*cos(60);}]
 \path 
  (0,0)  node(X) {$X$} 
  (90:R) node (Y1) {$Y_{\text{2-25}}$}
  (210:R) node (Y2) {$\P^2 \times \P^1$}
  (-30:R) node (Y3) {$V_7$}  
  (30:Rs) node(Z1) {$\P^3$} 
  (150:Rs) node(Z2) {$\P^1$} 
  (270:Rs) node(Z3) {$\P^2$};
 \path[commutative diagrams/.cd, every arrow, every label]
 (X) edge[swap] node {$f_1$} (Y1)
 (X) edge[swap] node {$f_2$} (Y2)
 (X) edge node {$f_3$} (Y3)
 (X) edge node {$\varphi_1$} (Z1)
 (X) edge[swap] node {$\varphi_2$} (Z2)
 (X) edge[swap] node {$\varphi_3$} (Z3)
 (Y1) edge node {$g_{11}$} (Z1)
 (Y1) edge[swap] node {$g_{12}$} (Z2)
 (Y2) edge node {$g_{22}=\pr_2$} (Z2)
 (Y2) edge[swap] node {$g_{23}=\pr_1$} (Z3)
 (Y3) edge node {$g_{33}$} (Z3)
 (Y3) edge[swap] node {$g_{31}$} (Z1);
\end{tikzpicture}
\]

\begin{proof}
Since $X$ has a conic bundle structure over $\P^2$ of type 2-34-vs-2-35, 
we get the above commutative diagram except for $f_1, g_{11}, g_{12}$. 
Since $\varphi_1 : X \to \P^3$ has a two-dimensional fibre, 
there exists a curve $C$ on $X$ contracted by $\varphi_1$ and $\varphi_2$. 
Thus (1) holds (Lemma \ref{l-exactly3}). 
Moreover, (b), (c), and (d) hold. 
Let $f_1 : X \to Y_1$ be the contraction of the remaining extremal ray. 
Note that $Y_1$ has two contractions 
$g_{11} : Y_1 \to \P^3$ and $g_{12}: Y_1 \to \P^1$. 
Since $\varphi_1$ is birational, $f_1$ is birational.

For the birational morphism $\varphi_1 : X \to \P^3$, 
we have $\Ex(\varphi_1)= D_X \cup E_X$ 
for $E_X := \Ex(f_3)$ and $D_X := f_3^{-1}D_{V_7}$, where  $D_{V_7} := \Ex(g_{31}) (\simeq \P^2)$. 
By $K_{V_7} =g_{31}^*K_{\P^3} +2D_{V_7}$, $-K_{V_7} \cdot B_3 = 14$, and 
$D_{V_7} \cdot B_3=1$ (Lemma \ref{l-P2-2-35}), we obtain 
$D_X \simeq \F_1$ and $-K_{\P^3} \cdot B_{\P^3} = 16$ for $B_{\P^3} :=g_{31}(B_3)$. 
Since $E_X$ is a $\P^1$-bundle over an elliptic curve by (c), 
we get $\Ex(f_1) = D_X$, and hence the extremal ray of $f_1$ must be of type $E_1$. 
In particular, $Y_1$ is a Fano threefold (Lemma \ref{l-nonFano-blowdown}) and $g_{11}$ is the blowup along 
the elliptic curve $B_{\P^3}$ of degree $4$. 
Then 
$Y_1$ is a Fano threefold of No.\ 2-25 (Subsection \ref{ss-table-pic2}). 
Thus we obtain the above commutative diagram. 
We get $(-K_{Y_{\text{2-25}}}) \cdot B_1 = 1$ by  
\[
28=(-K_X)^3  =(-K_{Y_{\text{2-25}}})^3 - 2(-K_{Y_{\text{2-25}}}) \cdot B_1 +2p_a(B_1)-2 
\]
\[
=32 - 2(-K_{Y_{\text{2-25}}}) \cdot B_1 +0 -2, 
\]
where the second equality follows from Lemma \ref{l-blowup-formula}. 
Hence (2) holds.  
Lemma \ref{l-ele-tf-K-relation} and Proposition \ref{p-pic2-Pic} imply 
\[
-2K_X \sim -f_2^*K_{\P^2 \times \P^1} - f_3^*K_{V_7} -\varphi_3^*B_{Z_3} 
\]
\[
\sim (2H_2 + 3H_3) + ( 2H_3 + 2H_1) - 3H_3  = 2(H_1+H_2+H_3). 
\]
Thus (3) holds. 

Let us show (4). 
Recall that 
$V_7$ is a Fano threefold of index $2$ and 
the blowup centre $B_3$ of $f_3$ is an elliptic curve with $-K_{V_7} \cdot B_3=14$ by (1). 
Let  $D$ be a Cartier divisor on $V_7$ satisfying  $-K_{V_7} \sim 2D$. 
It is enough to check that (i) and (ii) of  Lemma \ref{l-dP-fib-CI} holds. 
Lemma \ref{l-dP-fib-CI}(i) holds by 
$(-K_{V_7}) \cdot D^2 = (-K_{V_7})^3/4 =14 = (-K_{V_7}) \cdot B_3$. 
By (3) and $f_3^*(2D) \sim -f_3^*K_{V_7} \sim f_3^*(g_{31}^*\MO_{\P^3}(2) +g_{33}^*\MO_{\P^2}(2)) = 2H_1 + 2H_3$, 
we have $-K_X \sim H_1 +H_3 +  F \sim f_3^*D +F$ for a fibre $F$ of $\varphi_2 : X \to \P^1$. 
Hence $(f_3)_*F \sim -K_{V_7} -D \sim D$. 
Thus  Lemma \ref{l-dP-fib-CI}(ii) holds, which completes the proof of (4). 
\qedhere




\end{proof}

\begin{prop}[No.\ {\hyperref[table-3-12]{3-12}}]\label{p-pic3-12}
Let $X$ be a Fano threefold with $\rho(X)=3$ and $(-K_X)^3=28$. 
Assume $X$ has no conic bundle structure over $\P^2$ of type 2-34-vs-2-35. 
Then the following hold. 
\begin{enumerate}
\item $X$ has exactly three extremal rays. 
In what follows, we use Notation \ref{n-exactly3}. 
\item The contractions of the extremal faces are as in the following diagram. 
\begin{enumerate}
\item $f_1$ is of type $E_1$, $p_a(B_1)= 0$, $-K_{Y_1} \cdot B_1 = 4$. 
\item $f_2$ is of type $E_1$, $p_a(B_2)= 0$, $-K_{Y_2} \cdot B_2 = 12$. 
\item $f_3$ is of type $E_1$, $p_a(B_3)= 0$, $-K_{Y_3} \cdot B_3 = 12$. 
\item $\deg B_{Z_1} = 2$. 
\item $\varphi_2: X \to \P^3$ is a blowup along a disjoint union of 
a line and a rational cubic curve.  
\end{enumerate}
\item $-K_X \sim H_1 + H_2 + H_3$. 
\end{enumerate}
\end{prop}
\[
\begin{tikzpicture}[commutative diagrams/every diagram,
    declare function={R=3;Rs=R*cos(60);}]
 \path 
  (0,0)  node(X) {$X$} 
  (90:R) node (Y1) {$Y_{\text{2-27}}$}
  (210:R) node (Y2) {$Y_{\text{2-33}}$}
  (-30:R) node (Y3) {$\P^2 \times \P^1$}  
  (30:Rs) node(Z1) {$\P^2$} 
  (150:Rs) node(Z2) {$\P^3$} 
  (270:Rs) node(Z3) {$\P^1$};
 \path[commutative diagrams/.cd, every arrow, every label]
 (X) edge[swap] node {$f_1$} (Y1)
 (X) edge[swap] node {$f_2$} (Y2)
 (X) edge node {$f_3$} (Y3)
 (X) edge node {$\varphi_1$} (Z1)
 (X) edge[swap] node {$\varphi_2$} (Z2)
 (X) edge[swap] node {$\varphi_3$} (Z3)
 (Y1) edge node {$g_{11}$} (Z1)
 (Y1) edge[swap] node {$g_{12}$} (Z2)
 (Y2) edge node {$g_{22}$} (Z2)
 (Y2) edge[swap] node {$g_{23}$} (Z3)
 (Y3) edge node {$g_{33}=\pr_2$} (Z3)
 (Y3) edge[swap] node {$g_{31}=\pr_1$} (Z1);
\end{tikzpicture}
\]

\begin{proof}
Only (I) of  (\ref{n-pic3-vol}) holds. 
By Theorem \ref{t-ele-tr-P2}, $X$ has a conic bundle structure over $\P^2$ 
of type 2-27-vs-2-34. 
We then get the above commutative diagram except for $f_2, g_{22}, g_{23}$. 
By Lemma \ref{l-P2-2-27}, 
we have $D_1 \cap B_1 = \emptyset$ for $D_1 := \Ex(g_{12})$. 
Then we obtain the above commutative diagram except for $g_{23}$, 
and (a)-(e) hold. 

Lemma \ref{l-ele-tf-K-relation} and Proposition \ref{p-pic2-Pic} imply 
\[
-2K_X \sim -f_1^*K_{Y_{\text{2-27}}} - f_3^*K_{\P^2 \times \P^1} -\varphi_1^*B_1
\]
\[
\sim (H_1 +2H_2) + ( 3H_1 +2H_3) - 2H_1 = 2(H_1 + H_2 + H_3).  
\]
Thus (3) holds. 
Note that we have contractions $\psi : Y_{\text{2-33}} \to \P^1 =: \widetilde{Z}_3$ and $\widetilde{\varphi}_3 : X \xrightarrow{f_2} Y_{\text{2-33}} \xrightarrow{\psi} \P^1 = \widetilde{Z}_3$. 
Set $\widetilde{H}_3 := \widetilde{\varphi}_3^*\MO_{\P^3}(1)$. 
By Lemma \ref{l-K-disjoint-blowup} and Proposition \ref{p-pic2-Pic}, it holds that 
\[
-K_X \sim -f_1^*K_{Y_{\text{2-27}}} -f_2^*K_{Y_{\text{2-33}}} + \varphi_2^*K_{\P^3}
\]
\[
\sim (H_1 + 2H_2) + ( 3H_2 + \widetilde{H}_3) -4H_2 = H_1 + H_2 + \widetilde{H}_3. 
\]
Therefore, we get $H_3 \sim \widetilde{H}_3$, which implies $\varphi_3 = \widetilde{\varphi}_3$, (1), and (2). 
\qedhere



\end{proof}

\begin{prop}[No.\ {\hyperref[table-3-13]{3-13}}]\label{p-pic3-13}
Let $X$ be a Fano threefold with $\rho(X)=3$ and $(-K_X)^3=30$. 
Then the following hold. 
\begin{enumerate}
\item $X$ has exactly three extremal rays. 
In what follows, we use Notation \ref{n-exactly3}. 
\item The contractions of the extremal faces are as in the following diagram. 
\begin{enumerate}
\item $f_1$ is of type $E_1$, $p_a(B_1)= 0$, $-K_{Y_1} \cdot B_1 = 8$, 
and $B_1$ is of bidegree $(2, 2)$ with respect to $W \hookrightarrow \P^2 \times \P^2$. 
\item $f_2$ is of type $E_1$, $p_a(B_2)= 0$, $-K_{Y_2} \cdot B_2 = 8$ 
and $B_2$ is of bidegree $(2, 2)$ with respect to $W \hookrightarrow \P^2 \times \P^2$. 
\item $f_3$ is of type $E_1$, $p_a(B_3)= 0$, $-K_{Y_3} \cdot B_3 = 8$, 
and $B_3$ is of bidegree $(2, 2)$ with respect to $W \hookrightarrow \P^2 \times \P^2$. 
\item $\deg B_1 = 2$. 
\item $\deg B_2 = 2$. 
\item $\deg B_3 = 2$. 
\end{enumerate}
\item $-K_X \sim H_1 + H_2 + H_3$. 
\end{enumerate}
\end{prop}
\[
\begin{tikzpicture}[commutative diagrams/every diagram,
    declare function={R=3;Rs=R*cos(60);}]
 \path 
  (0,0)  node(X) {$X$} 
  (90:R) node (Y1) {$W$}
  (210:R) node (Y2) {$W$}
  (-30:R) node (Y3) {$W$}  
  (30:Rs) node(Z1) {$\P^2$} 
  (150:Rs) node(Z2) {$\P^2$} 
  (270:Rs) node(Z3) {$\P^2$};
 \path[commutative diagrams/.cd, every arrow, every label]
 (X) edge[swap] node {$f_1$} (Y1)
 (X) edge[swap] node {$f_2$} (Y2)
 (X) edge node {$f_3$} (Y3)
 (X) edge node {$\varphi_1$} (Z1)
 (X) edge[swap] node {$\varphi_2$} (Z2)
 (X) edge[swap] node {$\varphi_3$} (Z3)
 (Y1) edge node {$g_{11}$} (Z1)
 (Y1) edge[swap] node {$g_{12}$} (Z2)
 (Y2) edge node {$g_{22}$} (Z2)
 (Y2) edge[swap] node {$g_{23}$} (Z3)
 (Y3) edge node {$g_{33}$} (Z3)
 (Y3) edge[swap] node {$g_{31}$} (Z1);
\end{tikzpicture}
\]

\begin{proof}
Only (I) of  (\ref{n-pic3-vol}) holds. 
By Theorem \ref{t-ele-tr-P2}, $X$ has a conic bundle structure over $\P^2$ 
of type 2-32-vs-2-32. 
We then get the above commutative diagram except for $f_3, g_{31}, g_{33}$. 
For $f_2 : X \to W =: Y_2$, 
the blowup centre $B_2 \subset Y_2$ must be a subsection of $\varphi_3$, because (II) of  (\ref{n-pic3-vol}) does not hold. 
By Theorem \ref{t-ele-tr-P2}, we get another conic bundle structure over $\P^2$ 
of type 2-32-vs-2-32 consisting of $f_2, g_{23}, f_3, g_{33}$. 
Let $g_{34} : Y_3=W \to \P^2=:Z_4$ be the contraction of the extremal ray not corresponding to $g_{33}$. 
For the composition $\varphi_4 : X \xrightarrow{f_3} Y_3 =W \xrightarrow{g_{34}} Z_4 = \P^2$, 
we set  $H_4 := \varphi_4^*\MO_{\P^2}(1) = f_3^*g_{34}^*\MO_{\P^2}(1)$. 
In order to show (1), it is enough to prove that $H_1 \sim H_4$. 
By Lemma \ref{l-ele-tf-K-relation} and Proposition \ref{p-pic2-Pic}, we get 
\[
-2K_X \sim -f_1^*K_{Y_1} -f_2^*K_{Y_2} - \varphi_2^*B_{Z_2} \sim (2H_1+2H_2) +(2H_2+2H_3) -2H_2, 
\]
which implies $-K_X \sim H_1 + H_2 + H_3$. 
Similarly, we get $-K_X \sim H_2 + H_3+H_4$. Thus $H_1 \sim H_4$. 
Then (1) and (3) hold, and we get the above commutative diagram. 
By Theorem \ref{t-ele-tr-P2},  
(2) holds except for the assertion on the bidegree of each $B_i$. 
Note that $B_1$ is  of bidegree $(2, 2)$ by (d) and (e). 
Similarly, also $B_2$ and $B_3$ are of bidegree $(2. 2)$. Hence (2) holds. 
\qedhere




\end{proof}

\begin{lem}\label{l-P3-plane-curve-pt}
Take a point $P \in \P^3$ and a smooth curve $C$ on $\P^3$ such that 
$P \not\in C$, $\deg C \geq 2$, and $C$ is contained in a plane $V$ on $\P^3$. 
Assume that the blowup $X := \Bl_{P \amalg C}\,\P^3$ is Fano. 
Then $P \not\in V$. 
\end{lem}

\begin{proof}
Suppose $P \in V$. 
Let us derive a contradiction. 
Fix a line $L$ on $V =\P^2$ passing through $P$. 
Let $\sigma :X = \Bl_{P \amalg C}\,\P^3 \to \P^3$ be the induced blowup. 
Set $E_P$ and $E_C$ to be the $\sigma$-exceptional prime divisors lying over $P$ and $C$, respectively. 
For the proper transform $L_X$ of $L$ on $X$, it holds that 
\[
0> K_X \cdot L_X =(\sigma^*K_{\P^3} +E_C +2E_P) \cdot L_X \geq -4 +\deg C +2 \geq 0, 
\]
which is absurd. 
\end{proof}

\begin{prop}[No.\ {\hyperref[table-3-14]{3-14}}]\label{p-pic3-14}
Let $X$ be a Fano threefold with $\rho(X)=3$ and $(-K_X)^3=32$. 
Assume that $X$ has a conic bundle structure over $\P^2$ of type 2-35-vs-2-36. 
Then the following hold. 
\begin{enumerate}
\item $X$ has exactly four extremal rays. 
In what follows, we use Notation \ref{n-exactly3}. 
\item The contractions of the extremal faces are as in the following diagram. 
\begin{enumerate}
\item $f_1$ is of type $E_1$, $p_a(B_1)= 1$, $-K_{Y_1} \cdot B_1 = 12$. 
\item $f_2$ is of type $E_1$, $p_a(B_2)= 1$, $-K_{Y_2} \cdot B_2 = 15$. 
\item $f_3$ is of type $E_5$. 
\item $f_4$ is of type $E_2$.
\item $\varphi_1$ is a blowup along a disjoint union of a smooth plane cubic curve $C$ and a point $P$, where $P$ is not contained in the plane containing $C$.
\item $\deg B_{Z_2} = 3$. 
\item $Z_3$ is the cone over the Veronese surface, and 
$\varphi_3$ is a blowup of $Z_3$ along a disjoint union of the singular point 
and an elliptic curve $B_{Z_3}$ satisfying $-K_{Z_3} \cdot B_{Z_3} = 15$. 
\item $\varphi_4 : X \to Z_4$ is a birational morphism such that 
$\Ex(\varphi_4) = \Ex(f_3) \amalg \Ex(f_4)$. 
\end{enumerate}
\item $-K_X \sim H_1+H_3$. 
\end{enumerate}
\end{prop}
\[
\begin{tikzpicture}[commutative diagrams/every diagram,
    declare function={R=3.5;Rs=R*cos(45);}]
 \path 
  (0,0)  node(X) {$X$} 
  (45:R) node (Y1) {$V_7$}
  (135:R) node (Y2) {$\P_{\P^2}(\MO \oplus \MO(2))$}
  (-135:R) node (Y3) {$Y_3$}  
  (-45:R) node (Y4) {$Y_{\text{2-28}}$}  
  (0:Rs) node(Z1) {$\P^3$} 
  (90:Rs) node(Z2) {$\P^2$} 
  (180:Rs) node(Z3) {$Z_3$}
  (270:Rs) node(Z4) {$Z_4$};
 \path[commutative diagrams/.cd, every arrow, every label]
 (X) edge node {$f_1$} (Y1)
 (X) edge[swap] node {$f_2$} (Y2)
 (X) edge[swap] node {$f_3$} (Y3)
 (X) edge node {$f_4$} (Y4)
 (X) edge node {$\varphi_1$} (Z1)
 (X) edge[swap] node {$\varphi_2$} (Z2)
 (X) edge[swap] node {$\varphi_3$} (Z3)
 (X) edge node {$\varphi_4$} (Z4)
 (Y1) edge node {$g_{11}$} (Z1)
 (Y1) edge[swap] node {$g_{12}$} (Z2)
 (Y2) edge node {$g_{22}$} (Z2)
 (Y2) edge[swap] node {$g_{23}$} (Z3)
 (Y3) edge node {$g_{33}$} (Z3)
 (Y3) edge[swap] node {$g_{34}$} (Z4)
 (Y4) edge node {$g_{44}$} (Z4)
 (Y4) edge[swap] node {$g_{41}$} (Z1);
\end{tikzpicture}
\]

\begin{proof}
Since $X$ has a conic bundle structure over $\P^2$ of type 2-35-vs-2-36, 
we obtain the above commutative diagram except for the lower half of it, i.e., 
$g_{33}, f_3, \varphi_4, f_4, g_{41}, g_{34}, g_{44}$. 
Moreover, (a), (b), and (f) hold (Theorem \ref{t-ele-tr-P2}). 
Note that $Z_3$ is the cone over the Veronese surface, and the blowup centre $B_2$ of $f_2$ 
is disjoint from $D_{Y_2} := \Ex(g_{23})$ (Lemma \ref{l-P2-2-36}). 
In particular, we obtain the square diagram consisting of $f_2, g_{23}, f_3, g_{33}$ such that (c) and (g) hold. 

We now show that the blowup centre $B_1$ of $f_1$ is disjoint from $D_{Y_1} := \Ex(g_{11})$. 
By the first paragraph of the proof of Lemma \ref{l-P2-2-35}, 
we obtain $-K_{Y_1} \cdot B_{Y_1} = 4 \deg B_{Z_2} +2 D_{Y_1} \cdot B_1$. 
Then $12 \overset{{\rm (a)}}{=} -K_{Y_1} \cdot B_{Y_1} = 4 \deg B_{Z_2} +2 D_{Y_1} \cdot B_1 \overset{{\rm (f)}}{=} 4 \cdot 3 + 2D_{Y_1} \cdot B_1$, which implies $D_{Y_1} \cap B_1 = \emptyset$. 
Thus $\varphi_1$ is  a blowup along a disjoint union of a point $P$ and 
an elliptic curve $B_{\P^3} := g_{11}(B_{Y_1})$ of degree $3$. 
By the Riemann-Roch theorem for $B_{\P^3}$, 
$B_{\P^3}$ is contained in a plane on $\P^3$. 
Then the assertion (e) follows from Lemma \ref{l-P3-plane-curve-pt}. 
Hence we obtain the above commutative diagram except for $\varphi_4, g_{34}, g_{44}$. Moreover, (d) holds. 

Set $D^X_{Y_1}$ and $D^X_{Y_2}$ to be the proper transforms of $D_{Y_1}$ and $D_{Y_2}$ on $X$, respectively. 
Since $Y_3$ is singular and $Y_{\text{2-28}}$ is smooth, 
the extremal rays of $f_3$ and $f_4$ are different. 
As both $f_3(D^X_{Y_2})$ and $f_4(D^X_{Y_1})$ are points, 
we obtain $D^X_{Y_1} \cap D^X_{Y_2} = \emptyset$. 
Then there exists the birational contraction $\varphi_4 : X \to Z_4$ 
to a projective normal threefold $Z_4$ such that 
$\Ex(\varphi_4) = D^X_{Y_1} \amalg D^X_{Y_2}$ and 
both $\varphi_4(D^X_{Y_1})$ and $\varphi_4(D^X_{Y_2})$ are points (Lemma \ref{l-disjoint-cont}). 
Thus (1) and (2) hold and we get the above commutative diagram. 

Lemma \ref{l-ele-tf-K-relation} and Proposition \ref{p-pic2-Pic} imply  
\[
-2K_X \sim -f_1^*K_{Y_1} - f_2^*K_{Y_2} - \varphi_2^*B_{Z_2} \sim (2H_1+2H_2) +(H_2 + 2H_3) -3H_2 = 2H_1+2H_3. 
\]
Thus (3) holds. 
\qedhere 


\end{proof}

\begin{prop}[No.\ {\hyperref[table-3-15]{3-15}}]\label{p-pic3-15}
Let $X$ be a Fano threefold with $\rho(X)=3$ and $(-K_X)^3=32$. 
Assume that $X$ has no conic bundle structure over $\P^2$ of type 2-35-vs-2-36. 
Then the following hold. 
\begin{enumerate}
\item $X$ has exactly three extremal rays. 
In what follows, we use Notation \ref{n-exactly3}. 
\item The contractions of the extremal faces are as in the following diagram. 
\begin{enumerate}
\item $f_1$ is of type $E_1$, $p_a(B_1)= 0$, $-K_{Y_1} \cdot B_1 = 3$. 
\item $f_2$ is of type $E_1$, $p_a(B_2)= 0$, $-K_{Y_2} \cdot B_2 = 6$. 
\item $f_3$ is of type $E_1$, $p_a(B_3)= 0$, $-K_{Y_3} \cdot B_3 = 10$. 
\item $\varphi_2 : X \to Q$ is a blowup along a disjoint union of a line and a conic. 
\item $\deg B_{Z_3} = 2$. 
\end{enumerate}
\item $-K_X \sim H_1 + H_2 + H_3$. 
\end{enumerate}
\end{prop}
\[
\begin{tikzpicture}[commutative diagrams/every diagram,
    declare function={R=3;Rs=R*cos(60);}]
 \path 
  (0,0)  node(X) {$X$} 
  (90:R) node (Y1) {$Y_{\text{2-29}}$}
  (210:R) node (Y2) {$Y_{\text{2-31}}$}
  (-30:R) node (Y3) {$\P^2 \times \P^1$}  
  (30:Rs) node(Z1) {$\P^1$} 
  (150:Rs) node(Z2) {$Q$} 
  (270:Rs) node(Z3) {$\P^2$};
 \path[commutative diagrams/.cd, every arrow, every label]
 (X) edge[swap] node {$f_1$} (Y1)
 (X) edge[swap] node {$f_2$} (Y2)
 (X) edge node {$f_3$} (Y3)
 (X) edge node {$\varphi_1$} (Z1)
 (X) edge[swap] node {$\varphi_2$} (Z2)
 (X) edge[swap] node {$\varphi_3$} (Z3)
 (Y1) edge node {$g_{11}$} (Z1)
 (Y1) edge[swap] node {$g_{12}$} (Z2)
 (Y2) edge node {$g_{22}$} (Z2)
 (Y2) edge[swap] node {$g_{23}$} (Z3)
 (Y3) edge node {$g_{33}=\pr_1$} (Z3)
 (Y3) edge[swap] node {$g_{31}=\pr_2$} (Z1);
\end{tikzpicture}
\]
\begin{proof}
Only (I) of  (\ref{n-pic3-vol}) holds. 
By Theorem \ref{t-ele-tr-P2}, $X$ has a conic bundle structure over $\P^2$ 
of type 2-31-vs-2-34. 
We then get the above commutative diagram except for $f_1, g_{11}, g_{12}$. 
Moreover, (b), (c), and (e) hold. 
Since $D_{Y_{\text{2-31}}} := \Ex(g_{22})$ is disjoint from 
the blowup centre $B_2$ of $f_2 : X \to Y_{\text{2-31}}$ (Lemma \ref{l-P2-2-31}(3)), 
$\varphi_2 : X \to Q$ is a blowup along a disjoint union of two smooth curves. 
Since the blowup centre of $g_{22}$ is a line, 
we obtain the above commutative diagram except for $g_{11}$ such that (a)-(e) hold.

Lemma \ref{l-ele-tf-K-relation} and Proposition \ref{p-pic2-Pic} imply  
\[
-2K_X \sim -f_2^*K_{Y_{\text{2-31}}} - f_3^*K_{\P^1 \times \P^1} 
- \varphi_3^*B_{Z_3} 
\]
\[
\sim (2H_2 + H_3) +(3H_3 + 2H_1) -2H_3 = 2(H_1+H_2+H_3). 
\]
Thus (3) holds. 
For the induced contractions $\psi : Y_{\text{2-29}} \to \P^1$ and 
 $\widetilde{\varphi}_1 : X \xrightarrow{f_1} Y_{\text{2-29}} \xrightarrow{\psi} \P^1$, 
 we set $\widetilde{H}_1 := \widetilde{\varphi}_1^*\MO_{\P^1}(1)$. 
 Then it follows from 
 Lemma \ref{l-K-disjoint-blowup} and Proposition \ref{p-pic2-Pic} that 
\[
-K_X \sim -f_1^*K_{Y_{\text{2-29}}}-f_2^*K_{Y_{\text{2-31}}} +\varphi_2^*K_Q
\]
\[
\sim (\widetilde{H}_1 +2H_2) + (2H_2+H_3) -3H_2 = \widetilde{H}_1 + H_2 + H_3.  
\]
Hence we get $H_1 \sim \widetilde{H}_1$. 
Thus (1) and (2) hold and we get the above commutative diagram. 
\qedhere

    


\end{proof}

\begin{prop}[No.\ {\hyperref[table-3-16]{3-16}}]\label{p-pic3-16}
Let $X$ be a Fano threefold with $\rho(X)=3$ and $(-K_X)^3=34$. 
Then the following hold. 
\begin{enumerate}
\item $X$ has exactly three extremal rays. 
In what follows, we use Notation \ref{n-exactly3}. 
\item The contractions of the extremal faces are as in the following diagram. 
\begin{enumerate}
\item $f_1$ is of type $E_1$, $p_a(B_1)= 0$, $-K_{Y_1} \cdot B_1 = 1$. 
\item $f_2$ is of type $E_1$, $p_a(B_2)= 0$, $-K_{Y_2} \cdot B_2 = 6$. 
\item $f_3$ is of type $E_1$, $p_a(B_3)= 0$, $-K_{Y_3} \cdot B_3 = 10$, and 
the blowup centre $B_3$  of $f_3 : X \to V_7$ is the strict transform of a smooth cubic rational curve passing through the blowup centre of $g_{31}$. 
\item $\deg B_{Z_1} = 1$. 
\item $\deg B_{Z_2} = 2$. 
\end{enumerate}
\item $-K_X \sim H_1 + H_2 + H_3$. 
\end{enumerate}
\end{prop}
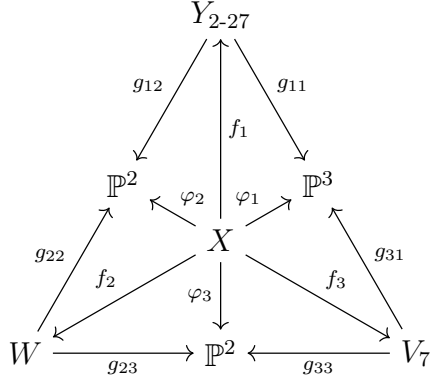
\[
\begin{tikzpicture}[commutative diagrams/every diagram,
    declare function={R=3;Rs=R*cos(60);}]
 \path 
  (0,0)  node(X) {$X$} 
  (90:R) node (Y1) {$Y_{\text{2-27}}$}
  (210:R) node (Y2) {$W$}
  (-30:R) node (Y3) {$V_7$}  
  (30:Rs) node(Z1) {$\P^3$} 
  (150:Rs) node(Z2) {$\P^2$} 
  (270:Rs) node(Z3) {$\P^2$};
 \path[commutative diagrams/.cd, every arrow, every label]
 (X) edge[swap] node {$f_1$} (Y1)
 (X) edge[swap] node {$f_2$} (Y2)
 (X) edge node {$f_3$} (Y3)
 (X) edge node {$\varphi_1$} (Z1)
 (X) edge[swap] node {$\varphi_2$} (Z2)
 (X) edge[swap] node {$\varphi_3$} (Z3)
 (Y1) edge node {$g_{11}$} (Z1)
 (Y1) edge[swap] node {$g_{12}$} (Z2)
 (Y2) edge node {$g_{22}$} (Z2)
 (Y2) edge[swap] node {$g_{23}$} (Z3)
 (Y3) edge node {$g_{33}$} (Z3)
 (Y3) edge[swap] node {$g_{31}$} (Z1);
\end{tikzpicture}
\]
\begin{proof}
First of all, 
we construct  the above commutative diagram except for $g_{31}$. 
Only (I) of  (\ref{n-pic3-vol}) holds, i.e., 
there exists a conic bundle structure $\varphi : X \to \P^2$. 
By Theorem \ref{t-ele-tr-P2}, 
$\varphi$ is of type 2-27-vs-2-32 or 2-32-vs-2-35. 
In any case, we get a contraction $f_2 : X \to W$ of type $E_1$ (recall that $W$ is of No.\ 2-32). 
Let $B_W$ be the blowup centre of $f_2$. 
Again by Theorem \ref{t-ele-tr-P2}, $-K_W \cdot B_W = 6$. 
Since (II) of  (\ref{n-pic3-vol}) does not hold, 
$B_W$ is a subsection for both contractions $g_{22} : W \to \P^2$ and $g_{23} : W \to \P^2$. 
By $-K_W \cdot B_W = 6$, $B_W$ is of bidegree $(1, 2)$ or $(2, 1)$. 
It follows from Theorem \ref{t-ele-tr-P2} that we get the above commutative diagram except for $g_{31}$. 
Moreover, (a), (b), (c)', (d), and (e) hold. 
\begin{enumerate}
\item[(c)'] $f_3$ is of type $E_1$, $p_a(B_3)= 0$, $-K_{Y_3} \cdot B_3 = 10$. 
\end{enumerate}

\medskip

Set $Y_4 := \P^3$ and let $g_{34} : V_7 \to Y_4 = \P^3$ be the contraction. 
Lemma \ref{l-ele-tf-K-relation}  and Proposition \ref{p-pic2-Pic} imply 
\[
-2K_X \sim  -f_1^*K_{Y_{\text{2-27}}} - f_2^*K_W -\varphi_{2}^*B_{Z_2} 
\sim (2H_1 +H_2) + (2H_2+2H_3) -H_2  = 2(H_1+H_2+H_3)
\]
\[
-2K_X \sim -f_2^*K_W - f_3^*K_{V_7} -\varphi_{3}^*B_{Z_3} 
\sim(2H_2+2H_3) + (2H_3+2H_4) -2H_3 = 2(H_2+H_3+H_4). 
\]
Hence 
$H_1 \sim H_4$. Thus (1) and (3) hold, and we get the above commutative diagram. 

By Lemma \ref{l-P2-2-35}, we have $D_{V_7} \cdot B_3 =1$ for $D_{V_7} := \Ex(g_{31})$. 
It follows from 
$K_{V_7} \sim g_{31}^*K_{\P^3} +2D_{V_7}$ and $-K_{Y_3} \cdot B_3 =10$ 
that $-K_{\P^3} \cdot B_{\P^3} =12$ for $B_{\P^3}:=g_{31}(B_3)$. 
Hence $B_{\P^3}$ is a smooth cubic rational curve passing through $g_{31}(\Ex(g_{31}))$. Thus (2) holds. 
\qedhere





\end{proof}

\begin{prop}[No.\ {\hyperref[table-3-17]{3-17}}]\label{p-pic3-17}
Let $X$ be a Fano threefold with $\rho(X)=3$ and $(-K_X)^3=36$. 
Assume that $X$ has a conic bundle structure. 
Then the following hold. 
\begin{enumerate}
\item $X$ has exactly three extremal rays. 
In what follows, we use Notation \ref{n-exactly3}. 
\item The contractions of the extremal faces are as in the following diagram. 
\begin{enumerate}
\item $f_1$ is of type $C_2$. 
\item $f_2$ is of type $E_1$, $p_a(B_2)= 0$, $-K_{Y_2} \cdot B_2 = 8$. 
\item $f_3$ is of type $E_1$, $p_a(B_3)= 0$, $-K_{Y_3} \cdot B_3 = 8$. 
\item $\deg B_{Z_3} = 2$
\end{enumerate}
\item $-K_X \sim H_1 + H_2 + 2H_3$. 
\item 
$\varphi_1 \times \varphi_2 \times \varphi_3 : X \to \P^1 \times \P^1 \times \P^2$ is a closed immersion and 
its image is of tridegree $(1, 1, 1)$. 
\end{enumerate}
\end{prop}
\[
\begin{tikzpicture}[commutative diagrams/every diagram,
    declare function={R=3;Rs=R*cos(60);}]
 \path 
  (0,0)  node(X) {$X$} 
  (90:R) node (Y1) {$\P^1_1 \times \P^1_2$}
  (210:R) node (Y2) {$\P^2 \times \P^1_1$}
  (-30:R) node (Y3) {$\P^2 \times \P^1_2$}  
  (30:Rs) node(Z1) {$\P^1_1$} 
  (150:Rs) node(Z2) {$\P^1_2$} 
  (270:Rs) node(Z3) {$\P^2$};
 \path[commutative diagrams/.cd, every arrow, every label]
 (X) edge[swap] node {$f_1$} (Y1)
 (X) edge[swap] node {$f_2$} (Y2)
 (X) edge node {$f_3$} (Y3)
 (X) edge node {$\varphi_1$} (Z1)
 (X) edge[swap] node {$\varphi_2$} (Z2)
 (X) edge[swap] node {$\varphi_3$} (Z3)
 (Y1) edge node {$g_{11}=\pr_1$} (Z1)
 (Y1) edge[swap] node {$g_{12}=\pr_2$} (Z2)
 (Y2) edge node {$g_{22}=\pr_2$} (Z2)
 (Y2) edge[swap] node {$g_{23}=\pr_1$} (Z3)
 (Y3) edge node {$g_{33}=\pr_1$} (Z3)
 (Y3) edge[swap] node {$g_{31}=\pr_2$} (Z1);
\end{tikzpicture}
\]

\begin{proof}
(I) or (IV) of  (\ref{n-pic3-vol}) holds. 
If (IV) holds, then $X$ had no birational contraction to $\P^2 \times \P^1$ (Lemma \ref{l-3-18}), which contradicts Theorem \ref{t-ele-tr-P2}. 
Thus (I) of  (\ref{n-pic3-vol}) holds. 
Then we get the above commutative diagram except for $f_1, g_{11}, g_{12}$ (Theorem \ref{t-ele-tr-P2}). 
Moreover, (b), (c), and (d) hold. 
Then (1) holds and we get the above commutative diagram (Lemma \ref{l-exactly3 2P^1}). 
Lemma \ref{l-ele-tf-K-relation} and Proposition \ref{p-pic2-Pic} imply 
\[
-2K_X \sim -f_2^*K_{\P^2 \times \P^1_1} -f_3^*K_{\P^2 \times \P^1_2} -\varphi^*_3B_{Z_3} 
\]
\[
\sim ( 2H_2 + 3H_3) + (3H_3 +2H_1) - 2H_3 = 2(H_1+H_2+2H_3). 
\]
Thus (3) holds. 

Let us show that $f_1$ is of type $C_2$, 
i.e., the bidegree $(d_1, d_2)$ of $\Delta_{f_1}$ is $(0, 0)$. 
We have $H_1 \cdot H_2 \cdot (2H_3) \overset{{\rm (3)}}{=} H_1 \cdot H_2 \cdot (-K_X) =2$ 
and $H_1 \cdot H_3^2 = 1$. 
Lemma \ref{l-P1P1-Delta-bideg} (applicable for $a_1 :=1, a_2 :=1$, and $D := 2H_3$) implies 
$(-K_X)^2 \cdot H_1 =  2 H_1 \cdot H_2 \cdot (2H_3) + H_1 \cdot (2H_3)^2 = 8$ and $d_2 = 8- (-K_X)^2 \cdot H_1 = 8=0$. 
By symmetry, $d_1=0$. Thus $f_1$ is of type $C_2$, and hence (2) holds. 
The assertion (4) follows from the same argument as that of Proposition \ref{p-pic3-3}(4).  
\qedhere

\end{proof}

\begin{prop}[No.\ {\hyperref[table-3-18]{3-18}}]\label{p-pic3-18}
Let $X$ be a Fano threefold with $\rho(X)=3$ and $(-K_X)^3=36$. 
Assume that $X$ has no conic bundle structure. 
Then the following hold. 
\begin{enumerate}
\item $X$ has exactly three extremal rays. 
In what follows, we use Notation \ref{n-exactly3}. 
\item The contractions of the extremal faces are as in the following diagram. 
\begin{enumerate}
\item $f_1$ is of type $E_1$, $p_a(B_1)= 0$, $-K_{Y_1} \cdot B_1 = 1$.  
\item $f_2$ is of type $E_1$, $p_a(B_2)= 0 $, $-K_{Y_2} \cdot B_2 = 4$. 
\item $f_3$ is of type $E_1$, $p_a(B_3)= 0$, $-K_{Y_3} \cdot B_3 = 8$. 
\item $\varphi_3 : X \to \P^3$ is a blowup of $\P^3$ along a disjoint union of a line and a conic. 
\end{enumerate}
\item $-K_X \sim H_1 + H_2 + H_3$. 
\end{enumerate}
\end{prop}
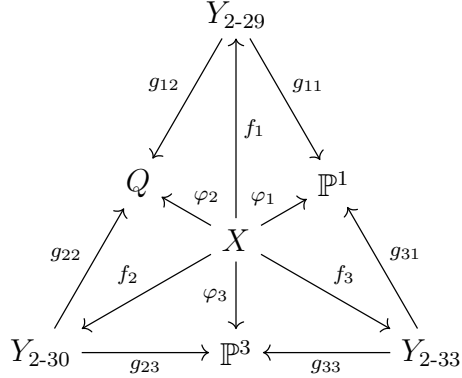
\[
\begin{tikzpicture}[commutative diagrams/every diagram,
    declare function={R=3;Rs=R*cos(60);}]
 \path 
  (0,0)  node(X) {$X$} 
  (90:R) node (Y1) {$Y_{\text{2-29}}$}
  (210:R) node (Y2) {$Y_{\text{2-30}}$}
  (-30:R) node (Y3) {$Y_{\text{2-33}}$}  
  (30:Rs) node(Z1) {$\P^1$} 
  (150:Rs) node(Z2) {$Q$} 
  (270:Rs) node(Z3) {$\P^3$};
 \path[commutative diagrams/.cd, every arrow, every label]
 (X) edge[swap] node {$f_1$} (Y1)
 (X) edge[swap] node {$f_2$} (Y2)
 (X) edge node {$f_3$} (Y3)
 (X) edge node {$\varphi_1$} (Z1)
 (X) edge[swap] node {$\varphi_2$} (Z2)
 (X) edge[swap] node {$\varphi_3$} (Z3)
 (Y1) edge node {$g_{11}$} (Z1)
 (Y1) edge[swap] node {$g_{12}$} (Z2)
 (Y2) edge node {$g_{22}$} (Z2)
 (Y2) edge[swap] node {$g_{23}$} (Z3)
 (Y3) edge node {$g_{33}$} (Z3)
 (Y3) edge[swap] node {$g_{31}$} (Z1);
\end{tikzpicture}
\]
\begin{proof}
Only (IV) of  (\ref{n-pic3-vol}) holds. 
By Lemma \ref{l-3-18}, (1) and (b)-(d) hold, we get the above commutative diagram, and $B_1 \simeq \P^1$.  
It follows from Lemma \ref{l-blowup-formula} that 
\[
36 = (-K_X)^3 = (-K_{Y_{\text{2-29}}})^3 -2 (-K_{Y_{\text{2-29}}}) \cdot B_1 +2p_a(B_1) -2 = 
40 -2 (-K_{Y_1}) \cdot B_1 +0 -2. 
\]
Hence $(-K_{Y_1}) \cdot B_1 =1$. Thus (2) holds. 
Lemma \ref{l-K-disjoint-blowup} and Proposition \ref{p-pic2-Pic} imply 
\[
-K_X \sim -f_2^*K_{Y_{\text{2-30}}} -f_3^*K_{Y_{\text{2-33}}} + \varphi_3^*K_{\P^3} 
\]
\[
\sim (H_2+2H_3) + (3H_3 + H_1) -4H_3 = H_1 + H_2 + H_3. 
\]
Thus (3) holds. 
\end{proof}

\begin{prop}[No.\ {\hyperref[table-3-19]{3-19}}]\label{p-pic3-19}
Let $X$ be a Fano threefold with $\rho(X)=3$ and $(-K_X)^3=38$. 
Assume that $X$ has a conic bundle structure over $\P^2$ of type 2-35-vs-2-35. 
Then the following hold. 
\begin{enumerate}
\item $X$ has exactly four extremal rays. 
In what follows, we use Notation \ref{n-exactly3}. 
\item The contractions of the extremal faces are as in the following diagram. 
\begin{enumerate}
\item $f_1$ is of type $E_2$. 
\item $f_2$ is of type $E_2$. 
\item $f_3$ is of type $E_1$, $p_a(B_3)= 0$, $-K_{Y_3} \cdot B_3 = 8$. 
\item $f_4$ is of type $E_1$, $p_a(B_4)= 0$, $-K_{Y_4} \cdot B_4 = 8$. 
\item $\varphi_1$ is a blowup along a disjoint union of a point and a conic. 
\item $\varphi_2$ is a blowup along a disjoint union of two points which are not collinear. 
\item $\varphi_3$ is a blowup along a disjoint union of a point and a conic. 
\item $\deg B_{Z_4}  = 2$. 
\end{enumerate}
\item $-K_X \sim H_1 + 2H_3 \sim H_2 +2H_4$. 
\end{enumerate}
\end{prop}
\[
\begin{tikzpicture}[commutative diagrams/every diagram,
    declare function={R=3;Rs=R*cos(45);}]
 \path 
  (0,0)  node(X) {$X$} 
  (45:R) node (Y1) {$Y_{\text{2-30}}$}
  (135:R) node (Y2) {$Y'_{\text{2-30}}$}
  (-135:R) node (Y3) {$V_7$}  
  (-45:R) node (Y4) {$V_7$}  
  (0:Rs) node(Z1) {$\P^3$} 
  (90:Rs) node(Z2) {$Q$} 
  (180:Rs) node(Z3) {$\P^3$}
  (270:Rs) node(Z4) {$\P^2$};
 \path[commutative diagrams/.cd, every arrow, every label]
 (X) edge node {$f_1$} (Y1)
 (X) edge[swap] node {$f_2$} (Y2)
 (X) edge[swap] node {$f_3$} (Y3)
 (X) edge node {$f_4$} (Y4)
 (X) edge node {$\varphi_1$} (Z1)
 (X) edge[swap] node {$\varphi_2$} (Z2)
 (X) edge[swap] node {$\varphi_3$} (Z3)
 (X) edge node {$\varphi_4$} (Z4)
 (Y1) edge node {$g_{11}$} (Z1)
 (Y1) edge[swap] node {$g_{12}$} (Z2)
 (Y2) edge node {$g_{22}$} (Z2)
 (Y2) edge[swap] node {$g_{23}$} (Z3)
 (Y3) edge node {$g_{33}$} (Z3)
 (Y3) edge[swap] node {$g_{34}$} (Z4)
 (Y4) edge node {$g_{44}$} (Z4)
 (Y4) edge[swap] node {$g_{41}$} (Z1);
\end{tikzpicture}
\]

\begin{proof}
By our assumption, 
we have an elementrary transform over $\P^2$ consisting 
of $f_3, g_{34}, f_4, g_{44}$ as in 
the above diagram. Moreover, 
(c), (d), and (h) hold (Theorem \ref{t-ele-tr-P2}). 
We then obtain contractions $\varphi_1, g_{41}, \varphi_3, g_{33}$. 
By Lemma \ref{l-P2-2-35}, the blowup centre $B_4$ of $f_4$ is disjoint from $D_{Y_4} := \Ex(g_{41})$. 
Thus  we obtain the square diagram consisting of $f_1, g_{11}, f_4, g_{41}$ 
such that (a) and (e) hold. 
By symmetry, we get the above commutative diagram except for $g_{22}$ 
satisfying  (a)-(e), (g), (h).

Proposition \ref{p-pic2-Pic} implies 
\begin{align*}   
-K_X &= -K_X -K_X +K_X\\
&\sim (-f_1^*K_{Y_{\text{2-30}}} -2\Ex (f_1)) + 
(-f_4^*K_{Y_4} -\Ex(f_4)) + (\varphi_1^*K_{\P^3} +2\Ex(f_1)+\Ex(f_4))\\
&= -f_1^*K_{Y_{\text{2-30}}}-f_4^*K_{Y_4} +\varphi_1^*K_{\P^3}\\ 
&\sim (2H_1 + H_2) + (2H_1 +2H_4) - 4H_1\\
&= H_2 + 2H_4. 
\end{align*}
Similarly, the following holds 
for the contraction $g'_{22} : Y'_{\text{2-30}} \to Q =:Y'_2$ and $H'_2 := f_2^*g'^*_{22}\MO_Q(1)$: 
\[
-K_X \sim H'_2 + 2H_4. 
\]
Thus $H_2 \sim H'_2$. 
Hence (1)  holds, and we obtain the above commutative diagram. 
Again by symmetry, we have $-K_X \sim H_1 + 2H_3$. Thus (3) holds.

It suffices to show (f). 
We see that  $\varphi_2 : X \to Q$ is a blowup along two points $P_1 \amalg P_2$. 
It suffices to show that $P$ and $Q$ are not collinear. 
Suppose that there exists a line $L$ on $Q$ passing through $P_1$ and $P_2$. 
Then we would get the following for the proper transform $L_X$ of $L$ on $X$: 
\[
K_X \cdot L_X =\varphi_2^*(K_Q + 2E_1 +2E_2) \cdot L_X = -3 + 2 + 2 >0, 
\]
which is absurd. 
\qedhere 



\end{proof}

\begin{prop}[No.\ {\hyperref[table-3-20]{3-20}}]\label{p-pic3-20}
Let $X$ be a Fano threefold with $\rho(X)=3$ and $(-K_X)^3=38$. 
Assume that $X$ has a conic bundle structure over $\P^2$ of type 2-31-vs-2-32. 
Then the following hold. 
\begin{enumerate}
\item $X$ has exactly three extremal rays. 
In what follows, we use Notation \ref{n-exactly3}. 
\item The contractions of the extremal faces are as in the following diagram. 
\begin{enumerate}
\item $f_1$ is of type $E_1$, $p_a(B_1)= 0$, $-K_{Y_1} \cdot B_1 = 4$. 
\item $f_2$ is of type $E_1$, $p_a(B_2)= 0$, $-K_{Y_2} \cdot B_2 = 3$. 
\item $f_3$ is of type $E_1$, $p_a(B_3)= 0$, $-K_{Y_3} \cdot B_3 = 3$. 
\item $\deg B_{Z_1} =1$. 
\item $\deg B_{Z_2} =1$. 
\item $\varphi_3 : X \to Q$ is a blowup along a disjoint of two lines. 
\end{enumerate}
\item $-K_X \sim H_1 + H_2 + H_3$. 
\end{enumerate}
\end{prop}
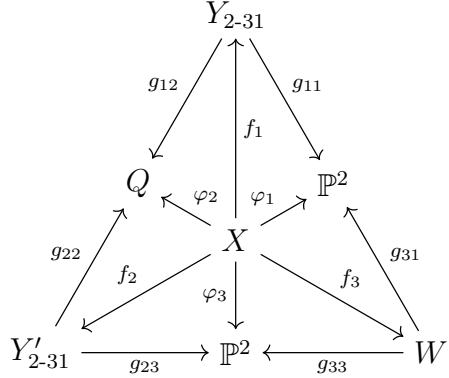
\[
\begin{tikzpicture}[commutative diagrams/every diagram,
    declare function={R=3;Rs=R*cos(60);}]
 \path 
  (0,0)  node(X) {$X$} 
  (90:R) node (Y1) {$W$}
  (210:R) node (Y2) {$Y_{\text{2-31}}$}
  (-30:R) node (Y3) {$Y'_{\text{2-31}}$}  
  (30:Rs) node(Z1) {$\P^2$} 
  (150:Rs) node(Z2) {$\P^2$} 
  (270:Rs) node(Z3) {$Q$};
 \path[commutative diagrams/.cd, every arrow, every label]
 (X) edge[swap] node {$f_1$} (Y1)
 (X) edge[swap] node {$f_2$} (Y2)
 (X) edge node {$f_3$} (Y3)
 (X) edge node {$\varphi_1$} (Z1)
 (X) edge[swap] node {$\varphi_2$} (Z2)
 (X) edge[swap] node {$\varphi_3$} (Z3)
 (Y1) edge node {$g_{11}$} (Z1)
 (Y1) edge[swap] node {$g_{12}$} (Z2)
 (Y2) edge node {$g_{22}$} (Z2)
 (Y2) edge[swap] node {$g_{23}$} (Z3)
 (Y3) edge node {$g_{33}$} (Z3)
 (Y3) edge[swap] node {$g_{31}$} (Z1);
\end{tikzpicture}
\]
\begin{proof}
Since $X$ has a conic bundle structure over $\P^2$ of type 2-31-vs-2-32, 
we obtain the above commutative diagram except for $f_3, g_{33}, g_{31}$ 
such that (a), (b), and (e) hold (Theorem \ref{t-ele-tr-P2}). 
By $-K_W \cdot B_1 = 4$ and $\deg B_{Z_2}=1$, 
$B_3$ is of bidegree $(1, 1)$. 
Hence $B_1$ is a regular subsection of $g_{11} : W \to \P^2$. 
We then obtain another conic bundle structure over $\P^2$ of type 2-31-vs-2-32 consisting of $f_1, g_{11}, f_3, g_{31}$ (Theorem \ref{t-ele-tr-P2}). 
Hence we get the above commutative diagram except for $g_{33}$. 
Lemma \ref{l-ele-tf-K-relation} and Proposition \ref{p-pic2-Pic} imply
\[
-2K_X \sim 
-f_1^*K_W -f_2^*K_{Y_{\text{2-31}}} +\varphi_2^*B_{Z_2} 
\]
\[
\sim (2H_1 +2H_2) + (H_2+2H_3) - H_2 = 2(H_1+H_2+H_3). 
\]
Thus (3) holds. 
Let $\widetilde{g}_{33} : Y'_{\text{2-31}} \to Q =:\widetilde{Y}_3$ be the contraction. 
For $\widetilde H_3 :=f_3^*\widetilde{g}^*_{33}\MO_Q(1)$, the same argument as above  implies 
$-K_X \sim H_1+H_2+\widetilde{H}_3$. 
Hence  $H_3 \sim \widetilde{H}_3$. 
Therefore, (1) and (a)-(e) hold, and we get the above commutative diagram. 
Finally, the blowup centre $B_2$ of $f_2 : X \to  Y_{\text{2-31}} $ is disjoint from $\Ex(g_{23})$ (Lemma \ref{l-P2-2-31}), 
and hence  (f) holds. 
\qedhere

\end{proof}

\begin{prop}[No.\ {\hyperref[table-3-21]{3-21}}]\label{p-pic3-21}
Let $X$ be a Fano threefold with $\rho(X)=3$ and $(-K_X)^3=38$. 
Assume that $X$ has no conic bundle structure over $\P^2$ of type 2-31-vs-2-32 
nor 2-35-vs-2-35. 
Then the following hold. 
\begin{enumerate}
\item $X$ has exactly three extremal rays. 
In what follows, we use Notation \ref{n-exactly3}. 
\item The contractions of the extremal faces are as in the following diagram. 
\begin{enumerate}
\item $f_1$ is of type $E_1$, $p_a(B_1)= 0$, $-K_{Y_1} \cdot B_1 = 7$, and $B_1$ is of bidegree $(1, 2)$.  
\item $f_2$ is of type $E_1$, $p_a(B_2)= 0$, $-K_{Y_2} \cdot B_2 = 0$. 
\item $f_3$ is of type $E_1$, $p_a(B_3)= 0$, $-K_{Y_3} \cdot B_3 = 0$.
\item $\deg B_{Z_2} = 1$. 
\item $\varphi_3$ is birational, 
$D:= \Ex(\varphi_3) = \Ex(f_2) = \Ex(f_3) \simeq \P^1 \times \P^1$, and $\varphi_3(D)$ is a point. 
\end{enumerate}
\item $-K_X \sim 2H_1 + 2H_2 + D$. 
\end{enumerate}
\end{prop}
\[
\begin{tikzpicture}[commutative diagrams/every diagram,
    declare function={R=3;Rs=R*cos(60);}]
 \path 
  (0,0)  node(X) {$X$} 
  (90:R) node (Y1) {$\P^2 \times \P^1$}
  (210:R) node (Y2) {$Y_{\text{non-Fano}}$}
  (-30:R) node (Y3) {$Y'_{\text{non-Fano}}$}  
  (30:Rs) node(Z1) {$\P^1$} 
  (150:Rs) node(Z2) {$\P^2$} 
  (270:Rs) node(Z3) {$Z_3$};
 \path[commutative diagrams/.cd, every arrow, every label]
 (X) edge[swap] node {$f_1$} (Y1)
 (X) edge[swap] node {$f_2$} (Y2)
 (X) edge node {$f_3$} (Y3)
 (X) edge node {$\varphi_1$} (Z1)
 (X) edge[swap] node {$\varphi_2$} (Z2)
 (X) edge[swap] node {$\varphi_3$} (Z3)
 (Y1) edge node {$g_{11}=\pr_2$} (Z1)
 (Y1) edge[swap] node {$g_{12}=\pr_1$} (Z2)
 (Y2) edge node {$g_{22}$} (Z2)
 (Y2) edge[swap] node {$g_{23}$} (Z3)
 (Y3) edge node {$g_{33}$} (Z3)
 (Y3) edge[swap] node {$g_{31}$} (Z1);
\end{tikzpicture}
\]

\begin{proof}
Only (I) of (\ref{n-pic3-vol}) holds. 
By Theorem \ref{t-ele-tr-P2} and our assumption, $X$ has a conic bundle structure over $\P^2$ of type 2-34-vs-non-Fano 
such that (a), (b), and (d) hold. 
Then we have the square diagram consisting of $f_1, g_{12}, f_2, g_{22}$. 
Since $f_2 : X \to Y_{\text{non-Fano}}$ is a contraction of type $E_1$ to a non-Fano threefold $Y_{\text{non-Fano}}$, 
we obtain the square diagram consisting of $f_2, g_{23}, f_3, g_{33}$ 
satisfying (c) and (e) (Proposition \ref{p-nonFano-flop}). 
Since $D=\Ex(\varphi_3) \simeq \P^1 \times \P^1$ is two-dimensional (Proposition \ref{p-nonFano-flop}), 
we can find a curve $C$ on $X$ contracted by $\varphi_1$ and $\varphi_3$. 
Hence (1) and (2) hold (Lemma \ref{l-exactly3}), and we obtain the above commutative diagram. 

Let us show (3). 
It follows from Lemma \ref{l-ele-tf-K-relation} that 
\[
-2K_X \sim -f_1^*K_{\P^2 \times \P^1} - f_2^*K_{Y_{\text{non-Fano}}} - \varphi_2^* B_{Z_2}. 
\]
This, together with $K_X \sim f_2^*K_{Y_{\text{non-Fano}}} +D$ 
and Proposition \ref{p-pic2-Pic}, implies 
\[
-K_X \sim  -f_1^*K_{\P^2 \times \P^1}+(K_X - f_2^*K_{Y_{\text{non-Fano}}})- \varphi_2^* B_{Z_2}
\]
\[
\sim (2H_1 + 3H_2)  +D - H_2 \sim 2H_1 +2H_2 +D. 
\]
Thus (3) holds. 
\qedhere



\end{proof}

\begin{prop}[No.\ {\hyperref[table-3-22]{3-22}}]\label{p-pic3-22}
Let $X$ be a Fano threefold with $\rho(X)=3$ and $(-K_X)^3=40$. 
Then the following hold. 
\begin{enumerate}
\item $X$ has exactly three extremal rays. 
In what follows, we use Notation \ref{n-exactly3}. 
\item The contractions of the extremal faces are as in the following diagram. 
\begin{enumerate}
\item $f_1$ is of type $E_1$, $p_a(B_1)= 0$, $-K_{Y_1} \cdot B_1 = 6$, 
and the blowup centre $B_1$ of $f_1$ is a conic on a plane $\P^2 \times \{ t\}$ for some closed point $t \in \P^1$. 
\item $f_2$ is of type $E_1$, $p_a(B_2)= 0$, $-K_{Y_2} \cdot B_2 = 10$. 
\item $f_3$ is of type $E_5$. 
\item $\deg B_{Z_2} =2$. 
\item $Z_3$ is the cone over the Veronese suraface $S \subset \P^5$ and 
$\varphi_3$ is a blowup along a disjoint union of the singular point and a 
smooth rational curve $C$ satisfying $-K_{Z_3} \cdot C= 10$. 
\end{enumerate}
\item $-K_X \sim H_1 + H_2 + H_3$. 
\end{enumerate}
\end{prop}
\[
\begin{tikzpicture}[commutative diagrams/every diagram,
    declare function={R=3;Rs=R*cos(60);}]
 \path 
  (0,0)  node(X) {$X$} 
  (90:R) node (Y1) {$\P^2 \times \P^1$}
  (210:R) node (Y2) {$\P_{\P^2}(\MO \oplus \MO(2))$}
  (-30:R) node (Y3) {$Y_3$}  
  (30:Rs) node(Z1) {$\P^1$} 
  (150:Rs) node(Z2) {$\P^2$} 
  (270:Rs) node(Z3) {$Z_3$};
 \path[commutative diagrams/.cd, every arrow, every label]
 (X) edge[swap] node {$f_1$} (Y1)
 (X) edge[swap] node {$f_2$} (Y2)
 (X) edge node {$f_3$} (Y3)
 (X) edge node {$\varphi_1$} (Z1)
 (X) edge[swap] node {$\varphi_2$} (Z2)
 (X) edge[swap] node {$\varphi_3$} (Z3)
 (Y1) edge node {$g_{11}=\pr_2$} (Z1)
 (Y1) edge[swap] node {$g_{12}=\pr_1$} (Z2)
 (Y2) edge node {$g_{22}$} (Z2)
 (Y2) edge[swap] node {$g_{23}$} (Z3)
 (Y3) edge node {$g_{33}$} (Z3)
 (Y3) edge[swap] node {$g_{31}$} (Z1);
\end{tikzpicture}
\]
\begin{proof}
Only (I) of  (\ref{n-pic3-vol}) holds. 
By Theorem \ref{t-ele-tr-P2}, $X$ has a conic bundle structure over $\P^2$ 
of type 2-34-vs-2-36. 
We then get the above commutative diagram except for $f_3, g_{31}, g_{33}$ such that 
(a)', (b), and (d) hold. 
\begin{enumerate}
\item[(a)'] $f_1$ is of type $E_1$, $p_a(B_1)= 0$, and $-K_{Y_1} \cdot B_1 = 6$. 
\end{enumerate}
Since $g_{23}$ has a two-dimensional fibre, 
so does $\varphi_3$. 
Then we can find a curve $C$ on $X$ contracted by $\varphi_1$ and $\varphi_3$. 
Hence (1) holds (Lemma \ref{l-exactly3}) and we obtain the above commutative diagram. 
Note that the blowup centre $B_2$ of $f_2$ is disjoint from $\Ex(g_{23})$ 
(Lemma \ref{l-P2-2-36}). 
Hence (c) and (e) hold. 
Lemma \ref{l-ele-tf-K-relation} and Proposition \ref{p-pic2-Pic} imply 
\[
-2K_X \sim -f_1^*K_{\P^2 \times \P^1} -f_2^*K_{\P_{\P^2}(\MO \oplus \MO(2))} -\varphi_2^*B_{Z_2} 
\]
\[
\sim (2H_1 + 3H_2) + (H_2 + 2H_3)-2H_2 = 2(H_1 + H_2 + H_3). 
\]
Thus (3) holds. 

It is enough to show (a). 
By the same argument as in the first paragraph of the proof of Lemma \ref{l-P2-2-34}, 
the following holds for a fibre $D:=\P^2 \times \{t\}$ of the projection $g_{11} : \P^2 \times \P^1 \to \P^1$ intersecting $B_1$: 
\[
2D \cdot B_1 = -K_{Y_1/\P^2} \cdot B_1 = -K_{Y_1} \cdot B_1 +K_{\P^2} \cdot B_{Z_2} \overset{{\rm (a)'(d)}}{=} 6-6=0. 
\]
Thus $B_1 \subset D$ and $B_1$ is a conic on the plane $D =\P^2 \times \{t\}$, 
because $B_{Z_2} = g_{12}(B_1)$ is a conic on $\P^2 =Z_2$. 
Thus (a) holds. 
\qedhere



\end{proof}

\begin{prop}[No.\ {\hyperref[table-3-23]{3-23}}]\label{p-pic3-23}
Let $X$ be a Fano threefold with $\rho(X)=3$ and $(-K_X)^3=42$. 
Assume that $X$ has a conic bundle structure over $\P^2$ of type 2-31-vs-2-35.  
Then the following hold. 
\begin{enumerate}
\item $X$ has exactly three extremal rays. 
In what follows, we use Notation \ref{n-exactly3}. 
\item The contractions of the extremal faces are as in the following diagram. 
\begin{enumerate}
\item $f_1$ is of type $E_1$, $p_a(B_1)= 0$, $-K_{Y_1} \cdot B_1 = 1$. 
\item $f_2$ is of type $E_1$, $p_a(B_2)= 0$, $-K_{Y_2} \cdot B_2 = 1$, and 
$B_2$ is a fibre of the induced $\P^1$-bundle $\Ex(g_{22}) \to g_{22}(\Ex(g_{22}))$. 
\item $f_3$ is of type $E_1$, $p_a(B_3)= 0$, $-K_{Y_3} \cdot B_3 = 6$, 
and $B_3$ is the proper transform of a conic passing through the blowup centre of $g_{31} : V_7 \to \P^3$. 
\item $\deg B_{Z_3} = 1$
\end{enumerate}
\item $-K_X \sim H_1 + H_2 + H_3$. 
\end{enumerate}
\end{prop}
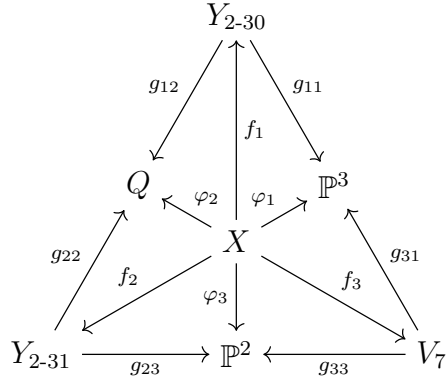
\[
\begin{tikzpicture}[commutative diagrams/every diagram,
    declare function={R=3;Rs=R*cos(60);}]
 \path 
  (0,0)  node(X) {$X$} 
  (90:R) node (Y1) {$Y_{\text{2-30}}$}
  (210:R) node (Y2) {$Y_{\text{2-31}}$}
  (-30:R) node (Y3) {$V_7$}  
  (30:Rs) node(Z1) {$\P^3$} 
  (150:Rs) node(Z2) {$Q$} 
  (270:Rs) node(Z3) {$\P^2$};
 \path[commutative diagrams/.cd, every arrow, every label]
 (X) edge[swap] node {$f_1$} (Y1)
 (X) edge[swap] node {$f_2$} (Y2)
 (X) edge node {$f_3$} (Y3)
 (X) edge node {$\varphi_1$} (Z1)
 (X) edge[swap] node {$\varphi_2$} (Z2)
 (X) edge[swap] node {$\varphi_3$} (Z3)
 (Y1) edge node {$g_{11}$} (Z1)
 (Y1) edge[swap] node {$g_{12}$} (Z2)
 (Y2) edge node {$g_{22}$} (Z2)
 (Y2) edge[swap] node {$g_{23}$} (Z3)
 (Y3) edge node {$g_{33}$} (Z3)
 (Y3) edge[swap] node {$g_{31}$} (Z1);
\end{tikzpicture}
\]
\begin{proof}
Since $X$ has a conic bundle structure over $\P^2$ of type 2-31-vs-2-35, 
there exists the above commutative diagram except for $f_1, g_{11}, g_{12}$ 
such that (b)', (c)', and (d) hold (Theorem \ref{t-ele-tr-P2}). 
\begin{enumerate}
\item[(b)'] $f_2$ is of type $E_1$, $p_a(B_2)= 0$, and 
$-K_{Y_2} \cdot B_2 = 1$. 
\item[(c)'] $f_3$ is of type $E_1$, $p_a(B_3)= 0$, and $-K_{Y_3} \cdot B_3 = 6$. 
\end{enumerate}
Then we get the above commutative diagram except for $g_{11}$ (Lemma \ref{l-P3-V7}). 
Let $\widetilde{g}_{11} : Y_{\text{2-30}} \to \P^3 =:\widetilde Y_1$ be the contraction and set $\widetilde H_1 := f_1^*\widetilde{g}_{11} \MO_{\P^3}(1)$. 
Lemma \ref{l-ele-tf-K-relation} and Proposition \ref{p-pic2-Pic} imply 
\[
-2K_X \sim -f_2^*K_{Y_{\text{2-31}}} - f_3^*K_{V_7}-\varphi_3^*K_{B_{Z_3}}
\]
\[
\sim (2H_2 + H_3) + (2H_3 +2H_1) - H_3 =2(H_1 + H_2 + H_3). 
\]
Thus (3) holds. 
In what follows, we set $E_{\psi} := \Ex(\psi)$ for a birational morphism $\psi$. 
Then (b) implies 
\[
K_X \sim f_2^*K_{Y_{\text{2-31}}} + E_{f_2} \sim f_2^*( g_{22}^*K_Q + E_{g_{22}}) +E_{f_2} \sim \varphi_2^* K_Q + E_{f_1} + 2E_{f_2}. 
\]
We have $K_X \sim f_1^*K_{Y_{\text{2-30}}} + E_{f_1}$. 
Hence it holds that 
\[
K_X \sim \varphi_2^* K_Q + E_{f_1} + 2E_{f_2} \sim \varphi_2^*K_Q + 
(K_X - f_1^*K_{Y_{\text{2-30}}}) +2(K_X -f_2^*K_{Y_{\text{2-31}}} ), 
\]
which implies 
\[
-2K_X \sim  - f_1^*K_{Y_{\text{2-30}}}-2f_2^*K_{Y_{\text{2-31}}}+ \varphi_2^*K_Q
\]
\[
\sim (2\widetilde{H}_1+H_2) +2 (2H_2 + H_3)-3H_2 = 2(\widetilde{H}_1 + H_2+H_3).
\]
Therefore, $H_1 \sim \widetilde{H}_1$. 
Then (1) holds, and we get the above commutative diagram.

It suffices to show (a)-(c). 
Lemma \ref{l-P2-2-31} and (b)' imply (b). 
Since the blowup centre $B_1$ of $f_1$ is the proper transform of the line 
$B_Q := g_{22}(\Ex(g_{22}))$ on $Q$ (Lemma \ref{l-P3-V7}), 
we obtain $p_a(B_1)=0$ and 
$K_{Y_{\text{2-30}}} \cdot B_1 = (g_{12}^*K_Q +2E_{g_{12}})\cdot B_1 
= K_Q \cdot B_Q + 2E_{g_{12}} \cdot B_1 =-3 +2 =-1$. 
Thus (a) holds. 
For the conic $B_{\P^3} := g_{11}(\Ex(g_{11}))$, 
$\varphi_1^{-1}(B_{\P^3}) = f_1^{-1}(g_{11}^{-1}(B_{\P^3}))$ is pure two-dimensional. 
Hence the blowup centre $B_3$ of $f_3 : X \to V_7$ coincides with the proper transform of the conic $B_{\P^3}$. 
Moreover, $B_{\P^3}$ passes through the blowup centre of $g_{31}$, 
because $-K_{\P^3} \cdot B_{\P^3} = 8 \neq 6 = -K_{Y_3} \cdot B_3$. 
Hence (c) holds. 
\qedhere


\end{proof}

\begin{prop}[No.\ {\hyperref[table-3-24]{3-24}}]\label{p-pic3-24}
Let $X$ be a Fano threefold with $\rho(X)=3$ and $(-K_X)^3=42$. 
Assume that $X$ has no conic bundle structure over $\P^2$ of type 2-31-vs-2-35.  Then the following hold. 
\begin{enumerate}
\item $X$ has exactly three extremal rays. 
In what follows, we use Notation \ref{n-exactly3}. 
\item The contractions of the extremal faces are as in the following diagram. 
\begin{enumerate}
\item $f_1$ is of type $C_2$. 
\item $f_2$ is of type $E_1$, $p_a(B_2)= 0$, $-K_{Y_2} \cdot B_2 = 2$, and $B_2$ is a fibre of the $\P^1$-bundle $g_{22}: W \to \P^2$.  
\item $f_3$ is of type $E_1$, $p_a(B_3)= 0$, $-K_{Y_3} \cdot B_3 = 5$. 
\item $\deg B_{Z_3}=1$. 
\end{enumerate}
\item $-K_X \sim H_1 + H_2 + 2H_3$. 
\end{enumerate}
\end{prop}
\[
\begin{tikzpicture}[commutative diagrams/every diagram,
    declare function={R=3;Rs=R*cos(60);}]
 \path 
  (0,0)  node(X) {$X$} 
  (90:R) node (Y1) {$\F_1$}
  (210:R) node (Y2) {$W$}
  (-30:R) node (Y3) {$\P^2 \times \P^1$}  
  (30:Rs) node(Z1) {$\P^1$} 
  (150:Rs) node(Z2) {$\P^2$} 
  (270:Rs) node(Z3) {$\P^2$};
 \path[commutative diagrams/.cd, every arrow, every label]
 (X) edge[swap] node {$f_1$} (Y1)
 (X) edge[swap] node {$f_2$} (Y2)
 (X) edge node {$f_3$} (Y3)
 (X) edge node {$\varphi_1$} (Z1)
 (X) edge[swap] node {$\varphi_2$} (Z2)
 (X) edge[swap] node {$\varphi_3$} (Z3)
 (Y1) edge node {$g_{11}$} (Z1)
 (Y1) edge[swap] node {$g_{12} =\tau$} (Z2)
 (Y2) edge node {$g_{22}$} (Z2)
 (Y2) edge[swap] node {$g_{23}$} (Z3)
 (Y3) edge node {$g_{33}=\pr_1$} (Z3)
 (Y3) edge[swap] node {$g_{31}=\pr_2$} (Z1);
\end{tikzpicture}
\]

\begin{proof}
Note that (I) or (II) of  (\ref{n-pic3-vol}) holds. 
In any case, 
Theorem \ref{t-ele-tr-P2} and Theorem \ref{t-F1-pic3} enable us to 
find a blowup $f_2 : X \to W$ along a smooth rational curve $B_2$ satisfying $-K_W \cdot B_2 =2$. 
Since the bidegree $(d_1, d_2)$ of $B_2$ satisfies $2d_1 +2d_2 = -K_W \cdot B_2 = 2$, 
we may assume that $d_1 =0$ and $d_2=1$. 
Then $B_2$ is a fibre of $g_{22} : W \hookrightarrow \P^2 \times \P^2 \xrightarrow{\pr_1} \P^2$ and a regular subsection of $g_{23}:  W \hookrightarrow \P^2 \times \P^2 \xrightarrow{\pr_2} \P^2$. 
Again by Theorem \ref{t-ele-tr-P2} and Theorem \ref{t-F1-pic3}, 
we obtain the above diagram except for $\varphi_1, g_{11}, g_{31}$ 
such that (a)-(d) hold. 
Let $g_{31} : \P^2 \times \P^1 \to \P^1$ be the contraction and take the composition $\varphi_1 : X \xrightarrow{f_3} \P^2 \times \P^1 \xrightarrow{g_{31}} \P^1$. 
Since the composition $\varphi_2 : X \xrightarrow{f_1}\F_1 \xrightarrow{g_{12}} \P^2$ has a two-dimensional fibre, 
there exists a curve $C$ on $X$ contracted by $\varphi_1$ and $\varphi_2$. 
Thus (1) holds (Lemma \ref{l-exactly3}) and we get the above commutative diagram. Moreover, (2) holds. 
Lemma \ref{l-ele-tf-K-relation} and Proposition \ref{p-pic2-Pic} imply 
\[
-2K_X \sim -f_2^*K_W -f_3^*K_{\P^2 \times \P^1} - \varphi_3^*B_{Z_3}
\]
\[
\sim (2H_2+2H_3) +(3H_3+2H_1) - H_3  =2(H_1 +H_2 +2H_3). 
\]
Thus (3) holds. 
\qedhere

%
\end{proof}

\begin{prop}[No.\ {\hyperref[table-3-25]{3-25}}]\label{p-pic3-25} 
Let $X$ be a Fano threefold with $\rho(X)=3$ and $(-K_X)^3=44$. 
Then the following hold. 
\begin{enumerate}
\item $X$ has exactly three extremal rays. 
In what follows, we use Notation \ref{n-exactly3}. 
\item The contractions of the extremal faces are as in the following diagram. 
\begin{enumerate}
\item $f_1$ is of type $C_2$. 
\item $f_2$ is of type $E_1$, $p_a(B_2)= 0$, $-K_{Y_2} \cdot B_2 = 4$. 
\item $f_3$ is of type $E_1$, $p_a(B_3)= 0$, $-K_{Y_3} \cdot B_3 = 4$. 
\item $\varphi_3 : X \to \P^3$ is a blowup of $\P^3$ along a disjoint union of two lines. 
\end{enumerate}
\item $-K_X \sim H_1 + H_2 + 2H_3$. 
\end{enumerate}
\end{prop}
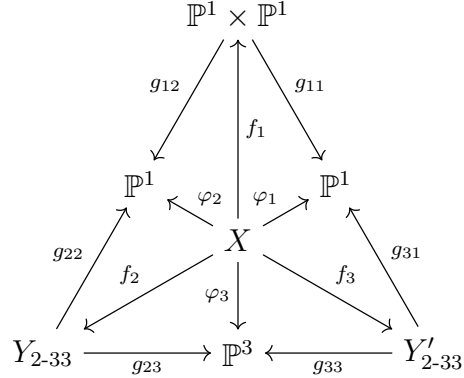
\[
\begin{tikzpicture}[commutative diagrams/every diagram,
    declare function={R=3;Rs=R*cos(60);}]
 \path 
  (0,0)  node(X) {$X$} 
  (90:R) node (Y1) {$\P^1_1 \times \P^1_2$}
  (210:R) node (Y2) {$Y_{\text{2-33}}$}
  (-30:R) node (Y3) {$Y'_{\text{2-33}}$}  
  (30:Rs) node(Z1) {$\P^1_1$} 
  (150:Rs) node(Z2) {$\P^1_2$} 
  (270:Rs) node(Z3) {$\P^3$};
 \path[commutative diagrams/.cd, every arrow, every label]
 (X) edge[swap] node {$f_1$} (Y1)
 (X) edge[swap] node {$f_2$} (Y2)
 (X) edge node {$f_3$} (Y3)
 (X) edge node {$\varphi_1$} (Z1)
 (X) edge[swap] node {$\varphi_2$} (Z2)
 (X) edge[swap] node {$\varphi_3$} (Z3)
 (Y1) edge node {$g_{11}=\pr_1$} (Z1)
 (Y1) edge[swap] node {$g_{12}=\pr_2$} (Z2)
 (Y2) edge node {$g_{22}$} (Z2)
 (Y2) edge[swap] node {$g_{23}$} (Z3)
 (Y3) edge node {$g_{33}$} (Z3)
 (Y3) edge[swap] node {$g_{31}$} (Z1);
\end{tikzpicture}
\]
\begin{proof}
Only (V) of  (\ref{n-pic3-vol}) holds.  
All the assertions except for (3) follow from Proposition \ref{p-V-2dP}. 
Lemma \ref{l-K-disjoint-blowup} and Proposition \ref{p-pic2-Pic} imply 
\[
-K_X \sim -f_2^*K_{Y_{\text{2-33}}}-f_3^*K_{Y'_{\text{2-33}}} + \varphi_3^*K_{\P^3} 
\sim (H_2 +3H_3) + (H_1+ 3H_3)-4H_3 = H_1 + H_2 +2H_3. 
\]
Thus (3) holds. 
\end{proof}

\begin{prop}[No.\ {\hyperref[table-3-26]{3-26}}]\label{p-pic3-26}
Let $X$ be a Fano threefold with $\rho(X)=3$ and $(-K_X)^3=46$. 
Then the following hold. 
\begin{enumerate}
\item $X$ has exactly three extremal rays. 
In what follows, we use Notation \ref{n-exactly3}. 
\item The contractions of the extremal faces are as in the following diagram. 
\begin{enumerate}
\item $f_1$ is of type $E_1$, $p_a(B_1)= 0$, $-K_{Y_1} \cdot B_1 = 3$. 
\item $f_2$ is of type $E_1$, $p_a(B_2)= 0$, $-K_{Y_2} \cdot B_2 = 4$. 
\item $f_3$ is of type $E_2$. 
\item $\deg B_{Z_2} =1$
\item $\varphi_3: X \to \P^3$ is a blowup along a disjoint union of a point and a line.  
\end{enumerate}
\item $-K_X \sim H_1 + 2H_2 + H_3$. 
\end{enumerate}
\end{prop}
\[
\begin{tikzpicture}[commutative diagrams/every diagram,
    declare function={R=3;Rs=R*cos(60);}]
 \path 
  (0,0)  node(X) {$X$} 
  (90:R) node (Y1) {$\P^2 \times \P^1$}
  (210:R) node (Y2) {$V_7$}
  (-30:R) node (Y3) {$Y_{\text{2-33}}$}  
  (30:Rs) node(Z1) {$\P^1$} 
  (150:Rs) node(Z2) {$\P^2$} 
  (270:Rs) node(Z3) {$\P^3$};
 \path[commutative diagrams/.cd, every arrow, every label]
 (X) edge[swap] node {$f_1$} (Y1)
 (X) edge[swap] node {$f_2$} (Y2)
 (X) edge node {$f_3$} (Y3)
 (X) edge node {$\varphi_1$} (Z1)
 (X) edge[swap] node {$\varphi_2$} (Z2)
 (X) edge[swap] node {$\varphi_3$} (Z3)
 (Y1) edge node {$g_{11}=\pr_2$} (Z1)
 (Y1) edge[swap] node {$g_{12}=\pr_1$} (Z2)
 (Y2) edge node {$g_{22}$} (Z2)
 (Y2) edge[swap] node {$g_{23}$} (Z3)
 (Y3) edge node {$g_{33}$} (Z3)
 (Y3) edge[swap] node {$g_{31}$} (Z1);
\end{tikzpicture}
\]
\begin{proof}
Only (I) of  (\ref{n-pic3-vol}) holds. 
By Theorem \ref{t-ele-tr-P2}, $X$ has a conic bundle structure over $\P^2$ of type 
2-34-vs-2-35. 
We then get the above commutative diagram except for $f_3, g_{31}, g_{33}$ 
such that (a), (b), and (d) hold. 
Moreover, (e) holds, because $D_{V_7} := \Ex(g_{23})$ is disjoint from 
the blowup centre $B_2$ of $f_2$ (Lemma \ref{l-P2-2-35}). 
Thus we obtain the above commutative diagram except for $g_{31}$. 
Moreover,  (c) holds. 

Lemma \ref{l-ele-tf-K-relation} and Proposition \ref{p-pic2-Pic} imply
\[
-2K_X \sim -f_1^*K_{\P^2 \times \P^1} - f_2^*K_{V_7} -\varphi_2^*B_{Z_2}
\]
\[
\sim (2H_1 + 3H_2) + (2H_2 + 2H_3) -H_2 = 2(H_1 +2H_2 + H_3). 
\]
Thus (3) holds. 
For the contraction $g_{34} : Y_{\text{2-33}} \to \P^1 =:Y_4$, the composition $\varphi_4 := g_{34} \circ f_3: X \to \P^1 =Y_4$, and $H_4 := \varphi_4^*\MO_{\P^1}(1)$, the following holds by (e): 
\[
-K_X \sim -f_2^*K_{V_7} -f_3^*K_{Y_{\text{2-33}}}+\varphi_3^*K_{\P^3}
\]
\[
\sim (2H_2 +2H_3) + (3H_3 + H_4) -4H_3 = 2H_2 +H_3 + H_4. 
\]
Therefore, we get $H_1 \sim H_4$, which imply (1) and (2). 
\qedhere



\end{proof}

\begin{prop}[No.\ {\hyperref[table-3-27]{3-27}}]\label{p-pic3-27}
Let $X$ be a Fano threefold with $\rho(X)=3$ and $(-K_X)^3=48$. 
Assume that $X$ has no conic bundle structure over $\F_1$. 
Then the following hold. 
\begin{enumerate}
\item $X$ has exactly three extremal rays. 
In what follows, we use Notation \ref{n-exactly3}. 
\item The contractions of the extremal faces are as in the following diagram. 
\begin{enumerate}
\item $f_1$ is of type $C_2$. 
\item $f_2$ is of type $C_2$. 
\item $f_3$ is of type $C_2$. 
\end{enumerate}
\item $-K_X \sim 2H_1 + 2H_2 + 2H_3$. 
\item 
$X \simeq \P^1 \times \P^1 \times \P^1$. 
\end{enumerate}
\end{prop}
\[
\begin{tikzpicture}[commutative diagrams/every diagram,
    declare function={R=3;Rs=R*cos(60);}]
 \path 
  (0,0)  node(X) {$X$} 
  (90:R) node (Y1) {$\P^1_1 \times \P^1_2$}
  (210:R) node (Y2) {$\P^1_2 \times \P^1_3$}
  (-30:R) node (Y3) {$\P^1_3 \times \P^1_1$}  
  (30:Rs) node(Z1) {$\P^1_1$} 
  (150:Rs) node(Z2) {$\P^1_2$} 
  (270:Rs) node(Z3) {$\P^1_3$};
 \path[commutative diagrams/.cd, every arrow, every label]
 (X) edge[swap] node {$f_1$} (Y1)
 (X) edge[swap] node {$f_2$} (Y2)
 (X) edge node {$f_3$} (Y3)
 (X) edge node {$\varphi_1$} (Z1)
 (X) edge[swap] node {$\varphi_2$} (Z2)
 (X) edge[swap] node {$\varphi_3$} (Z3)
 (Y1) edge node {$g_{11} =\pr_1$} (Z1)
 (Y1) edge[swap] node {$g_{12} =\pr_2$} (Z2)
 (Y2) edge node {$g_{22} =\pr_1$} (Z2)
 (Y2) edge[swap] node {$g_{23} =\pr_2$} (Z3)
 (Y3) edge node {$g_{33} =\pr_1$} (Z3)
 (Y3) edge[swap] node {$g_{31} =\pr_2$} (Z1);
\end{tikzpicture}
\]

\begin{proof}
Only (III) of  (\ref{n-pic3-vol}) holds.   
Then all the extremal rays are of type $C_2$ 
(Lemma \ref{l rho3 primitive type}). 
Hence   the assertions follow from \cite[Theorem 6.7]{ATIII}. 
\end{proof}

\begin{prop}[No.\ {\hyperref[table-3-28]{3-28}}]\label{p-pic3-28}
Let $X$ be a Fano threefold with $\rho(X)=3$ and $(-K_X)^3=48$. 
Assume that $X$ has a conic bundle structure over $\F_1$. 
Then the following hold. 
\begin{enumerate}
\item $X$ has exactly three extremal rays. 
In what follows, we use Notation \ref{n-exactly3}. 
\item The contractions of the extremal faces are as in the following diagram. 
\begin{enumerate}
\item $f_1$ is of type $C_2$.  
\item $f_2$ is of type $C_2$. 
\item $f_3$ is of type $E_1$, $p_a(B_3)= 0$, $-K_{Y_3} \cdot B_3 = 2$, and $B_3$ is a fibre of the projection $g_{33} : \P^2 \times \P^1 \to \P^2$.  
\end{enumerate}
\item $-K_X \sim 2H_1 + H_2 + 2H_3$. 
\item 
$X \simeq \F_1 \times \P^1$. 
\end{enumerate}
\end{prop}
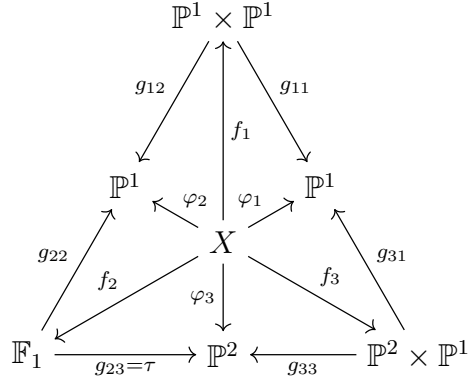
\[
\begin{tikzpicture}[commutative diagrams/every diagram,
    declare function={R=3;Rs=R*cos(60);}]
 \path 
  (0,0)  node(X) {$X$} 
  (90:R) node (Y1) {$\P^1_1 \times \P^1_2$}
  (210:R) node (Y2) {$\F_1$}
  (-30:R) node (Y3) {$\P^2 \times \P^1_1$}  
  (30:Rs) node(Z1) {$\P^1_1$} 
  (150:Rs) node(Z2) {$\P^1_2$} 
  (270:Rs) node(Z3) {$\P^2$};
 \path[commutative diagrams/.cd, every arrow, every label]
 (X) edge[swap] node {$f_1$} (Y1)
 (X) edge[swap] node {$f_2$} (Y2)
 (X) edge node {$f_3$} (Y3)
 (X) edge node {$\varphi_1$} (Z1)
 (X) edge[swap] node {$\varphi_2$} (Z2)
 (X) edge[swap] node {$\varphi_3$} (Z3)
 (Y1) edge node {$g_{11}=\pr_1$} (Z1)
 (Y1) edge[swap] node {$g_{12} =\pr_2$} (Z2)
 (Y2) edge node {$g_{22}$} (Z2)
 (Y2) edge[swap] node {$g_{23} =\tau$} (Z3)
 (Y3) edge node {$g_{33} =\pr_1$} (Z3)
 (Y3) edge[swap] node {$g_{31}=\pr_2$} (Z1);
\end{tikzpicture}
\]
\begin{proof}
Since $X$ has a conic bundle structure over $\F_1$, 
Theorem \ref{t-F1-pic3} implies  $X \simeq (\P^1 \times \P^2) \times_{\P^2} \F_1 \simeq \P^1 \times \F_1$, i.e., (4) holds. 
Moreover, (b) and (c) hold. 
In particular, we get the cartesian diagram consisting of $f_2, g_{23}, f_3, g_{33}$. 
Hence we obtain the above commutative diagram except for $f_1, g_{11}, g_{12}$. 
By Lemma \ref{l-exactly3 2P^1}, 
(1) holds (Lemma \ref{l-exactly3}) and 
we get the above commutative diagram. 

The extremal ray corresponding to $f_1$ is of type $C_2$ (i.e., (a) holds), because the scheme-theoretic fibre of $f_1$ over a closed point $(t_1, t_2) \in \P^1_1 \times \P^1_2$ is smooth: 
\[
(\P^1_1 \times \F_1) \times_{\P^1_1 \times \P^1_2} (t_1, t_2) 
\simeq \{t_1\} \times (\F_1 \times_{\P^1_2} \{ t_2\}) \simeq \P^1, 
\]
where $\P^1_1 := Z_1 = \P^1$ and $\P^1_2 := Z_2 = \P^1$. 
Thus we get (2). 
Lemma \ref{l-F1cart-KX} implies 
\[
-K_X \sim -f_3^*K_{\P^2 \times \P^1} -H_3 + H_2 
\]
\[
\sim (3H_3 + 2H_1) -H_3 + H_2 = 2H_1 +H_2 +2H_3. 
\]
Thus (3) holds. 
\end{proof}

\begin{prop}[No.\ {\hyperref[table-3-29]{3-29}}]\label{p-pic3-29}
Let $X$ be a Fano threefold with $\rho(X)=3$ and $(-K_X)^3=50$. 
Assume that $X$ has a conic bundle structure over $\P^2$. 
Then the following hold. 
\begin{enumerate}
\item $X$ has exactly three extremal rays. 
In what follows, we use Notation \ref{n-exactly3}. 
\item The contractions of the extremal faces are as in the following diagram. 
\begin{enumerate}
\item $f_1$ is of type $E_1$, $p_a(B_1)= 0$, $-K_{Y_1} \cdot B_1 = 2$, 
and $B_1$ is a line on $\Ex(g_{11}) \simeq \P^2$. 
\item $f_2$ is of type $E_1$, $p_a(B_2)= 0$, $-K_{Y_2} \cdot B_2 = 5$. 
\item $f_3$ is of type $E_5$. 
\item $\deg B_{Z_2} = 1$. 
\item $Z_3$ is the cone over the Veronese surface $S \subset \P^5$, and 
$\varphi_3$ is a blowup along a disjoint union of the singular point and a smooth rational curve $C$ 
satisfying $-K_{Z_3} \cdot C =5$. 
\end{enumerate}
\item $-K_X \sim H_1 + H_2 + H_3$. 
\end{enumerate}
\end{prop}
\[
\begin{tikzpicture}[commutative diagrams/every diagram,
    declare function={R=3;Rs=R*cos(60);}]
 \path 
  (0,0)  node(X) {$X$} 
  (90:R) node (Y1) {$V_7$}
  (210:R) node (Y2) {$\P_{\P^2}(\MO \oplus \MO(2))$}
  (-30:R) node (Y3) {$Y_3$}  
  (30:Rs) node(Z1) {$\P^3$} 
  (150:Rs) node(Z2) {$\P^2$} 
  (270:Rs) node(Z3) {$Z_3$};
 \path[commutative diagrams/.cd, every arrow, every label]
 (X) edge[swap] node {$f_1$} (Y1)
 (X) edge[swap] node {$f_2$} (Y2)
 (X) edge node {$f_3$} (Y3)
 (X) edge node {$\varphi_1$} (Z1)
 (X) edge[swap] node {$\varphi_2$} (Z2)
 (X) edge[swap] node {$\varphi_3$} (Z3)
 (Y1) edge node {$g_{11}$} (Z1)
 (Y1) edge[swap] node {$g_{12}$} (Z2)
 (Y2) edge node {$g_{22}$} (Z2)
 (Y2) edge[swap] node {$g_{23}$} (Z3)
 (Y3) edge node {$g_{33}$} (Z3)
 (Y3) edge[swap] node {$g_{31}$} (Z1);
\end{tikzpicture}
\]
\begin{proof}
By Theorem \ref{t-ele-tr-P2}, $X$ has a conic bundle structure over $\P^2$ of type 2-35-vs-2-36. 
We then get the above commutative diagram except for $f_3, g_{31}, g_{33}$ such that 
(a)', (b), and (d) hold. 
\begin{enumerate}
\item[(a)'] $f_1$ is of type $E_1$, $p_a(B_1)= 0$ and $-K_{Y_1} \cdot B_1 = 2$. 
\end{enumerate}
Then (e) holds by (b) and $D_{Y_2} \cap B_2  =\emptyset$ for $D_{Y_2} := \Ex(g_{23})$ (Lemma \ref{l-P2-2-36}). 
We then obtain the contraction $f_3 : X \to Y_3$ of type $E_5$, and we get the above commutative diagram except for $g_{31}$. 
In particular, (c) holds.

We now show that $B_1$ is a line on $D_{Y_1} := \Ex(g_{11}) \simeq \P^2$. 
By the first paragraph of the proof of Lemma \ref{l-P2-2-35}, 
we get $-K_{Y_1} \cdot B_1 = 4 \deg B_{Z_2} + 2 D_{Y_1} \cdot B_1$. 
This, together with $-K_{Y_1} \cdot B_1 \overset{{\rm (a)}}{=} 2$ and 
$\deg B_{Z_2} \overset{{\rm (d)}}{=} 1$, implies $D_{V_7} \cdot B_1 = -1$. 
Thus $B_1$ is a line on $D_{Y_1} \simeq \P^2$. 
Hence (a) holds.

Let $D_{Y_1}^X$ and $D^X_{Y_2}$ be the proper transforms of $D_{Y_1}$ and $D_{Y_2}$ on $X$, respectively. 
Let us show $D_{Y_1}^X = D_{Y_2}^X$. 
By the first (resp. second) paragraph of the proof of Lemma \ref{l-P2-2-35} 
(resp. Lemma \ref{l-P2-2-36}), we obtain 
\[
-K_{Y_1} \sim 2D_{Y_1} + g_{12}^*\MO_{\P^2}(4) 
\qquad (\text{resp.} -K_{Y_2} \sim 2D_{Y_2} + g_{22}^*\MO_{\P^2}(5)). 
\]
By (a) and $D_{Y_2} \cap B_2  =\emptyset$, we get $f_1^*D_{Y_1} = D_{Y_1}^X + \Ex(f_1)$ and $f_2^*D_{Y_2} = D_{Y_2}^X$, respectively. 
It holds that  
\[
-K_X +\Ex(f_1) \sim -f_1^*K_{Y_1} \sim f_1^*(2D_{Y_1} + g_{12}^*\MO_{\P^2}(4)) \sim 2D_{Y_1}^X + 2\Ex(f_1)+4H_2, 
\]
\[
-K_X + \Ex(f_2) \sim -f_2^*K_{Y_2} \sim f_2^*(2D_{Y_2} + g_{22}^*\MO_{\P^2}(5)) 
\sim 2D_{Y_2}^X +5 H_2. 
\]
Therefore, 
\[
2D_{Y_1}^X + \Ex(f_1)+4H_2 \sim -K_X \sim  2D_{Y_2}^X-\Ex(f_2) +5 H_2. 
\]
This, together with $H_2 \sim \varphi_2^*B_{Z_2} = \Ex(f_1)+\Ex(f_2)$, implies $D_{Y_1}^X \sim D_{Y_2}^X$. 
By $h^0(X, D_{Y_1}^X) = 1$, we get $D_{Y_1}^X = D_{Y_2}^X$.

By $D_{Y_1}^X = D_{Y_2}^X$, 
both $\varphi_1$ and $\varphi_3$ contract this prime divisor  to a point. 
In particular, $\varphi_1 \times \varphi_3 : X \to \P^3 \times Z_3$ is not a finite morphism. 
Thus (1) holds (Lemma \ref{l-exactly3}) and we get the above commutative diagram. 
Moreover, we get (2). 
Lemma \ref{l-ele-tf-K-relation} and Proposition \ref{p-pic2-Pic} imply 
\[
-2K_X \sim -f_1^*K_{V_7} -f_2^*K_{\P_{\P^2}(\MO \oplus \MO(2))} - \varphi_2^*B_{Z_2}
\]
\[
\sim (2H_1 + 2H_2) + (H_2 + 2H_3) - H_2 = 2(H_1 + H_2 + H_3). 
\]
Thus (3) holds. 
\qedhere

    

\end{proof}

\begin{prop}[No.\ {\hyperref[table-3-30]{3-30}}]\label{p-pic3-30}
Let $X$ be a Fano threefold with $\rho(X)=3$ and $(-K_X)^3=50$. 
Assume that $X$ has no conic bundle structure over $\P^2$. 
Then the following hold. 
\begin{enumerate}
\item $X$ has exactly three extremal rays. 
In what follows, we use Notation \ref{n-exactly3}. 
\item The contractions of the extremal faces are as in the following diagram. 
\begin{enumerate}
\item $f_1$ is of type $C_2$. 
\item $f_2$ is of type $E_1$, $p_a(B_2)= 0$, $-K_{Y_2} \cdot B_2 = 1$. 
\item $f_3$ is of type $E_1$, $p_a(B_3)= 0$, $-K_{Y_3} \cdot B_3 = 2$, 
and $B_3$ is the proper transform of a line passing through the blowup centre of $g_{33} : V_7 \to \P^3$. 
\end{enumerate}
\item $-K_X \sim H_1 + H_2 + 2H_3$. 
\end{enumerate}
\end{prop}
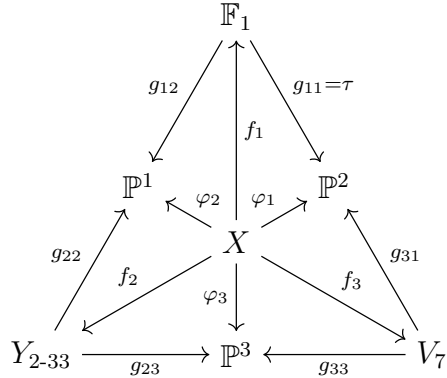
\[
\begin{tikzpicture}[commutative diagrams/every diagram,
    declare function={R=3;Rs=R*cos(60);}]
 \path 
  (0,0)  node(X) {$X$} 
  (90:R) node (Y1) {$\F_1$}
  (210:R) node (Y2) {$Y_{\text{2-33}}$}
  (-30:R) node (Y3) {$V_7$}  
  (30:Rs) node(Z1) {$\P^2$} 
  (150:Rs) node(Z2) {$\P^1$} 
  (270:Rs) node(Z3) {$\P^3$};
 \path[commutative diagrams/.cd, every arrow, every label]
 (X) edge[swap] node {$f_1$} (Y1)
 (X) edge[swap] node {$f_2$} (Y2)
 (X) edge node {$f_3$} (Y3)
 (X) edge node {$\varphi_1$} (Z1)
 (X) edge[swap] node {$\varphi_2$} (Z2)
 (X) edge[swap] node {$\varphi_3$} (Z3)
 (Y1) edge node {$g_{11} =\tau$} (Z1)
 (Y1) edge[swap] node {$g_{12}$} (Z2)
 (Y2) edge node {$g_{22}$} (Z2)
 (Y2) edge[swap] node {$g_{23}$} (Z3)
 (Y3) edge node {$g_{33}$} (Z3)
 (Y3) edge[swap] node {$g_{31}$} (Z1);
\end{tikzpicture}
\]
\begin{proof}
 Only (II) of  (\ref{n-pic3-vol}) holds. 
 By Theorem \ref{t-F1-pic3}, we obtain $X \simeq V_7 \times_{\P^2} \F_1$. 
Thus we get the above commutative diagram except for $f_2, g_{22}, g_{23}$ 
such that (a) and (c)' hold. 
\begin{enumerate}
\item[(c)'] $f_3$ is of type $E_1$, $p_a(B_3)= 0$, 
and $-K_{Y_3} \cdot B_3 = 2$. 
\end{enumerate}
Since $\varphi_3 : X \to \P^3$ has a two-dimensional fibre, 
we can find a curve $C$ on $X$ contracted by $\varphi_2$ and $\varphi_3$. 
Thus (1) holds (Lemma \ref{l-exactly3}). 
Let $f_2 : X \to Y_2$ be the contraction of the remaining extremal ray. 

Note that the blowup centre $B_3$ of $f_3 : X \to V_7$ is a fibre of the induced $\P^1$-bundle $g_{31} : V_7 \to \P^2$. 
Since $D_{V_7} :=\Ex(g_{33})$ is a section of $g_{31}$, we get $D_{V_7} \cdot B_3 = 1$ and $D_X  \simeq \F_1$ for the proper transform $D_X$ of $D_{V_7}$ on $X$. 
In particular, $B_{\P^3} := g_{33}(B_3)$ is a smooth rational curve. 
Then (c) holds, because 
\[
-2 = K_{V_7} \cdot B_3 = (g_{33}^*K_{\P^3} + 2D_{V_7}) \cdot B_3 
= K_{\P^3} \cdot B_{\P^3}+ 2. 
\]

We now show that $\Ex(f_2) =D_X$. 
Suppose $\Ex(f_2) \neq D_X$. 
By $\Ex(f_2) \subset \Ex(\varphi_3) = D_X \cup \Ex(f_3)$, 
we would get $\Ex(f_2) = \Ex(f_3)$. 
Since $\varphi_3(\Ex(f_3))$ 
 is a curve, also the image $f_2(\Ex(f_2))$ of $\Ex(f_2) = \Ex(f_3)$ to $Y_2$ is a curve. 
Then we can find a curve on $X$ contracted by $\Ex(f_2) \to f_2(\Ex(f_2))$ and $\Ex(f_3) \to f_3(\Ex(f_3))$, 
which contradicts the fact that $f_2$ and $f_3$ belong to distinct extremal rays.

By $\Ex(f_2) = D_X \simeq \F_1$, 
$f_2$ is of type $E_1$, $p_a(B_2)=0$, $Y_2$ is a Fano threefold (Lemma \ref{l-nonFano-blowdown}),  and $g_{23}(\Ex(g_{23})) =\varphi_3(\Ex(f_3)) = f_3(B_3) =  B_{\P^3}$.
Hence $Y_2$ is a Fano threefold of No.\ 2-33 (Subsection \ref{ss-table-pic2}). 
Moreover, 
it follows from Lemma \ref{l-blowup-formula} that 
\[
50 = (-K_X)^3  = (-K_{Y_{\text{2-33}}})^3 -2 (-K_{Y_{\text{2-33}}}) \cdot B_2 +2p_a(B_2) -2 
=54  - 2(-K_{Y_2}) \cdot B_2 -2, 
\]
i.e., $(-K_{Y_2}) \cdot B_2 =1$. Then (2) holds. 
Lemma \ref{l-F1cart-KX} and Lemma \ref{p-pic2-Pic} imply 
\[
-K_X \sim -f_3^*K_{V_7} - H_1 +H_2 \sim (2H_3+2H_1) -H_1 + H_2 = H_1 + H_2 + 2H_3. 
\]
Thus (3) holds. 
\end{proof}

\begin{prop}[No.\ {\hyperref[table-3-31]{3-31}}]\label{p-pic3-31}
Let $X$ be a Fano threefold with $\rho(X)=3$ and $(-K_X)^3=52$. 
Then the following hold. 
\begin{enumerate}
\item $X$ has exactly three extremal rays. 
In what follows, we use Notation \ref{n-exactly3}. 
\item The contractions of the extremal faces are as in the following diagram. 
\begin{enumerate}
\item $f_1$ is of type $C_2$. 
\item $f_2$ is of type $E_1$, $p_a(B_2)= 0$, $-K_{Y_2} \cdot B_2 = 0$. 
\item $f_3$ is of type $E_1$, $p_a(B_3)= 0$, $-K_{Y_3} \cdot B_3 = 0$. 
\item $D := \Ex(f_2) =\Ex(f_3) = \Ex(\varphi_3) \simeq \P^1 \times \P^1$, 
$\varphi_3(D)$ is a point, and $D$ is a section of $f_1$. 
\end{enumerate}
\item $-K_X \sim 3H_1 + 3H_2 + 2D$. 
\item 
$X \simeq \P_{\P^1 \times \P^1}( \MO_{\P^1 \times \P^1} \oplus \MO_{\P^1 \times \P^1}(1, 1))$. 
\end{enumerate}
\end{prop}
\[
\begin{tikzpicture}[commutative diagrams/every diagram,
    declare function={R=3;Rs=R*cos(60);}]
 \path 
  (0,0)  node(X) {$X$} 
  (90:R) node (Y1) {$\P^1_1 \times \P^1_2$}
  (210:R) node (Y2) {$Y_{\text{non-Fano}}$}
  (-30:R) node (Y3) {$Y'_{\text{non-Fano}}$}  
  (30:Rs) node(Z1) {$\P^1_1$} 
  (150:Rs) node(Z2) {$\P^1_2$} 
  (270:Rs) node(Z3) {$Z_3$};
 \path[commutative diagrams/.cd, every arrow, every label]
 (X) edge[swap] node {$f_1$} (Y1)
 (X) edge[swap] node {$f_2$} (Y2)
 (X) edge node {$f_3$} (Y3)
 (X) edge node {$\varphi_1$} (Z1)
 (X) edge[swap] node {$\varphi_2$} (Z2)
 (X) edge[swap] node {$\varphi_3$} (Z3)
 (Y1) edge node {$g_{11} =\pr_1$} (Z1)
 (Y1) edge[swap] node {$g_{12} =\pr_2$} (Z2)
 (Y2) edge node {$g_{22}$} (Z2)
 (Y2) edge[swap] node {$g_{23}$} (Z3)
 (Y3) edge node {$g_{33}$} (Z3)
 (Y3) edge[swap] node {$g_{31}$} (Z1);
\end{tikzpicture}
\]

\begin{proof}
 Only (III) of  (\ref{n-pic3-vol}) holds. 
Then $X$ has an extremal rays $R_1$ and $R_2$ such that 
$R_1$ is of type $C_2$ and $R_2$ is of type $E_1$ 
(Lemma \ref{l rho3 primitive type}). 
 By \cite[Lemma 6.10, Proposition 6.11]{ATIII}, 
 we get $X \simeq \P_{\P^1 \times \P^1}(\MO \oplus \MO(1, 1))$, 
 the induced $\P^1$-bundle $f_1 : X  \to \P^1 \times \P^1$, and 
a contraction $f_2: X \to Y_{\text{non-Fano}}$ of an extremal ray of type $E_1$ such that 
$D :=\Ex(f_2)$ is a section of $f_1$.  
In particular, (a) and (4) hold. 
As $X$ is primitive, $Y_{\text{non-Fano}}$ is actually non-Fano. 
By Lemma \ref{l-nonFano-blowdown}, Proposition \ref{p-nonFano-iff}, and Proposition \ref{p-nonFano-flop}, 
we obtain the above diagram except for $g_{31}$ such that (b), (c), and (d) hold. 
Since $\varphi_3(D)$ is a point, 
it follows from $f_1|_D :D \xrightarrow{\simeq} \P^1 \times \P^1$ that 
there exists a curve on $X$ contracted by $\varphi_1$ and $\varphi_3$. 
Thus (1) holds (Lemma \ref{l-exactly3}) and we get the above commutative diagram. 
Then (2) holds.

Since $D$ is a section of $f_1 : X \to \P^1 \times \P^1$, 
we can write $-K_X \equiv 2D +a_1 H_1 + a_2 H_2$  for some $a_1, a_2 \in \Q$. 
By $-K_X|_D \sim -D|_D \sim  \MO_{\P^1 \times \P^1}(1, 1)$ (Lemma \ref{l-nonFano-blowdown}),  we get $a_1 = a_2=3$, i.e., 
$-K_X \sim 3H_1 + 3H_2 + 2D$. Thus (3) holds. 
\end{proof}

\begin{dfn}\label{d-pic3}
Let $X$ be a Fano threefold with $\rho(X)=3$. 
We say that $X$ is {\em 3-xx} or of {\em No.\ 3-xx} 
if  $(-K_X)^3$, the types of the extremal rays, 
and the images of the contractions of the extremral rays 
satisfies the corresponding properties listed in Table \ref{table-pic3} in Subsection \ref{ss-table-pic3}. 
For example, the definitions of No.\ 3-1 and 3-6 are as follows. 
\begin{itemize}
\item  A Fano threefold $X$  is {\em 3-1} or  {\em of No.\ 3-1} if 
$\rho(X)=3$, $(-K_X)^3=12$, and 
there exist exactly three extremal rays $R_1, R_2, R_3$ such that  
all the contractions are of type $C_1$ and their images are $\P^1 \times \P^1$. 
\item  A Fano threefold $X$  is {\em 3-11} or  {\em of No.\ 3-11} if $\rho(X)=3$, $(-K_X)^3=28$, there exist exactly three extremal rays $R_1, R_2, R_3$, and the images of the contractions of $R_1, R_2, R_3$ are Fano threefolds of No.\ 2-25, 2-34, 2-35, respectively. 
\end{itemize}
\end{dfn}

\begin{thm}\label{t-pic3-main}
Let $X$ be a Fano threefold with $\rho(X)=3$. 
Then $X$ satisfies one and only one of the possibilities listed in Table \ref{table-pic3} in Subsection \ref{ss-table-pic3},  except for the column \lq\lq blowups". 
\end{thm}

\begin{proof}
The assertion follows from results in this subsection. 
For example, if $X$ is a Fano threefold with $\rho(X)=3$ and $(-K_X)^3=36$, 
then the assertion follows from Proposition \ref{p-pic3-17} and Proposition \ref{p-pic3-18}. 
\end{proof}

\begin{cor}\label{c-P3-disjoint}
Let $C_1$ and $C_2$ be smooth curves on $\P^3$ such that $C_1 \cap C_2 = \emptyset$. 
Assume that $\deg C_1 \leq \deg C_2$ and the blowup $X := \Bl_{C_1 \amalg C_2} \P^3$ is Fano. 
Then one of the following holds. 
\begin{enumerate}
\item $X$ is of No.\ 3-6, $C_1$ is a line, and $C_2$ is an elliptic curve of degree $4$. 
\item $X$ is of No.\ 3-12, $C_1$ is a line, and $C_2$ is a rational cubic curve. 
\item $X$ is of No.\ 3-18, $C_1$ is a line, and $C_2$ is a conic. 
\item $X$ is of No.\ 3-25 and both $C_1$  and $C_2$ are lines. 
\end{enumerate}
\end{cor}

\begin{proof}
The assertion follows from the classification obtained in this subsection. 
Indeed, except for the above cases (1)-(4), 
if there is a contraction $X \to \P^3$, 
then $X$ is of No.\ 3-11, 3-14, 3-19, 3-23, 3-26, 3-29, or 3-30. 
If $X$ is not 3-19, 
then there is a unique contraction $X \to \P^3$, and it factors through the blowup $V_7 \to \P^3$ (e.g., see Proposition \ref{p-pic3-11} for the case of No.\ 3-11). 
If $X$ is 3-19, then there exist exactly two  contractions $X \to \P^3$, 
and each of them factors through the blowup $V_7 \to \P^3$. 
\end{proof}

\begin{cor}\label{c-Q-disjoint}
Let $C_1$ and $C_2$ be smooth curves on $Q$ such that $C_1 \cap C_2 = \emptyset$. 
Assume that $\deg C_1 \leq \deg C_2$ and the blowup $X := \Bl_{C_1 \amalg C_2}\,Q$ is Fano. 
Then one of the following holds. 
\begin{enumerate}
\item $X$ is of No.\ 3-10 and both $C_1$ and $C_2$ are conics. 
\item $X$ is of No.\ 3-15, $C_1$ is a line, and $C_2$ is a conic. 
\item $X$ is of No.\ 3-20 and both $C_1$ and $C_2$ are lines. 
\end{enumerate}
\end{cor}

\begin{proof}
The assertion follows from the classification obtained in this subsection. 
Indeed, except for the above cases (1)-(4), 
if there is a contraction $X \to \P^3$, 
then $X$ is of No.\ 3-18, 3-19, or 3-23. 
For each case, there exists a unique contraction $X \to Q$ and it factors through 
the blowup $Y_{\text{2-30}} \to Q$ at a point. 
\end{proof}

\subsection{Fano conic bundles over $\P^1 \times \P^1$ ($\rho=3$)}\label{ss FCB P1P1 rho3}

For the classification for the case $\rho =4$, 
we shall need the following classification of Fano conic bundles over $\P^1 \times \P^1$ with $\rho=3$.

\begin{thm}\label{t-P1P1-pic3}
Set $S := \P^1 \times \P^1$  and let $f: X \to S = \P^1 \times \P^1$ be a Fano conic bundle. 
Assume that $\rho(X)=3$ and $d_1 \leq d_2$ for the bidegree $(d_1, d_2)$ of $\Delta_f$. 
Then one of the following holds. 
    \begin{center}
\begin{longtable}{cccccc}
No. & $(-K_X)^3$ &  bidegree of $\Delta_f$ \\ \hline
3-1 & $12$ &  $(4, 4)$\\ \hline
3-2 & $14$ & $(2, 5)$\\ \hline
3-3 & $18$ &  $(3, 3)$ \\ \hline
3-4 & $18$ &  $(2, 4)$\\ \hline
3-6 & $22$ &  $(2, 3)$\\ \hline
3-10 & $26$ & $(2, 2)$ \\ \hline
3-17 & $36$ & $(0, 0)$ \\ \hline
3-25 & $44$ & $(0, 0)$ \\ \hline
3-27 & $48$ &  $(0, 0)$ \\ \hline
3-28 & $48$ &   $(0, 0)$\\ \hline
3-31 & $52$ &  $(0, 0)$ \\ \hline
      \end{longtable}
  \end{center} 

\end{thm}

\begin{proof}
The assertion follows from Table \ref{table-pic3} in Subsection \ref{ss-table-pic3}, 
which is available by  Theorem \ref{t-pic3-main}. 
\end{proof}

\section{$\rho=4$}

The purpose of this subsection is to classify Fano threefolds of Picard number $4$. 
In what follows, we overview its proof and the contents of this section.

Let $X$ be a Fano threefold with $\rho(X)=4$. 
Since $X$ is imprimitive, there exists a blowup $f : X \to Y$ of a Fano threefold $Y$ along a smooth curve. 
By using the classification for the case of Picard number $3$, 
we prove the existence of a conic bundle $X\to S$ with 
$S \in \{\P^1 \times \P^1, \F_1\}$ (Subsection \ref{ss-pic4-structure}). 
Then we shall classify such conic bundles in 
Subsection \ref{ss-pic4-P1P1} ($S = \P^1 \times \P^1$) 
and Subsection \ref{ss-pic4-F1} ($S = \F_1$). 
In order to check the overlapping, 
we study the case when $\rho(X)=4$ and $(-K_X)^3 = 32$ in 
Subsection \ref{ss rho4 vol=32}. 
More concretely, this case is divided into  two subcases: 4-4 and 4-5, which are distinguished by whether 
there exists a smooth curve along which the blowup is Fano. 
In Subsection \ref{ss-pic4-classify},  we complete the classification of Fano threefolds of Picard number $4$.

\subsection{Existence of conic bundle structures}\label{ss-pic4-structure}

Let $X$ be a Fano threefold. 
In this subsection, we prove the following. 
\begin{itemize}
\item If $\rho(X) \geq 4$, then $X$ has a conic bundle structure (Theorem \ref{t-pic4-CB}).  
\item If $\rho(X) =4$, then there exists a conic bundle $X \to S$ with $S \in \{\P^1 \times \P^1, \F_1\}$ (Corollary \ref{c-pic4-CB}). 
\end{itemize}
We start with the following auxiliary result. 




\begin{lem}\label{l-pic4-3-18}
Let $Y$ be a Fano threefold of No.\ 3-18, i.e., $Y = \Bl_{L \amalg C}\,\P^3$ 
for a disjoint union of a line $L$ and a conic $C$ on $\P^3$. 
Let $\rho: Y= \Bl_{L \amalg C}\,\P^3 \to \P^3$ be the induced blowup. 
Take a smooth curve $\Gamma$ on $Y$. 
Assume that $X := \Bl_{\Gamma}\,Y$ is Fano. 
Then $\rho(\Gamma)$ is a point. Moreover, $(-K_X)^3 =32$. 
\end{lem}

\begin{proof}
Suppose that $\rho(\Gamma)$ is not a point. 
Then  $\Gamma$ is disjoint from $\Ex(\rho)$ (Lemma \ref{l-line-meeting}). 
Recall that the smallest linear subvariety $\langle C \rangle$ of $\P^3$ containing $C$ is a plane. 
Since $\langle C \rangle$ is an ample divisor on $\P^3$, 
it follows from $C \cap L = \emptyset$ and $C \cap \rho(\Gamma) = \emptyset$ that 
each of $\langle C \rangle \cap L$ and $\langle C \rangle \cap \rho(\Gamma)$ is non-empty and zero-dimensional. 
Hence we can find a line $Z$ on $\P^3$ such that  $Z \subset \langle C \rangle$, 
$Z \cap L \neq \emptyset$, and $Z \cap \rho(\Gamma) \neq \emptyset$.  
Recall that the induced birational morphism $h : X \to \P^3$ is the blowup along $L \amalg C \amalg \rho(\Gamma)$. 
For the proper transform $Z'$ of $Z$ on $X$, 
we get the following contradiction: 
\[
0< -K_X \cdot Z' = h^*\MO_{\P^3}(4) \cdot Z' -h^{-1}(L) \cdot Z' -h^{-1}(C) \cdot Z' -h^{-1}(\rho(\Gamma)) \cdot Z'  \leq 4 -1 -2 -1 =0, 
\]
It follows from $\Gamma \simeq \P^1$ and $-K_Y \cdot \Gamma=1$ 
that  $(-K_X)^3 = (-K_Y)^3 - 4 = 32$ (Lemma \ref{l-blowup-formula}). 
\end{proof}

\begin{thm}\label{t-pic4-CB}
Let $X$ be a Fano threefold with $\rho(X) \geq 4$. 
Then $X$ has a conic bundle structure. 
\end{thm}

\begin{proof}
By Proposition \ref{p-FCB-centre}, we may assume that $\rho(X) = 4$. 
Since $X$ is imprimitive, 
$X$ is obtained by taking a blowup $f: X \to Y$ of a Fano threefold $Y$ with $\rho(Y)=3$ along a smooth curve $\Gamma$ on $Y$. 
By Proposition \ref{p-FCB-centre} and Theorem \ref{t-pic3-structure}, 
the problem is reduced to the case when $Y$ is of No.\ 3-18, 
i.e., there exist a line $L$ and a conic $C$ on $\P^3$ such that $L \cap C =\emptyset$ and 
$Y \simeq \Bl_{L \amalg C} \P^3$. Let $\rho : Y \to \P^3$ be the induced blowup. 
By Lemma \ref{l-pic4-3-18}, $\rho(\Gamma)$ is a point. 
Hence  $\Gamma$ is either a fibre of $\rho^{-1}(L) \to L$ or a fibre of $\rho^{-1}(C) \to C$. 
In any case, we obtain birational morphisms $X \to Z \to V_7$, 
where $Z$ is a Fano threefold and 
each of $X \to Z$ and $Z \to V_7$ is of type $E_1$ (Lemma \ref{l-P3-V7}). 
Again by Proposition \ref{p-FCB-centre}, $X$ has a conic bundle structure. 
\qedhere



\end{proof}

\begin{cor}\label{c-pic4-CB}
Let $X$ be a Fano threefold with $\rho(X) = 4$. 
Then there exists  a conic bundle  $X \to S$ satisfying 
$S \in \{\P^1 \times \P^1, \F_1\}$. 
\end{cor}

\begin{proof}
By Theorem \ref{t-pic4-CB}, there is a conic bundle structure $f: X \to T$. 
If $T \simeq \P^1 \times \P^1$ or $T \simeq \F_1$, then we are done.

Let us show $T \not\simeq \P^2$. 
Suppose $T \simeq \P^2$. 
By $\rho(X) =4 >2 = \rho(T)+1$, there exist morphisms  
\[
f : X \xrightarrow{\sigma} Y \xrightarrow{g} T (\simeq \P^2), 
\]
where $g : Y \to T$ is a Fano conic bundle with $\rho(Y)=3$ and $\sigma : X \to Y$ is a blowup along a regular subsection $B_Y$ of $g$ (Proposition \ref{p-ele-tf}, Proposition \ref{p-ele-tf-1Fano}, Lemma \ref{l-FCB-pic-irre}). 
Again by $\rho(Y) =3 >2 = \rho(T)+1$, there exist morphisms  
\[
f : X \xrightarrow{\sigma} Y \xrightarrow{\tau} Z \xrightarrow{h} T (\simeq \P^2), 
\]
where $h : Z \to T$ is a Fano conic bundle with $\rho(Z)=2$ and $\tau : Y \to Z$ is a blowup along a regular subsection $B_Z$ of $h$. 
Then it holds that $\Delta_g \neq \emptyset$ and 
$g(B_Y) \cap \Delta_g = \emptyset$ (Proposition \ref{p-FCB-centre}). 
However, this is impossible, because $\Delta_g$ is ample.  

In what follows, we assume that $T$ is isomorphic to none of $\P^2, \P^1 \times \P^1, \F_1$. 
Then $X \simeq T \times \P^1$ (Proposition \ref{p-FCB-triv}), 
where $T$ is a smooth del Pezzo surface with  $K_T^2=7$. 
Then we have a blowup $\sigma : X = T \times \P^1 \to \F_1 \times \P^1$ along 
$\{ z \} \times \P^1$ for some closed point $z \in \F_1$. 
For the $\P^1$-bundle $\pi :\F_1 \to \P^1$, 
$\F_1 \times \P^1$ has a conic bundle structure 
$\pi \times {\rm id} : \F_1 \times \P^1 \to \P^1 \times \P^1$. 
Since $\{ z \} \times \P^1$ is a subsection of $\pi \times {\rm id}$, 
$X$ has a conic bundle structure over $\P^1 \times \P^1$ 
(Proposition \ref{p-FCB-centre}). 
\end{proof}

\subsection{Fano conic bunldes over $\P^1 \times \P^1$ ($\rho=4$)}\label{ss-pic4-P1P1}

The purpose of this subsection is to classify Fano conic bundles 
$f: X \to \P^1 \times \P^1$ with $\rho(X)=4$. 
We start with the following observation.

\begin{nasi}\label{nasi-pic4-P1P1}
Let $f: X \to \P^1 \times \P^1$ be a Fano conic bundle with $\rho(X)=4$. 
Then 
$X$ is obtained from 
another Fano conic bundle 
$g: Y \to \P^1 \times \P^1$ with $\rho(Y) =3$ by taking 
a blowup $X \to Y$ along a regular subsection 
$B_Y$ of $g$ (Proposition \ref{p-smaller-same-base}): 
\[
f : X \to Y \xrightarrow{g} \P^1 \times \P^1. 
\]
By Theorem \ref{t-P1P1-pic3}, $\Delta_g$ is either zero or ample. 
Note that  $\Delta_g$ is not ample (Proposition \ref{p-FCB-centre}). 
Hence it is enough to consider the case when $\Delta_g=0$, which is summarised in the following table (Theorem \ref{t-P1P1-pic3}): 



    \begin{center}
\begin{longtable}{ccp{7.5cm}cccc}
No. & $(-K_Y)^3$ & descriptions and extremal rays & conic bdl$/\P^2$  &   \\ \hline
3-17 & $36$ & 
$Y$ is a divisor on $\P^1 \times \P^1 \times \P^2$ of tridegree $(1, 1, 1)$
 && \\
&& $C_2: /\P^1 \times \P^1$ &  2-34-vs-2-34 & \\ 
 &  & $E_1:$ 2-34, $p_a(C)=0, -K \cdot C =8$  &  & \\
  &  & $E_1:$ 2-34, $p_a(C)=0, -K \cdot C =8$  &  & \\\hline
3-25 & $44$ & $Y$ is a blowup of $\P^3$ along a disjoint union of two lines 
& & \\
&&$C_2:$ $/\P^1 \times \P^1$  & none && \\
 &  & $E_1:$ 2-33, $p_a(C)=0, -K \cdot C =4$  &  &  & \\ 
 &  & $E_1:$ 2-33, $p_a(C)=0, -K \cdot C =4$  &  &  & \\ \hline
3-27 & $48$ & $Y = \P^1 \times \P^1 \times \P^1$ && \\
&&$C_2:$ $/\P^1 \times \P^1$ & none &   & \\ 
&  
& $C_2:$ $/\P^1 \times \P^1$  &&  & \\
 &  & $C_2:$ $/\P^1 \times \P^1$ & &  & \\ \hline
3-28 & $48$ & $Y = \F_1 \times \P^1$ & & \\
&& $C_2:$ $/\P^1 \times \P^1$  & none &  \\ 
 &  & $C_2:$ $/\F_1$ &  &  \\ 
 &  & $E_1:$ 2-34, $p_a(C)=0, -K \cdot C =2$  &  & \\ \hline
3-31 & $52$ & $Y =\P_{\P^1 \times \P^1}(\MO_{\P^1 \times \P^1} \oplus \MO_{\P^1 \times \P^1}(1, 1))$ && \\
&& $C_2: /\P^1 \times \P^1$   & none &\\  
&  & $E_1:$ non-Fano, $p_a(C)=0, -K \cdot C =0$  &  & \\
 &  & $E_1:$ non-Fano, $p_a(C)=0, -K \cdot C =0$  &  & \\ \hline
      \end{longtable}
  \end{center} 
Moreover, if the elementary transform $Y'$ of $X \to Y \to \P^1 \times \P^1$ is Fano, then also $Y'$ satisfies one of the above possibilities, and hence $(-K_{Y'})^3 \in \{36, 44, 48, 52\}$. 
\end{nasi}


\begin{nota}\label{n-ele-tr-P1P1}
Set $S:=\P^1 \times \P^1$ and let  $g: Y \to S = \P^1 \times \P^1$ be a Fano $\P^1$-bundle. 
Take a regular subsection $B_Y$ of $g$ and let $\sigma : X \to Y$ be the blowup along $B_Y$. 
Assume that $X$ is a Fano threefold. 
Let $Y'$ be the elementary transform of $f: X \xrightarrow{\sigma} Y \xrightarrow{g} S$: 
\[
\begin{tikzcd}
& X \arrow[ld, "\sigma"'] \arrow[rd, "\sigma'"] \arrow[dd, "f"]\\
Y \arrow[rd, "g"']& & Y' \arrow[ld, "g'"]\\
& S=\P^1 \times \P^1. 
\end{tikzcd}
\]
Set $B:= g(B_Y)$ and $B_{Y'} := \sigma'(\Ex(\sigma'))$, which implies $B_Y \xrightarrow{\simeq} B \xleftarrow{\simeq} B_{Y'}$. 
Let $(d_1, d_2)$ be the bidegree of $B$, i.e., 
$\MO_{\P^1 \times \P^1}(B) \simeq \MO_{\P^1 \times \P^1}(d_1, d_2)$.  
\end{nota}

\begin{rem}\label{r-P1P1-adjunction}
Let $B$ be a curve on $\P^1 \times \P^1$ of bidegree $(d_1, d_2)$. 
By the adjunction formula, we get 
\[
2p_a(B)-2 = (K_{\P^1 \times \P^1}+B) \cdot B = d_1(d_2-2) + d_2(d_1-2) = 2(d_1-1)(d_2-1) -2, 
\]
which implies $p_a(B) = (d_1-1)(d_2-1)$. 
In particular, $p_a(B)=0$ if and only if $d_1  =1$ or $d_2=1$. 
\end{rem}


\begin{lem}\label{l-P1P1-nonFano}
We use Notation \ref{n-ele-tr-P2}. 
Assume that $d_1 \leq d_2$ and $Y'$ is not Fano. 
Then the following holds. 
\begin{enumerate}
\item 
$B \simeq \P^1$ and $g^{-1}(B) \simeq \P^1 \times \P^1$. 
\item $d_1 = 1$ or $d_2=1$. 
\item If $d_1=1$, then 
$\deg B = (1, d_2), B^2=2d_2, -K_S \cdot B= 2d_2+2, -K_{Y/S} \cdot B_Y =2(B^2+1)= 4d_2 +2, 
-K_Y \cdot B_Y = 6d_2+4$. 
\end{enumerate}
\end{lem}

\begin{proof}
Proposition \ref{p-ele-tf-1Fano} implies $B \simeq \P^1, g^{-1}(B)=\P^1 \times \P^1$, and $-K_{Y/S}\cdot B_Y =2 (B^2 +1 )$. 
In particular, (1) holds. 
The assertion (2) follows from (1) and Remark \ref{r-P1P1-adjunction}. 
By direct computation, we see that (3) holds. 
\qedhere





\end{proof}

\begin{lem}\label{l-P1P1-3-27}
We use Notation \ref{n-ele-tr-P1P1}. 
Assume that $Y= \P^1_1 \times \P^1_2 \times \P^1_3$ (i.e., $Y$ is of No.\ 3-27) and 
$g: Y= \P^1_1 \times \P^1_2 \times \P^1_3 \to \P^1_1 \times \P^1_2 = S$ is the projection onto the first and second direct product factors. 
Let $(e_1, e_2, e_3)$ be the tridegree of $B_Y$. 
Assume $e_1 \geq e_2$. 
Then one of the following holds. 
\begin{enumerate}
    \item $(e_1, e_2, e_3) = (1, 0, 0)$, $Y'$ is a Fano threefold of No.\ 3-28, 
    $(-K_X)^3=42$, $\deg B = (0, 1)$, $p_a(B)=0$,  
$-K_Y \cdot B_Y = 2$, and $-K_{Y'} \cdot B_{Y'}= 2$ 
($X$ is 4-10). 
    \item $(e_1, e_2, e_3) = (1, 0, 1)$, $Y'$ is not a Fano threefold, 
    $(-K_X)^3=38$, $\deg B = (0, 1)$, $p_a(B)=0$,  
$-K_Y \cdot B_Y = 4$, and $-K_{Y'} \cdot B_{Y'}= 0$ ($X$ is 4-8). 
    \item $(e_1, e_2, e_3) = (1, 1, 3)$, $Y'$ is not a Fano threefold, 
    $(-K_X)^3=26$, $\deg B = (1, 1)$, $p_a(B)=0$,  
$-K_Y \cdot B_Y = 10$, and $-K_{Y'} \cdot B_{Y'}= 0$ ($X$ is 4-13). 
    \item  $(e_1, e_2, e_3) = (1, 1, 1)$, 
    $Y'$ is a Fano threefold of No.\ 3-25, 
    $(-K_X)^3=34$, $\deg B = (1, 1)$, $p_a(B)=0$,  
$-K_Y \cdot B_Y = 6$, and $-K_{Y'} \cdot B_{Y'}= 4$ 
($X$ is 4-6). 
    \item  
    $(e_1, e_2, e_3) = (2, 2, 2)$, 
    $Y'$ is a Fano threefold of No.\ 3-27 or 3-28, 
    $(-K_X)^3=24$, $\deg B = (2, 2)$, $p_a(B)=1$,  
$-K_Y \cdot B_Y = 12$, and $-K_{Y'} \cdot B_{Y'}= 12$ 
($X$ is 4-1). 
    \item  $(e_1, e_2, e_3) = (2, 1, 1)$, 
    $Y'$ is a Fano threefold of No.\ 3-28, 
    $(-K_X)^3=30$, $\deg B = (1, 2)$, $p_a(B)=0$,  
$-K_Y \cdot B_Y = 8$, and $-K_{Y'} \cdot B_{Y'}= 8$ 
($X$ is 4-3). 
    \item  $(e_1, e_2, e_3) = (3, 1, 1)$,
        $Y'$ is a Fano threefold of No.\ 3-31, 
    $(-K_X)^3=26$, $\deg B = (1, 3)$, $p_a(B)=0$,  
$-K_Y \cdot B_Y = 10$, and $-K_{Y'} \cdot B_{Y'}= 12$ 
($X$ is 4-13). 
    \item  $(e_1, e_2, e_3) = (1, 1, 0)$,        
    $Y'$ is a Fano threefold of No.\ 3-31, 
    $(-K_X)^3=38$, $\deg B = (1, 1)$, $p_a(B)=0$,  
$-K_Y \cdot B_Y = 4$, and $-K_{Y'} \cdot B_{Y'}= 6$ 
($X$ is 4-8). 
    \item $(e_1, e_2, e_3) = (1, 1, 2)$, 
    $Y'$ is a Fano threefold of No.\ 3-17, 
$(-K_X)^3=30$, $\deg B = (1, 1)$, $p_a(B)=0$,  
$-K_Y \cdot B_Y = 8$, and $-K_{Y'} \cdot B_{Y'}= 2$ 
($X$ is 4-3). 
\end{enumerate}
\end{lem}

\begin{proof}
By the isomorphism $g|_{B_Y} : B_Y \xrightarrow{\simeq} B$, 
the tridegree $(e_1, e_2, e_3)$ of $B_Y \subset \P^1_1 \times \P^1_2 \times \P^1_3$ 
satisfies $(e_2, e_1) = (d_1, d_2)$ 
for the bidegreee $(d_1, d_2)$ of $B \subset \P^1_1 \times \P^1_2$ (Remark \ref{r curve bidegree}). 
The following hold (Lemma \ref{l-blowup-formula}, Proposition \ref{p-ele-tf-numbers}(3)(4)): 
\[
B^2 = 2e_1e_2. 
\]
\[
-K_Y \cdot B_Y = 2(e_1+e_2+e_3). 
\]
\[
(-K_{Y'})^3 = (-K_Y)^3 -4 (-K_{Y/S}) \cdot B_Y +2 B^2 
=48 - 8e_3 +4e_1e_2. 
\]
\[
(-K_X)^3= (-K_Y)^3 -2(-K_Y) \cdot B_Y +2p_a(B) -2
=46 -4(e_1+e_2+e_3) +2p_a(B).
\]
\[
(-K_{Y'}) \cdot B_{Y'} = B^2 +2(-K_S) \cdot B -(-K_Y) \cdot B_Y
\]\[
=2e_1e_2 + 4(e_1+e_2) - 2(e_1+e_2+e_3) = 2e_1e_2 +2e_1+2e_2 -2e_3. 
\]

(1) Assume that  $B_Y$ is a fibre of some projection 
$Y =\P^1_1 \times \P^1_2 \times \P^1_3 \to \P^1_i \times \P^1_j$. 
Then  
it follows from $(e_1, e_2) \neq (0, 0)$ and $e_1 \geq e_2$ that 
$(i, j)=(2, 3)$ and $(e_1, e_2, e_3)=(1, 0, 0)$.  
In this case, the given  elemental transform is obtained by applying $(-) \times_k \P^1_1$ to 
the surface elementary transform 
$\P^1_2 \times \P^1_3 \leftarrow T \to \F_1$ over $\P^1_2$. 
In particular, we get $Y' \simeq \F_1 \times \P^1_1$, i.e., 
$Y'$ is a Fano threefolds of No.\ 3-28. 
Moreover,  $B^2 = 0, (-K_Y) \cdot B_Y =2, 
(-K_{Y'})^3 =48$, $(-K_X)^3 =46 -4(e_1+e_2+e_3) +2p_a(B) = 46-4 +0= 42,  
(-K_{Y'}) \cdot B_{Y'} =2$. 

\medskip

(2)-(9) In what follows, we assume that  $B_Y$ is a regular subsection for each projection $Y =\P^1_1 \times \P^1_2 \times \P^1_3 \to \P^1_i \times \P^1_j$. Let $B_{ij} \subset \P^1_i \times \P^1_j$ be the image of $B_Y$. 
By definition, we have $B = B_{12}$. 
By $B_Y \simeq B_{12} \simeq B_{23} \simeq B_{31}$, Remark \ref{r-P1P1-adjunction} implies 
\[
p_a(B_Y) = (e_1-1)(e_2-1) = (e_2-1)(e_3-1) = (e_3-1)(e_1-1). 
\]
In particular, $e_1 = e_2=e_3$ or two of $e_1, e_2, e_3$ are equal to $1$. 
Hence one of (i)-(iv) holds. 
\begin{enumerate}
\item[(i)] $e_1 = e_2 = e_3$. 
\item[(ii)] $e_1 \neq 1$ and $(e_2, e_3) = (1, 1)$. 
\item[(iii)] $e_2 \neq 1$ and $(e_1, e_3) = (1, 1)$. 
\item[(iv)] $e_3 \neq 1$ and $(e_1, e_2)=(1, 1)$. 
\end{enumerate}

\medskip

(2), (3) Assume that $Y'$ is not Fano. 
Lemma \ref{l-P1P1-nonFano}(2) implies $(e_1, e_2) \in \{ (1, 0), (e_1, 1)\}$. 
In order to treat these two cases simultaneously, we set 
\[
d:= \begin{cases}
0 \qquad (\text{if } (e_1, e_2)=(1, 0))\\
e_1 \qquad (\text{if } (e_1, e_2) =(e_1, 1)). 
\end{cases}
\]
Then the bidegree of $B$ is $(1, d)$ or $(d, 1)$. 
It follows from Proposition \ref{p-ele-tf-1Fano}(2) that 
$p_a(B) =0$, $-K_{Y'} \cdot B_{Y'}=0$,  and 
$2e_3 = -K_{Y/S} \cdot B_Y = 2(B^2 +1) =2(2d+1)$, which implies $e_3 = 2d+1$. 
Hence the tridegree of $B_Y$ is $(1, d, 2d+1)$ or $(d, 1, 2d+1)$. 
Since one of (i)-(iv) holds, we have that 
$d=0$ or $d=1$. 
\begin{enumerate}
\item[(2)] 
Assume $d=0$. 
Then 
$(e_1, e_2, e_3) = (1, 0, 1)$, $B^2 = 0, (-K_Y) \cdot B_Y =4$, 
$(-K_X)^3 =46 -4(e_1+e_2+e_3) +2p_a(B) = 46-8 +0= 38$.   
\item[(3)] 
Assume $d=1$. 
Then 
$(e_1, e_2, e_3) = (1, 1, 3)$, $B^2 = 2, (-K_Y) \cdot B_Y =10$, 
$(-K_X)^3 =46 -4(e_1+e_2+e_3) +2p_a(B) = 46-20+0= 26$.   
\end{enumerate}

\medskip

(4)-(9) In what follows, we assume that $Y'$ is Fano. 
In particular, it holds that $(-K_{Y'})^3 \in \{36, 44, 48, 52\}$ (\ref{nasi-pic4-P1P1}). 
We treat the four cases (i)-(iv) separately. 

\medskip


(i) Assume $e:= e_1 = e_2=e_3$. 
We have $\{36, 44, 48, 52\} \ni (-K_{Y'})^3 = 48-8e +4e^2 =4(e-1)^2 +44$. 
Then $(e, (-K_{Y'})^3) \in \{(1, 44), (2, 48)\}$. 
\begin{enumerate}
\item[(4)] If $(e, (-K_{Y'})^3)=(1, 44)$, then 
    $(e_1, e_2, e_3) = (1, 1, 1)$, 
    $Y'$ is a Fano threefold of No.\ 3-25, 
    $\deg B = (1, 1)$, $p_a(B)=0$,  
$-K_Y \cdot B_Y = 6$,   
$(-K_X)^3 =46 -4(e_1+e_2+e_3) +2p_a(B) = 46-12 +0= 34$,  
and $-K_{Y'} \cdot B_{Y'}=2e_1e_2 +2e_1+2e_2 -2e_3= 4$.  
\item[(5)] 
If $(e, (-K_{Y'})^3)=(2, 48)$, then 
    $(e_1, e_2, e_3) = (2, 2, 2)$, 
    $Y'$ is a Fano threefold of No.\ 3-27 or 3-28, 
    $\deg B = (2, 2)$, $p_a(B)=1$,  
$-K_Y \cdot B_Y = 12$, 
$(-K_X)^3 =46 -4(e_1+e_2+e_3) +2p_a(B) = 46-24 +2 =24$,  
and $-K_{Y'} \cdot B_{Y'}= 2e_1e_2 +2e_1+2e_2 -2e_3= 12$.  
\end{enumerate}

(ii) Assume $e_1 \neq 1$ and $(e_2, e_3) = (1, 1)$. 
By $e_1 \neq 1$ and $e_1 \geq e_2=1$, we obtain $e_1 \geq 2$. 
By $\{36, 44, 48, 52\} \ni (-K_{Y'})^3 = 48  -8e_3 + 4e_1e_2 = 40 +4e_1$, 
we get $(e_1, e_2, e_3, (-K_{Y'})^3) \in \{ (2, 1, 1, 48), (3, 1, 1, 52)\}$. 
\begin{enumerate}
\item[(6)] Assume $(e_1, e_2, e_3, (-K_{Y'})^3) = (2, 1, 1, 48)$. 
  Then  $Y'$ is a Fano threefold of No.\ 3-27 or 3-28, 
    $\deg B = (1, 2)$, $p_a(B)=0$,  
$-K_Y \cdot B_Y = 8$, 
    $(-K_X)^3=46 -4(e_1+e_2+e_3) +2p_a(B) = 30$, and $-K_{Y'} \cdot B_{Y'}= 2e_1e_2 +2e_1+2e_2 -2e_3= 8$. 
    Let us show that $Y'$ is of No.\ 3-28. 
    Fix $t \in \P^1_2$. 
    By $B \xrightarrow{\simeq} \P^1_2$,  the given elementary transform induces an elementary transform of ruled surfaces over $\P^1_1$ consisting of 
    the fibres $X_t, Y_t, Y'_t$ over $t \in \P^1_2$. 
    Note that $Y_t$ and $Y'_t$ are $\P^1$-bundles over $\P^1_1$ and $X_t \to Y_t$ and $X_t \to Y'_t$ is a blowup at a point. 
    Since $Y_t \simeq \P^1 \times \P^1$ and $Y'_t \in \{ \P^1 \times \P^1, \F_1\}$, 
    we get $Y'_t \simeq \F_1$, i.e., $Y'$ is 3-28. 
   \item[(7)] Assume $(e_1, e_2, e_3, (-K_{Y'})^3) = (3, 1, 1, 52)$. 
    Then $Y'$ is a Fano threefold of No.\ 3-31, 
    $\deg B = (1, 3)$, $p_a(B)=0$,  
$-K_Y \cdot B_Y = 10$,  
    $(-K_X)^3=46 -4(e_1+e_2+e_3) +2p_a(B) = 26$, 
    and $-K_{Y'} \cdot B_{Y'}= 2e_1e_2 +2e_1+2e_2 -2e_3= 6 + 6 + 2 -2 =12$. 
\end{enumerate}

(iii) 
Assume $e_2 \neq 1$ and $(e_1, e_3) = (1, 1)$. 
In this case, we get $e_2=0$ by $e_1 \geq e_2=1$ and $e_2 \neq 1$. 
Then $\{36, 44, 48, 52\} \ni (-K_{Y'})^3= 48  -8e_3 + 4e_1e_2  = 40$, which is absurd.

\medskip

(iv) 
Assume $(e_1, e_2)=(1, 1)$ and $e_3 \neq 1$. 
By $\{36, 44, 48, 52\} \ni 
(-K_{Y'})^3 = 48  -8e_3 + 4e_1e_2 = 52 -8e_3$, 
we obtain 
$(e_3, (-K_{Y'})^3) \in \{ (0, 52), (2, 36)\}$. 
\begin{enumerate}
\item[(8)] Assume $(e_1, e_2, e_3, (-K_{Y'})^3) = (1, 1, 0, 52)$. 
Then   $Y'$ is a Fano threefold of No.\ 3-31, 
    $\deg B = (1, 1)$, $p_a(B)=0$,  
$-K_Y \cdot B_Y = 4$, 
    $(-K_X)^3=46 -4(e_1+e_2+e_3) +2p_a(B) = 38$, and $-K_{Y'} \cdot B_{Y'}= 2e_1e_2 +2e_1+2e_2 -2e_3= 6$. 
   \item[(9)] Assume $(e_1, e_2, e_3, (-K_{Y'})^3) = (1, 1, 2, 36)$. 
Then    $Y'$ is a Fano threefold of No.\ 3-17, 
    $\deg B = (1, 1)$, $p_a(B)=0$,  
$-K_Y \cdot B_Y = 8$, 
    $(-K_X)^3=46 -4(e_1+e_2+e_3) +2p_a(B) = 46 - 16+0 = 30$, 
    and $-K_{Y'} \cdot B_{Y'}= 2e_1e_2 +2e_1+2e_2 -2e_3=  2$. 
\end{enumerate}
\end{proof}

\begin{lem}\label{l-P1P1-3-28}
We use Notation \ref{n-ele-tr-P1P1}. 
Assume that $Y$ is of No.\ 3-28. 
Let $Y= \F_1 \times \P^1_2 \to \F_1 \to \P^1_1$ be the induced contractions, 
where $\P^1_1 := \P^1$ and $\P^1_2 :=\P^1$. 
Set $B_{\F_1}$ to be the image of $B_Y$ to $\F_1$. 
Then one of the following holds. 
\begin{enumerate}
\item 
    $B_Y$ is a fibre of the first projection $Y = \F_1 \times \P^1 \to \F_1$, 
    $Y'$ is a Fano threefold of No.\ 3-27, 
    $(-K_X)^3=42$, $\deg B = (1, 0)$, $p_a(B)=0$,  
$-K_Y \cdot B_Y = 2$, and $-K_{Y'} \cdot B_{Y'}= 2$.  
($X$ is 4-10). 
\item 
$B_{\F_1}$ is the $(-1)$-curve on $\F_1$, 
   $Y'$ is a Fano threefold of No.\ 3-31, 
    $(-K_X)^3=44$, $\deg B = (0, 1)$, $p_a(B)=0$, 
$-K_Y \cdot B_Y = 1$, and $-K_{Y'} \cdot B_{Y'}= 3$.  
($X$ is 4-11). 
\item 
$B_{\F_1}$ is disjoint from the $(-1)$-curve on $\F_1$, 
$Y'$ is a Fano threefold of No.\ 3-27, 
    $(-K_X)^3=30$, $\deg B = (1, 2)$, $p_a(B)=0$, 
$-K_Y \cdot B_Y = 8$, and $-K_{Y'} \cdot B_{Y'}= 8$.  
($X$ is 4-3). 
\item 
$B_{\F_1}$ is disjoint from the  $(-1)$-curve on $\F_1$, 
$Y'$ is a Fano threefold of No.\ 3-25, 
    $(-K_X)^3=40$, $\deg B = (0, 1)$, $p_a(B)=0$, 
$-K_Y \cdot B_Y = 3$, and $-K_{Y'} \cdot B_{Y'}= 1$.  
($X$ is 4-9). 
\item 
$B_{\F_1}$ is disjoint from the $(-1)$-curve on $\F_1$, 
$Y'$ is a Fano threefold of No.\ 3-28, 
    $(-K_X)^3=36$, $\deg B = (1, 1)$, $p_a(B)=0$, 
$-K_Y \cdot B_Y = 5$, and $-K_{Y'} \cdot B_{Y'}= 5$.  
($X$ is 4-7). 
\item 
$B_{\F_1}$ is disjoint from the $(-1)$-curve on $\F_1$, 
$Y'$ is a Fano threefold of No.\ 3-31, 
    $(-K_X)^3=32$, $\deg B = (2, 1)$, $p_a(B)=0$, 
$-K_Y \cdot B_Y = 7$, and $-K_{Y'} \cdot B_{Y'}= 9$.  
($X$ is 4-5). 
\end{enumerate}
\end{lem}



\begin{proof}
Let $\pi : \F_1 \to \P^1_1$ be the contraction. 
For the $(-1)$-curve $\Gamma$ on $\F_1$, 
we have $-K_{\F_1} \sim 2\Gamma +\pi^*\MO_{\P^1_1}(3)$   
and $-K_Y \sim {\rm pr}_1^*(-K_{\F_1}) + {\rm pr}_2^*(-K_{\P^1_2}) 
\sim 3H_1 + 2 H_2 + 2\Gamma_Y$, where $\Gamma_Y := \pr_1^*\Gamma$ and $H_i$ denotes the pullback of $\MO_{\P^1_i}(1)$ on $Y$. 
Set 
$\gamma :=  \Gamma_Y \cdot B_Y$. 
We have that $-K_Y \cdot B_Y= 2d_1 + 3d_2 + 2 \gamma$ and 
$-K_{Y/S} \cdot B_Y= d_2 + 2\gamma$. 
The following hold  (Lemma \ref{l-blowup-formula}, Proposition \ref{p-ele-tf-numbers}(3)(4)): 
\[
(-K_{Y'})^3 = (-K_Y)^3 -4(-K_{Y/S}) \cdot B_Y + 2B^2 = 
48 - 4(d_2 + 2 \gamma)  +4 d_1d_2 = 48 -4d_2 -8 \gamma  +4d_1d_2. 
\]
\[
(-K_X)^3 = (-K_Y)^3 -2(-K_Y) \cdot B_Y +2p_a(B)-2 = 
 46  -2(-K_Y) \cdot B_Y +2p_a(B), 
\]
\[
(-K_{Y'}) \cdot B_{Y'} 
= B^2 +2(-K_S) \cdot B -(-K_Y) \cdot B_Y
\]
\[
= 2d_1d_2 +4(d_1+d_2) - (2d_1 + 3d_2 + 2 \gamma) 
=2d_1d_2 +2d_1 +d_2 -2 \gamma. 
\]
One of the following holds (Proposition \ref{p-FCB-centre}).  
\begin{enumerate}
    \item[(1)] $B_Y$ is a fibre of ${\rm pr}_1 : Y = \F_1 \times \P^1_2 \to \F_1$. 
    \item[(2)] $\pr_1|_{B_Y} : B_Y\xrightarrow{\simeq} B_{\F_1}$ and $B_{\F_1} = \Gamma$. 
    \item[(3)-(6)] $\pr_1|_{B_Y} : B_Y\xrightarrow{\simeq} B_{\F_1}$ and $B_{\F_1} \neq \Gamma $. 
\end{enumerate}

(1) Assume that $B_Y$ is a fibre of ${\rm pr}_1 : Y = \F_1 \times \P^1_2 \to \F_1$. 
In this case, $\gamma=0$ and $-K_Y \cdot B_Y =2$. 
By $2 = -K_Y \cdot B_Y= 2d_1 + 3d_2$, 
we obtain $(d_1, d_2) = (1, 0)$. 
Then $-K_{Y/S} \cdot B_Y= 0$, 
$(-K_{Y'})^3 =  48 -4d_2 -8 \gamma  +4d_1d_2 = 48$, 
$(-K_X)^3 = 42$, 
$(-K_{Y'}) \cdot B_{Y'} =2d_1d_2 +2d_1 +d_2 -2 \gamma = 2$. 
In this case, the given  elementary transform is obtained by applying the case change $(-) \times_k \P^1$ to 
the surface elementary transform between $\F_1$ and $\P^1 \times \P^1$. 
Thus $Y' \simeq \P^1 \times \P^1 \times \P^1$. 

(2) Assume that $B_Y\xrightarrow{\simeq} B_{\F_1}$ and $B_{\F_1} = \Gamma$. 
In this case, $\gamma =\Gamma^2 = -1$ and $p_a(B)=0$. 
Let us show that $Y'$ is Fano. 
If $Y'$ is non-Fano, then 
we would get 
$2d_1 + 3d_2 -2 =-K_Y \cdot B_Y= 6d_2 + 4$ (Lemma \ref{l-P1P1-nonFano}(3)), which implies $d_2 \in 2\Z$. 
Then Lemma \ref{l-P1P1-nonFano}(2) implies $d_1 = 1$. 
However, we obtain $3d_2 = 2d_1 + 3d_2 -2 = 6d_2 + 4$, which is absurd. 
Hence $Y'$ is Fano. Then (\ref{nasi-pic4-P1P1}) implies 
\[
\{ 36, 44, 48, 52 \} \ni (-K_{Y'})^3 = 56 + 4d_2(d_1-1). 
\]
In particular, we get $d_1=0$, which implies $d_2=1$ 
and $(-K_{Y'})^3 = 52$. 
Then $Y'$ is of No.\ 3-31, 
$p_a(B) = 0$, 
$-K_Y \cdot B_Y = 1$, 
$-K_{Y/S} \cdot B_Y= -1$, 
$(-K_X)^3 = 46 - 2 +0 =44$, 
$(-K_{Y'}) \cdot B_{Y'} =2d_1d_2 +2d_1 +d_2 -2 \gamma = 3$.

(3)-(6) 
Assume that $B_Y\xrightarrow{\simeq} B_{\F_1}$ and $B_{\F_1} \neq \Gamma $. 
In this case, $B_Y \cap \Gamma_Y = \emptyset$ and  $\gamma =0$ (Corollary \ref{c-FCB-(-1)-2}). 
For the birational contraction $\tau : \F_1 \to \P^2$ and $B_{\P^2} := \tau(B_{\F_1}) (\simeq B_{\F_1})$, 
we obtain $\deg B_{\P^2} = B_{\F_1} \cdot \tau^* \MO_{\P^2}(1) = B_{\F_1} \cdot \pi^* \MO_{\P^1_1}(1) = B_Y \cdot H_1 = d_2$ 
and $p_a(B_Y) = p_a(B_{\F_1}) = p_a(B_{\P^2}) = \frac{1}{2}(d_2-1)(d_2-2)$. 
This, together with $p_a(B_Y) = p_a(B) = (d_1-1)(d_2-1)$ (Remark \ref{r-P1P1-adjunction}),
 implies that either 
\begin{enumerate}
\item[(3)] $d_2 \neq 1$ and $2d_1 = d_2$, or 
\item[(4)-(6)]$d_2 =1$. 
\end{enumerate}

(3) Assume $2d_1 = d_2$. We obtain $-K_Y \cdot B_Y = 8d_1$. 
Lemma \ref{l-P1P1-nonFano} implies that $Y'$ is Fano. 
By (\ref{nasi-pic4-P1P1}),  we get 
\[
\{ 36, 44, 48, 52 \} \ni (-K_{Y'})^3 = 48 -8d_1 +8d_1^2 = 48 +8d_1(d_1-1) \in 8\Z. 
\]
Then $(d_1, d_2, (-K_{Y'})^3) = (1, 2, 48)$. 
Thus $p_a(B)=0, -K_Y \cdot B_Y = 8$, 
$-K_{Y/S} \cdot B_Y= 2$, 
$(-K_X)^3 = 46 - 16 +0 =30$, 
$(-K_{Y'}) \cdot B_{Y'} =2d_1d_2 +2d_1 +d_2 -2 \gamma = 8$. 
We see that $Y'$ is of No.\ 3-27 or 3-28. 
By the same argument as in Lemma \ref{l-P1P1-3-27}(6), $Y'$ is of No.\ 3-27. 

\medskip

(4)-(6) Assume $d_2=1$. We obtain $p_a(B)=0$, $-K_{Y/S} \cdot B_Y = 1$, 
$-K_Y \cdot B_Y = 3+2d_1 \not\in 2\Z$. Lemma \ref{l-P1P1-nonFano} implies that $Y'$ is Fano. 
Hence (\ref{nasi-pic4-P1P1}) implies 
\[
\{ 36, 44, 48, 52 \} \ni (-K_{Y'})^3 = 44 + 4d_1. 
\]
Then $(d_1, (-K_{Y'})^3) \in \{ (0, 44), (1, 48), (2, 52)\}$. 
We have $(-K_X)^3 = 40 -4d_1$. 

(4)
Assume $(d_1, d_2,  (-K_{Y'})^3) =  (0, 1, 44)$. 
Then $Y'$ is a Fano threefold  of No.\ 3-25, 
$-K_Y \cdot B_Y = 3$, 
$(-K_X)^3 = 40-0 =40$, 
$(-K_{Y'}) \cdot B_{Y'} = B^2 +2(-K_S) \cdot B -(-K_Y) \cdot B_Y = 0 + 4 - 3=1$. 

(5) 
Assume $(d_1, d_2,  (-K_{Y'})^3) =  (1, 1, 48)$. 
Then 
$-K_Y \cdot B_Y = 5$, 
$(-K_X)^3 = 40 - 4 =36$, 
$(-K_{Y'}) \cdot B_{Y'} = B^2 +2(-K_S) \cdot B -(-K_Y) \cdot B_Y = 2 + 8- 5 = 5$.  
Then $Y'$ is a Fano threefold, which is of No.\ 3-27 or 3-28. 
Note that $Y'$ is not of No.\ 3-27, as otherwise
the equality $(-K_X)^3 = 36$ would contradict Lemma \ref{l-P1P1-3-27}). 
Hence $Y'$ is of No.\ 3-28.

(6) 
Assume $(d_1, d_2,  (-K_{Y'})^3) =  (2, 1, 52)$. 
Then 
$Y'$ is a Fano threefold of No.\ 3-31, 
$-K_Y \cdot B_Y = 7$, 
$(-K_X)^3 = 40 - 8 =32$, 
$(-K_{Y'}) \cdot B_{Y'} = B^2 +2(-K_S) \cdot B -(-K_Y) \cdot B_Y = 4 + 12- 7 = 9$.  
\end{proof}

\begin{lem}\label{l-P1P1-3-17}
We use Notation \ref{n-ele-tr-P1P1}. 
Assume that $Y$ is of No.\ 3-17. 
Then $Y'$ is a Fano threefold of No.\ 3-27.  
\end{lem}

\begin{proof}
By Proposition \ref{p-pic3-17} or Subsection \ref{ss-table-pic3}, 
there is a conic bundle $h : Y \to \P^2$ with $\Delta_h \neq 0$. 
Hence the blowup centre $B_Y$ of $\sigma : X \to Y$ is a smooth fibre of $h$. 
In particular, $B_Y \simeq \P^1$, $-K_Y \cdot B_Y =2$, 
and $(-K_X)^3 = (-K_Y)^3 -6  = 30$. 
By using the closed immersion $Y \hookrightarrow Z_1 \times Z_2 \times Z_3 = \P^1 \times \P^1 \times \P^2$ 
as in Proposition \ref{p-pic3-17}, 
we have the tridegree $(d_1, d_2, d_3)$ of $B_Y$ for $d_i := H_i \cdot B_Y$, 
where $H_i$ denotes the pullback of the ample generator of $Z_i$ 
(note that $(d_2, d_1)$  coincides with the bidegree of $B$ in $\P^1 \times \P^1 = Z_1 \times Z_2$). 
Since $B_Y$ is a fibre of $Y \to Z_3$, we have $d_3 =0$. 
By Proposition \ref{p-pic3-17}(3), it holds that $2 = -K_Y \cdot B_Y = (H_1 +H_2 + 2H_3) \cdot B_Y = d_1 + d_2$. 
On the other hand, the inclusion 
$\P^1 \simeq B \subset S= \P^1 \times \P^1$ implies 
that $d_1=1$ or $d_2=1$ (Remark \ref{r-P1P1-adjunction}). 
Hence $(d_1, d_2) = (1, 1)$. 
Then $B^2 = 2$ and $-K_{Y/S} \cdot B_Y = -K_Y \cdot B_Y -  (-K_S) \cdot B = 2 -4 =-2$. 
Moreover, Proposition \ref{p-ele-tf-numbers} implies 
\[
(-K_{Y'})^3 = (-K_Y)^3 -4(-K_{Y/S}) \cdot B_Y + 2B^2 = 
36 - 4 \cdot (-2) +2 \cdot 2 = 48. 
\]
By Lemma \ref{l-P1P1-nonFano} and $-K_Y \cdot B_Y =2$, $Y'$ is Fano. 
It follows from Lemma \ref{l-P1P1-3-28} that $Y'$ is not of No.\ 3-28. 
Hence $Y'$ is of No.\ 3-27.  
\end{proof}

\begin{lem}\label{l-P1P1-3-31}
We use Notation \ref{n-ele-tr-P1P1}. 
Assume that $Y$ is of No.\ 3-31 and $(-K_{Y'})^3 \neq 48$. 
Then $Y'$ is a Fano threefold of No.\ 3-31, 
    $(-K_X)^3=28$, $\deg B = (2, 2)$, $p_a(B)=1$,  
$-K_Y \cdot B_Y = 12$, and $-K_{Y'} \cdot B_{Y'}= 12$.  
($X$ is 4-2).   
\end{lem}

\begin{proof}
First of all, we prove that $Y'$ is Fano. 
Recall that we have the $\P^1$-bundle structure $\pi : Y = \P_{\P^1 \times \P^1}(\MO_{\P^1 \times \P^1} \oplus \MO_{\P^1 \times \P^1}(1, 1)) \to \P^1 \times \P^1$ (Subsection \ref{ss-table-pic3}). 
The base change $g^{-1}(B) = \P_B(\MO_B \oplus (\MO_{\P^1 \times \P^1}(1, 1)|_B))$ is not isomorphic to $\P^1 \times \P^1$, as $\MO_{\P^1 \times \P^1}(1, 1)$ is ample. 
Thus $Y'$ is Fano (Lemma \ref{l-P1P1-nonFano}(1)).


Let $\varphi : Y = \P_{\P^1 \times \P^1}(\MO \oplus \MO(1, 1)) \to Z$ be the contraction of a $\pi$-section $D \simeq \P^1 \times \P^1$ to a point (Proposition \ref{p-pic3-31}). 
Recall that $-K_Y|_D \sim -D|_D \sim \MO_{\P^1 \times \P^1}(1, 1)$ (Lemma \ref{l-nonFano-blowdown}). 
For fibres $L_1, L_2$ of the first and second projections $\pr_1, \pr_2 : D =\P^1 \times \P^1 \to \P^1$, 
we obtain $-K_Y \cdot L_1 = -K_Y \cdot L_2=1$. 
Then it follows from Lemma \ref{l-line-meeting} that $D \cap B_Y = \emptyset$. 
Recall that $-K_Y \sim 3H_1 + 3H_2 + 2D$ (Proposition \ref{p-pic3-31}(3)). 
Thus $-K_Y \cdot B_Y = 3d_1 + 3d_2$ and $-K_{Y/S} \cdot B_Y = d_1+d_2$. 
Then Proposition \ref{p-ele-tf-numbers}(3)(4) implies 
\[
(-K_{Y'})^3 = (-K_Y)^3 -4(-K_{Y/S}) \cdot B_Y + 2B^2 
\]
\[= 
52 - 4(d_1+d_2) +2 \cdot 2 d_1d_2 = 48 +4(d_1 -1)(d_2-1). 
\]
\[
-K_{Y'} \cdot B_{Y'} = B^2 +2 (-K_S) \cdot B -(-K_Y) \cdot B_Y
\]
\[
= 2d_1d_2 +4(d_1+d_2) - 3(d_1+d_2) = 2d_1d_2 +(d_1+d_2). 
\]
We have $(d_1 -1) (d_2-1) \geq 0$, because  
$d_1d_2=0$  implies $(d_1, d_2) \in \{ (1, 0), (0, 1)\}$. 
Thus we obtain $(-K_{Y'})^3 \geq 48$. 
By $(-K_{Y'})^3 \in \{ 36, 44, 48, 52\}$ (\ref{nasi-pic4-P1P1}) and our assumption 
$(-K_{Y'})^3 \neq 48$, 
we obtain $(-K_{Y'})^3 = 52$, i.e.,  also $Y'$ is a Fano threefold of No.\ 3-31. 
Then  $(d_1, d_2)=(2, 2)$, $p_a(B)=1$,  
$-K_Y \cdot B_Y = 12$, $-K_{Y'} \cdot B_{Y'}= 12$, and 
$(-K_X)^3 = (-K_Y)^3 +2p_a(B_Y) -2 -2 (-K_Y) \cdot B_Y = 
52 +2-2 -24 = 28$ (Lemma \ref{l-blowup-formula}). 
\qedhere

\end{proof}

\begin{lem}\label{l-P1P1-3-25}
We use Notation \ref{n-ele-tr-P1P1}. 
Assume that $Y$ is No.\ 3-25, i.e., $Y = \Bl_{L_1 \amalg L_2} \P^3$  
for some lines $L_1$ and $L_2$ on $\P^3$ such that $L_1 \cap L_2 =\emptyset$. 
Let $\sigma: Y \to \P^3$ be the induced blowup. 
Then the following hold. 
\begin{enumerate}
\item 
One of the following holds. 
\begin{enumerate}
\item $B_Y$ is a one-dimensional fibre of $\sigma: Y \to \P^3$. In particular, $-K_Y \cdot B_Y =1$. 
\item $B_Y \cap \Ex(\sigma) = \emptyset$, $p_a(B)=0$ and  $-K_Y \cdot B_Y \in \{4,  12\}$.  
\item $B_Y \cap \Ex(\sigma) = \emptyset$, $p_a(B)=1$ and  $-K_Y \cdot B_Y  = 16$.  
\end{enumerate} 
\item 
$Y'$ is Fano. 
\item 
$Y'$ is not of No.\ 3-25. 
\end{enumerate}
\end{lem}

\begin{proof}
Let us show (1). 
Assume that $B_Y$ intersects $\Ex(\sigma)$. 
Then $B_Y$ must be a one-dimensional fibre of $\sigma$ (Lemma \ref{l-line-meeting}), i.e., (a) holds. 
Hence we may assume that $B_Y$ is disjoint from $\Ex(\sigma)$. 
For $B_{\P^3} := \sigma(B_Y)$, 
it follows from Proposition \ref{p-nonFano-iff} that 
$\Bl_{L_1 \amalg B_{\P^3}} \P^3$ is Fano. 
By Corollary \ref{c-P3-disjoint}. 
$B_Y$ is one of the following: (i) an elliptic curve of degree $4$, 
(ii) a cubic rational curve, 
(iii) a conic, and (iv) a line. 
The case (iii) contradicts Lemma \ref{l-pic4-3-18}. 
Hence one of (i), (ii), (iv) holds. 
Then (b) or (c) holds.  
This completes the proof of (1).


Let us show (2). 
Suppose that $Y'$ is not Fano. 
By symmetry, we may assume $d_1 \leq d_2$. 
If $d_1=1$, then 
we get $B \simeq \P^1$ and $-K_Y \cdot B_Y = 6d_2 +4$ (Lemma \ref{l-P1P1-nonFano}), which contradicts (1). 
Again by Lemma \ref{l-P1P1-nonFano}, we obtain $(d_1, d_2) =(0, 1)$. 
It follows from Proposition \ref{p-ele-tf-1Fano} that 
$-K_Y \cdot B_Y -2 = -K_Y \cdot B_Y -(-K_S) \cdot B= 2(B^2 +1) = 2$. 
Then (b) of (1) holds. 
Hence $B_{\P^3} = \sigma(B_Y)$ is a line. 
On the other hand, for the conic bundle structure $\pi :Y \to \P^1 \times \P^1$, 
each composite morphism $\pr_i \circ \pi : Y \to \P^1 \times \P^1 \to \P^1$ is a del Pezzo fibration. 
This del Pezzo fibration $\pr_i \circ \pi$ is induced by the pencil consisting of the planes containing $L_i$ after permuting the direct product factors of $\P^1 \times \P^1$ if necessary (because 
we have a factorisation 
$\pr_i \circ \pi : Y = \Bl_{L_1 \amalg L_2} \to \Bl_{L_i}\,\P^3 \to \P^1$ consisting of contractions by Proposition \ref{p-pic3-25}). 
In particular, $B_Y$ dominates each $\P^1$, i.e., $d_1 >0$ and $d_2 >0$. 
This is absurd. 
Thus (2) holds.


Let us show (3). 
Suppose that $Y'$ is a Fano threefold of No.\ 3-25. 
By $(-K_{Y'})^3 = (-K_Y)^3$ and $(-K_{Y'})^3 = (-K_Y)^3 -4(-K_{Y/S} \cdot B_Y) +2B^2$ (Proposition \ref{p-ele-tf-numbers}), we obtain $B^2 +2(-K_S) \cdot B = 2(-K_{Y}) \cdot B_Y$. 
It follows from $B^2 = 2d_1d_2$ and $(-K_S) \cdot B = 2(d_1+d_2)$ that 
\[
(d_1+2)(d_2+2) =  (-K_{Y}) \cdot B_Y +4. 
\]

(a) Assume that $B_Y$ intersects $\Ex(\sigma)$. 
By (1), we get $-K_Y \cdot B_Y = 1$. 
Hence $(d_1+2)(d_2+2) =  5$. 
There is no solution satisfying $(d_1, d_2) \in \Z_{\geq 0} \times \Z_{\geq 0}$.

(b) Assume that $B_Y \cap \Ex(\sigma) = \emptyset$ and $p_a(B)=0$. 
By (1), we get $-K_Y \cdot B_Y \in \{4, 12\}$. 
Then 
\[
(d_1+2)(d_2+2) =  (-K_{Y}) \cdot B_Y +4 \in \{ 8, 16\}. 
\]
By $d_1 \geq 0$, $d_2 \geq 0$, and $(d_1, d_2) \not\in \{ (0, 2), (2, 0)\}$,  
 there is no solution for $(d_1+2)(d_2+2)=8$. 
 Thus $(d_1+2)(d_2+2)=16$. 
 Then we obtain $(d_1, d_2)=(2, 2)$, and hence $B$ is an elliptic curve. 
However, this contradicts $p_a(B_Y)=0$. 

(c) Assume that $B_Y \cap \Ex(\sigma) = \emptyset$ and $p_a(B)\neq 0$. 
By (1), we get $p_a(B)=1$ and $-K_Y \cdot B_Y =16$. 
Then 
$(d_1+2)(d_2+2) = 20$. 
Then $(d_1, d_2) \in \{ (2, 3), (3, 2)\}$. 
In any case, we obtain $1 =p_a(B) = (d_1-1)(d_2-1) = 2$ (Remark \ref{r-P1P1-adjunction}), which is absurd. 
\qedhere 

\end{proof}

\begin{thm}\label{t-ele-tr-P1P1-pic4}
Let 
\[
\begin{tikzcd}
& X \arrow[ld, "\sigma"'] \arrow[rd, "\sigma'"] \arrow[dd, "f"]\\
Y \arrow[rd, "g"']& & Y' \arrow[ld, "g'"]\\
& S:=\P^1 \times \P^1 
\end{tikzcd}
\]
be an elementary transform of threefold conic bundles 
such that $\rho(X)=4$ and 
each of $f: X \to S=\P^1 \times \P^1$ and $g : Y \to S=\P^1 \times \P^1$ is a Fano conic bundle  (cf. Definition \ref{d-ele-tf}). 
Let $B_Y$ be the blowup centre of $\sigma: X \to Y$. 
Set $B:=g(B_Y)$. 
Assume that 
\begin{itemize}
    \item  $d_1 \leq d_2$ for the bidegree $(d_1, d_2)$ of $B$, and 
    \item $(-K_Y)^3 \leq (-K_{Y'})^3$ when also $Y'$ is Fano. 
\end{itemize} 
Then one of the following holds. 
    \begin{center}
\begin{longtable}{ccccccccc}
$X$ & $Y$ & $Y'$ & $(-K_X)^3$& $\deg B$ & $p_a(B)$ & $-K_Y \cdot B_Y$ & $-K_{Y'} \cdot B_{Y'}$ &\\ \hline
4-1 & 3-27 & 3-27 & $24$ & $(2, 2)$ & $1$ & $12$ & $12$\\ \hline
4-2 & 3-31 & 3-31 & $28$ & $(2, 2)$ & $1$ & $12$ & $12$\\ \hline
4-3 & 3-17 & 3-27 & $30$ & $(1, 1)$ & $0$ & $2$ & $8$\\ \hline
4-3 & 3-27 & 3-28 & $30$ & $(1, 2)$ & $0$ & $8$ & $8$\\ \hline
4-5 & 3-28 & 3-31 & $32$ & $(1, 2)$ & $0$ & $7$ & $9$\\ \hline
4-6 & 3-25 & 3-27 & $34$ & $(1, 1)$ & $0$ & $4$ & $6$\\ \hline
4-7 & 3-28 & 3-28 & $36$ & $(1, 1)$ & $0$ & $5$ & $5$\\ \hline
4-8 & 3-27 & 3-31 & $38$ & $(1, 1)$ & $0$ & $4$ & $6$\\ \hline
4-8 & 3-27 & non-Fano & $38$ & $(0, 1)$ & $0$ & $4$ & $0$\\ \hline
4-9 & 3-25 & 3-28 & $40$ & $(0, 1)$ & $0$ & $1$ & $3$\\ \hline
4-10 & 3-27 & 3-28 & $42$ & $(0, 1)$ & $0$ & $2$ & $2$\\ \hline
4-11 & 3-28 & 3-31 & $44$ & $(0, 1)$ & $0$ & $1$ & $3$\\ \hline
4-13 & 3-27 & 3-31 & $26$ & $(1, 3)$ & $0$ & $10$ & $12$\\ \hline
4-13 & 3-27 & non-Fano & $26$ & $(1, 1)$ & $0$ & $10$ & $0$\\ \hline
 \\
 \caption{Elementary transforms over $\P^1 \times \P^1$ ($\rho(X)=4$)}\label{table-ele-tr-P1P1}
      \end{longtable}
  \end{center} 
\end{thm}

We say that the above diagram is called an elemental transform over $\P^1 \times \P^1$ 
{\em of type 3-xx-vs-3-yy} if $Y$ is 3-xx and $Y'$ is 3-yy. 
In this case, we say that $X$ has a conic bundle structure over $\P^2$ {\em of type 3-xx-vs-3-yy}. 

\begin{proof}
Recall that $Y$ is of No.\ 3-17, 3-25, 3-27, 3-28, or 3-31 (\ref{nasi-pic4-P1P1}). 
The same conclusion holds for $Y'$ when $Y'$ is Fano. 
If one of $Y$ and $Y'$ is a Fano threefold of No.\ 3-27 (resp. 3-28), then 
the assertion  follows from 
Lemma \ref{l-P1P1-3-27} (resp. Lemma \ref{l-P1P1-3-28}). 
Here note that we have $(d_1, d_2) = (e_2, e_1)$ when we apply Lemma \ref{l-P1P1-3-27} (Remark \ref{r curve bidegree}). 
Then the case when $Y$ or $Y'$ is  Fano threefold of No.\  3-17 (resp. 3-31) is settled by Lemma \ref{l-P1P1-3-17} (resp. Lemma \ref{l-P1P1-3-31}). 
Therefore, we may assume that $Y$ is of No.\ 3-25, and hence we are done by Lemma \ref{l-P1P1-3-25}. 
\end{proof}

\subsection{Fano conic bunldes over $\mathbb F_1$ ($\rho=4$)}\label{ss-pic4-F1}

The purpose of this subsection is to classify Fano conic bundles 
$f: X \to \F_1$ with $\rho(X)=4$. 
Such a conic bunlde is obtained from a Fano conic bundle 
$g: Y \to \F_1$ with $\rho(Y) =3$ by taking 
a  blowup of $Y$ along a regular subsection $B_Y$ of $g$ (Proposition \ref{p-smaller-same-base}). 
Thus we work with the following situation.

\begin{nota}\label{n-ele-tr-F1-pic4}
Let $g: Y \to S$ be a Fano conic bundle over $S:=\F_1$ with $\rho(Y)=3$. 
Let $B_Y$ be a regular subsection of $g$ and let $\sigma : X \to Y$ be the blowup along $B_Y$. 
Assume that $X$ is a Fano threefold. 
Let $Y'$ be the elementary transform of $f: X \xrightarrow{\sigma} Y \xrightarrow{g} S$ (Definition \ref{d-ele-tf}): 
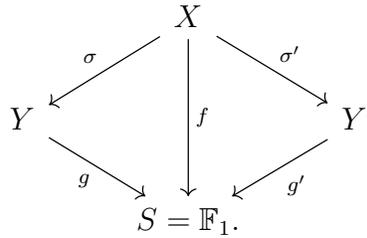
\[
\begin{tikzcd}
& X \arrow[ld, "\sigma"'] \arrow[rd, "\sigma'"] \arrow[dd, "f"]\\
Y \arrow[rd, "g"']& & Y' \arrow[ld, "g'"]\\
& S=\F_1. 
\end{tikzcd}
\]
Set $B:= g(B_Y)$ and $B_{Y'} := \sigma'(\Ex(\sigma'))$, which induces 
$B_Y \xrightarrow{\simeq} B \xleftarrow{\simeq} B_{Y'}$. 
Recall that we have $Y \simeq \widetilde Y \times_{\widetilde S} S$ 
for some Fano conic bundle $\widetilde g : \widetilde Y \to \widetilde S :=\P^2$ and 
the blowdown $\tau : S =\F_1 \to \widetilde S = \P^2$ of the $(-1)$-curve $\Gamma$ on $\F_1$ 
(Lemma \ref{l-F1-pic3}): 
\[
\begin{tikzcd}
Y \arrow[r, "\tau_Y"] \arrow[d, "g"] & \widetilde Y \arrow[d, "\widetilde g"] \\
S=\F_1 \arrow[r, "\tau"] & \widetilde S=\P^2. 
\end{tikzcd}
\]
The possibilities for $Y$ and $\widetilde Y$ are as follows (Theorem \ref{t-F1-pic3}): 

  \begin{center}
\begin{longtable}{cccccc}
$Y$ & $(-K_Y)^3$ & $\widetilde Y$ & $(-K_{\widetilde Y})^3$ & $\widetilde Y\to \widetilde S$ & $\deg \Delta_{\widetilde g}$\\ \hline
3-4 & $18$ & 2-18 & $24$ & $\widetilde Y \xrightarrow{2:1} \P^2\times \P^1 \xrightarrow{\pr_1} \P^2$ & $4$\\ \hline
3-8 & $24$ & 2-24 & $30$ & $\widetilde Y \hookrightarrow \P^2\times \P^2 \xrightarrow{\pr_1} \P^2, \widetilde Y \in |\MO(1, 2)|$ & $3$\\ \hline
3-24 & $42$ & 2-32 & $48$ & $\widetilde Y =W \hookrightarrow \P^2 \times \P^2 \xrightarrow{{\rm pr}_i} \P^2$ & $0$\\ \hline
3-28 & $48$ & 2-34 & $54$ & ${\rm pr}_1 : \widetilde Y=\P^2 \times \P^1 \to \P^2$ & $0$ \\ \hline
3-30 & $50$ & 2-35 & $56$ & $\widetilde Y=\P_{\P^2}(\MO_{\P^2} \oplus \MO_{\P^2}(1)) \to \P^2$ & $0$\\ \hline
      \end{longtable}
  \end{center} 
  If $Y'$ is Fano, then 
  the same conclusion holds for $Y'$, and hence $(-K_{Y'})^3 \in \{18, 24, 42, 48, 50\}$. 

\end{nota}

\begin{lem}\label{l-CC-2blowups}
Let $Z$ be a Fano threefold with $\rho(Z)=2$ 
such that  both extremal rays $R_1$ and $R_2$ of $\NE(Z)$ are of type $C$. 
Let $h_1 : Z \to \P^2$ and $h_2:Z \to \P^2$ be the contractions of $R_1$ and $R_2$, respectively. 
Let $\sigma': Y \to Z$ be a blowup along a smooth fibre $\Gamma_1$ of $h_1$ 
and let $\sigma : X \to Y$ be a blowup along a smooth curve on $Y$. 
Assume that $X$ is Fano. 
Then the composition $X \xrightarrow{\sigma} Y \xrightarrow{\sigma'} Z$ is 
a blowup along a disjoint union $\Gamma_1 \amalg \Gamma_2$, 
where $\Gamma_2$ is a smooth fibre of $h_2$. 
Moreover, both $h_1$ and $h_2$ are of type $C_2$. 
\end{lem}

\begin{proof}
For each $i \in \{1, 2\}$, let $g_i : Y \xrightarrow{\sigma'} Z \xrightarrow{h_i} \P^2$ be the composition. 
Since the blowup centre $\Gamma_1$ of $\sigma' : Y \to Z$ is a smooth fibre of $f_1 : Z \to \P^2$, 
$\Gamma_1$ is a regular subsection of $h_2 : Z \to \P^2$. 
Hence $h_2 : Z \to \P^2$ is of type $C_2$ (Proposition \ref{p-FCB-centre}(2)) and the discriminant divisor $\Delta_{g_2}$ on $\P^2$ is ample. 
Then the centre of the second blowup $\sigma : X \to Y$ must be a smooth fibre 
of the conic bundle $g_2 : Y \to \P^2$  (Proposition \ref{p-FCB-centre}(2)), 
which is disjoint from the exceptional divisor $\Ex(\sigma')$  (Proposition \ref{p-FCB-centre}(1)). 
By symmetry, also $h_1 : Z \to \P^2$ is of type $C_2$. 
\end{proof}

\begin{lem}\label{l-ele-tf-F1-3-4,8}
We use Notation \ref{n-ele-tr-F1-pic4}. 
Then $Y$ is of neither No.\ 3-4 nor 3-8. 
\end{lem}

\begin{proof}
Suppose that $Y$ is 3-4. 
Recall that $Y$
has a conic bundle structure $h : Y \to \P^1 \times \P^1$ 
such that the discriminant divisor $\Delta_h$  is ample (Subsection \ref{ss-table-pic3}). 
It follows from Proposition \ref{p-FCB-centre} that the blowup centre $B_Y$ of $\sigma : X \to Y$ is a smooth fibre of $h$ and we have a Fano conic bundle $X \to T$ with a smooth del Pezzo surface $T$ with $K_T^2 =7$. 
Then this conic bundle is trivial (Proposition \ref{p-FCB-triv}), i.e., $X \simeq T \times \P^1$. However, this would imply 
\[
 (-K_X)^3 = 42 >  18 = (-K_Y)^3, 
\]
which contradicts Lemma \ref{l-blowup-formula2}(2). 
Thus $Y$ is not 3-4.

Suppose that $Y$ is 3-8. 
Then $\widetilde Y$ is a Fano threefold of No.\ 2-24 (Notation \ref{n-ele-tr-F1-pic4}). 
Hence the extremal rays of $\widetilde Y$ are of type $C_1$ and $C_2$. 
This contradicts Lemma \ref{l-CC-2blowups}.  
Thus $Y$ is not 3-8. 
\end{proof}

\begin{lem}\label{l-ele-tf-F1-3-24}
We use Notation \ref{n-ele-tr-F1-pic4}. 
Assume that $Y$ is 3-24. 
Then $Y'$ is a Fano threefold  of No.\ 3-28,  $(-K_X)^3 = 36$, 
$B \in |\tau^*\MO_{\P^2}(1)|$, 
$p_a(B_Y)=0$, $(-K_Y) \cdot B_Y = 2$, $(-K_{Y'}) \cdot B_{Y'} =5$ ($X$ is 4-7). 
\end{lem}

\begin{proof} 
Recall that $Y = W \times_{\P^2} \F_1$ (Notation \ref{n-ele-tr-F1-pic4}). 
Let $\pi : W \to \P^2$ and $\pi' : W \to \P^2$ be the contractions of $W$. 
By the sequence $X \xrightarrow{\sigma} Y = W \times_{\P^2} \F_1 \xrightarrow{\pr_1} W$ of blowups, 
it follows from Lemma \ref{l-CC-2blowups} (applicable by setting $Z := W$) that $X \simeq \Bl_{B_W \amalg C_W} W$, 
where $B_W$ and $C_W$ are fibres of $\pi$ and $\pi'$ such that $B_W \cap C_W = \emptyset$. 
By symmetry, we may assume that $Y =\Bl_{C_W} W$ and 
$B_Y = \pr_1^{-1}(B_W) (\xrightarrow{\simeq} B_W)$. 
In particular, $p_a(B_Y)=p_a(B_W) =0$ and $-K_Y \cdot B_Y  = -K_W \cdot B_W = 2$. 
For the blowup $\tau : \F_1 \to \P^2$ at $\pi'(C_W)$, 
we obtain the following diagram in which each square is cartesian: 
\[
\begin{tikzcd}
\Bl_{B_W} W \arrow[d, "b_{B_W}"'] & X = \Bl_{B_W \amalg C_W} W \arrow[l, "b_{C_W}"'] \arrow[d, "\sigma = b_{B_Y}"] \\
W \arrow[d, "\pi'"'] & Y =W \times_{\P^2}\F_1 \arrow[l, "\pr_1 = b_{C_W}"'] \arrow[d, "g=\pr_2"] \\
\P^2  & \F_1, \arrow[l, "\tau"']
\end{tikzcd}
\]
where $b_{\Gamma}$ denotes the blowup along $\Gamma$. 
Since $B_Y$ is dijoint from $\Ex(b_{C_W} : W \to Y)$, 
we see that its image $B$ on $\F_1$ is disjoint from $\Ex(\tau)$. 
For $B_{\P^2} := \tau(B)$, it holds that $B_{\P^2}$ is a line on $\P^2$ (Proposition \ref{p-pic3-24}), 
and hence we get $B \in |\tau^*\MO_{\P^2}(1)|$ and $B^2=1$. 
We obtain $ (-K_{Y/S}) \cdot B_Y = -K_Y \cdot B - (-K_S) \cdot B = 2 -3 = -1$. 
Proposition \ref{p-ele-tf-numbers} and Lemma \ref{l-blowup-formula} imply  $(-K_{Y'})^3  = (-K_Y)^3 -4 (-K_{Y/S}) \cdot B_Y +2B^2 
= 42 - 4 \cdot (-1) + 2 \cdot 1 = 48$, 
$(-K_X)^3= (-K_Y)^3 -2(-K_Y) \cdot B_Y +2p_a(B) -2
=42 -2 \cdot 2 +0-2 =36$, and 
$(-K_{Y'}) \cdot B_{Y'} = B^2 +2(-K_S) \cdot B -(-K_Y) \cdot B_Y 
=  1 + 2 \cdot 3- 2 =5$. 
Note that $Y'$ is Fano, because $-K_{Y/S} \cdot B_Y = -1 \not\in 2\Z$ 
(Proposition \ref{p-ele-tf-1Fano}). Then $Y'$ is of No.\ 3-28 (Notation \ref{n-ele-tr-F1-pic4}). 
\end{proof}


\begin{lem}\label{l-ele-tf-F1-3-28}
We use Notation \ref{n-ele-tr-F1-pic4}. 
Assume that $Y$ is a Fano threefold of No.\ 3-28 and $Y'$ is Fano. 
Then one of the following holds. 
\begin{enumerate}
\item 
$Y'$ is a Fano threefold of No.\ 3-28, $(-K_X)^3 =30$, $B \in |\tau^*\MO_{\P^2}(2)|$, 
$-K_Y \cdot B_Y =8$, 
$-K_{Y'} \cdot B_{Y'} = 8$ ($X$ is 4-3). 
\item 
$Y'$ is a Fano threefold of No.\ 3-30, $(-K_X)^3 =40$, $B \in |\tau^*\MO_{\P^2}(1)|$, 
$-K_Y \cdot B_Y =3$, 
$-K_{Y'} \cdot B_{Y'} = 4$ ($X$ is 4-9). 
\item 
$Y'$ is a Fano threefold of No.\ 3-24, $(-K_X)^3 =36$, $B \in |\tau^*\MO_{\P^2}(1)|$, 
$-K_Y \cdot B_Y =5$, 
$-K_{Y'} \cdot B_{Y'} = 2$ ($X$ is 4-7). 
\end{enumerate}
\end{lem}

\begin{proof}
We have the induced contractions $Y = \F_1 \times  \P^1_2\xrightarrow{\pr_1} \F_1 \xrightarrow{\pi} \P^1_1$, 
where $\P^1_1 := \P^1$ and $\P^1_2 :=\P^1$.

We now prove that $B_Y$ dominates $\P_1^1$. 
Otherwise,  $B \subset \F_1$  would be a fibre of $\pi : \F_1 \to \P^1_1$. 
This contradicts the fact that 
the $(-1)$-curve $\Gamma$ on $\F_1$ 
is disjoint from or a connected component of 
the discriminant divisor $\Delta_f$ (Corollary \ref{c-FCB-(-1)-2}). 

Since $B_Y$ dominates $\P^1_1$, 
$B_Y$ is a subsection of the conic bundle $Y \to \P^1_1 \times \P^1_2$ (cf.\ Proposition \ref{p-pic3-28}). 
Then one of (1)-(6) of Lemma \ref{l-P1P1-3-28} holds. 
Set $B_{\P^1 \times \P^1}$ to be the image of $B_Y$ to $\P^1_1 \times \P^1_2$.  
Since $B_Y$ is a subsection of $Y =\F_1 \times \P^1 \to \F_1$, 
 Lemma \ref{l-P1P1-3-28}(1) does not occur. 

 \medskip

Assume Lemma \ref{l-P1P1-3-28}(2), 
i.e., $B$ is the $(-1)$-curve, $(-K_X)^3 =44$, and 
$-K_Y \cdot B_Y =1$. 
Then $-K_{Y/S} \cdot B_Y = -K_Y \cdot B_Y -(-K_S) \cdot B  =0$, 
$(-K_{Y'})^3 = (-K_Y)^3 -4 (-K_{Y/S}) \cdot B_Y +2B^2 
= 48 -0 -2= 46$ (Proposition \ref{p-ele-tf-numbers}). 
This contradicts $(-K_{Y'})^3 \in \{ 18, 24, 42, 48, 50\}$ (Notation \ref{n-ele-tr-F1-pic4}). 
 

 \medskip

In what follows, we treat the case when  $B \cap \Gamma = \emptyset$. 
In this case, we have $B \in |\tau^*\MO_{\P^2}(d_2)|$ 
for the bidegree $(d_1, d_2)$ of $B_{\P^1 \times \P^1}$, because 
\[
d_2 = \pr_1^*\MO_{\P^1}(1) \cdot B_{\P^1 \times \P^1} 
= \pr_1^*\pi^*\MO_{\P^1}(1) \cdot B_Y = \pi^*\MO_{\P^1}(1) \cdot B 
= \tau^*\MO_{\P^2}(1) \cdot B. 
\]

\medskip

(1) 
Assume Lemma \ref{l-P1P1-3-28}(3), 
i.e., $B \cap \Gamma = \emptyset$, $(-K_X)^3 =30$, 
$-K_Y \cdot B_Y =8$, and $B_{\P^1 \times \P^1}$ is of bidegree $(1, 2)$. 
In particular,  $B \in |\tau^*\MO_{\P^2}(2)|$. 
Then $-K_{Y/S} \cdot B_Y = -K_Y \cdot B_Y -(-K_S) \cdot B  =8 - (-K_{\P^2}) \cdot \MO_{\P^2}(2) =2$, 
$(-K_{Y'})^3 = (-K_Y)^3 -4 (-K_{Y/S}) \cdot B_Y +2B^2 = 48 - 8 +8 =48$, and 
$(-K_{Y'}) \cdot B_{Y'} = B^2 +2(-K_S) \cdot B -(-K_Y) \cdot B_Y 
=4 + 12 - 8=8$ 
(Proposition \ref{p-ele-tf-numbers}). 
Then $Y'$  is of No.\ 3-28 (Notation \ref{n-ele-tr-F1-pic4}). 

 \medskip

(2) 
Assume Lemma \ref{l-P1P1-3-28}(4), 
i.e., $B \cap \Gamma = \emptyset$, $(-K_X)^3 =40$, 
$-K_Y \cdot B_Y =3$, and $B_{\P^1 \times \P^1}$ is of bidegree $(0, 1)$. 
In particular,  $B \in |\tau^*\MO_{\P^2}(1)|$. 
Then $-K_{Y/S} \cdot B_Y = -K_Y \cdot B_Y -(-K_S) \cdot B_Y  =3 - (-K_{\P^2}) \cdot \MO_{\P^2}(1) =0$, 
$(-K_{Y'})^3 = (-K_Y)^3 -4 (-K_{Y/S}) \cdot B_Y +2B^2 = 48 -0 +2 =50$, and  $(-K_{Y'}) \cdot B_{Y'} = B^2 +2(-K_S) \cdot B -(-K_Y) \cdot B_Y 
=1 + 6 - 3=4$ 
(Proposition \ref{p-ele-tf-numbers}). 
Then $Y'$ is of No.\ 3-30 (Notation \ref{n-ele-tr-F1-pic4}).

 \medskip

(3) Assume Lemma \ref{l-P1P1-3-28}(5), 
i.e., $B \cap \Gamma = \emptyset$, $(-K_X)^3 =36$, 
$-K_Y \cdot B_Y =5$, and $B_{\P^1 \times \P^1}$ is of bidegree $(1, 1)$. 
In particular,  $B \in |\tau^*\MO_{\P^2}(1)|$. 
Then $-K_{Y/S} \cdot B_Y = -K_Y \cdot B_Y -(-K_S) \cdot B_Y  =5 - (-K_{\P^2}) \cdot \MO_{\P^2}(1) =2$, 
$(-K_{Y'})^3 = (-K_Y)^3 -4 (-K_{Y/S}) \cdot B_Y +2B^2 = 48 -8 +2 =42$, and 
$(-K_{Y'}) \cdot B_{Y'} = B^2 +2(-K_S) \cdot B -(-K_Y) \cdot B_Y 
=1 + 6 - 5=2$ 
(Proposition \ref{p-ele-tf-numbers}). 
Then $Y'$ is of No.\ 3-24  (Notation \ref{n-ele-tr-F1-pic4}).

 \medskip

Assume Lemma \ref{l-P1P1-3-28}(6), 
i.e., $B \cap \Gamma = \emptyset$, $(-K_X)^3 =32$, 
$-K_Y \cdot B_Y =7$, and $B_{\P^1 \times \P^1}$ is of bidegree $(2, 1)$. 
In particular,  $B \in |\tau^*\MO_{\P^2}(1)|$. 
Then $-K_{Y/S} \cdot B_Y = -K_Y \cdot B_Y -(-K_S) \cdot B_Y  =7 - (-K_{\P^2}) \cdot \MO_{\P^2}(1) =4$, 
$(-K_{Y'})^3 = (-K_Y)^3 -4 (-K_{Y/S}) \cdot B_Y +2B^2 = 48 - 16 +2 =34$  
(Proposition \ref{p-ele-tf-numbers}). 
However, this contradicts 
$(-K_{Y'})^3 \in \{ 18, 24, 42, 48, 50\}$  (Notation \ref{n-ele-tr-F1-pic4}). 
\qedhere

\end{proof}

\begin{lem}\label{l-ele-tf-F1-3-30}
We use Notation \ref{n-ele-tr-F1-pic4}. 
Assume that both $Y$ and $Y'$ are Fano threefolds of No.\ 3-30.  
Then it holds that 
\[
X \simeq \widetilde X \times_{\P^2} \F_1
\]
for a Fano threefold $\widetilde X$ of No.\ 3-19, $(-K_X)^3 =32$, $B \in |\tau^*\MO_{\P^2}(2)|$,  
$p_a(B)=0$, 
$-K_Y \cdot B_Y =8$, and 
$-K_{Y'} \cdot B_{Y'} = 8$ ($X$ is 4-4). 
\end{lem}

\begin{proof}
We have $Y \simeq Y' \simeq V_7 \times_{\P^2} \F_1$ (Proposition \ref{p-pic3-30}). 
Suppose that $B = \Gamma$. 
By $(-K_Y)^3 =(-K_{Y'})^3$ and $(-K_{Y'})^3 = (-K_Y)^3 -4 (-K_{Y/S}) \cdot B_Y +2B^2$ (Proposition \ref{p-ele-tf-numbers}), 
we obtain $-2 =2B^2 = 4 (-K_{Y/S}) \cdot B_Y \in 4\Z$, which is absurd. 

Thus the curve $B \subset \F_1$ is disjoint from the $(-1)$-curve $\Gamma$ (Corollary \ref{c-FCB-(-1)-2}). 
The elementary transform as in Notation \ref{n-ele-tr-F1-pic4} 
is obtained by applying the base change $(-) \times_{\P^2} \F_1$ 
to an elementary transform over $\P^2$ of type 2-35-vs-2-35 (cf.\ the table in Notation \ref{n-ele-tr-F1-pic4}). 
Thus $X \simeq \widetilde X \times_{\P^2} \F_1$ for a Fano threefold $\widetilde X$ of No.\ 3-19 (Theorem \ref{t-ele-tr-P2}). 
In particular, $(-K_X)^3 = (-K_{\widetilde X})^3 - 6 = 38-6=32$.  
The remaining assertions $B \in |\tau^*\MO_{\P^2}(2)|$,  $p_a(B)=0$, 
$-K_Y \cdot B_Y =8$, and 
$-K_{Y'} \cdot B_{Y'} = 8$ follow from the corresponding results 
for Fano threefolds 
of No.\ 3-19 in Theorem \ref{t-ele-tr-P2}. 
\end{proof}

\begin{lem}\label{l-ele-tf-F1-nonFano}
We use Notation \ref{n-ele-tr-F1-pic4}. 
Assume that $Y'$ is not Fano. 
Then one of the following holds. 
\begin{enumerate}
\item 
$(-K_X)^3 =44$, $Y$ is of No.\ 3-28, $B$ is the $(-1)$-curve $\Gamma$ on $\F_1$, 
$\Delta_f =\Gamma$, $-K_Y \cdot B_Y =1, -K_{Y'}\cdot B_{Y'} =0$ ($X$ is 4-11). 
\item 
$(-K_X)^3 =46$, $Y$ is of No.\ 3-30, $B$ is the $(-1)$-curve $\Gamma$ on $\F_1$, 
$\Delta_f = \Gamma$, $-K_Y \cdot B_Y =1, -K_{Y'} \cdot B_{Y'} =0$ ($X$ is 4-12). 
\item 
$X \simeq \widetilde X \times_{\P^2} \F_1$ for a Fano threefold $\widetilde X$ of No.\ 3-21, $(-K_X)^3 =32$, 
$Y$ is of No.\ 3-28, $B \in |\tau^*\MO_{\P^2}(1)|$, 
$p_a(B)=0$, $\Delta_f =B$, 
$-K_Y \cdot B_Y = 7,  -K_{Y'} \cdot B_{Y'} =0$ ($X$ is 4-5). 
\end{enumerate}
\end{lem}

\begin{proof}
For $D := g^{-1}(B) \subset Y$, the following hold (Proposition \ref{p-ele-tf-1Fano}): 
\[
B \simeq \P^1, \qquad 
-K_{Y'} \cdot B_{Y'} =0, \qquad 
-K_{Y/S} \cdot B_Y = 2(B^2+1), \qquad D \simeq \P^1 \times \P^1. 
\]
Note that $Y$ is either 3-28 or 3-30, because 
$Y$ is of No.\ 3-4, 3-8, 3-24, 3-28, or 3-30 (Notation \ref{n-ele-tr-F1-pic4}), and 
the case of No.\ 3-4, 3-8 (resp. 3-24) is excluded by 
Lemma \ref{l-ele-tf-F1-3-4,8} (resp. Lemma \ref{l-ele-tf-F1-3-24}). 

\medskip

(1), (2) 
We first treat the case when $B$ is the $(-1)$-curve $\Gamma$ on $\F_1$. 
In this case, $B^2 =-1$ and $-K_Y \cdot B_Y = -K_S \cdot B =1$. 
Then $(-K_X)^2 = (-K_Y)^3 - 4$ (Lemma \ref{l-blowup-formula}). 
\begin{enumerate}
    \item If $Y$ is 3-28, 
then $(-K_X)^3 = (-K_Y)^3 -4 = 44$. 
    \item If $Y$ is 3-30, 
then $(-K_X)^3 = (-K_Y)^3 -4 = 46$.  
\end{enumerate}


\medskip

(3) 
In what follows, we assume that $B$ is not the $(-1)$-curve $\Gamma$ on $\F_1$. 
Then $B \cap \Gamma = \emptyset$ (Corollary \ref{c-FCB-(-1)-2}). 
We have $X \simeq \Bl_{B_{\widetilde Y} \amalg F}\,\widetilde Y$, 
where $\rho : Y \to \widetilde Y$ is the blowup along a smooth fibre $F$ of 
a Fano conic bundle $\widetilde Y \to \P^2$ and $B_{\widetilde Y} := \rho(B_Y)$. 
Set $\widetilde X := \Bl_{B_{\widetilde Y}} \widetilde Y$. 
Since $X$ and $\widetilde Y$ are Fano, so is $\widetilde X$ 
(Corollary \ref{c-disjoint-blowup}). 
Let $\widetilde Y'$ be the elementary  transform of $\widetilde X \to \widetilde Y \to \P^2$. 
Then the elementary transform over $\F_1$ as in 
Notation \ref{n-ele-tr-F1-pic4} is obtained by applying the base change $(-) \times_{\P^2} \F_1$ to  
the elementary transform over $\widetilde S := \P^2$ consisting of $\widetilde X, \widetilde Y, \widetilde Y'$. 
The relation $-K_{Y/S} \cdot B_Y = 2(B^2+1)$ implies that 
\[
-K_{\widetilde Y/\widetilde S} \cdot B_{\widetilde Y} 
= 
-K_{Y/S} \cdot B_Y = 2(B^2+1)  =2(\widetilde B^2+1), 
\]
where $\widetilde B \subset \widetilde S = \P^2$ 
denotes the image of $B \subset S = \F_1$. 
Since $Y$ is 3-28 or 3-30, $\widetilde Y$ is 2-34 or 2-35, respectively 
(Theorem \ref{t-F1-pic3}). 
For $d:= \deg \widetilde B$, we get 
\[
-K_{\widetilde Y} \cdot B_{\widetilde Y} 
= -K_{\widetilde S} \cdot \widetilde{B} +2(\widetilde{B}^2+1)  
=2d^2 +3d +2. 
\]
We now compare Table \ref{table-ele-tr-P2} in Theorem \ref{t-ele-tr-P2} and the following list: 
\begin{itemize}
\item $d=1$: $-K_{\widetilde Y} \cdot B_{\widetilde Y} = 2d^2 +3d+2 = 7$. 
\item $d=2$: $-K_{\widetilde Y} \cdot B_{\widetilde Y} =2d^2 +3d+2 = 16$.
\item $d=3$: $-K_{\widetilde Y} \cdot B_{\widetilde Y} =2d^2 +3d+2 = 29$.
\item $d=4$: $-K_{\widetilde Y} \cdot B_{\widetilde Y} =2d^2 +3d+2 = 46$.
\end{itemize}
Then $\widetilde Y$ is 2-34, $\widetilde Y'$ is non-Fano, and 
$\widetilde X$ is either 3-5 or 3-21. 
In particular, $Y$ is 3-28 (Theorem \ref{t-F1-pic3}). 
If $\widetilde X$ is 3-21, then (3) holds by Theorem \ref{t-ele-tr-P2} and 
$(-K_X)^3 = (-K_{\widetilde X})^3 -6 =32$. 

It is enough to show that $\widetilde X$ is not 3-5. 
Suppose that $\widetilde X$ is 3-5. 
Then $(-K_{\widetilde X})^3  = 20$. 
Since $X$ is obtained by the blowup of $\widetilde X$ along a smooth fibre of $\widetilde X \to \P^2$, we get $(-K_X)^3 = (-K_{\widetilde X})^3 -6 = 14$. 
On the other hand, $Y$ is 3-28. 
Since $Y$ has a conic bundle structure over $\P^1 \times \P^1$ (Theorem \ref{t-P1P1-pic3}), 
so does $X$  (as otherwise, $X \simeq \P^1 \times T$ by 
Proposition \ref{p-FCB-centre} 
and 
Proposition \ref{p-FCB-triv}, which would imply $(-K_X)^3 = 42$). 
However, there exists no Fano threefold $X$ 
such that $\rho(X)=4$, $(-K_X)^3 =14$, and $X$ has a conic bundle structure over $\P^1 \times \P^1$ (Theorem \ref{t-ele-tr-P1P1-pic4}). 
\end{proof}

\begin{thm}\label{t-ele-tr-F1-pic4} 
Let 
\[
\begin{tikzcd}
& X \arrow[ld, "\sigma"'] \arrow[rd, "\sigma'"] \arrow[dd, "f"]\\
Y \arrow[rd, "g"']& & Y' \arrow[ld, "g'"]\\
& S:=\F_1 
\end{tikzcd}
\]
be an elementary transform of threefold conic bundles 
such that $\rho(X)=4$ and 
each of 
$f: X \to S=\F_1$ and $g : Y \to S=\F_1$ is a  Fano conic bundle (cf. Definition \ref{d-ele-tf}). 
Let $B_Y$ be the blowup centre of $\sigma: X \to Y$. 
Set $B:=g(B_Y)$. 
Assume that $(-K_Y)^3 \leq (-K_{Y'})^3$ when also $Y'$ is Fano. 
Then one of the following holds, where $\Gamma$ denotes the $(-1)$-curve on $\F_1$ and let $\tau : \F_1 \to \P^2$ be the contraction. 
    \begin{center}
\begin{longtable}{ccccccccc}
$X$ & $Y$ & $Y'$ & $(-K_X)^3$& $\Delta_f$ & $p_a(B)$ & $-K_Y \cdot B_Y$ & $-K_{Y'} \cdot B_{Y'}$ &\\ \hline
4-3 & 3-28 & 3-28 & $30$ & $\tau^*\MO_{\P^2}(2)$ & $0$ & $8$ & $8$\\ \hline
4-4 & 3-30 & 3-30 & $32$ & $\tau^*\MO_{\P^2}(2)$ & $0$ & $8$ &  $8$\\ \hline
4-5 & 3-28 & non-Fano & $32$ & $\tau^*\MO_{\P^2}(1)$ & $0$ & $7$ & $0$\\ \hline
4-7 & 3-24 & 3-28 & $36$ & $\tau^*\MO_{\P^2}(1)$ & $0$ & $2$ & $5$\\ \hline
4-9 & 3-28 & 3-30 & $40$ & $\tau^*\MO_{\P^2}(1)$ & $0$ & $3$ & $4$\\ \hline
4-11 & 3-28 & non-Fano& $44$ & $\Gamma$  & $0$  & $1$ & $0$\\ \hline
4-12 & 3-30 & non-Fano & $46$ & $\Gamma$  & $0$  & $1$ & $0$\\ \hline
 \\
 \caption{Elementary transforms over $\F_1$ ($\rho(X)=4$)}\label{table-ele-tr-F1-rho4}
      \end{longtable}
  \end{center} 
\end{thm} 

\begin{proof}
If $Y'$ is not Fano, then the assertion follows from 
Lemma \ref{l-ele-tf-F1-nonFano}. 
Hence we may assume that both $Y$ and $Y'$ are Fano. 
Then each of $Y$ and $Y'$ is 3-24, 3-28, or 3-30 (Notation \ref{n-ele-tr-F1-pic4}, Lemma \ref{l-ele-tf-F1-3-4,8}). 
If $Y$ is 3-24 (resp. 3-28), 
then the assertion follows from 
Lemma \ref{l-ele-tf-F1-3-24} 
(resp. Lemma \ref{l-ele-tf-F1-3-28}). 
Hence the remaining case is when each of $Y$ and $Y'$ is a Fano threefold of No.\ 3-30, 
which is settled by Lemma \ref{l-ele-tf-F1-3-30}. 
\end{proof}

\begin{thm}\label{t-pic4-fibre-blowup}
Let $\widetilde{f} : \widetilde{X} \to \widetilde{S}$ be a Fano conic bundle with $\rho(\widetilde{X})=3$. 
Let $C$ be a smooth fibre of $\widetilde f$ such that 
$X$ is Fano for the blowup $\sigma: X \to \widetilde X$ along $C$. 
Then one of the following holds. 
   \begin{center}
\begin{longtable}{ccccccccc}
$X$ & $\widetilde X$  & $(-K_X)^3$& $(-K_{\widetilde X})^3$ & $\widetilde S$ &  $\Delta_{\widetilde f}$ &\\ \hline
4-3 & 3-17  & $30$ & $36$ & $\P^2$ &   $\deg =2$\\ \hline
4-4 & 3-19  & $32$ & $38$ &$\P^2$  &   $\deg=2$\\ \hline
4-5 & 3-21  & $32$ & $38$ & $\P^2$  & $\deg =1$\\ \hline
4-7 & 3-24  & $36$ & $42$ & $\P^2$  & $\deg =1$\\ \hline
4-9 & 3-26  & $40$ & $46$ & $\P^2$  &  $\deg =1$\\ \hline
4-10 & 3-27  & $42$ & $48$ & $\P^1 \times \P^1$  &  $\emptyset$ \\ \hline
4-10 & 3-28  & $42$ & $48$ & $\F_1$  &  $\emptyset$ \\ \hline
 \\
 \caption{Fano conic bundles obtained by fibre blowups ($\rho(X)=4$)}\label{table-ele-tr-F1-rho4-2}
      \end{longtable}
  \end{center} 
\end{thm}

\begin{proof}
By $\rho(\widetilde S) < \rho(\widetilde X)=3$, 
we get $\widetilde S \in \{ \P^2, \P^1 \times \P^1, \F_1\}$.

\medskip

We now treat the case when $\widetilde S  =\P^1 \times \P^1$ or $\widetilde S = \F_1$. 
In this case, we get  $X \simeq S \times \P^1$ 
for a smooth projective del Pezzo surface 
with $K_S^2 =7$ (Proposition \ref{p-FCB-triv}). 
Then $(-K_{\widetilde X})^3 = (-K_X)^3 +6 = 48$. 
Hence $\widetilde X\simeq \P^1 \times \P^1 \times \P^1$ (No.\ 3-27) 
or $\widetilde X\simeq \F_1 \times \P^1$ (No.\ 3-28) by the classification list (Subsection \ref{ss-table-pic3}).  
If $\wt{X} \simeq \P^1 \times \P^1 \times \P^1$, 
then $\wt{S} = \P^1 \times \P^1$ by (Proposition \ref{p-pic3-27}). 
Assume that $\wt{X} \simeq \F_1 \times \P^1$. 
By Proposition \ref{p-pic3-28}, 
$\Delta_{\wt{f}} = \emptyset$ and $\wt{S} \in \{ \F_1, \P^1 \times \P^1\}$. 
Suppose  $\wt{S} = \P^1 \times \P^1$. 
Since $C$ is contracted by $\wt{f} : \wt{X} \to \wt{S}$ and this is a contraction of an extremal ray, 
$C$ is not contracted by a contraciton of any other extremal ray, 
and hence $C$ is not contracted by the projection $\pr_1 : \F_1 \times \P^1 \to \F^1$ (Proposition \ref{p-pic3-28}). 
Then $C$ is a regular subsection of $\pr_1 : \F_1 \times \P^1 \to \F^1$. Moreover, the image $C_{\F_1}$ of $C$ on $\F_1$ is a fibre over the $\P^1$-bundle $\pi : \F_1 \to \P^1$, 
because $C$ is contracted by $\pi \circ \pr_1$. 
However, this contradicts Corollary \ref{c-FCB-(-1)-2}. 
This completes the proof for the case when  $\widetilde S  =\P^1 \times \P^1$ or $\widetilde S = \F_1$. 

\medskip




In what follows, we assume that $\widetilde S =\P^2$. 
By $\rho(\widetilde{X}) =3 > 2 = \rho(\widetilde{S})+1$ and 
Lemma \ref{l-FCB-pic-irre}, 
we may apply Proposition \ref{p-ele-tf}, and hence  we obtain the following right square diagram which is an elementary transform of 
threefold conic bundles over $\widetilde{S} = \P^2$. 
Let $\tau : S \to \widetilde S$ be the   blowdown of the $(-1)$-curve $\Gamma$ on $S=\F_1$. 
By applying the base change $(-) \times_S S'$ to the right square, 
we get the left square  which is an elementary transform of threefold conic bundles over $S := \F_1$.  
\[
\begin{tikzcd}
& X \arrow[ld, "\sigma"'] \arrow[rd, "\sigma'"] \arrow[dd, "f"] &&&& \widetilde X\arrow[ld, "\widetilde \sigma"'] \arrow[rd, "\widetilde \sigma'"] \arrow[dd, "\widetilde f"]\\
Y \arrow[rd, "g"']& & Y' \arrow[ld, "g'"] && \widetilde Y \arrow[rd, "\widetilde g"']& & \widetilde Y' \arrow[ld, "\widetilde g'"]\\
& S=\F_1 \arrow[rrrr, "\tau"] &&&& \widetilde{S}=\P^2
\end{tikzcd}
\]
After permuting $\widetilde{Y}$ and $\widetilde{Y'}$ if necessary, 
we may assume that 
\begin{itemize}
    \item $Y$ is Fano (Proposition \ref{p-ele-tf-1Fano}), and 
    \item $(-K_Y)^3 \geq (-K_{Y'})^3$ when also $Y'$ is Fano. 
\end{itemize}

The left square satisfies the assumption and the conclusion of 
Theorem \ref{t-ele-tr-F1-pic4}. 
Moreover, we have $\Delta_f \neq \Gamma$. 
We have $(-K_{\widetilde X})^3 = (-K_{X})^3 +6$. 
Then we are done by comparing 
Table \ref{table-ele-tr-P2} in Theorem \ref{t-ele-tr-P2}
and 
Table \ref{table-ele-tr-F1-rho4} in Theorem \ref{t-ele-tr-F1-pic4}. 
Indeed, there are exactly five possibilities for the triple 
$((-K_X)^3, -K_Y \cdot B_Y, -K_{Y'} \cdot B_{Y'})$ by 
Table \ref{table-ele-tr-F1-rho4} in Theorem \ref{t-ele-tr-F1-pic4}. 
For example, assume that   
 $((-K_X)^3, -K_Y \cdot B_Y, -K_{Y'} \cdot B_{Y'}) = (32, 7, 0)$. 
By construction, we get  $((-K_{\wt{X}})^3, -K_{\wt{Y}} \cdot B_{\wt{Y}}, -K_{\wt{Y}'} \cdot B_{\wt{Y}'}) = (38, 7, 0)$ 
for the image $B_{\wt{Y}}$ (resp. $B_{\wt{Y}'}$) of $B_Y$ (resp. $B_{Y'}$) on $Y$ (resp. $Y'$). 
By Table \ref{table-ele-tr-P2} in Theorem \ref{t-ele-tr-P2}, 
we see that $\wt{X}$ is a Fano threefold of No.\ 3-21. 
We are done for the case when  $((-K_X)^3, -K_Y \cdot B_Y, -K_{Y'} \cdot B_{Y'}) = (32, 7, 0)$. 
The same argument is applicable for the remaining four cases. 
\qedhere


\end{proof}

\subsection{Case when $\rho(X)=4$ and $(-K_X)^3 = 32$}\label{ss rho4 vol=32}

In order to distinguish the cases of No.\ 4-4 and 4-5, 
we shall use the following proposition. 
Note that the Fano threefold $X$ appearing below is of No.\ 5-1. 

\begin{prop}\label{p-5-1-ample}
Let $C$ be a conic on $Q$. 
Fix three closed points $P_1, P_2, P_3 \in C$. 
Set $Y := \Bl_C\,Q$. 
For the induced blowup $\rho : Y := \Bl_C\,Q \to Q$, 
let $B_1, B_2, B_3$ be the fibres of $\rho$ over $P_1, P_2, P_3$, respectively. 
Set $X := \Bl_{B_1 \amalg B_2 \amalg B_3}\,Y$ and we have the induced birational morphisms: 
\[
\widetilde \sigma : X =  \Bl_{B_1 \amalg B_2 \amalg B_3}\,Y \xrightarrow{\sigma} Y=\Bl_C\,Q \xrightarrow{\rho} Q. 
\]
Then $|-K_X|$ is base point free and $-K_X$ is ample. 
\end{prop}

\setcounter{step}{0}
\begin{proof}
Set $E_Y :=\Ex(\rho)$ and $E_X := \sigma_*^{-1} E_Y$, i.e., 
$E_X$ is the proper transform of $E_Y$ on $X$. 
Let $D_1, D_2, D_3$ be the $\sigma$-exceptional prime divisors lying over $B_1, B_2, B_3$, respectively. 
The following hold: 
\[
\sigma^*E_Y = E_X +D_1 + D_2 +D_3,\qquad K_Y \sim \rho^*K_Q +E_Y, 
\]
\[
K_X \sim \sigma^* K_Y +D_1+D_2+ D_3 
\sim \sigma^*(\rho^*K_Q +E_Y) +D_1+D_2+ D_3 
\]
\[
= \sigma^*\rho^*K_Q +E_X + 2D_1 +2D_2 +2D_3 
\sim -3H_Q +E_X + 2D_1 +2D_2 +2D_3,  
\]
where $H_Q := \sigma^*\rho^*\MO_Q(1)$. 

For the closed embeddings $C \subset Q \subset \P^4$, 
the linear subvariety $\langle C \rangle$ generated by $C$ in $\P^4$ is a plane. 
We have the scheme-theoretic intersection $Q \cap \langle C \rangle = C$. 
As we can write $\langle C \rangle = H \cap H'$ for some hyperplanes 
$H$ and $H'$ on $\P^4$, 
we get $T \cap T'  = C$ for $T := H|_Q \in |\MO_Q(1)|$ and 
$T' :=H'|_Q \in |\MO_Q(1)|$. 
Note that $T$ is a quadric surface smooth around $C$, and hence $T$ is a normal prime divisor. 
Let $T_Y$ be the strict transform of $T$ on $Y$. 
Then $\rho^*T =T_Y + E_Y$ and $H_Q \sim \sigma^*\rho^*T = \sigma^*T_Y +E_X + D_1 + D_2 + D_3$. 
It follows from $T \cap T' = C$ that 
$|T_Y|$ is base point free and $T_Y \xrightarrow{\simeq} T$.


\begin{step}\label{s1-5-1-ample}
$\Bs\,|-K_X| \subset E_X$. 
\end{step}

\begin{proof}[Proof of Step \ref{s1-5-1-ample}]
It holds that 
\begin{eqnarray*}
-K_X 
&\sim& 3H_Q - E_X -2D_1 -2D_2 -2D_3\\
&\sim& H_Q +2( \sigma^*T_Y +E_X + D_1 + D_2 + D_3)  - E_X -2D_1 -2D_2 -2D_3\\
&=& H_Q + 2\sigma^*T_Y +E_X. 
\end{eqnarray*}
Since $|H_Q|$ and $|\sigma^*T_Y|$ are base point free, 
we obtain $\Bs\,|-K_X| \subset E_X$. 
This completes the proof of Step \ref{s1-5-1-ample}. 
\end{proof}

\begin{step}\label{s2-5-1-ample}
The restriction map 
\[
H^0(X, \MO_X(-K_X)) \to H^0(E_X, \MO_X(-K_X)|_{E_X})
\]
is surjective. 
\end{step}

\begin{proof}[Proof of Step \ref{s2-5-1-ample}]
Consider the following exact sequence: 
\[
0 \to \MO_X(-K_X-E_X) \to \MO_X(-K_X) \to \MO_X(-K_X)|_{E_X} \to 0. 
\]
It suffices to show $H^1(X, \MO_X(-K_X-E_X))=0$. 
By the proof of Step \ref{s1-5-1-ample}, we have that 
\[
-K_X -E_X \sim H_Q + 2\sigma^*T_Y \sim \sigma^*( \rho^*\MO_Q(1) + 2T_Y).  
\]
Consider the following Leray spectral sequence:  
\[
E_2^{i, j} =H^i(Y, R^j\sigma_*\MO_X(-K_X-E_X))
\Rightarrow H^{i+j}(X, \MO_X(-K_X-E_X)) =E^{i+j}. 
\]
Note that $-K_X-E_X-K_X \sim \sigma^*( \rho^*\MO_Q(1) + 2T_Y)-K_X$. 
Hence $-K_X-E_X-K_X$ is $\sigma$-ample. 
Since the relative Kodaira vanishing theorem holds for $\sigma$ (cf.\ \cite[Theorem 0.5]{Tan15}), 
we obtain $R^j\sigma_*\MO_X(-K_X-E_X)=0$ for every $j>0$. 
Hence 
\[
H^i(X, \MO_X(-K_X-E_X)) =E^i \simeq  E_2^{i, 0} \simeq  
H^i(Y, \rho^*\MO_Q(1) \otimes \MO_Y( 2T_Y)). 
\]
Recall that $T_Y (\simeq T)$ is a quadric surface. 
We have the following exact sequence for all integers $i>0, s \geq 0$: 
\[
H^i(Y, \rho^*\MO_Q(1)  \otimes \MO_Y(sT_Y))
\to
H^i(Y,  \rho^*\MO_Q(1)  \otimes \MO_Y((s+1)T_Y))
\]
\[
\to 
H^i(T_Y, (\rho^*\MO_Q(1)  \otimes \MO_Y((s+1)T_Y))|_{T_Y}) \overset{(\star)}{=} 0, 
\]
where the equality $(\star)$ follows from the fact that $T_Y$ is toric and $(\rho^*\MO_Q(1) \otimes  \MO_Y(s+1)T_Y)|_{T_Y}$ is nef (recall that $|T_Y|$ is base point free). 
Then it is enough to show that $H^i(Y,  \rho^*\MO_Q(1))=0$ for every $i>0$. 
This follows from 
\[
H^i(Y,  \rho^*\MO_Q(1)) 
\overset{(\star\star)}{\simeq}
H^i(Q, \MO_Q(1))=0, 
\]
where the isomorphism $(\star\star)$ can be checked by 
using the following Leray spectral sequence: 
\[
H^i(Q, R^j\rho_*\rho^*\MO_Q(1))
\Rightarrow H^{i+j}(Y, \rho^*\MO_Q(1)). 
\]
This completes the proof of Step \ref{s2-5-1-ample}. 
\end{proof}

\begin{step}\label{s3-5-1-ample}
$-K_X|_{E_X} \sim (\sigma|_{E_X})^*(T_Y|_{E_Y} + B_1)$ and $|-K_X|$ is base point free, 
where $\sigma|_{E_X} : E_X \to E_Y$ denotes the induced isomorphism. 
\end{step}

\begin{proof}[Proof of Step \ref{s3-5-1-ample}]
By identifying $E_X$ and $E_Y$ via the induced isomorphism $\sigma|_{E_X} : E_X \xrightarrow{\simeq} E_Y$, we obtain 
\[
-K_X|_{E_X} \sim (2H_Q + \sigma^*T_Y -D_1-D_2-D_3)|_{E_X} 
\sim (\rho^*\MO_Q(2) +T_Y)|_{E_Y} -B_1-B_2-B_3. 
\]
We have the $\P^1$-bundle structure $\rho|_{E_Y} : E_Y \to C \simeq \P^1$. 
All of $B_1, B_2, B_3$ are fibres of $\rho|_{E_Y}$. 
Note that $\MO_Q(2) \cdot C = 4$. 
Therefore, 
\[
-K_X|_{E_X} \sim (\rho^*\MO_Q(2) +T_Y)|_{E_Y} -B_1-B_2-B_3 \sim T_Y|_{E_Y} +B_1, 
\]
which means $-K_X|_{E_X} \sim (\sigma|_{E_X})^*(T_Y|_{E_Y} + B_1)$. 
Then $|-K_X|_{E_X}|$ is base point free. 
This, together with Step \ref{s1-5-1-ample} and Step \ref{s2-5-1-ample}, 
implies that $|-K_X|$ is base point free. 
This completes the proof of Step \ref{s3-5-1-ample}. 
\end{proof}

\begin{step}\label{s4-5-1-ample}
$-K_X$ is ample. 
\end{step}

\begin{proof}[Proof of Step \ref{s4-5-1-ample}]
Fix a curve $\Gamma$ on $X$. 
Since $|-K_X|$ is base point free (Step \ref{s3-5-1-ample}), 
it is enough to show that $-K_X \cdot \Gamma >0$.  
We treat the following three cases separately. 
\begin{enumerate}
\item[(i)] $\Gamma \not\subset E_X \cup D_1 \cup D_2 \cup D_3$. 
\item[(ii)] $\Gamma \subset E_X$. 
\item[(iii)] $\Gamma \subset D_1 \cup D_2 \cup D_3$. 
\end{enumerate}

(i) Assume $\Gamma \not\subset E_X \cup D_1 \cup D_2 \cup D_3 $. 
Recall that $E_X \cup D_1 \cup D_2 \cup D_3 =\Ex(\widetilde \sigma: X \to Q)$. 
Then 
\[
-K_X \cdot \Gamma = ( H_Q +2\sigma^*T_Y +E_X) \cdot \Gamma \geq H_Q \cdot \Gamma >0. 
\]

(ii) Assume $\Gamma \subset E_X$. 
By Step \ref{s3-5-1-ample}, we obtain $-K_X|_{E_X} \sim  (\sigma|_{E_X})^*(T_Y|_{E_Y} + B_1)$ for $\sigma|_{E_X} : E_X \xrightarrow{\simeq} E_Y$. 
In what follows, we identify $E_X$ and $E_Y$ via this isomorphism. 
Note that $T_Y|_{E_Y}$ is a nonzero effective divisor on $E_Y$ such that $\rho(T_Y|_{E_Y}) = C$. 
If $\Gamma$ is a fibre of $\rho_{E_X} : E_X \to C$, 
then $-K_X \cdot \Gamma = (-K_X|_{E_X}) \cdot \Gamma =
(T_Y|_{E_Y} + B_1) \cdot \Gamma = (T_Y|_{E_Y}) \cdot \Gamma >0$. 
If $\Gamma$ is not a fibre of $\rho_{E_X} : E_X \to C$, then 
$-K_X \cdot \Gamma = (-K_X|_{E_X}) \cdot \Gamma =
(T_Y|_{E_Y} + B_1) \cdot \Gamma \overset{(*)}{\geq}  B_1 \cdot \Gamma >0$, 
where the inequality $(*)$ follows from the fact that $T_Y$ is nef. 

(iii) 
Assume $\Gamma \subset D_1 \cup D_2 \cup D_3$. 
By symmetry, the problem is reduced to the case when $\Gamma \subset D_1$. 
Moreover, we may assume, by the case (ii), that $\Gamma \not\subset E_X$. 
If $\Gamma$ is a fibre of $D_1 \to \sigma(D_1)$, then $-K_X \cdot \Gamma =1$. 
Hence we may assume that $\sigma(\Gamma)=B_1$. 
It holds that 
$-K_X \cdot \Gamma =(H_Q+2\sigma^*T_Y +E_X) \cdot \Gamma \geq 2 \sigma^*T_Y \cdot \Gamma$. 
We obtain 
\[
\sigma^*T_Y \cdot \Gamma  = T_Y \cdot \sigma_*(\Gamma) = T_Y \cdot (nB_1) = n (T_Y|_{E_Y}) \cdot B_1 >0
\]
for some integer $n>0$, 
where the inequality $(T_Y|_{E_Y}) \cdot B_1 >0$ 
follows from the fact that $B_1$ is a fibre of the $\P^1$-bundle 
$E_Y \to C$ and 
$T_Y|_{E_Y}$ is a nonzero effective divisor on $E_Y$ 
which dominates $C$. 
This completes the proof of Step \ref{s4-5-1-ample}. 
\qedhere



\end{proof}
Step \ref{s3-5-1-ample} and Step \ref{s4-5-1-ample} complete  the proof of Proposition \ref{p-5-1-ample}. 
\end{proof}

\subsection{Classification  ($\rho=4$)}\label{ss-pic4-classify}


\begin{nasi}\label{n-pic4-vol}
Let $X$ be a Fano threefold with $\rho(X)=4$. 
By Corollary \ref{c-pic4-CB}, one of (I) and (II) holds. 
\begin{enumerate}
\item[(I)] $X$ has a conic bundle structure over $\P^1 \times \P^1$. 
In this case, the following holds (Theorem \ref{t-ele-tr-P1P1-pic4}):   
\[
(-K_X)^3 \in \{24, 26, 28, 30, 32, 34, 36, 38, 40, 42, 44\}. 
\]
\item[(II)]  $X$ has a conic bundle structure over $\F_1$. 
In this case, the following holds (Theorem \ref{t-ele-tr-F1-pic4}):   
\[
(-K_X)^3 \in \{30, 32, 36, 40, 44, 46\}. 
\]
\end{enumerate}
\end{nasi}

\begin{nasi}\label{n-pic4-blowup}
Let $X$ be a Fano threefold with $\rho(X)=4$, 
let $Y$ be a Fano threefold with $\rho(Y)=3$, and 
let $\sigma : X \to Y$ be a blowup along a smooth curve $B_Y$ on $Y$. 
By Proposition \ref{t-pic3-structure} and Lemma \ref{l-FCB-P^2-1}, one of (A)-(D) holds. 
\begin{enumerate}
\item[(A)] There is a conic bundle $g: Y \to S = \P^1 \times \P^1$ 
and $B_Y$ is a regular subsection of $g$. 
In this case, the following holds (Theorem \ref{t-ele-tr-P1P1-pic4}):   
\[
(-K_X)^3 \in \{24, 26, 28, 30, 32, 34, 36, 38, 40, 42, 44\}. 
\]
\item[(B)] There is a conic bundle $g: Y \to S = \F_1$ 
and $B_Y$ is a regular subsection of $g$. 
In this case, the following holds (Theorem \ref{t-ele-tr-F1-pic4}):   
\[
(-K_X)^3 \in \{30, 32, 36, 40, 44, 46\}. 
\]
\item[(C)] There is a conic bundle $g: Y \to S$ 
and $B_Y$ is a smooth fibre of $g$. 
In this case, the following holds (Theorem \ref{t-pic4-fibre-blowup}):  
\[
(-K_X)^3 \in \{ 30, 32, 36, 40, 42 \}. 
\]
\item[(D)] $Y$ is a Fano threefold of No.\ 3-18. 
In this case, $(-K_X)^3  =32$ (Lemma \ref{l-pic4-3-18}).  
\end{enumerate}
\end{nasi}



\begin{lem}\label{l-pic4-over-3-28}
Let $X$ be a Fano threefold with $\rho(X)=4$. 
Assume that there exists a blowup $\sigma : X \to Y := \F_1 \times \P^1$ along a smooth curve $B_Y$ on $Y = \F_1 \times \P^1$. 
Let $\pi : Y  = \F_1 \times \P^1 \to \F_1$ and $\pi' : \F_1 \times \P^1 \to \P^1 \times \P^1$ be the contractions (cf.\ Proposition \ref{p-pic3-28}).
Then one of the following holds. 
\begin{enumerate}
\item 
$B_Y$ is a regular subsection of each of the conic bundles $\pi$ and $\pi'$. 
Moreover, 
$X$ has a conic bundle structure over $\P^1 \times \P^1$ and $X$ has a conic bundle structure over $\F_1$. 
    \item 
    $X \simeq S \times \P^1$, where $S$ is a smooth del Pezzo surface with $K_S^2 =7$. Moreover, $(-K_X)^3 = 42$. 
\end{enumerate}
\end{lem}

\begin{proof}
Recall that one of the following holds (Proposition \ref{p-FCB-centre}). 
\begin{enumerate}
\item[(1)'] $B_Y$ is a regular subsection of each of the conic bundles 
$\pi$ and $\pi'$. 
\item[(2)'] $B_Y$ is a fibre of one of $\pi$ and $\pi'$. 
\end{enumerate}
Again by Proposition \ref{p-FCB-centre}, (1)' implies (1). 
Assume (2)'. 
Then (2) holds by Proposition \ref{p-FCB-centre} and 
Proposition \ref{p-FCB-triv}. 
\end{proof}

\begin{dfn}\label{d bdown pic4}
Let $X$ be a Fano threefold with $\rho(X)=4$. 
We define the finite set 
\[
\Blowdown(X) \subset \{ \text{3-1, 3-2, 3-3, ..., 3-31}\}
\]
by the following condition: $\text{3-xx} \in \Blowdown(X)$ if and only if 
there exist a Fano threefold $Y_{\text{3-xx}}$ of No.\ 3-xx and a smooth curve $C$ on $Y_{\text{3-xx}}$ such that $X$ is isomorphic to the blowup of $Y_{\text{3-xx}}$ along $C$. 
\end{dfn}

\begin{prop}[No.\ {\hyperref[table-4-1]{4-1}}]\label{p-pic4-1}
Let $X$ be a Fano threefold with $\rho(X)=4$ and $(-K_X)^3=24$. 
Then the following hold. 
\begin{enumerate}
\item $\Blowdown(X) = \{ \text{3-27}\}$. 
\item There is a conic bundle structure $f:X \to \P^1 \times \P^1$ of type 3-27-vs-3-27 such that $\deg \Delta_f = (2, 2)$. 
\item 
$X$ is isomorphic to 
a divisor on $\P^1 \times \P^1 \times \P^1 \times \P^1$ of multi-degree $(1, 1, 1, 1)$, i.e., a member of $|\MO_{\P^1 \times \P^1 \times \P^1 \times \P^1}(1, 1, 1, 1)|$. 
\end{enumerate}
\end{prop}

\begin{proof}
Only (A) of (\ref{n-pic4-blowup}) holds. 
By Proposition \ref{t-ele-tr-P1P1-pic4}, (1) and (2) hold.

Let us show (3). 
 Take an elementary transform over $\P^1 \times \P^1$ as in (2) (cf. Theorem \ref{t-ele-tr-P1P1-pic4}): 
\[
\begin{tikzcd}
& X \arrow[ld, "\sigma"'] \arrow[rd, "\sigma'"] \arrow[dd, "f"]\\
Y=\P^1_1 \times \P^1_2 \times \P^1_3 \arrow[rd, "g"']& & Y' = \P^1_2 \times \P^1_3 \times \P^1_4 \arrow[ld, "g'"]\\
& S:=\P^1_2 \times \P^1_3 
\end{tikzcd}
\]
For each $i \in \{1, 2, 3, 4\}$, we have the contraction $\varphi_i : X \to \P^1_i$. 
Let $\varphi := \varphi_1 \times \varphi_2 \times \varphi_3 \times \varphi_4 : X \to \P^1_1 \times \P^1_2 \times \P^1_3 \times \P^1_4$ be the induced morphism. 
Since 
the extremal ray $R_{\sigma}$ of $\sigma = \varphi_1 \times \varphi_2 \times \varphi_3$ 
is not contained in the extremal face $F_{\varphi_4}$ of $\varphi_4 : X \to \P^1_4$, 
we obtain $R_{\sigma} \cap F_{\varphi_4} = \{0\}$, i.e., 
$\varphi : X \to \P^1_1 \times \P^1_2 \times \P^1_3 \times \P^1_4$ is a finite morphism. 
For $X' := \varphi(X)$, let $\psi : X \to X'$ be the induced finite surjective morphism. 
For 
\[
\varphi_i : X \xrightarrow{\psi} X' \xrightarrow{\varphi'_i} \P^1_i, 
\]
we set $H_i := \varphi_i^*\MO_{\P^1}(1)$ and $H'_i :=  \varphi_i'^*\MO_{\P^1}(1)$. 
Lemma \ref{l-ele-tf-K-relation} implies 
\[
-2K_X \sim -\sigma^*K_Y -\sigma'^*K_{Y'} -2f^*\Delta_f 
\sim 2(H_1+H_2+H_3) + 2(H_2+H_3+H_4) -2(H_2 +H_3). 
\]
Hence  $-K_X \sim H_1 + H_2 + H_3 + H_4$. 
Since we have the following factorisation: 
\[
\sigma: X \xrightarrow{\psi} X' \hookrightarrow \P^1_1 \times \P^1_2 \times \P^1_3 \times \P^1_4 \to 
\P^1_1 \times \P^1_2 \times \P^1_3, 
\]
$\psi : X \to X'$ is birational.

Let $(d_1, d_2, d_3, d_4)$ be the multi-degree of $X'$ in $\P^1_1 \times \P^1_2 \times \P^1_3 \times \P^1_4$. 
If $d_i=0$ for some $i \in \{1, 2, 3, 4\}$, 
then $X'$ is contained in $\P^1 \times \P^1 \times \P^1 \times \{t\} \simeq \P^1 \times \P^1 \times \P^1$, and hence $X \simeq X' \simeq \P^1 \times \P^1 \times \P^1$, which contradicts $\rho(X)=4$. 
Hence $d_1 >0, d_2 >0, d_3>0, d_4>0$. 
We have that 
\[
24 = (-K_X)^3 =(H_1 +H_2 + H_3 + H_4)^3 
\]
\[
= 
6(H_2 \cdot H_3 \cdot H_4+ H_1 \cdot H_3 \cdot H_4 
 + H_1 \cdot H_2 \cdot H_4
+H_1 \cdot H_2 \cdot H_3) = 6(d_1 + d_2 + d_3+d_4). 
\]
Hence $d_1 = d_2 =d_3 = d_4 =1$. 
For the conductor divsior  $C \subset X$ of the normalisation $\psi : X \to X'$, we get $\MO_X(K_X+C) \simeq \psi^*\omega_{X'}$. 
By the adjunction formula, we have $\omega_{X'}^{-1} \simeq \MO_{X'}(H'_1 + H'_2+H'_3+H'_4)$. 
This, together with  $-K_X \sim H_1 + H_2 + H_3 + H_4$, 
implies $C \sim 0$, and hence $C =0$. 
Therefore, $X \xrightarrow{\psi, \simeq} X' \subset \P^1_1 \times \P^1_2 \times \P^1_3 \times \P^1_4$. Thus (3) holds. 
\end{proof}

\begin{prop}[No.\ {\hyperref[table-4-2]{4-2}}]\label{p-pic4-2}
Let $X$ be a Fano threefold with $\rho(X)=4$ and $(-K_X)^3=28$. 
Then the following hold. 
\begin{enumerate}
\item $\Blowdown(X) = \{ \text{3-31}\}$. 
\item There is a conic bundle structure $f:X \to \P^1 \times \P^1$ of type 3-31-vs-3-31 such that $\deg \Delta_f = (2, 2)$. 
\item 
$X \simeq \Bl_{C}\,Y_{\text{3-31}}$, where 
$Y_{\text{3-31}} := \P_{\P^1 \times \P^1}(\MO \oplus \MO(1, 1))$, 
$C$ is a smooth curve on a section $T$ of 
the $\P^1$-bundle 
$\pi : Y_{\text{3-31}} = \P_{\P^1 \times \P^1}(\MO \oplus \MO(1, 1)) \to \P^1 \times \P^1$, $\pi(C)$ is an elliptic curve of bidegree $(2, 2)$, 
and $T$ is disjoint from the section $S$ of $\pi$ such that $\MO_{Y_{\text{3-31}}}(-S)|_S$ is ample. 
\end{enumerate}
\end{prop}

\begin{proof}
Only (A) of (\ref{n-pic4-blowup}) holds. 
By Proposition \ref{t-ele-tr-P1P1-pic4}, (1) and (2) hold.

Let us show (3). Take an elementary transform over $\P^1 \times \P^1$ as in (2) (cf. Proposition \ref{t-ele-tr-P1P1-pic4}): 
\[
\begin{tikzcd}
& X \arrow[ld, "\sigma"'] \arrow[rd, "\sigma'"] \arrow[dd, "f"]\\
Y:=\P_{\P^1 \times \P^1}(\MO \oplus \MO(1, 1)) \arrow[rd, "g"']& & Y':=\P_{\P^1 \times \P^1}(\MO \oplus \MO(1, 1)) \arrow[ld, "g'"]\\
& \P^1 \times \P^1. 
\end{tikzcd}
\]
Set $D$ (resp. $D'$) to be the section of $g$ (resp. $g'$) such that  $-D|_{D}$  (resp. $-D'|_{D'}$) is ample. 
Let $D_X$ and $D'_X$ be the proper transforms of $D$ and $D'$ on $X$, respectively. 

By Proposition \ref{t-ele-tr-P1P1-pic4}, it is enough to show that 
\begin{enumerate}
\item[(i)]the section $\sigma(D'_X)$ of $g$ contains the blowup centre $B_Y$ of $\sigma: X \to Y$, and 
\item[(ii)] $D \cap \sigma(D'_X) = \emptyset$. 
\end{enumerate}
We have $D \simeq \P^1 \times \P^1$ and 
$K_Y|_D \sim D|_D \sim \MO_{\P^1 \times \P^1}(-1, -1)$ 
(Lemma \ref{l-nonFano-blowdown}, Proposition \ref{p-pic3-31}). 
For a fibre $\zeta_i$ of each projection $\pr_i : D \xrightarrow{\simeq} \P^1 \times \P^1 \xrightarrow{\pr_i} \P^1$, we have that $-K_Y \cdot \zeta_i  = (-K_Y|_D) \cdot \zeta_i =1$. 
Hence the blowup centre $B_Y$ of $\sigma : X \to Y$ is disjoint from $D$ (Lemma \ref{l-line-meeting}). 
By symmetry, the blowup centre $B_{Y'}$ of $\sigma' : X \to Y'$ is disjoint from $D'$. 
For the irreducible decomposition $f^{-1}(\Delta_f) = F_1 \cup F_2$, 
$D_X$ and $D'_X$ are sections of $f: X \to \P^1 \times \P^1$ 
each of which 
intersects one and only one of $F_1$ and $F_2$. 
Therefore, the section $\sigma(D'_X)$ of $g: Y \to \P^1 \times \P^1$ contains $B_Y$. 
Thus (i) holds. 

Let us show (ii). 
By construction, we have $D_X \neq D'_X$. 
Proposition \ref{p-2ray} and 
Proposition \ref{p-nonFano-flop} enable us to find a birational contraction $\varphi : X \to Z$ 
(resp. $\varphi' : X \to Z'$) such that 
$Z$ (resp. $Z'$) is a projective normal threefold, 
$\Ex(\varphi) = D_X$ (resp. $\Ex(\varphi') = D'_X$), and 
$\varphi(\Ex(\varphi))$ (resp. $\varphi'(\Ex(\varphi'))$) is a point. 
By $D'_X \simeq \P^1 \times \P^1$, we get $D_X \cap D'_X =\emptyset$ (as otherwise,  for an ample Cartier divisors $A_Z$ on $Z$, 
$(\varphi^*A_Z)|_{D'_X}$ would be nef and big but not ample 
by $(\varphi^*A_Z)|_{D_X \cap D'_X} \equiv 0$, 
which contradicts the fact that every nef and big divisor on $D'_X (\simeq \P^1 \times \P^1)$ is ample).
Thus  (ii) holds. 
This completes the proof of  (3).  
\end{proof}

\begin{prop}[No.\ {\hyperref[table-4-3]{4-3}}]\label{p-pic4-3}
Let $X$ be a Fano threefold with $\rho(X)=4$ and $(-K_X)^3=30$. 
Then the following hold. 
\begin{enumerate}
\item $\Blowdown(X) = \{ \text{3-17}, \text{3-27}, \text{3-28}\}$. 
\item 
\begin{enumerate}
\item  There is a conic bundle structure $f_1:X \to \P^1 \times \P^1$ of type 3-17-vs-3-27 such that $\deg \Delta_{f_1} = (1, 1)$. 
\item  There is a conic bundle structure $f_2:X \to \P^1 \times \P^1$ of type 3-27-vs-3-28 such that $\deg \Delta_{f_2} = (1, 1)$. 
\item  There is a conic bundle structure $f_3:X \to \F_1$ of type 3-28-vs-3-28 such that 
$\Delta_{f_3} \in |\tau^*\MO_{\P^2}(2)|$. Moreover, $X \simeq Y_{\text{3-17}} \times_{\P^2} \F_1$. 
\end{enumerate}
\item 
$X \simeq \Bl_C\,(\P^1 \times \P^1 \times \P^1)$, 
where $C$ is a smooth curve on $\P^1 \times \P^1 \times \P^1$ of tridegree $(1, 1, 2)$. 
\end{enumerate}
\end{prop}

\begin{proof}
First of all, we prove (i) and (ii) below. 
\begin{enumerate}
\item[(i)] (a) $\Leftrightarrow$ (b) $\Leftrightarrow$ (3). 
\item[(ii)] (a), (b), and (c) hold. 
\end{enumerate}

(i) If (3) holds, then (a) and (b) hold by taking suitable projections $\P^1 \times \P^1 \times \P^1 \to \P^1 \times \P^1$ (Theorem \ref{t-ele-tr-P1P1-pic4}). 
Conversely, assume that (a) or (b) holds. 
By Theorem \ref{t-ele-tr-P1P1-pic4}, 
there is a blowup $f: X \to \P^1 \times \P^1 \times \P^1=Y$ along a smooth curve $B_Y$ 
such that $-K_Y \cdot B_Y = 8$ and $(d_1, d_2) \in \{ (1, 1), (1, 2), (2, 1)\}$ for the tridegree $(d_1, d_2, d_3)$ of $B_Y$. 
By $2(d_1+d_2+d_3) =-K_Y \cdot B_Y = 8$, 
we obtain $(d_1, d_2, d_3) \in \{ (1, 1, 2), (1, 2, 1), (2, 1, 1)\}$. 
Hence we get the implications (a) $\Rightarrow$ (3) and (b) $\Rightarrow$ (3). 
Thus (i) holds.

(ii) 
Recall that (I) or (II) of (\ref{n-pic4-vol}) holds. 
In any case, 
there exists a blowup $X \to \F_1 \times \P^1$ along a smooth curve on $\F_1 \times \P^1$ (Theorem \ref{t-ele-tr-P1P1-pic4}, Theorem \ref{t-ele-tr-F1-pic4}). 
Then $X$ has conic bundle structures over $\P^1 \times \P^1$ and $\F_1$ (Lemma \ref{l-pic4-over-3-28}(1)). 
It follows from Theorem \ref{t-ele-tr-P1P1-pic4} (resp. Theorem \ref{t-ele-tr-F1-pic4} and Theorem \ref{t-pic4-fibre-blowup}) 
that (a) or (b) (resp. (c)) holds. 
Then (i) implies that both (a) and (b) hold. 
This completes the proof of (ii).


By (i) and (ii), we see that (2) and (3) hold. 
In particular, $\Blowdown(X) \supset \{ \text{3-17}, \text{3-27}, \text{3-28}\}$. 
By (\ref{n-pic4-blowup}), 
the opposite inclusion 
$\Blowdown(X) \subset \{ \text{3-17}, \text{3-27}, \text{3-28}\}$ 
follows from 
Theorem \ref{t-ele-tr-P1P1-pic4}, Theorem \ref{t-ele-tr-F1-pic4}, 
and Theorem \ref{t-pic4-fibre-blowup}. 
Thus (1) holds. 
\end{proof}

\begin{lem}\label{l-pic4-vol32}
Let $X$ be a Fano threefold with $\rho(X)=4$ and $(-K_X)^3=32$. 
Then
\[
\Blowdown(X) \subset  \{ \text{3-18}, \text{3-19}, \text{3-30}\} 
\cup \{ \text{3-21}, \text{3-28}, \text{3-31}\}. 
\]
\end{lem}

\begin{proof}
Let $\sigma : X \to Y$ be a blowup along a smooth curve $B_Y$ on a Fano threefold $Y$. 
Then one of (A)-(D) in (\ref{n-pic4-blowup}) holds. 
Hence the assertion follows from 
Theorem \ref{t-ele-tr-P1P1-pic4}, Theorem \ref{t-ele-tr-F1-pic4}, and Theorem \ref{t-pic4-fibre-blowup}. 
\end{proof}

\begin{prop}[No.\ {\hyperref[table-4-4]{4-4}}]\label{p-pic4-4}
Let $X$ be a Fano threefold with $\rho(X)=4$ and $(-K_X)^3=32$. 
Assume that $\Blowdown(X)$ contains one of 
$\text{3-18}, \text{3-19}, \text{3-30}$. 
Then the following hold. 
\begin{enumerate}
\item $\Blowdown(X) = \{ \text{3-18}, \text{3-19}, \text{3-30}\}$. 
\item  There is a conic bundle structure $f:X \to \F_1$ of type 3-30-vs-3-30 such that $\deg \Delta_f \in |\tau^*\MO_{\P^2}(2)|$. Moreover, $X \simeq Y_{\text{3-19}} \times_{\P^2} \F_1$. 
\item 
$X \simeq \Bl_{B_1 \amalg B_2} Y_{\text{2-29}}$, 
where $Y_{\text{2-29}} = \Bl_C\,Q$ for a conic $C$ on $Q$, 
and $B_1$ and $B_2$ are distinct one-dimensional fibres 
of the induced blowup $\Bl_C\,Q \to Q$. 
\item 
There exists a smooth curve $\Gamma$ on $X$ such that the blowup of $X$ along $\Gamma$ is Fano. 
\end{enumerate}
\end{prop}

The proof of (1) will be given in Lemma \ref{l-4-4-4-5-(1)}. 

\begin{proof}[Proof except for (1)]
We shall prove the following (i)-(iii).  
\begin{enumerate}
\item[(i)] $\text{3-19} \in \Blowdown(X) \Leftrightarrow \text{3-30} \in \Blowdown(X)$. 
\item[(ii)] $\text{3-18} \in \Blowdown(X) \Rightarrow \text{3-30} \in \Blowdown(X)$. 
\item[(iii)] $\text{3-19} \in \Blowdown(X) \Rightarrow$ (3) $\Rightarrow \text{3-18} \in \Blowdown(X)$. 
\end{enumerate}
We now finish the proof by assuming (i)-(iii). 
Our assumption and (i)-(iii) imply (3) and 
$\Blowdown(X) \supset \{ \text{3-18}, \text{3-19}, \text{3-30}\}$. 
Then 
(2) holds by Theorem \ref{t-ele-tr-P1P1-pic4}, Theorem \ref{t-ele-tr-F1-pic4}, and Theorem \ref{t-pic4-fibre-blowup}. 
Note that (4) follows from (3) and Proposition \ref{p-5-1-ample}. 

\medskip

It is enough to show (i)-(iii). 
Let us show (i). 
If $\text{3-19} \in \Blowdown(X)$ or $\text{3-30} \in \Blowdown(X)$, 
then (B) or (C) of (\ref{n-pic4-blowup}) 
holds by Theorem \ref{t-ele-tr-P1P1-pic4}. 
Hence (i) follows by 
comparing $\Delta_f$ in Theorem \ref{t-ele-tr-F1-pic4} and $\Delta_{\wt{f}}$ in Theorem \ref{t-pic4-fibre-blowup}.

Let us show (ii). 
Assume  $\text{3-18} \in \Blowdown(X)$. 
Let $\sigma:X \to Y:= Y_{\text{3-18}}$ be a blowup along a smooth curve $B_Y$. 
Note that $(-K_Y)^3-(-K_X)^3  = 36-32= 4$. 
Recall that we have $Y_{\text{3-18}} = \Bl_{L\amalg C}\,\P^3$ for 
a line $L$ and a conic $C$ on $\P^3$ which are mutually disjoint. 
For the induced blowup $\rho : Y = Y_{\text{3-18}} \to \P^3$, 
$B_Y$ intersects $\Ex(\rho)$, as otherwise we would get 
$-K_Y \cdot B_Y = -\rho^*K_{\P^3} \cdot B_Y \in 4\Z$, 
which contradicts Lemma \ref{l-blowup-formula2}(3). 
Then $B_Y$ is a one-dimensional fibre of $\rho : Y = Y_{\text{3-18}} \to \P^3$ 
(Lemma \ref{l-line-meeting}), 
and hence $-K_Y \cdot B_Y = 1$. 
If $B_Y$ is lying over the conic $C$, then 
we get a blowup $X \to Y_{\text{3-23}}$ along the inverse image of the line (Lemma \ref{l-P3-V7}, Proposition \ref{p-pic3-23}), 
and hence  $\text{3-23} \in \Blowdown(X)$, which contradicts Lemma \ref{l-pic4-vol32}. 
Thus $B_Y$ is lying over $L$. 
Then, for the blowup $Z$ of $\Bl_L\,\P^3$ along the image of $B_Y$ (which is a one-dimensional fibre of the blowup $\Bl_L\,\P^3 \to \P^3$), 
$Z$ is a Fano threefold (Proposition \ref{p-nonFano-iff}). 
It follows from $(-K_Z)^3 = 50$ and Lemma \ref{l-pic4-vol32} that 
$Z$ is of No.\ 3-30, which implies $\text{3-30} \in \Blowdown(X)$. 
Thus (ii) holds.

Let us show (iii). 
The implication  (3) $\Rightarrow \text{3-18} \in \Blowdown(X)$ 
follows from the fact that a Fano threefold of No.\ 3-18 
is isomorphic to $\Bl_{B_1}\,Y_{\text{2-29}}$ as in (3) (Proposition \ref{p-pic3-18}). 
Assume 
$\text{3-19} \in \Blowdown(X)$. 
Let us show (3). 
By Theorem \ref{t-ele-tr-P1P1-pic4}, Theorem \ref{t-ele-tr-F1-pic4}, and Theorem \ref{t-pic4-fibre-blowup}, 
we have $X \xrightarrow{\sigma} Y=Y_{\text{3-19}} \xrightarrow{\rho} Q$, where $X = Y_{\text{3-19}} \times_{\P^2} \F_1$, $
\sigma : X = Y_{\text{3-19}} \times_{\P^2} \F_1 \to Y_{\text{3-19}}$ is the first projection, and 
$\rho : Y_{\text{3-19}} \to Q$ 
is the blowup along along two points $P_1$ and $P_2$ (Proposition \ref{p-pic3-19}). 
In particular, $\sigma: X \to Y=Y_{\text{3-19}}$ is a blowup along a smooth rational curve $B_Y$ with $-K_Y \cdot B_Y =2$. 
Let $E_1, E_2 \subset Y_{\text{3-19}}$ be the $\rho$-exceptional prime divisors  
lying over $P_1, P_2$, respectively. 
Then $E_1$ and $E_2$ dominates $\P^2$ (Proposition \ref{p-pic3-19}). 
We have 
\[
-K_Y \sim -\rho^*K_Q - 2E_1 -2E_2. 
\]
By Proposition \ref{p-pic3-19}(3), 
we get $2 = (-K_Y) \cdot B_Y = \rho^*\MO_Q(1) \cdot B_Y$. 
Then $(2E_1 + 2E_2) \cdot B_Y = (-\rho^*K_Q -(-K_Y)) \cdot  B_Y = 6-2=4$. 
This, together with $E_1 \cdot  B_Y>0$ and $E_2 \cdot  B_Y >0$, 
implies $E_1 \cdot  B_Y =E_2 \cdot B_Y =1$. 
Then we get $-\rho^*K_Q \cdot B_Y= (-K_Y  +2E_1 +2E_2) \cdot B_Y =2+2+2=6$, 
which implies that $\rho(B_Y)$ is a conic. 
Applying Lemma \ref{l-P3-V7} twice, (3) holds.  
This completes the proof of (iii). 
\qedhere






\end{proof}

\begin{prop}[No.\ {\hyperref[table-4-5]{4-5}}]\label{p-pic4-5}
Let $X$ be a Fano threefold with $\rho(X)=4$ and $(-K_X)^3=32$. 
Assume that $\Blowdown(X)$ contains one of 
$\text{3-21}, \text{3-28}, \text{3-31}$. 
Then the following hold. 
\begin{enumerate}
\item $\Blowdown(X) = \{ \text{3-21}, \text{3-28}, \text{3-31}\}$. 
\item 
\begin{enumerate}
\item  There is a conic bundle structure $f_1:X \to \P^1 \times \P^1$ of type 3-28-vs-3-31 such that $\deg \Delta_{f_1} = (1, 2)$. 
\item  There is a conic bundle structure $f_2:X \to \F_1$ of type 3-28-vs-non-Fano such that 
$\Delta_{f_2} \in |\tau^*\MO_{\P^2}(1)|$. Moreover, $X \simeq Y_{\text{3-21}} \times_{\P^2} \F_1$. 
\end{enumerate}
\item 
$X \simeq \Bl_{C_1 \amalg C_2} (\P^2 \times \P^1)$, 
where $C_1$ and $C_2$ are smooth curves of bidegree $(0, 1)$ and $(1, 2)$, respectively. 
\item 
There exists no smooth curve $\Gamma$ on $X$ such that the blowup of $X$ along $\Gamma$ is Fano. 
\end{enumerate}
\end{prop}

The proof of (1) will be given in Lemma \ref{l-4-4-4-5-(1)}. 

\begin{proof}[Proof except for (1)]
By Theorem \ref{t-ele-tr-P1P1-pic4}, Theorem \ref{t-ele-tr-F1-pic4}, and Theorem \ref{t-pic4-fibre-blowup}, 
we get $\text{3-28} \in \Blowdown(X)$, i.e., 
there is a blowup $X \to Y := \F_1 \times \P^1$ along a smooth curve $B_Y$ on $Y$. 
For the $\P^1$-bundles $\pi: \F_1 \times \P^1 \to \F_1$ and 
$\pi' : \F_1 \times \P^1 \to \P^1 \times \P^1$, 
it follows from Theorem \ref{t-pic4-fibre-blowup} that 
$B_Y$ is a regular subsection of each of $\pi$ and $\pi'$. 
Then Theorem \ref{t-ele-tr-P1P1-pic4} (resp. Theorem \ref{t-ele-tr-F1-pic4}) 
implies (a) (resp. (b)), and hence (2) holds. 
Then (a) implies  (4). 

Let us show (3). 
We have 
\[
X \simeq \Bl_{C_2'} (\F_1 \times \P^1) \simeq \Bl_{C_1 \amalg C_2} (\P^2 \times \P^1), 
\]
where 
\begin{itemize}
\item $C_1$ is a fibre of $\pr_1 : \P^2 \times \P^1 \to \P^2$, 
\item $C_2$ is a smooth curve on $\P^2 \times \P^1$ with $C_1 \cap C_2 = \emptyset$, and 
\item 
the inverse image $C'_2$ of $C_2$ on $\F_1 \times \P^1$ satisfies 
$-K_{\F_1 \times \P^1} \cdot C'_2 = 7$ and 
$(C'_2)_{\F_1} \in |\tau^*\MO_{\P^2}(1)|$ for the image 
$(C'_2)_{\F_1}$ of $C'_2$ on $\F_1$
(Theorem \ref{t-ele-tr-F1-pic4}). 
\end{itemize}
Hence $C_1$ is of bidegree $(0, 1)$. 
Set $\deg C_2 := (d_1, d_2)$ to be the bidegree of $C_2 \subset \P^2 \times \P^1$. 
Since $C'_2$  is a regular subsection of $\P^2 \times \P^1$, 
it follows from $(C'_2)_{\F_1} \in |\tau^*\MO_{\P^2}(1)|$  that 
$d_1 = \tau^*\MO_{\P^2}(1) \cdot (C'_2)_{\F_1} =1$. 
We then get $d_2=2$ by 
\[
7 = -K_{\F_1 \times \P^1} \cdot C'_2 
= -K_{\P^2 \times \P^1} \cdot C_2 = 3d_1 + 2d_2 = 3 + 2d_2. 
\]
Thus (3) holds. 
\qedhere

\end{proof}

\begin{lem}\label{l-4-4-4-5-(1)}
Proposition \ref{p-pic4-4}(1) and Proposition \ref{p-pic4-5}(1) hold. 

\end{lem}

\begin{proof} 
Let $X$ be a Fano threefold with $\rho(X)=4$ and $(-K_X)^3=32$. 
Let us show Proposition \ref{p-pic4-4}(1). 
Assume that $\Blowdown(X)$ contains one of 
$\text{3-18}, \text{3-19}, \text{3-30}$. 
By Proposition \ref{p-pic4-4}(2)(3), 
we obtain $\{ \text{3-18}, \text{3-19}, \text{3-30}\} \subset \Blowdown(X)$. 
If $\{ \text{3-18}, \text{3-19}, \text{3-30}\} \subsetneq \Blowdown(X)$, 
then Lemma \ref{l-pic4-vol32} implies that one of \text{3-21}, \text{3-28}, \text{3-31} is contained in $\Blowdown(X)$. 
Then both Proposition \ref{p-pic4-4}(4)  and Proposition \ref{p-pic4-5}(4) hold, 
which is absurd. 
Hence Proposition \ref{p-pic4-4}(1) holds. 
The proof of  Proposition \ref{p-pic4-5}(1) is similar. 
\qedhere

\end{proof}

\begin{prop}[No.\ {\hyperref[table-4-6]{4-6}}]\label{p-pic4-6}
Let $X$ be a Fano threefold with $\rho(X)=4$ and $(-K_X)^3=34$. 
Then the following hold. 
\begin{enumerate}
\item $\Blowdown(X) = \{ \text{3-25}, \text{3-27}\}$. 
\item 
There is a conic bundle structure $f:X \to \P^1 \times \P^1$ of type 3-25-vs-3-27 such that 
$\deg \Delta_f  =(1, 1)$. 
\item 
$X \simeq \Bl_C\,(\P^1 \times \P^1 \times \P^1)$, 
where $C$ is a smooth curve on $\P^1 \times \P^1 \times \P^1$ of tridegree $(1, 1, 1)$. 
\item 
$X \simeq \Bl_{L_1 \amalg L_2 \amalg L_3}\,\P^3$ for a disjoint union of three lines $L_1, L_2, L_3$ on $\P^3$. 
\end{enumerate}
\end{prop}

\begin{proof}
Only (A) of (\ref{n-pic4-blowup}) holds. 
Theorem \ref{t-ele-tr-P1P1-pic4} implies  (1) and (2). 
By (2), we have contractions $\sigma : X \to Y :=Y_{\text{3-25}} $ and $\sigma' : X \to Y' := \P^1 \times \P^1 \times \P^1$. 
It follows from Theorem \ref{t-ele-tr-P1P1-pic4} that 
the blowup centre $B_{Y'}$ of $\sigma': X \to 
Y' = \P^1 \times \P^1 \times \P^1$ satisfies $-K_{Y'} \cdot B_{Y'} = 6$ and $(d_1, d_2) = (1, 1)$ 
 for the tridegree $(d_1, d_2, d_3)$ of $B_{Y'}$. 
By $2(d_1+d_2+d_3) =-K_{Y'} \cdot B_{Y'} = 6$, 
we obtain $(d_1, d_2, d_3) =(1, 1, 1)$. 
Thus (3) holds. 

Let us show (4). 
Recall that $Y_{\text{3-25}} = \Bl_{L_1 \amalg L_2}\,\P^3$ 
for mutually disjoint lines $L_1$ and $L_2$ (Subsection \ref{ss-table-pic3}). 
Let $\rho : Y_{\text{3-25}} = \Bl_{L_1 \amalg L_2}\,\P^3 \to \P^3$ be the blowup along $L_1 \amalg L_2$. 
By Theorem \ref{t-ele-tr-P1P1-pic4}, we have $-K_Y \cdot B_Y = 4$ 
for the blowup centre $B_Y \subset Y_{\text{3-25}}$ of $\sigma : X \to Y_{\text{3-25}}$. 
It follows from Lemma \ref{l-line-meeting} that 
$B_Y \cap \Ex(\rho) = \emptyset$, as otherwise we would get $-K_Y \cdot B_Y = 1$. 
By $-K_Y \cdot B_Y = 4$, its image $L_3 := \rho(B_Y)$ is a line on $\P^3$, 
which is disjoint from $L_1$ and $L_2$. 
Therefore, $X \simeq \Bl_{L_1 \amalg L_2 \amalg L_3}\,\P^3$. 
Thus (4) holds. 
\end{proof}

\begin{prop}[No.\ {\hyperref[table-4-7]{4-7}}]\label{p-pic4-7}
Let $X$ be a Fano threefold with $\rho(X)=4$ and $(-K_X)^3=36$. 
Then the following hold. 
\begin{enumerate}
\item $\Blowdown(X) = \{ \text{3-24}, \text{3-28}\}$. 
\item 
\begin{enumerate}
\item  There is a conic bundle structure $f_1:X \to \P^1 \times \P^1$ of type 3-28-vs-3-28 such that $\deg \Delta_{f_1} = (1, 1)$. 
\item  There is a conic bundle structure $f_2:X \to \F_1$ of type 3-24-vs-3-28 such that 
$\Delta_{f_2} \in |\tau^*\MO_{\P^2}(1)|$. Moreover, $X \simeq Y_{\text{3-24}} \times_{\P^2} \F_1$. 
\end{enumerate}
\item 
$X \simeq \Bl_{C_1 \amalg C_2}\,W$, where $C_1$ and $C_2$ are mutually disjoint smooth curves of bidegree $(1, 0)$ and $(0, 1)$, respectively. 
\end{enumerate}
\end{prop}

\begin{proof}
First of all, we prove that (I) and (II) of (\ref{n-pic4-vol}) hold. 
Recall that one of (I) and (II) holds. 
If (I) (resp. (II)) holds, then 
it follows from Theorem \ref{t-ele-tr-P1P1-pic4} (resp. Theorem \ref{t-ele-tr-F1-pic4}) that 
there is a blowup $X \to \F_1 \times \P^1$ 
along a smooth curve. 
Then Lemma \ref{l-pic4-over-3-28} implies that (I) and (II) of (\ref{n-pic4-vol}) hold.

By Theorem \ref{t-ele-tr-P1P1-pic4} 
(resp. Theorem \ref{t-ele-tr-F1-pic4} and Theorem \ref{t-pic4-fibre-blowup}), 
we obtain (a) (resp. (b)). 
Thus (2) holds and $\Blowdown(X) \supset \{ \text{3-24}, \text{3-28}\}$. 
The opposite inclusion 
 $\Blowdown(X) \subset \{ \text{3-24}, \text{3-28}\}$ 
 holds by (\ref{n-pic4-blowup}), 
Theorem \ref{t-ele-tr-P1P1-pic4}, Theorem \ref{t-ele-tr-F1-pic4}, Theorem \ref{t-pic4-fibre-blowup}. 
Thus (1) holds. 
Lemma \ref{l-CC-2blowups} implies (3), 
because a Fano threefold of No.\ 3-24 is a blowup of $W$ along a smooth fibre of a $\P^1$-bundle $W \to \P^2$ (Subsection \ref{ss-table-pic3}). 
\end{proof}

\begin{prop}[No.\ {\hyperref[table-4-8]{4-8}}]\label{p-pic4-8}
Let $X$ be a Fano threefold with $\rho(X)=4$ and $(-K_X)^3=38$. 
Then the following hold. 
\begin{enumerate}
\item $\Blowdown(X) = \{ \text{3-27}, \text{3-31}\}$. 
\item 
\begin{enumerate}
\item  There is a conic bundle structure $f_1:X \to \P^1 \times \P^1$ of type 3-27-vs-3-31 such that $\deg \Delta_{f_1} = (1, 1)$. 
\item  There is a conic bundle structure $f_1:X \to \P^1 \times \P^1$ of type 3-27-vs-non-Fano such that $\deg \Delta_{f_2} = (0, 1)$. 
\end{enumerate}
\item 
$X \simeq \Bl_C\,\P^1 \times \P^1 \times \P^1$, 
where $C$ is a smooth curve on $\P^1 \times \P^1 \times \P^1$ of tridegree $(0, 1, 1)$. 
\end{enumerate}
\end{prop}

\begin{proof}
Only (A) of (\ref{n-pic4-blowup}) holds. 
By Theorem \ref{t-ele-tr-P1P1-pic4}, 
there exists a blowup $f: X \to \P^1 \times \P^1 \times \P^1=Y$ along a smooth curve $B_Y$ 
such that $-K_Y \cdot B_Y = 4$ and $(d_1, d_2) \in \{ (1, 1), (0, 1)\}$ for the tridegree $(d_1, d_2, d_3)$ of $B_Y$. 
By $2(d_1+d_2+d_3) =-K_Y \cdot B_Y = 4$, 
we obtain $(d_1, d_2, d_3) \in \{ (1, 1, 0), (0, 1, 1)\}$. 
Thus (3) holds. 
Then (3) implies (2) by taking suitable projections $\P^1 \times \P^1 \times \P^1 \to \P^1 \times \P^1$ (Theorem \ref{t-ele-tr-P1P1-pic4}). 
In particular, we obtain $\Blowdown(X) \supset \{ \text{3-27}, \text{3-31}\}$. 
Since only (A) of (\ref{n-pic4-blowup}) holds, the opposite inclusion $\Blowdown(X) \subset \{ \text{3-27}, \text{3-31}\}$ follows from Theorem \ref{t-ele-tr-P1P1-pic4}. 
Thus (1) holds. 
\end{proof}

\begin{prop}[No.\ {\hyperref[table-4-9]{4-9}}]\label{p-pic4-9}
Let $X$ be a Fano threefold with $\rho(X)=4$ and $(-K_X)^3=40$. 
Then the following hold. 
\begin{enumerate}
\item $\Blowdown(X) = \{ \text{3-25, 3-26, 3-28, 3-30}\}$. 
\item 
\begin{enumerate}
\item  There is a conic bundle structure $f_1:X \to \P^1 \times \P^1$ of type 3-25-vs-3-28 such that $\deg \Delta_{f_1} = (0, 1)$. 
\item  There is a conic bundle structure $f_2:X \to \F_1$ of type 3-28-vs-3-30 such that 
$\Delta_{f_2} \in |\tau^*\MO_{\P^2}(1)|$. Moreover, $X \simeq Y_{\text{3-26}} \times_{\P^2} \F_1$. 
\end{enumerate}
\item 
$X \simeq \Bl_C Y_{\text{3-25}}$, where 
$Y_{\text{3-25}} := \Bl_{L_1 \amalg L_2}\,\P^3$ for 
mutually disjoint lines $L_1$ and $L_2$ on $\P^3$ 
and $C$ is a one-dimensional fibre of the induced blowup $\rho : Y_{\text{3-25}} = \Bl_{L_1 \amalg L_2} \P^3 \to \P^3$. 
\end{enumerate}
\end{prop}

\begin{proof}
First of all, we prove that (I) and (II) of (\ref{n-pic4-vol}) hold. 
Recall that one of (I) and (II) holds. 
If (I) (resp. (II)) holds, then 
it follows from Theorem \ref{t-ele-tr-P1P1-pic4} (resp. Theorem \ref{t-ele-tr-F1-pic4}) that 
there exists a blowup $X \to \F_1 \times \P^1$ along a smooth curve. 
Then Lemma \ref{l-pic4-over-3-28} implies that (I) and (II) of (\ref{n-pic4-vol}) hold. 

Then Theorem \ref{t-ele-tr-P1P1-pic4} 
(resp. Theorem \ref{t-ele-tr-F1-pic4} and Theorem \ref{t-pic4-fibre-blowup}) implies (a) (resp. (b)). 
Thus (2) holds. We then get  $\Blowdown(X) \supset \{ \text{3-25, 3-26, 3-28, 3-30}\}$. 
Since (D) of (\ref{n-pic4-blowup}) does not hold, 
the opposite inclusion $\Blowdown(X) \subset \{ \text{3-25, 3-26, 3-28, 3-30}\}$ follows from Theorem \ref{t-ele-tr-P1P1-pic4}, Theorem \ref{t-ele-tr-F1-pic4}, and 
Theorem \ref{t-pic4-fibre-blowup}. 
Thus (1) holds. 

It suffices to  show (3). 
By (a) and Theorem \ref{t-ele-tr-P1P1-pic4}, 
there is a blowup $\sigma: X \to Y := Y_{\text{3-25}}$ along a smooth curve $B_Y$ with $-K_Y \cdot B_Y =1$. 
It is enough to prove  that $B_Y$ is a one-dimensional fibre of 
the induced blowup $\rho : Y_{\text{3-25}}  = \Bl_{L_1 \amalg L_2} \P^3 \to \P^3$. 
If $B_Y$ is disjoint from $\Ex(\rho)$, then we would get the following contradiction: $-1 =K_{Y_{\text{3-25}}} \cdot B_Y = (\rho^*K_{\P^3} + \Ex(\rho)) \cdot B_Y = \rho^*K_{\P^3} \cdot B_Y \in 4\Z$. 
Hence $B_Y \cap \Ex(\rho) \neq \emptyset$. 
Then Lemma \ref{l-line-meeting} implies that $B_Y$ is a one-dimensional fibre of $\rho$. 
Thus (3) holds. 
\end{proof}

\begin{prop}[No.\ {\hyperref[table-4-10]{4-10}}]\label{p-pic4-10}
Let $X$ be a Fano threefold with $\rho(X)=4$ and $(-K_X)^3=42$. 
Then the following hold. 
\begin{enumerate}
\item $\Blowdown(X) = \{ \text{3-27, 3-28}\}$. 
\item 
There is a conic bundle structure $f:X \to \P^1 \times \P^1$ of type 3-27-vs-3-28 such that $\deg \Delta_{f} = (0, 1)$. 
\item 
$X \simeq S \times \P^1$, where $S$ is a smooth del Pezzo surface with $K_S^2 =7$. 
\end{enumerate}
\end{prop}

\begin{proof}
Only (I) of (\ref{n-pic4-vol}) holds. 
By Lemma \ref{l-P1P1-3-27} and Theorem \ref{t-ele-tr-P1P1-pic4}, (2) and (3) hold. 
In particular, we obtain $\Blowdown(X) \supset \{ \text{3-27, 3-28}\}$. 
Since only (A) and (C) of (\ref{n-pic4-blowup})  holds, 
the opposite inclusion $\Blowdown(X) \subset \{ \text{3-27, 3-28}\}$ follows from 
Theorem \ref{t-ele-tr-P1P1-pic4} and Theorem \ref{t-pic4-fibre-blowup}. 
\end{proof}

\begin{prop}[No.\ {\hyperref[table-4-11]{4-11}}]\label{p-pic4-11}
Let $X$ be a Fano threefold with $\rho(X)=4$ and $(-K_X)^3=44$. 
Then the following hold. 
\begin{enumerate}
\item $\Blowdown(X) = \{ \text{3-28, 3-31}\}$. 
\item 
\begin{enumerate}
\item  There is a conic bundle structure $f_1:X \to \P^1 \times \P^1$ of type 3-28-vs-3-31 such that $\deg \Delta_{f_1} = (0, 1)$. 
\item  There is a conic bundle structure $f_2:X \to \F_1$ of type 3-28-vs-non-Fano such that 
$\Delta_{f_2}$ is equal to the $(-1)$-curve. 
\end{enumerate}
\item 
$X \simeq \Bl_C\,(\F_1 \times \P^1)$ for $C = \Gamma \times \{t\}$, 
where $\Gamma$ is the $(-1)$-curve on $\F_1$ and $t \in \P^1$ is a closed point. 
\end{enumerate}
\end{prop}

\begin{proof}
First of all, we prove that (I) and (II) of (\ref{n-pic4-vol}) hold. 
Recall that one of (I) and (II) holds. 
If (I) (resp. (II)) holds, then 
it follows from Theorem \ref{t-ele-tr-P1P1-pic4} (resp. Theorem \ref{t-ele-tr-F1-pic4}) that 
there exists a blowup $X \to \F_1 \times \P^1$ along a smooth curve. 
Then Lemma \ref{l-pic4-over-3-28} implies that (I) and (II) of (\ref{n-pic4-vol}) hold.

By  Theorem \ref{t-ele-tr-P1P1-pic4} and Theorem \ref{t-ele-tr-F1-pic4}, we get (a) and (b), respectively. 
Hence (2) holds. 
We then get  $\Blowdown(X) \supset \{ \text{3-28, 3-31}\}$. 
Since none of (C) nor (D) of (\ref{n-pic4-blowup}) holds, 
the opposite inclusion $\Blowdown(X) \subset \{ \text{3-28, 3-31}\}$ follows from 
Theorem \ref{t-ele-tr-P1P1-pic4} and Theorem \ref{t-ele-tr-F1-pic4}. 
Thus (1) holds. 

It suffices to show (3). 
By Theorem \ref{t-ele-tr-F1-pic4}, 
there is a blowup $\sigma : X \to Y :=\F_1 \times \P^1$ along a smooth curve $B_Y$ 
on $Y = \F_1 \times \P^1$ such that 
$-K_Y \cdot B_Y = 1$, 
$B_Y$ is a regular subsection of $\pr_1 : Y = \F_1 \times \P^1 \to \F_1$, and the image $\Gamma := \pr_1(B_Y) \subset \F_1$ is the $(-1)$-curve on $\F_1$.   
We have  $\Gamma_Y := \pr_1^{-1}(\Gamma) = \Gamma \times \P^1 \simeq \P^1 \times \P^1$. 
Let $(d_1, d_2)$ be the bidegree of $B_Y$ in $\Gamma_Y = \Gamma \times \P^1$. 
By $\pr_1|_{B_Y} : B_Y \xrightarrow{\simeq} \Gamma$, we get $d_2  =1$. 
Since $-K_Y|_{\Gamma_Y}$ is ample, we obtain 
$(-K_Y|_{\Gamma_Y}) \cdot \MO_{\Gamma_Y}(1, 0)>0$ and 
$(-K_Y|_{\Gamma_Y}) \cdot \MO_{\Gamma_Y}(0, 1)>0$. 
Then $1 = -K_Y \cdot B_Y = (-K_Y|_{\Gamma_Y}) \cdot B_Y$ implies $d_1 =0$. 
Thus $B_Y$ is a divisor on $\Gamma \times \P^1 \simeq \P^1 \times \P^1$ of bidegree $(0, 1)$, i.e., $B_Y = \Gamma \times \{t\}$ for some $t \in \P^1$. 
Thus (3) holds. 
\end{proof}

\begin{prop}[No.\ {\hyperref[table-4-12]{4-12}}]\label{p-pic4-12}
Let $X$ be a Fano threefold with $\rho(X)=4$ and $(-K_X)^3=46$. 
Then the following hold. 
\begin{enumerate}
\item $\Blowdown(X) = \{ \text{3-30}\}$. 
\item There is a conic bundle structure $f:X \to \F_1$ of type 3-30-vs-non-Fano such that 
$\Delta_{f}$ is equal to the $(-1)$-curve. 
\item 
$X \simeq \Bl_{C_1 \amalg C_2}\,Y_{\text{2-33}}$, where 
$Y_{\text{2-33}} = \Bl_L\,\P^3$ for a line $L$, and 
$C_1$ and $C_2$ are mutually distinct  one-dimensional fibres of 
the induced blowup $Y_{\text{2-33}} = \Bl_L \P^3 \to \P^3$. 
\end{enumerate}
\end{prop}

\begin{proof}
Only (B) of (\ref{n-pic4-blowup}) holds. 
By Theorem \ref{t-ele-tr-F1-pic4}, (1) and (2) hold.

It suffices to show (3). 
By (2) and Theorem \ref{t-ele-tr-F1-pic4}, 
there is a blowup $\sigma: X \to Y := Y_{\text{3-30}}$ 
along a smooth curve $B_Y$ with $-K_Y \cdot B_Y = 1$. 
Recall that $Y_{\text{3-30}} = \Bl_{C_1}\,Y_{\text{2-33}}$ 
for a one-dimensional fibre $C_1$ of the induced blowup $\rho : Y_{\text{2-33}} = \Bl_L\,\P^3 \to \P^3$ 
(Lemma \ref{l-P3-V7}, Proposition \ref{p-pic3-30}), 
where $L$ is a line on $\P^3$. 
Let $\widetilde{\rho} : Y_{\text{3-30}}  \xrightarrow{\sigma'} Y_{\text{2-33}} \xrightarrow{\rho} \P^3$ be the induced birational morphism. 
Then $\Ex(\widetilde{\rho}) \cap B_Y \neq \emptyset$, as otherwise we would get a contradiction: 
$-1 = K_Y \cdot B_Y = \widetilde{\rho}^*K_{\P^3} \cdot B_Y \in 4\Z$. 
For $D := \Ex(\sigma')$, 
it is enough to show that $D \cap B_Y = \emptyset$ (Lemma \ref{l-line-meeting}). 
Suppose $D \cap B_Y \neq \emptyset$. 
Lemma \ref{l-line-meeting} implies that $B_Y$ is a fibre of $D \to \sigma'(D)$. 
Again by Lemma \ref{l-line-meeting}, it is enough to show  $-K_{Y} \cdot C =1$ for the intersection $C := D \cap E_{\text{3-30}}$, 
where $E_{\text{3-30}}$ denotes the proper transform of $E := \Ex(\rho)$. 
This follows from the following: 
\[
-K_{Y} \cdot C = (-K_Y)|_{E_{\text{3-30}}} \cdot C= (-K_Y)|_{E_{\text{3-30}}} \cdot  \zeta_Y = (-K_{Y_{\text{2-33}}})|_E \cdot \sigma'(\zeta_Y) =1, 
\]
where $\zeta_Y$ denotes a fibre of $E_{\text{3-30}} \xrightarrow{\simeq}E \to \rho(E)$ other than $C$. 
\end{proof}

\begin{prop}[No.\ {\hyperref[table-4-13]{4-13}}]\label{p-pic4-13}
Let $X$ be a Fano threefold with $\rho(X)=4$ and $(-K_X)^3=26$. 
Then the following hold. 
\begin{enumerate}
\item $\Blowdown(X) = \{ \text{3-27, 3-31}\}$. 
\item 
\begin{enumerate}
\item  There is a conic bundle structure $f_1:X \to \P^1 \times \P^1$ of type 3-27-vs-3-31 such that $\deg \Delta_{f_1} = (1, 3)$. 
\item  There is a conic bundle structure $f_2:X \to \P^1 \times \P^1$ of type 3-27-vs-non-Fano such that $\deg \Delta_{f_2} = (1, 1)$. 
\end{enumerate}
\item 
$X \simeq \Bl_C\,(\P^1 \times \P^1 \times \P^1)$, 
where $C$ is a smooth curve on $\P^1 \times \P^1 \times \P^1$ of tridegree $(1, 1, 3)$. 
\end{enumerate}
\end{prop}

\begin{proof}
Only (A) of (\ref{n-pic4-blowup}) holds. 
By Theorem \ref{t-ele-tr-P1P1-pic4}, 
there is a blowup $f: X \to  Y := \P^1 \times \P^1 \times \P^1$ along a smooth curve $B_Y$ 
such that $-K_Y \cdot B_Y = 10$ and $(d_1, d_2) =(1, 1)$ for the tridegree $(d_1, d_2, d_3)$ of $B_Y \subset  \P^1 \times \P^1 \times \P^1$. 
By $2(d_1+d_2+d_3) =-K_Y \cdot B_Y = 10$, 
we obtain $(d_1, d_2, d_3) = (1, 1, 3)$. 
Hence  (3) holds. 
Taking suitable projections $\P^1 \times \P^1 \times \P^1 \to \P^1 \times \P^1$, 
we obtain (2) by Theorem \ref{t-ele-tr-P1P1-pic4}. 
In particular, we get $\Blowdown(X) \supset \{ \text{3-27, 3-31}\}$. 
Since only (A) of \ref{n-pic4-blowup} holds, 
the opposite inclusion $\Blowdown(X) \subset \{ \text{3-27, 3-31}\}$ follows from Theorem \ref{t-ele-tr-P1P1-pic4}. 
\end{proof}

\begin{dfn}\label{d-pic4}
Let $X$ be a Fano threefold with $\rho(X)=4$. 
We say that $X$ is {\em 4-xx} or of {\em No.\ 4-xx} 
if $(-K_X)^3$ and $\Blowdown(X)$ satisfies the corresponding properties listed in Table \ref{table-pic4} in Subsection \ref{ss-table-pic4}. 
For example, the definitions of No.\ 4-1 and 4-4 are as follows. 
\begin{itemize}
\item A Fano threefold $X$ is  {\em 4-1} or  of {\em No.\ 4-1} if 
$\rho(X)=4$, $(-K_X)^3=24$, and $\Blowdown(X) = \{ \text{3-27}\}$. 
\item  A Fano threefold $X$ is   {\em 4-4} or  {\em of No.\ 4-4} if $\rho(X)=4$, $(-K_X)^3=32$, and 
$\Blowdown(X) = \{ \text{3-18, 3-19, 3-30}\}$. 
\end{itemize}
\end{dfn}

\begin{thm}\label{t-pic4-main}
Let $X$ be a Fano threefold with $\rho(X)=4$. 
Then $X$ satisfies one and only one of the possibilities listed in Table \ref{table-pic4} in Subsection \ref{ss-table-pic4},  except for the column \lq\lq blowups". 
\end{thm}

\begin{proof}
The assertion follows from results in this subsection. 
For example, if $X$ is a Fano threefold with $\rho(X)=4$ and $(-K_X)^3=32$, 
then the assertion follows from Proposition \ref{p-pic4-4} and Proposition \ref{p-pic4-5}. 
\end{proof}

\section{$\rho \geq 5$}\label{s-pic5}


The purpose of this subsection is to classify Fano threefolds of Picard number $\geq 5$. 
The main part is the case of  Picard number $5$ (Theorem \ref{t-pic5-main}). 
Given a Fano threefold $X$ with $\rho(X)=5$ and $X \not\simeq S \times \P^1$, 
a key step is to 
prove $X \simeq \Bl_{C \amalg C'}\,Z$ 
for some Fano threefold $Z$ and 
mutually disjoint smooth curves $C$ and $C'$ on $Z$ 
(Lemma \ref{l-pic5-F1}). 
Let us introduce $\Blowdown(X)$ as before. 

\begin{dfn}
Let $X$ be a Fano threefold with $\rho(X)=5$. 
We define the finite set 
\[
\Blowdown(X) \subset \{ \text{4-1, 4-2, 4-3, ..., 4-13}\}
\]
by the following condition: $\text{4-xx} \in \Blowdown(X)$ if and only if 
there exist a Fano threefold $Y_{\text{4-xx}}$ of No.\ 4-xx and a smooth curve $C$ on $Y_{\text{4-xx}}$ such that $X$ is isomorphic to the blowup of $Y_{\text{4-xx}}$ along $C$. 
\end{dfn}


\begin{prop}\label{p-dP-prod-R}
Let $T$ be a smooth del Pezzo surface with $K_T^2 \leq 7$ and set $X := T \times \P^1$. 
Take an extremal ray $R$ of $\NE(X)$.  
Then one and only one of the following holds. 
\begin{enumerate}
\item The contraction of $R$ is the first projection $\pr_1 : X =T \times \P^1 \to T$. 
\item There exists a $(-1)$-curve $\ell_T$ on $T$ such that 
the contraction of $R$ is 
the induced birational morphism $\psi \times \P^1 : X = T\times \P^1 \to T' \times \P^1$ 
for the blowdown $\psi : T \to T'$ of $\ell_T$. 
\end{enumerate}
\end{prop}

\begin{proof}
Fix an extremal rational curve $\ell$ on $X$ satisfying $R = \R_{\geq 0} [\ell]$. 
Let $p : X = T \times \P^1 \to T$ and $q : X \to \P^1$ be the first and second projections, respectively. 
We have 
\[
-K_X \sim -p^*K_T - q^*K_{\P^1}. 
\]
Since $-p^*K_T$ and $-q^*K_{\P^1}$ are nef, there are the following three cases. 
\begin{enumerate}
\renewcommand{\labelenumi}{(\roman{enumi})}
\item $(p^*K_T) \cdot \ell =0$. 
\item $(q^*K_{\P^1}) \cdot \ell =0$. 
\item $(-p^*K_T) \cdot \ell >0$ and $(-q^*K_{\P^1}) \cdot C>0$. 
\end{enumerate}

(i) Assume $(p^*K_T) \cdot \ell =0$, i.e., $p(\ell)$ is a point for the first projection 
$p : X= T \times \P^1 \to T$. 
In this case, (1) holds. 

\medskip

(ii) Assume $(q^*K_{\P^1}) \cdot \ell =0$, i.e., $q(\ell)$ is a point for 
the second projection $q: X= T \times \P^1 \to \P^1$. 
For the point $s :=q(\ell) \in \P^1$, we have $\ell \subset T \times \{ s\}$. 
By $K_T^2 \leq 7$, $\ell$ is a $(-1)$-curve on $T \times  \{s\} (\simeq T)$. 
Let $\ell_T$ be the corresponding $(-1)$-curve on $T$, i.e., 
it holds that $\ell = \ell_T \times \{s\}$. 
Let $\psi : T \to T'$ be the contraction of $\ell_T$. 
Let us prove that $\psi \times \P^1 : X =T \times \P^1 \to T' \times \P^1$ 
is the contraction of $R$. 
Recall that $\ell = \ell_T \times \{s\}$. 
Fix another closed point $s' \in \P^1$ and set $\ell' := \ell_T \times \{s'\}$. 
It suffices to show that $\ell \equiv \ell'$, i.e., $D \cdot \ell = D \cdot \ell'$ for a Cartier divisor  $D$ on $X$. 
By $\Pic\, X \simeq \Pic\,T \times \Pic\,\P^1$, we have $D \sim p^*D_T + q^*D_{\P^1}$ 
for some Cartier divisors $D_T$ on $T$ and $D_{\P^1}$ on $\P^1$. 
It holds that 
$D \cdot \ell = (p^*D_T + q^*D_{\P^1}) \cdot \ell = p^*D_T \cdot \ell =D_T \cdot \ell_T$. 
Similarly, we obtain $D \cdot \ell' = D_T \cdot \ell_T$. Thus (2) holds. 

\medskip

(iii) 
Assume $-p^*K_T \cdot \ell >0$ and $-q^*K_{\P^1} \cdot \ell >0$. 
It suffices to derive a contradiction. 
In this case, we obtain $-K_X \cdot \ell = (-p^*K_T  -q^*K_{\P^1})\cdot \ell \geq 2$. 
Then $R$ is of type $E_2$ or $C_2$, because  the length $-K_X \cdot \ell$ of $R$ is $\geq 2$. 
In particular, $-K_X \cdot \ell=2$, and hence 
$-p^*K_T \cdot \ell =-q^*K_{\P^1} \cdot \ell=1$. 

Suppose that $R$ is of type $E_2$. 
Let $\sigma: X \to Y$ be the contraction of $R$. 
Set $E := \Ex(\sigma) \simeq \P^2$. 
Since the composition $\ell \hookrightarrow X =T \times \P^1 \xrightarrow{\pr_2} \P^1$ is surjective, 
so is $E \hookrightarrow X =T \times \P^1 \xrightarrow{\pr_2} \P^1$ (note that $\ell \subset E$). 
This is absurd, because $E (\simeq \P^2)$  dominates no curve.

Hence $R$ is of type $C_2$. 
Then its contraction $p' : X \to T'$ is trivial (Proposition \ref{p-FCB-triv}), i.e., $X \simeq T' \times \P^1$ and $f$ is the projection. Let $q' : X \xrightarrow{\simeq} T' \times \P^1 \to \P^1$ be the second projection. 
Let $F_q$ and $F_{q'}$ be the extremal faces of $q$ and $q'$, respectively. 
If 
$F_q = F_{q'}$, then we obtain  
\[
2 = (-K_X) \cdot \ell =(-p'^*K_{T'} -q'^*K_{\P^1}) \cdot \ell = 
-p'^*K_{T'}\cdot \ell  -q^*K_{\P^1}\cdot \ell = 0 +1, 
\]
which is a contradiction. 

Hence it is enough show $F_q = F_{q'}$. 
Fix a closed point $s \in \P^1$ and a $(-1)$-curve $m$ on $T \times \{s\}$. 
By $-p^*K_T -q^*K_{\P^1} \sim -K_X \sim -p'^*K_{T'} -q'^*K_{\P^1}$, 
it holds that  
\[
1 = (-p^*K_T -q^*K_{\P^1}) \cdot m = -K_X \cdot m = (-p'^*K_{T'} -q'^*K_{\P^1}) \cdot m. 
\]
By 
$-p'^*K_{T'} \cdot m \in \Z_{\geq 0}$ and $-q'^*K_{\P^1} \cdot m \in \Z_{\geq 0}$, we get $-p'^*K_{T'} \cdot m=0$ or $-q'^*K_{\P^1} \cdot m=0$. 
Note that $m$ is not a fibre of $p' : X=T' \times \P^1 \to T'$, because 
a fibre $\zeta'$ of $p'$ satisfies $-K_X \cdot \zeta' =2$. 
Hence we obtain $-p'^*K_{T'} \cdot m =1$ and $-q'^*K_{\P^1} \cdot m =0$. 
Then $m$ is contained in a fibre of $q' : X = T' \times \P^1 \to \P^1$. 
Therefore, we get $F_q \subset F_{q'}$, 
because $F_q$ is generated by the $(-1)$-curves $m$ on $T \times \{s\}$ (recall that $K_T^2 \leq 7$). 
We then get a factorisation $q' : X \xrightarrow{q} \P^1 \xrightarrow{\theta} \P^1$. 
Since $\theta : \P^1 \to \P^1$ is a contraction, 
$\theta$ is automatically an isomorphism. 
Therefore, we obtain $F_q = F_{q'}$, 
\end{proof}

\begin{lem}\label{l-pic5-blowdown}
Let $Y$ be a Fano threefold with $\rho(X)=4$ and let $\sigma: X \to Y$ be a blowup along a smooth curve $B_Y$ such that $X$ is Fano. 
Then the following hold. 
\begin{enumerate}
\item The No.\ of $Y$ is one of $\text{4-4, 4-9, 4-10, 4-11, 4-12}$. In particular, 
\[
\Blowdown(X) \subset \{ \text{4-4, 4-9, 4-10, 4-11, 4-12} \}. 
\]
\item 
$\text{4-10} \in \Blowdown(X)$ if and only if $X \simeq T \times \P^1$ for a smooth del Pezzo surface $T$  with $K_T^2 =6$. 
\item 
If $Y$ is not 4-10, then 
$Y$ has a conic bundle structure $g : Y \to \F_1$ and $B_Y$ is a regular subsection of $g$. 
\end{enumerate}
\end{lem}

\begin{proof}
Let us show (2). 
Assume that $Y$ is 4-10. 
Then $Y \simeq \P^1 \times U$ for a smooth del Pezzo surface $U$ with $K_U^2 =7$ 
(Subsection \ref{ss-table-pic4}). 
It follows from  Proposition \ref{p-FCB-centre} and Proposition \ref{p-FCB-triv} 
that $X \simeq \P^1 \times T$ for a smooth projective surface $T$. 
Since $T$ is a del Pezzo surface with $\rho(T) = \rho(X) -1 =4$ (Proposition \ref{p-FCB-dP}), we get $K_T^2 =6$. 
Conversely, if  $X \simeq T \times \P^1$, then 
$\text{4-10} \in \Blowdown(X)$ by Proposition \ref{p-dP-prod-R}. 
This completes the proof of (2).

Let us show (1) and (3). 
In what follows, we assume that $Y$ is not 4-10. 
By Proposition \ref{p-FCB-centre} and Proposition \ref{p-FCB-triv}, 
there exists no conic bundle structure $h : Y \to \P^1 \times \P^1$ 
such that $\Delta_h$ is ample. 
By the classification list (Subsection \ref{ss-table-pic4}), 
$Y$ is one of 4-4, 4-9, 4-11, 4-12. 
Thus (1) holds. 
Again by  the classification list (Subsection \ref{ss-table-pic4}), 
there is a Fano conic bundle $g: Y \to \F_1$. 
Then $B_Y$ is a regular subsection of $g$ (Proposition \ref{p-FCB-triv}). 
Thus (3) holds. 
\end{proof}

\begin{lem}\label{l-pic5-F1}
Let $Y$ be a Fano threefold with $\rho(Y)=4$. 
Let $\sigma : X \to Y$ be a blowup along a smooth curve $B_Y$. 
Assume that $X$ is Fano and $Y$ is not of No.\ 4-10. 
Then there exists a commutative  diagram 
\begin{equation}
\begin{tikzcd}
    & X \simeq \Bl_{C \amalg C'}\,Z \arrow[ld, "\sigma"'] \arrow[rd, "\sigma'"] \arrow[dd, "\varphi"]\\
    Y\simeq \Bl_C\,Z \arrow[rd, "\rho"] \arrow[rdd, bend right, "g"'] & & Y':=\Bl_{C'}\,Z \arrow[ld, "\rho'"'] \arrow[ldd, bend left, "g'"]\\
    & Z \arrow[d, "h"]\\
    & \F_1
\end{tikzcd}
\end{equation}
such that 
\begin{enumerate}
\item $g, g', h$ are Fano conic bundles, 
\item $\rho : Y \xrightarrow{\simeq} \Bl_C\,Z \to Z$ is a blowup along a regular subsection $C$ of $h$, 
\item $\rho' :  Y' = \Bl_C\,Z \to Z$  is a blowup along a regular subsection $C'$ of $h$, 
\item $C \cap C' = \emptyset$, and   
\item $\sigma'$ and $\varphi$ are the induced blowups. 
\end{enumerate}
In particular, we obtain $X \simeq Y \times_Z Y'$. 
Moreover, one of the following holds. 
\begin{enumerate}
\item[(A)] $(Y, Y', Z)$ is (4-4, 4-12, 3-30) or  (4-12, 4-4, 3-30). 
In this case, $(-K_X)^3 = 28$ ($X$ is 5-1). 
\item[(B)] $(Y, Y', Z)$ is (4-9, 4-11, 3-28) or  (4-11, 4-9, 3-28). 
In this case, $(-K_X)^3 = 36$ ($X$ is 5-2). 
\item[(C)] $(Y, Y', Z)$ is (4-9, 4-12, 3-30) or  (4-12, 4-9, 3-30). 
In this case, $(-K_X)^3 = 36$ ($X$ is 5-2). 
\end{enumerate}
\end{lem}




    \begin{center}
\begin{longtable}{cclcccc}
No.\ & $(-K_Y)^3$ & \hspace{5mm}conic bundles /$\F_1$  & $\Delta$  &  \\ \hline
4-4 & $32$ 
    &  3-30-vs-3-30  
    & $\tau^*\MO_{\P^2}(2)$ & 
    &   \\ \hline
4-9 & $40$ 
     & 3-28-vs-3-30  
     &   $\tau^*\MO_{\P^2}(1)$& \\ \hline
4-11 & $44$ 
    & 3-28-vs-non-Fano   & $(-1)$ & \\
      \hline
4-12 & $46$ 
    & 3-30-vs-non-Fano  & $(-1)$ & 
    &   \\ \hline
      \end{longtable}
  \end{center} 

\begin{proof}
Since $Y$ is not 4-10, 
$Y$ is one of 4-4, 4-9, 4-11, 4-12 (Lemma \ref{l-pic5-blowdown}). 
Again by Lemma \ref{l-pic5-blowdown}, 
$Y$ has a conic bundle structure $g : Y \to \F_1$ and $B_Y$ is a regular subsection of $g$. 
Moreover, we get 
\[
g : Y \xrightarrow{\rho} Z \xrightarrow{h} \F_1, 
\]
where $h: Z \to \F_1$ is a Fano conic bundle with $\rho(Z)=3$ and 
$\rho$ is a blowup along a regular subsection of $h$ (Proposition \ref{p-smaller-same-base}). 
The possibilities for $(Y, Z)$ are  as above (Subsection \ref{ss-table-pic4}), i.e., one of 
$(\text{4-4, 3-30}), (\text{4-9, 3-28}), (\text{4-9, 3-30}), (\text{4-11, 3-28}), (\text{4-12, 3-30})$. 

In what follows, we only treat the case when $(Y, Z)$ is $(\text{4-4, 3-30})$, 
as the proofs for all the  cases are identical. 
Set $Y_{\text{4-4}} := Y$ and $Z_{\text{3-30}} :=Z$. 
In this case, $\Delta_g \in |\tau^*\MO_{\P^2}(2)|$, and hence $\Delta_g \cap \Gamma =\emptyset$ (Corollary \ref{c-FCB-(-1)-2}), 
where $\Gamma$ denotes the $(-1)$-curve on $\F_1$. 
Again by Corollary \ref{c-FCB-(-1)-2}, 
the blowup centre $B_Y$ of $\sigma : X \to Y$ satisfies $g(B_Y)=\Gamma$. 
Hence we obtain $X \simeq Y_{\text{4-4}} \times_{Z_{\text{3-30}}} Y'$ 
for $Y' := \Bl_{C'}\,Z_{\text{3-30}}$ and $C' := \rho(B_Y)$. 
Then $Y'$ is Fano (Corollary \ref{c-disjoint-blowup}), which must be of No.\ 4-12 by the above table. 
Thus (1)-(5) hold. 

We have that $(Y, Y', Z)$ is (4-4, 4-12, 3-30). 
Then it follows from Lemma \ref{l-blowup-formula} that
\[
(-K_X)^3 -(-K_Z)^3  = ( (-K_Y)^3 - (-K_Z)^3 ) + ((-K_{Y'})^3 - (-K_Z)^3), 
\]
which implies $(-K_X)^3 = (-K_Y)^3 +(-K_{Y'})^3 -(-K_Z)^3 = 32 +46-50=28$. 
Thus (A) holds. 
\qedhere





\end{proof}

\begin{prop}[No.\ {\hyperref[table-5-1]{5-1}}]\label{p-pic5-1}
Let $X$ be a Fano threefold with $\rho(X)=5$ and $(-K_X)^3 = 28$. 
Then the following hold. 
\begin{enumerate}
\item $\Blowdown(X) = \{ \text{4-4, 4-12}\}$. 
\item $X \simeq Y_{\text{4-4}} \times_{Z_{\text{3-30}}} Y_{\text{4-12}}$, 
where 
each of $Y_{\text{4-4}} \to Z_{\text{3-30}}$ and $Y_{\text{4-12}} \to Z_{\text{3-30}}$ is a blowup of $Z_{\text{3-30}}$ 
along a regular subsection of a Fano conic bundle $Z_{\text{3-30}} \to \F_1$. 
\item 
$X \simeq \Bl_{B_1 \amalg B_2 \amalg B_3}\,Y_{\text{2-29}}$, 
where $Y_{\text{2-29}} := \Bl_C\,Q$ for a conic $C$ and 
$B_1, B_2, B_3$ are mutually distinct  one-dimensional fibres of 
the induced blowup $\rho: Y_{\text{2-29}} := \Bl_C\,Q \to Q$. 
\end{enumerate}
\end{prop}

\begin{proof}
By $(-K_X)^3 =28$, we obtain $\text{4-10} \not\in \Blowdown(X)$ (Lemma \ref{l-pic5-blowdown}). 
Only (A) of Lemma \ref{l-pic5-F1} holds. 
Thus (1) and (2) hold (Lemma \ref{l-pic5-blowdown}(3), Lemma \ref{l-pic5-F1}). 

Let us show (3). 
We have a blowup  $\sigma : X \to Y = Y_{\text{4-4}}$ along a smooth curve $B_Y$. 
By $(-K_Y)^3 - (-K_X)^3 = 32- 28 =4$, we get 
$(p_a(B_Y), -K_Y \cdot B_Y) \in \{ (0, 1), (1, 2), (2, 3)\}$ (Lemma \ref{l-blowup-formula2}). 
Recall that $Y_{\text{4-4}} =\Bl_{B_1 \amalg B_2}\,Z_{\text{2-29}}$ and 
$Z_{\text{2-29}} = \Bl_C\,Q$, where $C$ is a conic on $Q$ and $B_1$ and $B_2$ 
are mutually distinct one-dimensional fibres of the induced blowup 
$\psi : \Bl_C\,Q \to Q$ (Proposition \ref{p-pic4-4}): 
\[
X \xrightarrow{\sigma} Y =Y_{\text{4-4}}=\Bl_{B_1 \amalg B_2}\,Z_{\text{2-29}} \xrightarrow{\varphi} 
Z= Z_{\text{2-29}} =\Bl_C\,Q
 \xrightarrow{\psi} Q. 
\]
Set $E_Z :=\Ex(\psi)$ and $E_Y := \varphi_*^{-1}E_Z$, which is the strict transform of $E_Z$ on $Y$. 
Let $D_1$ and $D_2$ be the $\varphi$-exceptional prime divisors lying over $B_1$ and $B_2$, respectively. 
It is enough to show (i) and (ii) below.  
\begin{enumerate}
\item[(i)] $B_Y$ is disjoint from $D_1 \amalg D_2$. 
\item[(ii)] $B_Z := \varphi(B_Y)$ is a fibre of $E_Z \to C$. 
\end{enumerate}
Indeed, (i) and (ii) implies that $B_1, B_2, B_3 := B_Z$ are mutually distinct fibres of $E_Z \to C$, 
and hence (3) holds.

(i) 
Suppose that $B_Y$ intersects $D_1 \amalg D_2$. 
By symmetry, we may assume that $B_Y \cap D_1 \neq \emptyset$. 
It follows from Lemma \ref{l-line-meeting} that 
$B_Y$ is a fibre of $D_1 \to \varphi(D_1) = B_1$. 
Then $B_Y$ properly intersects $D_1 \cap E_Y =: \zeta_1$. 
In order to derive a contradiction, it is enough to prove $-K_Y \cdot \zeta_1 =1$ 
(Lemma \ref{l-line-meeting}). 
This follows from 
\[
-K_Y \cdot \zeta_1 = (-K_Y|_{E_Y}) \cdot \zeta_1 = (-K_Y|_{E_Y}) \cdot \xi_{E_Y} = (-K_Z|_{E_Z}) \cdot \xi_{E_Z} = -K_Z \cdot \xi_{E_Z} = 1, 
\]
where $\xi_{E_Y}$ and $\xi_{E_Z}$ denote general fibres of 
$E_Y \xrightarrow{\simeq} E_Z \to C$ and $E_Z \to C$, respectively. 

(ii) 
We now show that $B_Y$ intersects $\Ex(\psi \circ \varphi : Y \to Q)$. 
Otherwise, 
for the image $B_Q (\simeq B_Y) \subset Q$ of $B_Y$, 
we would get $\{1, 2, 3\} \ni -K_Y \cdot B_Y = -K_Q\cdot B_Q \in 3\Z$, 
which implies $(p_a(B_Y), -K_Y \cdot B_Y) = (2, 3)$. 
However, the image $B_Q (\simeq B_Y)$ of $B_Y$ is a line on $Q$ by $-K_Q \cdot B_Q=1$. This contradicts $p_a(B_Q) =p_a(B_Y)=2$. 

Thus $B_Y$ intersects $\Ex(\psi \circ \varphi : Y \to Q) = E_Y \cup D_1 \cup D_2$. 
By (i) and Lemma \ref{l-line-meeting}, $B_Y$ is a fibre of $E_Y \to C$ which is disjoint from $D_1 \cup D_2$. 
Then $B_Z = \varphi(B_Y)$  is a fibre of  $E_Z \to C$ which is different from $B_1$ and $B_2$. 
Thus (ii) holds. 
\end{proof}





\begin{prop}[No.\ {\hyperref[table-5-2]{5-2}}]\label{p-pic5-2}
Let $X$ be a Fano threefold with $\rho(X)=5$ and $(-K_X)^3 = 36$. 
Assume $\text{4-10} \not\in \Blowdown(X)$. 
Then the following hold. 
\begin{enumerate}
\item $\Blowdown(X) = \{ \text{4-9, 4-11, 4-12}\}$. 
\item $X \simeq Y_{\text{4-9}} \times_{Z_{\text{3-28}}} Y_{\text{4-11}}$, 
where 
each of $Y_{\text{4-9}} \to Z_{\text{3-28}}$ and $Y_{\text{4-11}} \to Z_{\text{3-28}}$ is a blowup of $Z_{\text{3-28}}$ 
along a regular subsection of a Fano conic bundle $Z_{\text{3-28}} \to \F_1$. 
\item $X \simeq Y_{\text{4-9}} \times_{Z_{\text{3-30}}} Y_{\text{4-12}}$, 
where 
each of $Y_{\text{4-9}} \to Z_{\text{3-30}}$ and $Y_{\text{4-12}} \to Z_{\text{3-30}}$ is a blowup of $Z_{\text{3-30}}$ 
along a regular subsection of a Fano conic bundle $Z_{\text{3-30}} \to \F_1$. 
\item 
$X \simeq \Bl_{B \amalg B'}\, Y_{\text{3-25}}$, 
where $Y_{\text{3-25}} := \Bl_{L_1 \amalg L_2}\,\P^3$ 
for a mutually disjoint lines $L_1$ and $L_2$, and 
both $B$ and $B'$ are one-dimensional fibres of 
the induced blowup $\rho: Y_{\text{3-25}} =\Bl_{L_1 \amalg L_2}\,\P^3 \to \P^3$ which are lying over $L_1$. 
\end{enumerate}
\end{prop}

\begin{proof}
First of all, we show that both (2) and (3) hold. 
By $\text{4-10} \not\in \Blowdown(X)$, (2) or (3) holds (Lemma \ref{l-pic5-F1}). 
In any case, we get a blowup $\sigma: X \to Y:= Y_{\text{4-9}}$ along a regular subsection $B_Y$ 
of a conic bundle structure $g: Y_{\text{4-9}} \to \F_1$ such that $B_{\F_1} := g(B_Y)$ is the $(-1)$-curve on $\F_1$ (Theorem \ref{t-ele-tr-F1-pic4}). 
By $p_a(B_Y)=0$ and $(-K_Y)^3  -(-K_X)^3= 40-36 =4$, 
it follows from Lemma \ref{l-blowup-formula} that $-K_Y \cdot B_Y = 1$. 
By $\rho(Y_{\text{4-9}}) > \rho(\F_1)+1$ and Theorem \ref{t-ele-tr-F1-pic4}, we have an elementary transform over $\F_1$ of type 3-28-vs-3-30 such that $\Delta_g$ is disjoint from the $(-1)$-curve  $B_{\F_1}$ on $\F_1$: 
\[
\begin{tikzcd}
& Y_{\text{4-9}} \arrow[ld, "\tau"'] \arrow[rd, "\tau'"] \arrow[dd, "g"]\\
Z_{\text{3-28}} \arrow[rd, "h"']& & Z_{\text{3-30}} \arrow[ld, "h'"]\\
& \F_1.
\end{tikzcd}
\]
In particular, the blowup centre $C_{\text{3-28}}$ (resp. $C_{\text{3-30}}$) 
of $\tau: Y_{\text{4-9}} \to Z_{\text{3-28}}$ 
(resp. $\tau': Y_{\text{4-9}} \to Z_{\text{3-30}}$) 
is disjoint from the inverse image of $B_{\F_1}$. 
Then $\tau$ and $\tau'$ satisfy (B) and (C) of Lemma \ref{l-pic5-F1}, respectively. 
Therefore, both (2) and (3) hold. 

Let us show (1). 
By (2), (3), and Lemma \ref{l-pic5-blowdown}, we get 
\[
\{ \text{4-9, 4-11, 4-12}\} \subset \Blowdown(X) \subset 
\{ \text{4-4, 4-9, 4-10, 4-11, 4-12}\}. 
\]
If $\text{4-4} \in \Blowdown(X)$, then $(-K_X)^3=28$ (Lemma \ref{l-pic5-F1}), which is absurd. 
This, together with  $\text{4-10} \not\in \Blowdown(X)$, implies (1).

Let us show (4). 
Recall that $\sigma: X \to Y= Y_{\text{4-9}}$ is a blowup along a smooth curve $B_Y$ satisfying $-K_Y \cdot B_Y = 1$. 
By Proposition \ref{p-pic4-9}, 
we obtain $Y_{\text{4-9}} = \Bl_B\,Z_{\text{3-25}}$, 
where $Z_{\text{3-25}} := \Bl_{L_1 \amalg L_2}\,\P^3$ for a mutually disjoint lines $L_1$ and $L_2$ on $\P^3$ 
and $B$ is a one-dimensional fibre of the induced blowup $\psi : Z_{\text{3-25}} := \Bl_{L_1 \amalg L_2}\,\P^3 \to \P^3$ lying over $L_1$. 
We get the following induced blowups: 
\[
X \xrightarrow{\sigma} Y = Y_{\text{4-9}} = \Bl_C\,Z_{\text{3-25}} 
\xrightarrow{\varphi} Z=Z_{\text{3-25}} = \Bl_{L_1 \amalg L_2}\,\P^3 \xrightarrow{\psi} \P^3. 
\]
It follows from  $-K_Y \cdot B_Y =1$ that $B_Y$ intersects $\Ex(\varphi \circ \psi)$, as otherwise we would get 
$-1 = K_Y \cdot B_Y = (\varphi\circ \psi)^*K_{\P^3} \cdot B_Y \in 4\Z$. 
By the same argument as in the proof of Proposition \ref{p-pic5-1} (especially, (i) and (ii) in the proof),  
\begin{itemize}
\item $B_Y$ is disjoint from $\Ex(\varphi)$ and 
\item $B' :=\varphi(B_Y)$ is a one-dimensional fibre of $\psi : Z = \Bl_{L_1 \amalg L_2}\,\P^3 \to \P^3$ which is disjoint from $B$. 
\end{itemize}
In particular, $X \simeq \Bl_{B \amalg B'}\,Z_{\text{3-25}}$ and $\psi(B') \in L_1 \amalg L_2$. 
It suffices to show that $\psi(B') \in L_1$. 
Suppose $\psi(B') \in L_2$. 
For each $i \in \{1, 2\}$, 
let $E_i^Z$ be the $\psi$-exceptional prime divisor lying over $L_i$, and 
let $E_i^X$ be the strict transform of $E_i^Z$ on $X$.  
Let $D$ and $D'$ be the $(\varphi \circ \sigma)$-exceptional prime divisors lying over $B$ and $B'$, respectively. 
Then $K_Z \sim \psi^*K_{\P^3} + E_1^Z + E_2^Z$ and 
\[
K_X \sim \sigma^*\varphi^*K_Z +D +D'
\sim  \sigma^*\varphi^*\psi^*K_{\P^3} +2D +2D' +E_1^X +E_2^X. 
\]
For the line $L$ on $\P^3$ passing through the points $\psi(B)$ and $\psi(B')$, 
its proper transform $L_X$ on $X$ intersects both $D_1$ and $D_2$. 
Hence we get the following contradiction: 
\[
0> K_X \cdot L_X = (\sigma^*\varphi^*\psi^*K_{\P^3} +2D_1 +2D_2 +E_1^X +E_2^X) \cdot L_X 
\geq -4 +2 +2+0+0 =0. 
\]
Thus (4) holds. 
\end{proof}

\begin{prop}[No.\ {\hyperref[table-5-3]{5-3}}]\label{p-pic5-3}
Let $X$ be a Fano threefold with $\rho(X)=5$ and $(-K_X)^3 = 36$. 
Assume $\text{4-10} \in \Blowdown(X)$. 
Then the following hold. 
\begin{enumerate}
\item $\Blowdown(X) = \{ \text{4-10}\}$. 
\item 
$X \simeq T \times \P^1$ for a smooth del Pezzo surface $T$ with $K_T^2 =6$. 
\end{enumerate}
\end{prop}

\begin{proof}
By $\text{4-10} \in \Blowdown(X)$, (2) holds (Lemma \ref{l-pic5-blowdown}). 
The assertion  (1) follows from Proposition \ref{p-dP-prod-R}. 
\end{proof}

\begin{dfn}\label{d-pic5}
\begin{enumerate}
\item A Fano threefold $X$ is {\em 5-1} or of {\em No.\ 5-1} if 
$\rho(X)=5$, $(-K_X)^3 =28$, and $\Blowdown(X) = \{ \text{4-4, 4-12}\}$. 
\item A Fano threefold $X$ is {\em 5-2} or of {\em No.\ 5-2} if 
$\rho(X)=5$, $(-K_X)^3 =36$, and $\Blowdown(X) = \{ \text{4-9, 4-11, 4-12}\}$. 
\item A Fano threefold $X$ is {\em 5-3} or of {\em No.\ 5-3} if 
$\rho(X)=5$, $(-K_X)^3 =36$, and $\Blowdown(X) = \{ \text{4-10}\}$. 
\end{enumerate}
\end{dfn}

\begin{thm}\label{t-pic5-main}
Let $X$ be a Fano threefold with $\rho(X)=5$. 
Then $X$ satisfies one and only one of the possibilities listed in Table \ref{table-pic5} in Subsection \ref{ss-table-pic5},  except for the column \lq\lq blowups". 
\end{thm}

\begin{proof}
Since $X$ is imprimitive, 
we get $\Blowdown(X) \neq \emptyset$. 
If $\text{4-10} \in \Blowdown(X)$
 then $X$ is 5-3 and the assertion holds (Proposition \ref{p-pic5-3}). 
Assume that $\text{4-10} \not\in \Blowdown(X)$. 
Then Lemma \ref{l-pic5-F1} implies $(-K_X)^3 \in \{ 28, 36\}$. 
If $(-K_X)^3 =28$ (resp. $(-K_X)^3 =36$), 
then 
it follows from Proposition \ref{p-pic5-1} (resp. Proposition \ref{p-pic5-2}) 
that $X$ is 5-1 (resp. 5-2) and the assertion holds. 
\end{proof}

\begin{thm}\label{t-pic6}
Let $X$ be a Fano threefold with $\rho(X) \geq 6$. 
Then the following hold. 
\begin{enumerate}
\item Let $X \to Y$ be a blowup of a Fano threefold $Y$ along a smooth curve.  
Then $Y \simeq T \times \P^1$ for a smooth del Pezzo surface $T$. 
\item $X \simeq S \times \P^1$ for a smooth del Pezzo surface $S$. 
\end{enumerate}
\end{thm}

\begin{proof}
In what follows, we only treat the case when $\rho(X)=6$ 
(see Remark \ref{r t-pic6} for a rigorous proof). 

Let us show (1). 
Recall that $Y$ is 5-1, 5-2, or 5-3 (Theorem \ref{t-pic5-main}). 
It suffices prove that $Y$ is of No.\ 5-3. 
Suppose that $Y$ is of No.\ 5-1 or 5-2. 
By Proposition \ref{p-pic5-1} 
and Proposition \ref{p-pic5-2}, 
 there exists a conic bundle $g: Y \to \F_1$ such that 
$\Delta_g = \Gamma \amalg C$, where $\Gamma$ is the $(-1)$-curve 
on $\F_1$ and $C$ is a smooth curve on $\F_1$ disjoint from $\Gamma$. 
Since $-\Gamma +  m C$ is ample for a large integer $m \gg 0$, 
any curve $B$ on $\F_1$ intersects $\Delta_g$. 
Therefore, $B_Y$ must be a smooth fibre of $g : Y \to \F_1$ (Proposition \ref{p-FCB-centre}). 
However, this contradicts Proposition \ref{p-FCB-triv}, because the base change 
$X = X \times_{\F_1} T \to T$ is a non-trivial Fano conic bundle by $\Delta_g \neq \emptyset$, where $T \to \F_1$ is a blowup at a point. 
Thus (1) holds. 

Let us show (2). 
Since $X$ is imprimitive, 
there exists a blowup $X \to Y$ of a Fano threefold $Y$ along a smooth curve $B_Y$. 
By (1), we get  $Y = T \times \P^1$ for a smooth del Pezzo surface $T$. 
By Proposition \ref{p-FCB-centre} and Proposition \ref{p-FCB-triv}, $B_Y$ is a smooth fibre of the projection $Y = T \times \P^1 \to T$. 
Therefore, $X \simeq S \times \P^1$ for a smooth del Pezzo surface $S$. 
Thus (2) holds. 
\end{proof}

\begin{rem}\label{r t-pic6}
We can give a  rigorous proof of Theorem \ref{t-pic6} as follows. 
For an integer $\rho \geq 6$, consider the assertions $(1)_{\rho}$ and $(2)_{\rho}$ that are the assertions (1) and (2) of Theorem \ref{t-pic6} for the case when $\rho(X) = \rho$. 
In the above proof, we treat the case when $\rho(X)=6$, i.e., 
$(1)_6$ and $(2)_6$ hold. 
Moreover, the implication 
$(1)_{\rho} \Rightarrow (2)_{\rho}$ holds by the same argument as the above proof of $(2)_6$. 
The implication $(2)_{\rho} \Rightarrow (1)_{\rho+1}$ is clear. 
By induction, 
$(1)_{\rho}$ and $(2)_{\rho}$ hold for every integer $\rho \geq 6$. 
\end{rem}



\section{Classification tables}\label{s-table}

In this section, we provide the classification tables for Fano threefolds in characteristic $p>0$.  
In addition to terminologies given in Subsection \ref{ss-notation}, 
we now summarise notation used in this section. 

  \begin{enumerate}
\item We say that $f : X \to Y$ is a {\em split double cover} if 
$f$ is a finite surjective morphism of projective normal varieties such that $\MO_Y \to f_*\MO_X$ splits as an $\MO_Y$-module homomorphism and the induced field extension $K(X) \supset K(Y)$ is of degree two. For a split double cover $f : X \to Y$, we set $\mathcal L := (f_*\MO_X/\MO_Y)^{-1}$, which is an invertible sheaf on $Y$ \cite[Remark 2.2]{ATIII}. 
Note that all the double covers appearing in the following tables are split by 
\cite[Lemma 2.5]{ATIII} ($\rho(X) \neq 2$) and 
\cite[Table 5 in Section 9]{ATIII} ($\rho(X)=2$). 
\item 
In what follows, the centre of every blowup is assumed to be smooth whenever it is a curve. 
\item 
Let $f: X \to S$ be a Fano conic bundle. 
\begin{enumerate}
\item If $X$ is not of No.\ 2-24 nor 3-10, 
then $f$ is generically smooth (Theorem \ref{t-wild-cb}). 
In this case, $\Delta$ denotes the discriminant divisor  $\Delta_f$ of $f$. 
\item If $X$ is of No.\ 2-24 or 3-10, 
then 
$\Delta$ denotes the discriminant bundle $\Delta_f^{{\rm bdl}}$ of $f$. 
Recall that $\Delta_f^{{\rm bdl}} \simeq \MO_S(\Delta_f)$ when $f$ is generically smooth. 
\end{enumerate}
\item 
If an extremal ray is of type $D$, 
then $X_t$ denotes a fibre of its contraction $X \to \P^1$, where $t$ is  a closed point of $\P^1$. 
In particular, $(-K_X)^2 \cdot X_t$ coincides with $(-K_{X_K})^2$ for the generic fibre $X_K$. 
\item 
If an extremal ray is of type $E$, 
then $f: X \to Y$ denotes its contraction. 
Moreover, if $f$ is of type $E_1$, then $C$ denotes its blowup centre, which is a smooth curve on $Y$. 
\item $\tau : \F_1 \to \P^2$ denotes the blowdown of the $(-1)$-curve on $\F_1$. 
\item For an integer $d$ satisfying $1 \leq d \leq 7$, $S_d$ is a smooth del Pezzo surface with $K_{S_d}^2 = d$. 
\item 
Let $X$ be a Fano threefold. 
$\Blowdown(X)$ is the finite set defined by the following condition: 
$\text{x-yz} \in \Blowdown(X)$ if and only if 
there exist a Fano threefold $Y$ of No.\ x-yz and a smooth curve $C$ on $Y$ such that 
$X \simeq \Bl_C\,Y$. 
In the following tables, the column \lq\lq blowdowns" gives $\Blowdown(X)$. 
For example, if $X$ is a Fano threefold of No.\ 4-3 (resp. 3-1), then 
$\Blowdown(X) = \{ \text{3-17, 3-27, 3-28}\}$ 
(resp. $\Blowdown(X) = \emptyset$). 
\item 
Let $X$ be a Fano threefold. 
$\Blowup(X)$ is the finite set defined by the following condition: 
$\text{x-yz} \in \Blowup(X)$ if and only if 
there exists a smooth curve $C$ on $X$ such that 
$\Bl_C\,X$ is a Fano threefold of No.\ x-yz. 
In the following tables, 
$\Blowup(X)$ is contained in the list in the column \lq\lq blowups". 
For example, if $X$ is a Fano threefold of No.\ 2-24, then 
we have $\Blowup(X) \subset \{ \text{3-8}\}$. 
Note that this inclusion is not necessarily an equality (e.g., if $p=2$ and $X := \{ x_0y_0^2 +x_1y_1^2 + x_2y_2^2 =0\} \subset \P^2 \times \P^2$, then 
$X$ is a Fano threefold of No.\ 2-24 satisfying $\Blowup(X) = \emptyset$). 
\end{enumerate}

\begin{rem}\label{r column blowup}
Except for the column \lq\lq blowups",
the assertions in all the following tables have already established. 
The column \lq\lq blowups" is confirmed by comparing with another table. 
For example, the column  \lq\lq blowups" of No.\ 2-25 states \lq\lq 3-6, 3-11", 
because only 3-6 and 3-11 include  \lq\lq $E_1$: 2-25" in 
the $\rho=3$ table: Table \ref{table-pic3} (i.e., $\text{2-25} \in \Blowdown(X)$ if and only if $X$ is 3-6 or 3-11). 
\end{rem}

\subsection{$\rho=1$}\label{ss-table-pic1}

Let $X$ be a Fano threefold with $\rho(X)=1$. 
Let $r_X$ be the index of $X$. 
When $r_X=1$, we define $g$ by $(-K_X)^3  =2g -2$. 
Up to isomorphisms, one and only one of the following possibilities listed in Table \ref{table-pic1} occurs by 
\begin{itemize}
    \item \cite[Theorem 2.18 and Theorem 2.23]{TanI} ($r_X \geq 2$), 
    \item \cite[Theorem 1.1]{TanI} ($r_X =1$ and $|-K_X|$ is not very ample), and 
\item \cite[Theorem 1.1 and Proposition 2.8]{TanII} ($r_X =1$ and $|-K_X|$ is very ample). 
\end{itemize}

    \begin{center}
\begin{longtable}{cccp{10cm}ccc}
      \caption{$\rho(X)=1$}\label{table-pic1}\\
$r_X$ & $(-K_X)^3$ & $g$ & descriptions   &  &  \\ \hline
$1$ & $2$ & $2$ & $f: X \to \P^3$ is a split double cover with 
$\mathcal L \simeq \MO_{\P^3}(3)$.\\
\hline
$1$ & $4$ & $3$ & $X$ is a hypersurface in $\P^4$ of degree $4$, or 
$f: X \to Q$ is a split double cover with $\mathcal L \simeq \MO_{Q}(2)$.\\
\hline
$1$ & $6$ & $4$ & $X \subset \P^5$ is a complete intersection of a quadric hypersurface and a cubic hypersurface.  \\
\hline
$1$ & $8$ & $5$ & $X \subset \P^6$ is a complete intersection of three quadric hypersurfaces.  \\
\hline
$1$ & $10$ & $6$ &  \\
\hline
$1$ & $12$ & $7$ &  \\
\hline
$1$ & $14$ & $8$ &  \\
\hline
$1$ & $16$ & $9$ & \\
\hline
$1$ & $18$ & $10$ &  \\
\hline
$1$ & $22$ & $12$ &  \\
\hline
$2$ & $8 \cdot 1$ &  &  $X =V_1$, which is a weighted hypersurface in $\P(1, 1, 1, 2, 3)$ of degree $6$. \\
\hline
$2$ & $8 \cdot 2$ &  & $X =V_2$, which is a weighted hypersurface in $\P(1, 1, 1, 1, 2)$ of degree $4$. \\
\hline
$2$ & $8 \cdot 3$ &  & $X =V_3$, which is a cubic hypersurface in $\P^4$.\\
\hline
$2$ & $8 \cdot 4$ &  & $X =V_4$, which is a complete intersection of two quadric hypersurfaces in $\P^5$.\\
\hline
$2$ & $8 \cdot 5$ &  & $X =V_5$, which is an intersection of ${\rm Gr}(2, 5) \subset \P^9$ 
and a linear subvariety $L$ in $\P^9$ of codimension $3$.\\
\hline
$3$ & $54$ &  & $X =Q$. \\
\hline
$4$ & $64$ &  & $X = \P^3$. \\
\hline
      \end{longtable}
  \end{center}

\subsection{$\rho=2$}\label{ss-table-pic2}

The definition of 2-xx is given in Definition \ref{d-pic2}.  
For  a Fano threefold $X$ with $\rho(X)=2$, 
 one and only of the following possibilities listed in Table \ref{table-pic2} 
occurs up to isomorphisms. 
 Except for the column \lq\lq blowups", this is proven in \cite[Section 9]{ATIII} (cf. Subsection \ref{ss-pic2-III}). 
For the case when $\rho(X)=2$ (Table \ref{table-pic2}), we determine the types of the extremal rays. For example, if $X$ is a Fano threefold of No.\ 2-1, then one of the extremal rays is of type $D_1$ and the other is of type $E_1$.

  \begin{center}
\begin{longtable}{ccp{10cm}c}
      \caption{$\rho(X)=2$}\label{table-pic2}\\
No. & $(-K_X)^3$ & descriptions and extremal rays  & blowups\\ \hline
2-1 & $4$ & $D_1: (-K_X)^2 \cdot X_t=1$ & none\\ 
&  & 
$E_1:$ blowup of $V_1$ along an elliptic curve of degree $1$  which is a complete intersection of two members of $|-\frac{1}{2}K_{V_1}|$ & \\ \hline
2-2 & $6$ & $X$ is a split double cover of $\mathbb{P}^2\times \mathbb{P}^1$ with $\mathcal L \simeq \MO(2, 1)$. &  none\\ 
 &  & $C_1: \deg \Delta=8$, $X \xrightarrow{{\rm 2:1}} \P^2 \times \P^1 \xrightarrow{{\rm pr}_1} \P^2$  & \\ 
&  & $D_1: (-K_X)^2 \cdot X_t = 2$& \\ \hline
2-3 & $8$ & $D_1: (-K_X)^2 \cdot X_t=2$ &  none\\ 
 &  & $E_1:$ blowup of $V_2$ along an elliptic curve of degree $2$  which is a complete intersection of two members of $|-\frac{1}{2}K_{V_2}|$ & \\ \hline
2-4 & $10$ & $D_1: (-K_X)^2 \cdot X_t=3$ & none\\ 
 &  & $E_1:$ blowup of $\P^3$ along a curve of genus $10$ and degree $9$  which is a complete intersection of two cubic surfaces & \\ \hline
2-5 & $12$ & $D_1: (-K_X)^2 \cdot X_t =3$ &  none\\ 
&  & $E_1:$ blowup of $V_3$ along an elliptic curve of degree $3$  which is a complete intersection of two members of $|-\frac{1}{2}K_{V_3}|$ & \\ \hline
2-6 & $12$ & 
$X$ is a  divisor on $\mathbb{P}^2\times \mathbb{P}^2$ of bidegree $(2,2)$, or 
a split double cover of $W$ with $\mathcal L^{\otimes 2} \simeq \omega_W^{-1}$. &  none\\
 &  & $C_1: \deg \Delta=6$, $X \to \P^2 \times \P^2 \xrightarrow{{\rm pr_1}} \P^2$  & \\ 
&  &  $C_1: \deg \Delta=6$, $X \to \P^2 \times \P^2 \xrightarrow{{\rm pr_2}} \P^2$ & \\ \hline
2-7 & $14$ & $D_1: (-K_X)^2 \cdot X_t=4$ &  none\\ 
 &  & $E_1:$ blowup of $Q$ along a curve of genus $5$ and degree $8$ which is a complete intersection of two members of $|\MO_Q(2)|$ & \\ \hline
2-8 & $14$ &  $X$ is a split double cover of $V_7$ 
 with $\mathcal L^{\otimes 2} \simeq \omega_{V_7}^{-1}$ & none \\ 
 &  &  $C_1: \deg \Delta = 6$, $X \xrightarrow{{\rm 2:1}} 
 V_7 = \mathbb{P}_{\P^2}(\MO_{\mathbb{P}^2}\oplus \MO_{\mathbb{P}^2}(1)) \xrightarrow{{\rm pr}} \P^2$  & \\ 
&  & $E_3\,{\rm or}\,E_4$& \\ \hline
2-9 & $16$ & $C_1: \deg \Delta =5$ &  none\\ 
 &  & $E_1:$ blowup of $\P^3$ along a curve of genus $5$ and degree $7$ & \\ \hline
2-10 & $16$ & $D_1:$ $(-K_X)^2 \cdot X_t=4$ &  none\\ 
 &  & $E_1:$ blowup of $V_4$ along an elliptic curve of degree $4$  which is a complete intersection of two members of $|-\frac{1}{2}K_{V_4}|$ & \\ \hline
2-11 & $18$ & $C_1:$ $\deg \Delta =5$ &  none\\ 
 &  & $E_1:$ blowup of $V_3$ along a line & \\ \hline
2-12 & $20$ & $E_1:$ blowup of $\P^3$ along a curve of genus $3$ and degree $6$ &  none\\ 
 &  & $E_1:$ blowup of $\P^3$ along a curve of genus $3$ and degree $6$ & \\ \hline
2-13 & $20$ & $C_1:$ $\deg \Delta =4$ &  none\\ 
 &  & $E_1:$ blowup of $Q$ along a curve of genus $2$ and degree $6$ & \\ \hline
2-14 & $20$ & $D_1:$ $(-K_X)^2 \cdot X_t = 5$ &  none\\ 
 &  & $E_1:$ blowup of $V_5$ along an elliptic curve of degree $5$  which is a complete intersection of two members of $|-\frac{1}{2}K_{V_5}|$ & \\ \hline
2-15 & $22$ & $E_1:$ blowup of $\P^3$ along a curve of genus $4$ and degree $6$ &  none\\ 
&  &$E_3\,{\rm or}\,E_4$ & \\ \hline
2-16 & $22$ & $C_1: \deg \Delta =4$ &  none\\ 
 &  & $E_1:$ blowup of $V_4$ along a conic & \\ \hline
2-17 & $24$ & $E_1:$ blowup of $\P^3$ along an elliptic curve of degree $5$ &  none\\
 &  & $E_1:$ blowup of $Q$ along an elliptic curve of degree $5$ & \\ \hline
2-18 & $24$ & $X$ is a split double cover of $\mathbb{P}^2\times \mathbb{P}^1$ 
      with $\mathcal L \simeq \MO_{\P^2 \times \P^1}(1, 1)$. & 3-4\\
 &  & $C_1: \deg \Delta =4$, $X \xrightarrow{{\rm 2:1}} \P^2 \times \P^1 \xrightarrow{\pr_1} \P^2$  & \\ 
&  & $D_2: (-K_X)^2 \cdot X_t = 8$ & \\ \hline
2-19 & $26$ & $E_1:$ blowup of $\P^3$ along a curve of genus $2$ and degree $5$ &  none\\ 
 &  & $E_1:$ blowup of $V_4$ along a line & \\ \hline
2-20 & $26$ & $C_1: \deg \Delta=3$ & none \\ 
&  & $E_1:$  blowup of $V_5$ along a cubic rational curve & \\ \hline
2-21 & $28$ & 
$E_1:$ blowup of $Q$ along a rational curve of degree $4$ &  none\\
  &  & 
  $E_1:$ blowup of $Q$ along a rational curve of degree $4$ & \\ \hline
2-22 & $30$ & $E_1:$ blowup of $\P^3$ along a rational curve of degree $4$ &  none\\ 
&  & $E_1:$ blowup of $V_5$ along a conic & \\ \hline
2-23 & $30$ & 
$E_1:$ blowup of $Q$ along an elliptic curve of degree $4$ &  none\\ 
&  & $E_3\,{\rm or}\,E_4$ &\\ \hline
2-24 & $30$ & $X$ is a  divisor on $\mathbb{P}^2\times \mathbb{P}^2$ of bidegree $(1,2)$. & 3-8\\
 &  & $C_1: \deg \Delta = 3$ & \\ 
&  & $C_2$ & \\ \hline
2-25 & $32$ & $D_2: (-K_X)^2 \cdot X_t =8$ & 3-6, 3-11\\ 
 &  & $E_1:$ blowup of $\P^3$ along an elliptic curve of degree $4$  which is a complete intersection of two quadric  surfaces & \\ \hline
2-26 & $34$ & $E_1:$ blowup of $Q$ along a cubic rational curve  &  none\\ 
& & $E_1:$ blowup of $V_5$ along a line  & \\ \hline
2-27 & $38$ &$C_2$ &  3-12, 3-16\\ 
 &  & $E_1:$ blowup of $\P^3$ along a cubic rational curve & \\ \hline
2-28 & $40$ & 
$E_1:$ blowup of $\P^3$ along an elliptic curve of degree $3$ &  none\\ 
&  &$E_5$ & \\ \hline
2-29 & $40$ & $D_2: (-K_X)^2 \cdot X_t =8$ & 3-10,  3-15, 3-18\\ 
 &  & $E_1:$ blowup of $Q$ along a conic  which is a complete intersection of two members of $|\MO_Q(1)|$ & \\ \hline
2-30 & $46$ & $E_1:$ blowup of $\P^3$ along a conic & 3-18,  3-23\\ 
&  & $E_2:$ blowup of $Q$ at a point & \\ \hline
2-31 & $46$ & $C_2$ & 3-15, 3-20, 3-23 \\ 
&  & $E_1:$ blowup of $Q$ along a line & \\ \hline
2-32 & $48$ & $X$ is a divisor $W$ on $\P^2 \times \P^2$ of bidegree $(1, 1)$ & 3-7, 3-13, 3-16,\\
 &  & $C_2: W \hookrightarrow \P^2 \times \P^2 \xrightarrow{\pr_1} \P^2$ & 3-20, 3-24 \\ 
 &  & $C_2: W \hookrightarrow \P^2 \times \P^2 \xrightarrow{\pr_1} \P^2$ & \\ \hline
2-33 & $54$ & $D_3: (-K_X)^2 \cdot X_t =9$ & 3-6, 3-12, 3-18\\ 
 &  & $E_1:$ blowup of $\P^3$ along a line &3-25, 3-30 \\ \hline
2-34 & $54$ & $X=\P^2 \times \P^1$. & 3-3, 3-5, 3-7,   \\
 & &  $C_2:$ the projection $\P^2 \times \P^1 \to \P^2$ & 3-8, 3-11,3-12, \\ 
 &  & $D_3:$ the projection $\P^2 \times \P^1 \to \P^1$ &   3-15, 3-17, 3-21,\\ 
 & & &  3-22, 3-24, \\
  & & & 3-26, 3-28\\ \hline
2-35 & $56$ & $X=V_7 =\P_{\P^2}(\MO_{\P^2} \oplus \MO_{\P^2}(1))$.&  3-11, 3-14, 3-16,\\
 & & $C_2:$ the projection $\mathbb{P}_{\P^2} (\MO_{\mathbb{P}^2}\oplus \MO_{\mathbb{P}^2}(1)) \to \P^2$ & 
  3-19, 3-23, 3-26,\\ 
 &  & $E_2:$ blowup of $\P^3$ at a point & 3-29, 3-30\\ \hline
2-36 & $62$  & 
$X=\mathbb{P}_{\P^2}(\MO_{\mathbb{P}^2}\oplus \MO_{\mathbb{P}^2}(2))$. & 3-9, 3-14, \\
& & 
$C_2:$ the projection $\mathbb{P}_{\P^2}(\MO_{\mathbb{P}^2}\oplus \MO_{\mathbb{P}^2}(2)) \to \P^2$  & 3-22, 3-29\\ 
&  &  $E_5:$ blowup at the singular point of the cone over the Veronese surface & \\ \hline
      \end{longtable}
  \end{center} 

\subsection{$\rho=3$}\label{ss-table-pic3}

The definition of 3-xx is given in Definition \ref{d-pic3}.  
For  a Fano threefold $X$ with $\rho(X)=3$, 
 one and only one of the following possibilities listed in Table \ref{table-pic3}  occurs up to isomorphisms (Theorem \ref{t-pic3-main}, Remark \ref{r column blowup}). 
 \begin{enumerate}
 \item For the case when $\rho(X)=3$ (Table \ref{table-pic3}), we determine 
 the number of extremal rays and their types. 
For example, if $\rho(X)=3$ and $X$ is none of  3-9, 3-14, 
nor 3-19 (resp. one of 3-9, 3-14, 
and 3-19), then $X$ has exactly three (resp. four) extremal rays. 
 \item 
For an extremal ray $R$ of $X$, 
$X \to Y$ denotes its contraction and we use the following terminologies. 
\begin{itemize}
\item If $R$ is of type $C$, then 
\lq\lq $/S$" means that 
the target $Y$ of the contraction of $R$ is isomorphic to $S$. 
\item If $R$ is of type $E_1$ or $E_2$, then 
\lq\lq 2-xx"  (resp. \lq\lq non-Fano")  means  that $Y$ is a Fano threefold of No.\ 2-xx (resp. non-Fano). Moreover, if $R$ is of type $E_1$, then $C$ denotes the blowup centre of $X \to Y$. 
\end{itemize}
\item As for the  column \lq\lq conic bdl/$\P^2$", 
we use the following terminologies. 
\begin{itemize}
\item \lq\lq None" means that $X$ has no conic bundle structure over $\P^2$. 
\item \lq\lq 2-xx-vs-2-yy" means that $X$ has a conic bundle structure over $\P^2$ 
of type 2-xx-vs-2-yy. 
\end{itemize}
For example, a Fano threefold of No.\ 3-3 has exactly one conic bundle structure and it is of type 2-34-vs-2-34 (Proposition \ref{p-pic3-3}). 
A Fano threefold of No.\ 3-20 has exactly two conic bundle structures and 
both of them are of type 2-31-vs-2-32 (Proposition \ref{p-pic3-20}). 
 \end{enumerate}



    \begin{center}
\begin{longtable}{ccp{7.5cm}cp{1.8cm}ccc}
      \caption{$\rho(X)=3$}\label{table-pic3}\\
No. & $(-K_X)^3$ & descriptions and extremal rays & conic bdl$/\P^2$  & blowups   \\ \hline
\hyperref[p-pic3-1]{3-1}\label{table-3-1}
& $12$ & 
$X$ is a split double cover of $\P^1 \times \P^1 \times \P^1$ 
with $\mathcal L   \simeq \MO_{\P^1 \times \P^1 \times \P^1}(1, 1, 1)$.  \\
&&    $C_1:$  $/\P^1 \times \P^1$, $\deg \Delta = (4, 4)$  & none & none 
    \\
 &  & $C_1:$ $/\P^1 \times \P^1$, $\deg \Delta = (4, 4)$   &  &  &  \\ 
 &  & $C_1:$ $/\P^1 \times \P^1$, $\deg \Delta = (4, 4)$   &  &  & \\ \hline
\hyperref[p-pic3-2]{3-2}\label{table-3-2} 
& $14$ &
$X$ is a member of $|\MO_P(2)\otimes\pi^*\MO_{\mathbb{P}^1\times \mathbb{P}^1}(2,3)|$ on the $\mathbb{P}^2$-bundle $\pi\colon P=\mathbb{P}(\MO_{\mathbb{P}^1\times \mathbb{P}^1}\oplus \MO_{\mathbb{P}^1\times \mathbb{P}^1}(-1,-1)^{\oplus 2})\to \mathbb{P}^1\times \mathbb{P}^1$. \\
 &   & $C_1:$ $/\P^1 \times \P^1$, $\deg \Delta = (2, 5)$ &none & none\\
 &  & $E_1:$ non-Fano, $p_a(C)=0, -K_Y \cdot C =0$  & & \\ 
 &  & $E_1:$ non-Fano, $p_a(C)=0, -K_Y \cdot C =0$ &  & \\ \hline
\hyperref[p-pic3-3]{3-3}\label{table-3-3}  & $18$ & $X$ is a divisor on $\P^1 \times \P^1 \times \P^2$ of tridegree $(1, 1, 2)$.\\
& & $C_1: /\P^1 \times \P^1, \deg \Delta = (3, 3)$  &  2-34-vs-2-34 & none\\ 
 &  & $E_1:$ 2-34, $p_a(C)=3, -K_Y \cdot C =20$  & & \\
  &  & $E_1:$ 2-34, $p_a(C)=3, -K_Y \cdot C =20$  &  & \\\hline
\hyperref[p-pic3-4]{3-4}\label{table-3-4}  & $18$ & 
$X$ is a blowup of $Y_{\text{2-18}}$ along a smooth fibre of the contraction $Y_{\text{2-18}} \to \P^2$.
\\
&    & $C_1:$ $/\P^1 \times \P^1$, $\deg \Delta =(2, 4)$ & none & none\\ 
 & &  $C_1:$ $/\F_1$, $\Delta \in |\tau^*\MO_{\P^2}(4)|$ &  & \\
 & &  $E_1:$ 2-18, $p_a(C)=0, -K_Y \cdot C =2$ &  & \\ \hline
\hyperref[p-pic3-5]{3-5}\label{table-3-5}  & $20$ & 
$X$ is a blowup of $\P^2 \times \P^1$ along a smooth curve $C$ 
of bidegree $(2, 5)$ such that $C \hookrightarrow \P^2 \times \P^1 \xrightarrow{\pr_1} \P^2$ is a closed immersion.\\
&& $E_1:$ 2-34, $p_a(C)=0, -K_Y \cdot C =16$   & 2-34-vs-non-Fano & none
\\
  &  & $E_1:$ non-Fano, $p_a(C)=0, -K_Y \cdot C =0$   &   & \\
    &  & $E_1:$ non-Fano, $p_a(C)=0, -K_Y \cdot C =0$  &  & \\\hline
\hyperref[p-pic3-6]{3-6}\label{table-3-6}  & $22$ &  $X$ is a blowup of $\P^3$ along a disjoint union of a line and an elliptic curve of degree four.&  \\
&&$C_1: /\P^1 \times \P^1, \deg \Delta =(3, 2)$ & none & none \\ 
&  & $E_1:$ 2-25, $p_a(C)=1, -K_Y \cdot C =16$ &  && \\
 &  & $E_1:$ 2-33, $p_a(C)=0, -K_Y \cdot C =4$ &  & &\\\hline
\hyperref[p-pic3-7]{3-7}\label{table-3-7}  & $24$ & 
$X$ is a blowup of $W$ along an elliptic curve which is a complete intersection of two members of $|-\frac{1}{2}K_W|$. & &
\\
&& $E_1:$ 2-32, $p_a(C)=1, -K_Y \cdot C =12$  &   2-32-vs-2-34 &  none\\
  &  & $E_1:$ 2-34, $p_a(C)=1, -K_Y \cdot C =15$  &  2-32-vs-2-34  & \\
    &  & $E_1:$ 2-34, $p_a(C)=1, -K_Y \cdot C =15$   &  & \\\hline
\hyperref[p-pic3-8]{3-8}\label{table-3-8}  & $24$ & 
$X$ is 
a divisor on $\F_1 \times \P^2$ 
which is a member of $|\pr_1^*\tau^*\MO_{\P^2}(1) \otimes \pr_2^*\MO_{\P^2}(2)|$. && \\
& 
    & $C_1:$ $/\F_1$, $\Delta \in |\tau^*\MO_{\P^2}(3)|$ & 2-24-vs-2-34&  none \\ 
 & &  $E_1:$ 2-24, $p_a(C)=0, -K_Y \cdot C =2$ &  & \\
 & &  $E_1:$ 2-34, $p_a(C)=0, -K_Y \cdot C =14$ &  & \\ \hline
\hyperref[p-pic3-9]{3-9}\label{table-3-9} 
& $26$ & $X$ is a blowup of $\P_{\P^2}(\MO \oplus \MO(2))$ along a smooth curve $C$ on a section $S$ of the $\P^1$-bundle $\pi : \P_{\P^2}(\MO \oplus \MO(2)) \to \P^2$ such that $\pi(C)$ is a quartic curve. & & \\
&&$E_1:$ 2-36, $p_a(C)=3, -K_Y \cdot C =20$ & 2-36-vs-2-36  & none\\ 
 &  & $E_1:$ 2-36, $p_a(C)=3, -K_Y \cdot C =20$   &  & \\
  &  & $E_5$   &  & \\
   &  & $E_5$   &  & \\ \hline
\hyperref[p-pic3-10]{3-10}\label{table-3-10} & $26$ & 
$X$ is a blowup of $Q$ along a disjoint union of two conics. \\
&&$C_1: /\P^1 \times \P^1, \deg \Delta =(2, 2)$ & none &  none & \\ 
 &  & $E_1:$ 2-29, $p_a(C)=0, -K_Y \cdot C =6$  &  &  & \\
 &  & $E_1:$ 2-29, $p_a(C)=0, -K_Y \cdot C =6$  &  &   & \\\hline
\hyperref[p-pic3-11]{3-11}\label{table-3-11} & $28$ & 
$X$ is a  blowup of $V_7$ along an elliptic curve which is a complete intersection of two members of $|-\frac{1}{2}K_{V_7}|$. \\
&& $E_1:$ 2-25, $p_a(C)=0, -K_Y \cdot C =1$  & 2-34-vs-2-35  & none \\ 
&  & $E_1:$ 2-34, $p_a(C)=1, -K_Y \cdot C =13$   &  & \\
&  & $E_1:$ 2-35, $p_a(C)=1, -K_Y \cdot C =14$   &  & \\\hline
\hyperref[p-pic3-12]{3-12}\label{table-3-12}
& $28$ &  $X$ is a blowup of $\P^3$ along a disjoint union of 
a line and a rational cubic curve.  \\
&&$E_1:$ 2-27, $p_a(C)=0, -K_Y \cdot C =4$  & 2-27-vs-2-34  & none \\ 
  &  & $E_1:$ 2-33, $p_a(C)=0, -K_Y \cdot C =12$   &  & \\
    &  & $E_1:$ 2-34, $p_a(C)=0, -K_Y \cdot C =12$    &  & \\\hline
\hyperref[p-pic3-13]{3-13}\label{table-3-13} & $30$ & 
$X$ is a blowup of $W$ along a curve $C$ of bidegree $(2, 2)$ 
such that $W \hookrightarrow \P^2 \times \P^2 \xrightarrow{\pr_i} \P^2$ is a closed immersion 
for each $i\in \{1, 2\}$.\\
&&$E_1:$ 2-32, $p_a(C)=0, -K_Y \cdot C =8$  & 2-32-vs-2-32  & none \\
&  &$E_1:$ 2-32, $p_a(C)=0, -K_Y \cdot C =8$   & 2-32-vs-2-32 & \\
&  & $E_1:$ 2-32, $p_a(C)=0, -K_Y \cdot C =8$   & 2-32-vs-2-32 & \\ \hline
\hyperref[p-pic3-14]{3-14}\label{table-3-14}
& $32$ & $X$ is a blowup of $\P^3$ along a disjoint union of a point $P$ and a plane cubic curve $C$, where $P$ is not contained in the plane containing $C$.\\
&& $E_1:$ 2-35, $p_a(C)=1, -K_Y \cdot C =12$  & 2-35-vs-2-36  & none \\ 
 &  & $E_1:$ 2-36, $p_a(C)=1, -K_Y \cdot C =15$   &  & \\
  &  & $E_2:$ 2-28   &  & \\
   &  & $E_5$   &  & \\\hline
\hyperref[p-pic3-15]{3-15}\label{table-3-15}
& $32$ & $X$ is a blowup of $Q$ along a disjoint union of a line and a conic. \\
&& $E_1:$ 2-29, $p_a(C)=0, -K_Y \cdot C =3$ & 2-31-vs-2-34  & none \\ 
 &  & $E_1:$ 2-31, $p_a(C)=0, -K_Y \cdot C =6$ &  & \\
  &  & $E_1:$ 2-34, $p_a(C)=0, -K_Y \cdot C =10$   &  & \\\hline
\hyperref[p-pic3-16]{3-16}\label{table-3-16} & $34$ & 
$X$ is a blowup of $V_7$ along  the strict transform of a smooth cubic rational curve passing through the blowup centre of $V_7 \to \P^3$. \\
&&$E_1:$ 2-27, $p_a(C)=0, -K_Y \cdot C =1$ & 2-27-vs-2-32  & none \\ 
 &  & $E_1:$ 2-32, $p_a(C)=0, -K_Y \cdot C =6$    & 2-32-vs-2-35  &\\ 
 &  & $E_1:$ 2-35, $p_a(C)=0, -K_Y \cdot C =10$  &  & \\ \hline
\hyperref[p-pic3-17]{3-17}\label{table-3-17} & $36$ & 
$X$ is a divisor on $\P^1 \times \P^1 \times \P^2$ of tridegree $(1, 1, 1)$.
 && \\
&& $C_2: /\P^1 \times \P^1$ &  2-34-vs-2-34 & 4-3\\ 
 &  & $E_1:$ 2-34, $p_a(C)=0, -K_Y \cdot C =8$  &  & \\
  &  & $E_1:$ 2-34, $p_a(C)=0, -K_Y \cdot C =8$  &  & \\\hline
\hyperref[p-pic3-18]{3-18}\label{table-3-18} & $36$ & 
$X$ is a blowup of $\P^3$ along a disjoint union of a line and a conic.&& \\
&&$E_1:$ 2-29, $p_a(C)=0, -K_Y \cdot C =1$  & none  & 4-4&  \\ 
  &  & $E_1:$ 2-30, $p_a(C)=0, -K_Y \cdot C =4$  &  & \\
    &  & $E_1:$ 2-33, $p_a(C)=0, -K_Y \cdot C =8$  &  & \\\hline
\hyperref[p-pic3-19]{3-19}\label{table-3-19} & $38$ & 
$X$ is a blowup of $Q$ along a disjoint union of two points which are not collinear. && \\
&& $E_1:$ 2-35, $p_a(C)=0, -K_Y \cdot C =8$  & 2-35-vs-2-35  & 4-4\\
 &  & $E_1:$ 2-35, $p_a(C)=0, -K_Y \cdot C =8$   &  & \\
  &  &$E_2:$ 2-30   &  & \\
   &  & $E_2:$ 2-30   &  & \\\hline
\hyperref[p-pic3-20]{3-20}\label{table-3-20} & $38$ & $X$ is a blowup of $Q$ along a disjoint union of two lines. \\
&& $E_1:$ 2-31, $p_a(C)=0, -K_Y \cdot C =3$   & 2-31-vs-2-32  & none\\ 
&  &$E_1:$ 2-31, $p_a(C)=0, -K_Y \cdot C =3$   & 2-31-vs-2-32  & \\
&  & $E_1:$ 2-32, $p_a(C)=0, -K_Y \cdot C =4$ &  & \\
\hline
\hyperref[p-pic3-21]{3-21}\label{table-3-21} & $38$ & $X$ is a blowup of $\P^2 \times \P^1$ along a curve of bidegree $(1, 2)$.&& \\
&& $E_1:$ 2-34, $p_a(C)=0, -K_Y \cdot C =7$  & 2-34-vs-non-Fano  & 4-5 \\ 
 &  & $E_1:$ non-Fano, $p_a(C)=0, -K_Y \cdot C =0$   &  & \\
  &  &$E_1:$ non-Fano, $p_a(C)=0, -K_Y \cdot C =0$   &  & \\\hline
\hyperref[p-pic3-22]{3-22}\label{table-3-22} & $40$ & $X$ is a blowup of $\P^2 \times \P^1$ along a conic on a plane $\P^2 \times \{ t\}$ for some closed point $t \in \P^1$. \\
&&$E_1:$ 2-34, $p_a(C)=0, -K_Y \cdot C =6$  &  2-34-vs-2-36 & none \\ 
 &  & $E_1:$ 2-36, $p_a(C)=0, -K_Y \cdot C =10$  &  & \\
  &  & $E_5$   &  & \\\hline
\hyperref[p-pic3-23]{3-23}\label{table-3-23} & $42$ & 
$X$ is a blowup of $V_7$ along  the strict transform of a conic passing through the blowup centre of $V_7 \to \P^3$.  \\
&& $E_1:$ 2-30, $p_a(C)=0, -K_Y \cdot C =1$  & 2-31-vs-2-35  & none \\ 
 &  & $E_1:$ 2-31, $p_a(C)=0, -K_Y \cdot C =1$  &  & \\
  &  & $E_1:$ 2-35, $p_a(C)=0, -K_Y \cdot C =6$   &  & \\\hline
\hyperref[p-pic3-24]{3-24}\label{table-3-24} & $42$ & $X =W \times_{\P^2} \F_1$ for a contraction $W \to \P^2$ and the blowdown $\tau : \F_1 \to \P^2$. && \\
&&$C_2:$ $/\F_1$   & 2-32-vs-2-34 & 4-7\\ 
 &  & $E_1:$ 2-32, $p_a(C)=0, -K_Y \cdot C =2$  &  & \\ 
 &  & $E_1:$ 2-34, $p_a(C)=0, -K_Y \cdot C =5$  &  & \\ \hline
\hyperref[p-pic3-25]{3-25}\label{table-3-25} & $44$ & $X$ is a blowup of $\P^3$ along a disjoint union of two lines. 
& & \\
&&$C_2:$ $/\P^1 \times \P^1$  & none & 4-6,  4-9 & \\
 &  & $E_1:$ 2-33, $p_a(C)=0, -K_Y \cdot C =4$  &  &  & \\ 
 &  & $E_1:$ 2-33, $p_a(C)=0, -K_Y \cdot C =4$  &  &  & \\ \hline
\hyperref[p-pic3-26]{3-26}\label{table-3-26} & $46$ & 
$X$ is a blowup of $\P^3$ along a disjoint union of a point and a line. && 
\\
&&$E_1:$ 2-34, $p_a(C)=0, -K_Y \cdot C =3$  & 2-34-vs-2-35  & 4-9\\ 
 &  & $E_1:$ 2-35, $p_a(C)=0, -K_Y \cdot C =4$   &  & \\
  &  & $E_2:$ 2-33   &  & \\\hline
  \hyperref[p-pic3-27]{3-27}\label{table-3-27}& $48$ & $X = \P^1 \times \P^1 \times \P^1$. && \\
&&$C_2:$ $/\P^1 \times \P^1$ & none &  4-1, 4-3,  & \\ 
&  
& $C_2:$ $/\P^1 \times \P^1$  && 4-6, 4-8,  & \\
 &  & $C_2:$ $/\P^1 \times \P^1$ & & 4-10, 4-13 & \\ \hline
\hyperref[p-pic3-28]{3-28}\label{table-3-28} & $48$ & $X = \F_1 \times \P^1$. & & \\
&& $C_2:$ $/\P^1 \times \P^1$  & none & 4-3, 4-5, \\ 
 &  & $C_2:$ $/\F_1$ &  & 4-7, 4-9, \\ 
 &  & $E_1:$ 2-34, $p_a(C)=0, -K_Y \cdot C =2$  &  & 4-10, 4-11\\ \hline
\hyperref[p-pic3-29]{3-29}\label{table-3-29} & $50$ & $X$ is a blowup of $V_7$ along a line on the exceptional divisor of 
the blowup $V_7 \to \P^3$.\\
&& $E_1$: 2-35, $p_a(C)=0, -K_Y \cdot C =2$  & 2-35-vs-2-36  & none\\ 
&  & $E_1:$ 2-36, $p_a(C)=0, -K_Y \cdot C =5$   &  & \\
&  & $E_5$   &  & \\\hline
\hyperref[p-pic3-30]{3-30}\label{table-3-30} & $50$ & $X$ is a blowup of $V_7$ along  the strict transform of a line passing through the blowup centre of $V_7 \to \P^3$. &&
\\
&& $C_2:$ $/\F_1$  &  none &4-4, 4-9, \\ 
 &  & $E_1:$ 2-33, $p_a(C)=0, -K_Y \cdot C =1$  &  & 4-12 \\
 &  & $E_1:$ 2-35, $p_a(C)=0, -K_Y \cdot C =2$  &  & \\\hline
\hyperref[p-pic3-31]{3-31}\label{table-3-31} & $52$ & $X =\P_{\P^1 \times \P^1}(\MO_{\P^1 \times \P^1} \oplus \MO_{\P^1 \times \P^1}(1, 1))$. && \\
&& $C_2: /\P^1 \times \P^1$   & none &4-2, 4-5, \\  
&  & $E_1:$ non-Fano, $p_a(C)=0, -K_Y \cdot C =0$  &  & 4-8, 4-11, \\
 &  & $E_1:$ non-Fano, $p_a(C)=0, -K_Y \cdot C =0$  &  & 4-13\\ \hline
      \end{longtable}
  \end{center}

\subsection{$\rho=4$}\label{ss-table-pic4}

The definition of 4-xx is given in Definition \ref{d-pic4}. 
For  a Fano threefold $X$ with $\rho(X)=4$, 
one and only one of the following possibilities listed in Table \ref{table-pic4}  occurs up to isomorphisms (Theorem \ref{t-pic4-main}, Remark \ref{r column blowup}). 

As for the  column \lq\lq descriptions and conic bundles", 
we use the following terminologies. 
\begin{enumerate}
\item \lq\lq 3-xx-vs-3-yy /$\P^1 \times \P^1$" 
(resp. \lq\lq 3-xx-vs-3-yy /$\F_1$" )
means that $X$ has a conic bundle structure over $\P^1 \times \P^1$ (resp. $\F_1$)  
of type 2-xx-vs-2-yy. 
Moreover, the column \lq\lq $\Delta$" gives the bidegree (resp. the linear equivalence) of its discriminant divisor. 
\item \lq\lq $(\text{3-zz}) \times_{\P^2} \F_1$" means that $X \simeq Y \times_{\P^2} \F_1$, 
where $Y$ is a Fano threefold of No.\ 3-zz, 
$\pi : Y \to \P^2$ is a conic bundle, and 
$\F_1 \to \P^2$ is a blowup at a point outside $\Delta_{\pi}$. 
\end{enumerate}
For example, if $X$ is a Fano threefold $X$ of No.\ 4-5, then 
there exist conic bundles 
$f : X \to \P^1 \times \P^1$
and $g: X \to \F_1$ such that 
$f$ is of type 3-28-vs-3-31, $\Delta_f$ is of bidegree $(1, 2)$, 
$g$ is of type 3-24-vs-non-Fano, and $\Delta_g \sim \tau^*\MO_{\P^2}(1)$. 


    \begin{center}
\begin{longtable}{ccp{7.5cm}ccp{1.9cm}}
      \caption{$\rho(X)=4$}\label{table-pic4}\\
No. & $(-K_X)^3$ & descriptions and conic bundles  & $\Delta$  & blowups & blowdowns \\ \hline
\hyperref[p-pic4-1]{4-1}\label{table-4-1}
& $24$ & $X$ is a divisor on $\P^1 \times \P^1 \times \P^1 \times \P^1$ of multi-degree $(1, 1, 1, 1)$. & \\
&    & 3-27-vs-3-27 $/\P^1 \times \P^1$ & $(2, 2)$ &  none & 3-27
   \\ \hline
\hyperref[p-pic4-2]{4-2}\label{table-4-2} & $28$  &
$X$ is a blowup along an elliptic curve $C$ on a section $T$ of 
the $\P^1$-bundle 
$\pi : Y_{\text{3-31}} = \P_{\P^1 \times \P^1}(\MO \oplus \MO(1, 1)) \to \P^1 \times \P^1$ 
such that $\pi(C)$ is of bidegree $(2, 2)$. &
\\
&  & 3-31-vs-3-31 $/\P^1 \times \P^1$  &  $(2, 2)$&  none & 3-31
   \\ \hline
\hyperref[p-pic4-3]{4-3}\label{table-4-3} & $30$ & $X$ is a blowup along a curve on $\P^1 \times \P^1 \times \P^1$ of tridegree $(1, 1, 2)$.
&& \\
&    & 3-17-vs-3-27 $/\P^1 \times \P^1$ & $(1, 1)$ &  
none & 3-17, 3-27,   \\
   &  & 3-27-vs-3-28 $/\P^1 \times \P^1$   &   $(1, 1)$&  & 3-28 \\ 
      &  & 3-28-vs-3-28 $/\F_1$, $(3\text{-}17)\times_{\P^2} \F_1$  &   $\tau^*\MO_{\P^2}(2)$& \\ \hline
\hyperref[p-pic4-4]{4-4}\label{table-4-4} & $32$ &
$X$ is a blowup of $Y_{\text{2-29}}$ along $B_1 \amalg B_2$, 
where 
$B_1$ and $B_2$ are mutually distinct fibres of the blowup 
$Y_{\text{2-29}} \to Q$ along a conic. & \\
& 
    &  3-30-vs-3-30 $/\F_1$, $(3\text{-}19)\times_{\P^2} \F_1$ & $\tau^*\MO_{\P^2}(2)$&5-1 &
3-18, 3-19,   \\ 
&&&&& 3-30 \\
\hline
\hyperref[p-pic4-5]{4-5}\label{table-4-5} & $32$ & 
$X$ is a blowup of $\P^2 \times \P^1$ along $C_1 \amalg C_2$, 
where $C_1$ and $C_2$ are mutually disjoint smooth curves of bidegree $(0, 1)$ and $(1, 2)$. 
& \\
&
    & 3-28-vs-3-31 $/\P^1 \times \P^1$  & $(1, 2)$ &none & 3-21, 3-28,
      \\ 
     &  & 3-28-vs-non-Fano $/\F_1$, $(3\text{-}21)\times_{\P^2} \F_1$   &  $\tau^*\MO_{\P^2}(1)$ && 3-31\\
     \hline
\hyperref[p-pic4-6]{4-6}\label{table-4-6} & $34$ & 
$X$ is a blowup of $\P^1 \times \P^1 \times \P^1$ along a curve of tridegree $(1, 1, 1)$.& \\
&
    & 3-25-vs-3-27 $/\P^1 \times \P^1$ &  $(1, 1)$
    & none & 3-25, 3-27   \\ \hline
\hyperref[p-pic4-7]{4-7}\label{table-4-7} & $36$ & 
$X$ is a blowup of $W$ along $C_1 \amalg C_2$, 
where $C_1$ and $C_2$ are mutually disjoint  curves of bidegree $(1, 0)$ and $(0, 1)$&. \\
&
    & 3-28-vs-3-28 $/\P^1 \times \P^1$ &  $(1, 1)$& none & 3-24, 3-28  \\ 
    &  &  3-24-vs-3-28 $/\F_1$, $(3\text{-}24)\times_{\P^2} \F_1$  &  $\tau^*\MO_{\P^2}(1)$ & \\
      \hline
\hyperref[p-pic4-8]{4-8}\label{table-4-8} & $38$ & 
$X$ is a blowup of $\P^1 \times \P^1 \times \P^1$ along a curve of tridegree $(0, 1, 1)$. & \\
&
    & 3-27-vs-3-31 $/\P^1 \times \P^1$ &  $(1, 1)$&none & 3-27, 3-31  \\ 
    & & 3-27-vs-non-Fano $/\P^1 \times \P^1$ &  $(0, 1)$&  
    &   \\ \hline
\hyperref[p-pic4-9]{4-9}\label{table-4-9} & $40$ &
$X$ is a blowup of $Y_{\text{3-25}}$ along $C$, where 
$\rho : Y_{\text{3-25}} \to \P^3$ is a blowup along a disjoint union $L_1 \amalg L_2$ of two lines $L_1$ and $L_2$, and $C$ is a one-dimensional fibre of $\rho$.
& 
\\
&
    & 3-25-vs-3-28 $/\P^1 \times \P^1$  & $(0, 1)$  &5-2 & 3-25, 3-26,   \\ 
    &  & 3-28-vs-3-30 $/\F_1$, $(3\text{-}26)\times_{\P^2} \F_1$   &   $\tau^*\MO_{\P^2}(1)$&&  3-28, 3-30 \\ \hline
\hyperref[p-pic4-10]{4-10}\label{table-4-10} & $42$ & $X = S_7 \times \P^1$. 
&& \\
&
    & 3-27-vs-3-28 $/\P^1 \times \P^1$  &  $(0, 1)$ & 
    5-3 & 3-27, 3-28   \\ \hline
\hyperref[p-pic4-11]{4-11}\label{table-4-11} & $44$ & 
$X$ is a blowup of $\F_1 \times \P^1$ along $C = \Gamma \times \{t\}$, 
where $\Gamma$ is the $(-1)$-curve on $\F_1$ and $t$ is a point on $\P^1$. 
&&\\
&
    & 3-28-vs-3-31 $/\P^1 \times \P^1$  &  $(0, 1)$& 5-2 & 3-28, 3-31  \\ 
      &  & 3-28-vs-non-Fano $/\F_1$  & $(-1)$-curve & \\
      \hline
\hyperref[p-pic4-12]{4-12}\label{table-4-12} & $46$ & 
$X$ is a blowup of $Y_{\text{2-33}}$ along $C_1 \amalg C_2$, where 
$\rho : Y_{\text{2-33}} \to \P^3$ is a blowup along a line $L$, and 
$C_1$ and $C_2$ are mutually distinct one-dimensional fibres of $\rho$. 
&& \\
&
    & 3-30-vs-non-Fano $/\F_1$ & $(-1)$-curve & 
5-1, 5-2& 3-30   \\ \hline
\hyperref[p-pic4-13]{4-13}\label{table-4-13} & $26$ & 
$X$ is a blowup of $\P^1 \times \P^1 \times \P^1$ along a curve of tridegree $(1, 1, 3)$. && \\
&
    & 3-27-vs-3-31 $/\P^1 \times \P^1$ &  $(1, 3)$& none & 3-27, 3-31 \\
    & & 3-27-vs-non-Fano $/\P^1 \times \P^1$ &  $(1, 1)$&  
    &   \\ \hline
      \end{longtable}
  \end{center}

\subsection{$\rho=5$}\label{ss-table-pic5}

The definition of 5-xx is given in Definition \ref{d-pic5}.  
    For  a Fano threefold $X$ with $\rho(X)=5$, 
one and only one of the following possibilities listed in Table \ref{table-pic5}  occurs up to isomorphisms (Theorem \ref{t-pic5-main}, Remark \ref{r column blowup}). 
 
    \begin{center}
\begin{longtable}{ccp{7.5cm}cp{1.9cm}}
      \caption{$\rho(X)=5$}\label{table-pic5}\\
No. & $(-K_X)^3$ & descriptions and conic bundles   & blowups & blowdowns \\ \hline
\hyperref[p-pic5-1]{5-1}\label{table-5-1} & $28$ & $X$ is a blowup of $Y_{\text{2-29}}$ 
along $B_1 \amalg B_2 \amalg B_3$, 
where $\rho : Y_{\text{2-29}} \to Q$ is a blowup along a conic and 
$B_1, B_2, B_3$ are mutually distinct one-dimensional fibres of $\rho$.   & none & 4-4, 4-12 \\ \hline
\hyperref[p-pic5-2]{5-2}\label{table-5-2}  & $36$  &
$X$ is a blowup of $Y_{\text{3-25}}$ 
along $B \amalg B'$, 
where $\rho : Y_{\text{3-25}} \to \P^3$ 
is a blowup along a disjoint union $L_1 \amalg L_2$ of 
 lines $L_1$ and $L_2$, and 
both $B$ and $B'$ are mutually distinct one-dimensional fibres of $\rho$  which are lying over $L_1$.  & none & 4-9, 4-11, 4-12
\\ \hline
\hyperref[p-pic5-3]{5-3}\label{table-5-3}  & $36$ & $X = S_6 \times \P^1$. 
& 6-1& 4-10\\
\hline
      \end{longtable}
  \end{center}

\subsection{$\rho \geq 6$}\label{ss-table-pic6}

For  a Fano threefold $X$ with $\rho(X)\geq 6$, 
one and only one of the following possibilities listed in Table \ref{table-pic6}  occurs up to isomorphisms (Theorem \ref{t-pic6}, Remark \ref{r column blowup}). 

    
    \begin{center}
\begin{longtable}{ccp{4.5cm}cp{1.9cm}}
      \caption{$\rho(X)\geq 6$}\label{table-pic6}\\
No. & $(-K_X)^3$ & descriptions   & blowups & blowdowns \\ \hline
6-1 & $30$ & $X = S_5 \times \P^1$. 
& 7-1& 5-3\\
\hline
7-1 & $24$ & $X = S_4 \times \P^1$.
& 8-1& 6-1\\
\hline
8-1 & $18$ & $X = S_3 \times \P^1$. 
& 9-1& 7-1\\
\hline
9-1 & $12$ & $X = S_2 \times \P^1$. 
& 10-1& 8-1\\
\hline
10-1 & $6$ & $X = S_1 \times \P^1$. 
& none & 9-1\\
\hline
      \end{longtable}
  \end{center}

\bibliographystyle{skalpha}
\bibliography{reference.bib}

\end{document}